\documentstyle[fleqn]{article}

\newtheorem{lem}{Lemma}[section]
\newtheorem{thm}[lem]{Theorem}
\newtheorem{cor}[lem]{Corollary}

\newtheorem{prop}[lem]{Proposition}
\newtheorem{fact}[lem]{Fact}
\newtheorem{definition}[lem]{Definition}

\newtheorem{note}[lem]{Note}
\newtheorem{notation}[lem]{Notation}
\newtheorem{example}[lem]{Example}
\newtheorem{rem}[lem]{Remark}
\newtheorem{remark}[lem]{Remark}

\newtheorem{lemnot}[lem]{Lemma-Notation}
\newtheorem{sta}[lem]{}
\newcommand{\proof}{\noindent \it Proof \hbox{       }}
\newcommand{\claim}{{\noindent \bf Claim }}
\def\prf{\begin{proof}}
\def\eprf{\qed \end{proof}}
\def\bY{{\bar{Y}}}
\def\by{{\bar{y}}}

 \def\lm{{\lambda}}
 \def\lemref#1{{Lemma \ref{#1}}}
\font\msbmten msbm10
\def \Bbb#1{\hbox{\msbmten #1}}
 \def\cS{{\cal S}}
 
\def\Aa{\Bbb A}
\def\Pp {\Bbb P}
\def\Qq {\Bbb Q}
\def\Zz {\Bbb Z}
\def\Ff {\Bbb F}
\def\Cc{\Bbb C}  \def\Rr{\Bbb R}
\def\Nn{\Bbb N}

  \def\etl{\acute{e}t}

\def\tS{\tilde{S}}

\def\g{\gamma}

\def\si{\sigma}

\def\iso{\simeq}

 \def\lbl#1{ \label{#1}}
\def\proof{ { {\noindent \it Proof} \hbox{       } } }

\def\tensor{{\otimes}}
\def\meet{\cap}
\def\union{\cup}

\def\Sum{\Sigma}

\def\set#1{\{#1\}}
\def\abs#1{ \left| #1 \right| }
\def\doub#1{   \abs{\abs{#1}}   }

\def\lb{{\hfill \break}}
\def\qed{\hfill $\Box$ \break}

\def\pr{\rm pr}
\def\U{{\cal U}}
\def\spe{{\rm Spec^\si \ }}
\def\spec{{\rm Spec \ }}
\def\proj{{\rm Proj \ }}
\def\projs{{\rm Proj^\si \ }}

\def\invv{{\rm inv}}
\def\inv{{^{-1}}}

\def\Zsi{{ \Zz} {\raise 1pt \hbox{$[\si,\si^{-1}]$}               }  }
\def\Nsi{{  \Nn} {\raise 1pt \hbox{ $[\si]  $     }     }                   }

 \def\hsp{\,}

\def\codim#1#2{  {\rm codim}_{#1} ({#2})  }
 \def\cdim{{\rm codim}}
\def\pr{{\rm pr \ }}
\def\half{\frac{1}{2}}

\def\implies{ \Rightarrow }

 \def\alg{{\rm alg \hsp}}
 \def\mod{{\rm mod \hsp }}
 \def\dim{ {\rm dim \hsp } }
 
 \def\pr{{\rm pr \hsp }}
 
 \def\fixs{{\rm Fix } ^\si \hsp}
 \def\Implies{ $\Rightarrow$ }

 \def\cO{{\cal O}}
 \def\spe{{\rm Spec^\si \hsp }}
 \def\spec{{\rm Spec \hsp }}
 \def\proj{{\rm Proj \hsp }}
 \def\projs{{\rm Proj^\si \hsp }}

 \def\dimk{{\rm dim}_k}

 \def\Nsi{{\Nn[\si]}}
 
 \def\<{\begin}
 \def\>{\end}
 \def\st { \begin{sta} \ \end{sta}}
 \def\deff#1{\begin{def} \ \lbl{#1} \end{def}}
 \def\sk{\smallskip}
  
 \def\dimk{{\rm dim}_k}

 \def\mult{{\rm Mult}}
 
 \def\tV{{\tilde V}}
 
\def\tilde{\widetilde}

\def\Ak{\breve{A}_k} \def\kt{{k[\breve{t}]_\si}}  \def\Ft{{F[\breve{t}]_\si}} 
\def\hT{{\breve{h_t}}} 
\def\Lt{{\breve{L_t}}} 
\def\htt{{h_t}^*}
 \def\val{{\, \rm val \,}}
  \def\Vv{{\underline V}}

\def\tC{{\tilde{C}}}
  \def\bK{{\bar{K}}}

 \def\Xvs{   {X_{\rightarrow 0}  }  }

\def\size#1{{|#1|}}
\def\dd{{\underline{d}}}

\def\res{{\rm res \,}}
\def\M{{\cal M}}
\def\OK{{\cal O}_K}
\def\Oo{{\cal O}}
\def\vrk{{rk_{\val}}}  
\def\ba{{\bar{a}}}  
\def\Z-si{{\Zz_\si}}  
\def\Q-si{{\Qq_\si}}
\def\G{{\Gamma}} 
 \def\tM{{\tilde{M}}}
   
\def\tdim{t.dim}
 
 \def\onto{{\to >}}
 
 \def\srad{   \sqrt[\si]}
  \def\rwm{ {\rm WM_{rad}}}

\def\fy{{Y^\psi}}

\def\m{{\setminus}}

\def\RES{\, {\rm RES} \,}
\def\fn{{\Phi_{\rm fn}}} \def\dimV{{\rm dim}_V }
 \newtheorem{lemdef}[lem]{Lemma and Definition }

\def\hat{\widehat}

\def\B{{\mathcal B}}
\def\CV{{\cal C}(V)}   \def\CVf{  {\cal C}_f(V)   }  
\def\CVsep{{\CV}_{\rm sep}}
 \def\ZCV{\Zz \CV}   \def\ZCVf{\Zz \CVf} 
\def\NCV{{\ZCV}_{\rm rr}}

\def\uS{{\underline{S}}}
\def\uPhi{{\underline{\Phi}}}
\def\reff#1{    \ref{#1} }
\def\lemm#1{ \begin{lem} \lbl{#1} }
\title {The Elementary Theory of the Frobenius Automorphisms}

\author{Ehud Hrushovski   \footnote{The work was carried out at the Department of Mathematics, Hebrew University, Jerusalem,  
first posted in Feb. 2004.  It was supported in part  by the Israel Science Foundation grants no. 244/03, 1048/07.  Current address:  Mathematical Institute, Oxford, hrushovski@maths.ox.ac.uk}}

\begin{document}

\maketitle   

\abstract

We lay down   elements of a geometry based on difference
equations.  Various  constructions of
algebraic geometry are shown to have meaningful analogs:  dimensions,
blowing-up, moving lemmas.  

Analogy aside, 
the geometry of difference equations has  two quite different functorial connections with   algebraic geometry.  On the one hand, a difference scheme is locally
 a scheme with an automorphism;  the scheme is usually not of finite type, but is determined by certain finite  algebro-geometric data.    On the other hand, for each prime power $p^m$, one has
a functor into algebraic schemes over $\Ff_p$, where the structure
endomorphism becomes Frobenius.    

Transformal zero-cycles have a rich structure in the new geometry.  
 In particular, the Frobenius reduction
functors show that they encapsulate data described in classical cases
by zeta or $L$-functions.   A theory of rational and algebraic equivalence
of $0$-cycles is initiated, via a study of the transformal analog of discrete
valuation rings.

Our main result is the determination of 
 the elementary theory of the class of Frobenius
difference fields (algebraically closed fields of
characteristic $p>0$, enriched with a symbol for  $x \mapsto x^{p^m}$.  
This theory coincides with the model companion ACFA of  the theory of 
difference fields.   In particular,  a sentence $P$
holds  in $(\Ff_p^a, x \mapsto x^p)$ for almost all primes $p$ if and only if it is true
  for a generic automorphism of a field of characteristic $0$; i.e. true in 
$(L,\si)$ for a co-meager set of $\si \in Aut(L)$, where 
$L$ is the field of algebraic functions in denumerably  many variables over $\Qq$. 

  The proof requires a twisted version of the Lang-Weil estimates, related
to Deligne's conjecture but less precise and valid more generally.  It is proved on the
basis of the preceding work on difference geometry.

Some applications are given, in particular to finite simple groups, and
to the Jacobi bound for difference equations.

\newpage
\tableofcontents
\newpage
\begin{section}{Introduction}

\<{subsection}{The main theorem}
     A classical result of algebraic model theory 
(\cite{ax},\cite{fried-jarden})
shows that all elementary statements about finite fields are determined
by three
facts from geometry. The fundamental fact is
 Weil's Riemann Hypothesis for curves, entering via the Lang-Weil
estimates \cite{LW}.
The auxiliary statements are the Cebotarev density theorem,
and the cyclic nature of the Galois group. All three facts are proved,
within algebraic geometry, number theory and Galois theory respectively,
by viewing a finite field as the fixed field of the Frobenius
automorphism.
One might therefore hope that the
pertinent geometry
could be used directly to derive the full elementary theory
of the Frobenius maps. The elementary theory of finite fields
 would follow as a special case.

 This hope turns out to be correct, except that the Lang-Weil estimates 
do not
suffice, and need to be replaced by a   more general principle.
Let $k$
be an algebraically closed field.
 If $V$ is a variety over $k$ and $\si$ is an automorphism
of $k$, we denote by $V^\sigma$ the variety obtained
from $V$ by applying $\sigma$ to the defining parameters.
 (In terms of schemes, if $f: V \to Spec(K) $
is a scheme over $Spec(K)$, then $V^\si$ is the same scheme, but with 
$f$ replaced
by $f \circ \sigma^{-1}$. ) $\phi_q$ denotes the Frobenius automorphism 
$x \mapsto x^q$.
 $\Phi_q \subset (X \times X^{\phi_q})$ denotes the graph of $\phi_q$, 
viewed as a subvariety.

\begin{thm}\lbl{1} Let $X$ be an affine variety over $k$, and let
 $S \subset (X \times X^{\phi_q})$
be an irreducible subvariety. Assume $\dim(S)=\dim(X)=d$,  
  the maps $S \to X$, $S \to X'$ are dominant, and one is quasi-finite.
Let $a = [S:X]/[S:X']_{insep}$.  
  Then there exists a constant $c$ depending only on the degrees of
  $X,S$ such that for $q>c$, 
$$|S(k) \meet \Phi_q(k)| = aq^d + e \hbox{ with } |e| \leq  c q^{d-\half}$$ \end{thm}

The numbers $[S:X]$,$[S:X']_{insep}$ refer to the degree (respectively
purely inseparable degree) of the field extensions $k(S)/k(X)$,
$k(S)/k(X')$.   In particular if $S$ is the graph of a separable morphism
$X \to X$, we have $a=1$.  The degrees of $X,S$ referred to can be
made explicit  as follows. If $X \subset 
\Aa^m$, let
$\bar S$ be the closure of $S \subset (\Aa^m \times \Aa^m) \subset \Pp^m \times \Pp^m$;
then $c$ depends only
on $m$ and the bidegrees of $\bar S$.
 In particular, if $q$ is large compared to these degrees,
the theorem implies that
$S(k) \meet \Phi_q(k) \neq \emptyset$.

  We will give   other versions
of the statement, better suited to the inherent uniformity of the 
situation.  Below (\ref{1}A) we state it  in the language of difference algebra.  A  algebro-geometric version  will be given 
later (\S \ref{geostat}.)

Observe that $X$ need not be defined over the fixed field of
$\phi_q$ (or indeed over any finite field.)  If $X$ is defined
over $GF(q)$, then $X=X^{\phi_q}$, and the diagonal 
(or the graph any other dominant morphism $X \to X$) becomes 
a possible choice of $S$.  In the case of the diagonal, one obtains the Lang-Weil
estimates.  

When $X$ and $S$ descend  to a finite field, 
  the projection $S \to X$ is proper, and $S \to X'$ is 
quasi-finite, Theorem \ref{1}  follows from Deligne's conjecture    (\cite{fujiwara},\cite{pink})  together with his theorem on 
eigenvalues of Frobenius.  
  In general these last two assumptions 
  (\cite{pink} 7.1.1) cannot be simultaneously obtained in our context, as far as I can see.   

Since $X$ is affine, an arbitrary proper divisor can be removed,
preserving the same estimate.  We thus find that
the set of points of $S \meet (X \times X^{\phi_q})$ is 
``asymptotically Zariski dense''.  Let us state a special case of this
(when $S$ is a morphism, it was used by Borisov and Sapir \cite{borisov-sapir} for a 
group theoretic application.)

\<{cor}\lbl{1c}  Let $X$ be an affine variety over $k=GF(q)$, and let
 $S \subset X^2$
be an irreducible subvariety over $K^a$. Assume    
 the two projections $S \to X$  are dominant.     Then for any
proper subvariety $W$ of $X$, for large enough $m$,
there exist $x \in X(k^a)$ with $(x,\phi_q^m(x) ) \in S$ and $x \notin W$.  \end{cor}

($k^a$ denotes the algebraic closure of $k$.)

\<{paragraph}{Difference algebraic statement}

Here is the same theorem in the language of difference algebra,
that will be used in most of this paper, followed by a completely
algebraic corollary; we do not know any purely algebraic proof of 
this corollary.  

  A {\em difference field} is a field $K$ with a distinguished 
endomorphism $\si$.  $K$ is {\em inversive} if $\si$ is
an automorphism.  A subring of $K$, closed under $\si$, 
is called a {\em difference domain}.    
A difference ring $R$ is {\em simple} if every difference ring homomorphism is
injective or zero on $R$. 

If $q$ is a power of a prime $p$, let $K_q$ be the difference field
consisting of an
algebraically closed field of characteristic $p$,
together with the $q$-Frobenius automorphism $\phi_q(x)=x^q$.
$K_q$  will be called a {\em Frobenius difference field}.  

\smallskip \noindent {\bf Theorem {\ref{1}A}} {\em
Let $D \subset R$ be finitely generated  
difference domains.   Assume there exists an embedding of $D$ into 
an algebraically closed, inversive difference field $L$, such that $R 
\tensor_D L$
is a difference domain, of transcendence degree $d$ over $L$.  Then 
there exist
$b,b' \in \Nn$, $0<c \in \Qq$, and $d \in D$ such that  for
 any     prime power $q \geq b$, and any homomorphism of difference
rings $h: D \to K_q$ with $h(d) \neq 0$, 
$h$ extends to a homomorphism $\bar h: R \to K_q$.
Moreover,
the number of different $\bar h$ is $cq^d +e$, where  $e \leq b' q^{d-1/2}$.
                }

\begin{cor} Let $D$ be a finitely generated  difference ring.  Assume
$D$ is simple, and has no zero divisors.  Then $D$ is a finite field, with a power of Frobenius.  \end{cor} 

The assumption that $D$ is a domain can be weakened to the
condition:    $ab =0$ implies $a\si(b)=0$.  A difference ring
satisfying this condition will be called  {\em well-mixed}.


\>{paragraph}

\begin{paragraph}{Model theoretic consequences}
With this in hand, the methods of Ax easily generalize to yield the 
first order
theory of the Frobenius difference fields.
Let $T_\infty$ be the set of all first-order sentences $\theta$ such 
that
for all sufficiently large $q$, $K_q \models \theta$.

\begin{thm}\lbl{dec}    $T_\infty$ is decidable. It coincides with
the model companion ``ACFA" of difference fields. \end{thm}

Here is a more precise presentation of the theorem (restricted
to primes), that does not explicitly mention the axioms of ACFA.    Fix
a countable algebraically closed field $L$ of infinite transcendence
degree over $\Qq$.  Let $G = Aut(L)$; for $\si \in Aut(\Qq^a)$,
let $G_\si = \{\tau \in G: \tau | \Qq^a = \si \}$.  $G_\si$ is a Polish
space, and it makes sense to talk of meager and co-meager sets;
we will mean ``almost every'' in this sense.

{\em Let $\phi$ be a first order sentence in the language of difference rings.  The following numbers are equal:

1)  The Dirichlet density of the set of primes $p$ such that 
$K_p \models \phi$.

2)  The Haar measure of the set of $\si \in Aut(\Qq^a)$, such that, for almost every $\tau \in G_\si$,   $(L,\tau) \models \phi$.
}

The theory ACFA is described and studied in \cite{CH}.
The axioms state that the field is algebraically closed, and:

{\em Let $V$ be an absolutely irreducible variety over $K$,
and let $S$ be an irreducible subvariety of $V \times V^\sigma$ 
projecting dominantly onto $V$ and
onto $V^\sigma$. Then there exists $a \in V(K)$ with $(a,\si(a)) \in 
S$.}

With different axioms, the existence of a model companion for
difference fields was discovered earlier by Angus Macintyre and Lou
Van den Dries, precisely in connection with the search for the theory
of the Frobenius.   See \cite{macintyre}.
\
In positive characteristic, we can also conclude:

\begin{cor}\lbl{1cor}  Let $F = GF(p)^{\alg}$.  
For almost every automorphism $\si$ of $F$ (in the sense of
Baire category,)
$(F,\si) \models ACFA$.  \end{cor}

I.e. the set of
exceptional automorphisms is meager.  

A similar result holds, much more trivially, if $GF(p)^{\alg}$
is replaced by an algebraically closed field of infinitely countable
transcendence degree over the prime field.  However, an example of
Cherlin and Jarden shows that a generic automorphism of $\bar{\Qq}$
does not yield a model of ACFA, nor is this true
for any other field of finite transcendence degree.
\

We can also consider the theory $T_{\rm Frobenius}$ consisting of those
sentences true in {\em all} $K_q$, or in all $K_p$ with $p$ prime.
    The completions of this theory are just
the completions of $ACFA$ together with the complete theories of
the individual
$K_q$. With some extra care, we can also conclude:

\begin{thm}\lbl{4} The theory of all $K_p$, $p$ prime, is decidable.   
\end{thm}

A similar result holds for the theory of all $K_q$, $q$ an arbitrary  
prime power; or
with $q$ ranging over all powers of a given prime; etc.

\end{paragraph}

\>{subsection}

\<{subsection}{Some applications}
\begin{paragraph}{Finite simple groups}

The theory of definable
groups has been worked out for $ACFA$, and has some suggestive 
corollaries
for finite simple groups. This work has not yet been published;  but an 
important
special case (groups over the fixed field) contains the main ideas, see
 \cite{HP}.  Here we give a preview.

We say that a family of
finite groups indexed by prime powers $q$ is {\em uniformly definable} 
if there exist
 first order formulas
$\phi,\psi$ such that $\phi$ defines a finite subset of each $K_q$, 
$\psi$ defines
a group operation on $\phi$, and the family consists of these groups for 
the various $K_q$.
Examples of such families include $G_n(q)$, where $n$ is fixed, and 
$G_n(q)$ is
one of the families of finite simple groups (e.g. $PSL_n(q)$). All but 
two of
these families are already definable over finite fields.   However, the 
Ree
and Suzuki families are not. Instead they are defined by the formula:
$$\si(x) = \phi(x), x \in G$$
where $G$ is a certain algebraic group, and $\phi$ is an algebraic 
automorphism
whose square is the Frobenius $\phi_2$ or $\phi_3$.

From Theorem \ref{dec}, we obtain immediately:

\begin{thm} For each fixed $n$, the first order theory of each of the 
classes
$G_n(q)$ of finite simple groups is decidable. \end{thm}

The remaining results use the theory of groups of finite S1-rank.

\begin{thm} Every uniformly definable family of finite simple groups is (up to a 
finite set)
contained in a finite union of families $G_n(q)$ \end{thm}

More precisely, for each $q$, there is a definable isomorphism in $K_q$ 
between
the group $\phi(K_q)$, and some $G_n(q)$; only finitely many values of 
$n$ occur,
and finitely many definable isomorphisms. The theorem is proved without 
using the
classification of the finite simple groups. Note that it puts the Ree
and Suzuki groups into a natural general context.  (See also 
\cite{bombieri}
in this connection.)

\begin{thm}\label{finite-simple} Let $G_n(q)$ be one of the families of 
finite simple groups. Then
there exists an integer $r=r(n)$ such that for any
$q$, every nontrivial conjugacy class of $G_n(q)$
generates $G_n(q)$ in at most $r$ steps.
\end{thm}

This is a special case of more general results on generation of 
subgroups
by definable subsets.  One can indeed take $r(n) = 2 \dim H_n$, where 
$H_n$
is the associated algebraic group.

\end{paragraph}

\<{paragraph}{Jacobi's bound for difference equations}  
A linear difference or differential equation of order $h$, it is well-known, has
a solution set whose dimension is at most $h$.  The same is true for a nonlinear
equation, once the appropriate definitions of dimension are made.  
 Jacobi, in \cite{jacobi}, proposed
a generalization to systems of $n$ differential equations in $n$ variables.  The
statement is still conjectural today; an analogous conjecture for difference equations
was formulated by Cohn.  We prove this in \S \ref{jacobiS}.  
Our method of proof illustrates theorem \ref{1}:  we show that our analog
of the Lefschetz principle translates Cohn's conjecture to a very strong, but
known, form of Bezout's theorem.  
\>{paragraph}

See \S \ref{proofs} for proofs and some other applications.

\>{subsection}
\<{subsection}{Difference algebraic geometry}

Take a system $X$ of difference equations:  to fix ideas,
take the equations of the unitary group, $\si(u) u^t = 1$, where
$u$ is an $n \times n$ -matrix variable; or in one variable, 
take $\si(x) = x^2+c$, for some constant $c$.  Consider four
 approaches to the study of $X$.

It may be viewed as part of the category of objects defined by
difference equations; i.e. part of an autonomous geometry of 
difference equations.  See the discussion and examples in  \cite{H-MM}.  Here $X$ is identified, in the first instance,
with the set of solutions of the defining equations in a difference
field.  The
rough geography of the  category of definable sets and maps is 
worked out in \cite{CH}, \cite{CHP}.    A great deal
remains to be done, for example with respect to topology.  The 
definition of difference schemes in Part I can be viewed as a small step 
in a related direction.

Alternatively, it is possible to view a difference equation in terms
of the algebro-geometry data defining it.    One can often represent the
system as a correspondence:  a pair $(V,S)$, where $V$ is an algebraic variety
over a difference field $k$, $S$ a subvariety of $V \times V^\si$; the
equation is meant to be $X = \{x: (x,\si(x)) \in S\}$.  See \ref{direct-ritt-0}.  
The great strength of this approach is that results of algebraic 
geometry become directly accessible, especially when
$V$ is be smooth and complete, and the fibers of $S \to V, S \to V^\si$
of equal dimensions.  It is generally not  possible to meet these
desiderata.   Moreover, given a system $(V,S)$
in this form, simple difference-theoretic questions (such as irreducibility)
are not geometrically evident; and simple operations require
radical changes in $V$.  Nevertheless, we are often obliged
to take this approach, in order to use techniques  
that are not available directly for difference varieties or schemes;
especially the cohomological results used in \S \ref{virtual}.  We will
in fact stay with this formalism for as long as possible, particularly in \S \ref{geopre}, even when a shorter  difference-theoretic treatment is visible; cf. the proof of the moving lemmas in \S \ref{moving}.   
One obtains more explicit statements, and a contrast with those
points that seem to really require a difference-theoretic treatment.

The third and most classical approach is that of dynamics.  See
Gromov's \cite{SAG} for a combination of dynamics
with algebraic geometry in a wider context, replacing our single
endomorphism with a finitely generated group; it also  contains, in  a different way than ours, reductions to finite objects.   Here
the space corresponding to $(V,S)$ is 
$$Y = \{(a_0,a_1,\ldots): (a_n,a_{n+1}) \in S \}$$  This is a pro-algebraic
variety; we can topologize it by the coarsest topology making
the projections to algebraic varieties continuous, if the latter
are given the discrete topology.   Note that the points of $X$ and
$Y$ are quite different; in particular it makes sense to talk of points
of $Y$ over a field (rather than a difference field.)

  For the time being
the work described here barely touches upon dynamical issues;  though 
just under the surface, the proofs in \cite{CHP} and in \S \ref{virtual}, the 
definition of transformal degree, and other points do strongly involve iteration.   It is also true that by showing density of Frobenius
difference fields among all difference fields, we  show density
of periodic points (at least if one is willing to change characteristics); raising hopes for future contact.

A fourth point of view, introduced in \S \ref{frobred}, is to view the symbol
$\si$ as standing for a variable Frobenius.  For each prime power $q=p^m$ we obtain
an algebraic variety (or scheme) $X_q$, by interpreting $\si$ as
the Frobenius $x \mapsto x^q$ over a field of characteristic $p$. 
The main theorem 
can be read as saying that this approach is equivalent to the others.

It can be quite curious to see information, hidden in $(V,S)$, 
released  by the two apparently disparate processes of iteration,
and intersection with Frobenius.   
Even simple dimension-theoretic facts can sometimes be seen
more readily  in this way.  Cf. \ref{illustrate-dirdim}, or
the proof of the Jacobi bound for difference equations in \S \ref{jacobiS}
(where I know no other way.)

Here is a the simplest example of a definition and a number from the various 
points of view.  The following are essentially equivalent:
\<{itemize}
\item $X$ has transformal dimension $0$, or equivalently finite 
SU-rank in the sense of \cite{CH}, or finite total order (\S \ref{dimensions});
\item The map $S \to V$ is quasi-finite;
\item The space $Y$ is locally compact;
\item Almost all $X_q$ are finite.
\>{itemize}
If this is the case, and $X$ is irreducible, the following numbers are essentially equal:
\<{itemize}
 \item 
the degree of the correspondence $S$ (adjusted, in positive characteristic, by
the purely inseparable degrees) 
\item the entropy of $Y$;
\item  the  number $c$ such
that $X_q$ has asymptotic size $cd^q$, where $d = \dim(V)$.
\>{itemize}

The associated invariant of $X$ needs to be defined in terms
of difference schemes, rather than varieties, and will be described in \S
\ref{redseq}.

''Essentially" indicates that these equivalences are not  
tautologies, and require certain nondegeneracy conditions; for instance it suffices
for 
$S \to V$ to be generically quasi-finite, if the exceptional infinite fibers
do not meet $X$.  The various technicalities will be taken up in \S 2- 6.

\>{subsection}  
\<{subsection} {Description of the paper}

\S 2 - 6 introduce  a theory of difference schemes.  As in algebraic
geometry, they are obtained by gluing together basic schemes of
the form $\spe A$, with $A$ a difference ring.  The $\spe A$ are
more general than difference varieties in that $A$ may have nonzero nilpotent elements $a$,
both in the usual sense ($a^2=0$), and in the transformal sense, such as $\si(a)=0$.   

\S 2 develops the basic properties of a well-mixed difference ring.   Any difference ring
without zero-divisors,  as well as any difference ring whose structure endomorphism
is Frobenius, is well-mixed.  Well-mixed rings admit a reasonable
theory of localization; in particular,  Proposition \ref{wm-localization}
shows that a difference scheme based on well-mixed rings makes sense intrinsically, without including the generating difference rings explicitly in the data.   

In \S \ref{dimensions}, the dimension and degree of a difference scheme is defined.  In fact
there are two dimensions: transformal dimension, and when that vanishes,
total dimension.  The transformal affine line has transformal dimension
one, and contains many difference subschemes of transformal
dimension zero and finite total dimension, notably the ``short'' affine line
over the fixed field, of total dimension one.
Here too, for  schemes based on well-mixed rings we obtain a 
smoother theory. We improve
a theorem of Ritt and Raudenbush by proving   Noetherianity
of finitely generated difference rings with respect to radical, well-mixed
ideals, in finite total dimension.  We do not know whether the same theorem holds for
all well-mixed ideals.   (Ritt-Raudenbush proved it for
perfect ideals; cf. \cite{cohn}.)

\S \ref{trans-mult} is concerned with the kernel of the endomorphism
$\si$ on the local rings, or with transformal nilpotents; it is
necessary to show how their effect on total dimension is distributed
across a difference scheme.  \S \ref{dirpres} deals with the presentation
of a difference scheme as a correspondence $(V,S)$, and
how in simple case one may test for e.g. irreduciblity. 

Further basic material, not used in the main line of the paper,
is relegated to \S \ref{complements}.  In \S \ref{transdimdeg},
the transformal dimension of the difference scheme corresponding
to $(V,S)$ is determined in general; I do not know a proof
that does not go through the main theorem of the paper.    \S \ref{projective} looks at  transformal projective space.  \S \ref{blowup}
introduces two transformal versions of the blowing-up construction.  They  are shown to essentially coincide; as one (based on the short line,
over the fixed field) has no Frobenius analog, this may be geometrically
interesting.   

The idea of a moving family of varieties, and of a limit of such a family,
is characteristic of geometry.  It is poorly understood model-theoretically, and has not so far been studied in difference algebraic geometry.  In 
algebraic geometry, it is possible to give a completely synthetic treatment:
given a family $V_t$ of varieties over the affine $t$-line (except over
$0$), take the closure and let the limit be the fiber $V_0$.  This works
in part since an irreducible family over the line is automatically flat.  In 
difference algebraic geometry, this still works over the short line, but not over the 
transformal affine line.  It is possible to use blowing ups to remedy the situation, but the use of valuation theory appears to be much clearer.

Transformal valued fields are   studied  in   \S \ref{td1}.   They
differ from the valued difference fields of \cite{scanlon}, primarily in 
that the endomorphism does not fix the value group:  on the
contrary, it acts as a rapidly increasing map.  We are not concerned
with the first order theory here, but
rather with the structure of certain analogues of discrete valuation
rings.  Their value group is not $\Zz$ but $\Zz[\si]$, the
polynomial ring over $\Zz$ in a variable $\si$, with $\si > \Zz$.
A satisfactory structure theory can be obtained, especially
for the completions.  Extensions
of such fields  of finite transcendence degree can be described
in terms of inertial and totally ramified extensions, and the completion 
process; there are no other immediate extensions.  As a result,
the sum of the inertial transcendence degree extension and the
ramification dimension is a good invariant; cf. \ref{vrk-dec}.  This
valuative dimension will be used in \S \ref{transvalan}.

Given a family $X_t$ of difference varieties, we define in \S \ref{transspec}
 the specialization
$X_{\to 0}$, a difference subvariety of the ``naive'' special fiber $X_0$.
The total dimension of $X_{\to 0}$ is at most that of the generic
fiber $X_t$; this need not be the case for $X_0$.
   Together with
the material in \S \ref{transvalan}, this lays the basis
for a theory of rational and algebraic equivalence of transformal
cycles.  (It will be made explicit in a future work.)

  The main discovery allowing the description of specializations
in terms of transformal valuation rings is this:  definable sets
 of finite total dimension are analyzable in terms of the residue field.
The simplest case is of  a definable set $X_t$ in $K$, contained
entirely in the valuation ring, and such that the residue map is injective
on $X_t$.   The equation $\si(x) = x$ for   is immediately 
seen to have this property.    In general, for an equation of transformal
dimension zero, $X_t$ is {\em analyzable} over the residue
field  (\S \ref{fintdra}.)  This model-theoretic notion is explained in
\S \ref{anresdim}.   It means that sets of finite total dimension
over the generic fiber can actually be viewed as belonging
to the residue field; but in a somewhat sophisticated sense, that
will be treated separately.  Here we will content ourselves with
the numerical consequence for Frobenius fields, bounding the number of points  specializing to a given point of $X_0$ in terms of the valuative
dimension. With this relation in mind, 
 \S \ref{dirgen} contains a somewhat technical result, needed in
the main estimates,
bounding the valuative dimension in terms of directly visible data,
when a difference scheme is presented geometrically as $(V,S)$.

The proof of the main theorem begins in \S \ref{geopre}.  In this
section, it is reduced to an intersection problem on a smooth,
complete algebraic variety $V$:  one has correspondence
$S \leq V \times V^\si$, and wants to estimate the number
of points of $S \meet \Phi_q$, where $\Phi_q$ is a graph
of Frobenius; or rather those points lying outside a proper
subvariety $W$ of $V$.  The main difficulty is that, away
from $W$, the projection $S \to V$ need not be quasi-finite.
For instance, $S$ may contain a ``square'' $C \times C$,
and then any Frobenius will meet $S$ in an infinite set.   
A moving lemma in \S \ref{moving} shows that, while this can perhaps not be avoided,
when viewed as a cycle $S$ can be moved to another, $S'$,
whose components do not raise a similar difficulty.

In \S \ref{virtual} we obtain an estimate for the intersection number 
$S \cdot \Phi_q$ in the
sense of intersection theory.  See \ref{cycle-t} for a quick proof in the case of projective space, and \ref{weil} for a proof for curves:   indeed
Weil's proof, using positivity  in the intersection product
on a surface, works in our case too.  

  For general varieties, in the absence of the standard
conjectures, we use the cohomological representation and
Deligne's theorem.  
According to the Lefschetz fixed point formula, the question amounts to a question on
the eigenvalues of the composed correspondence $  \Phi^{-1} \circ S$,
acting on the cohomology $H^*(V)$ of $V$.  When $V= V^\si$ and $S,\Phi$
commute, one can apply Deligne's theorem on the eigenvalues
of $\Phi$.  But in general, $\Phi$ induces a map between two
distinct spaces
$H^*(V)$ and $H^*(V^\si)$; so one cannot speak of eigenvalues.  
 Some trickery is therefore  needed ,
in order to apply the results of Deligne.  It is in these manipulations   that
the precision of Deligne's theorem is needed,  rather than just the
Lang-Weil estimates.  It may suffice to know
that every eigenvalue on $H^i$, $i< 2n$, has absolute value $\leq q^{2n-1/2}$;
but because of
possible cancellations in the trace (when e.g.  all $m$'th roots of unity are eigenvalues for some $m$),
this still seems beyond Lang-Weil.

  It would be interesting to sharpen the statement
of the theorem, so as to give a precise rather than asymptotic cohomological account of the size of Frobenius specializations of zero-dimensional difference schemes.    
 

Having found the intersection number, we still do not know the 
number of (isolated) points of the intersection $S \meet \Phi_q (\setminus W)$; this is due to
 the possible existence of infinite components, mentioned above.
 Cf. \cite{fulton}, and Kleiman's essay in \cite{seidenberg}. \S \ref{equivalence} is 
devoted
to estimating the ``equivalence" of these components.  A  purely
geometric proof using the theory of \cite{fulton}
appears to be difficult;  the main trouble
lies in telling apart the 0-dimensional distinguished subvarieties of
the intersection, from the true isolated points (not embedded in larger
components) that we actually wish to count.  Instead, we use
the methods of specialization of difference schemes, developed
above.  By the moving lemma, there exists a better-behaved cycle
$S_t$ on $V \times V^\si$, specializing to $S$.  The intersection 
numbers corresponding to $S_t,S$ are the same; and we can 
assume that the equivalence   has been bounded for $S_t$.  The 
problem is to bound the number of points of $S_t \meet \Phi_q$
specializing into $W$; as well as the number (with multiplicities) of points of $S_0 \meet \Phi_q$ to which no point of $S_t \meet \Phi_q$ specializes.
The inverse image of $W$ under the residue map is shown
to be analyzable over the residue field, of dimension smaller than
$\dim(V)$, proving this point. 

This use of intersection theory does not flow smoothly,
since the difference varieties do not really wish to remain restricted
to $V \times V^\si$.  They are often more naturally represented
by subvarieties of,  say, $V \times V^\si \times V^{\si^2}$.
Even the decomposition into irreducible components of distinct
transformal dimensions cannot be carried out in $V \times V^\si$; e.g.
 above, the closure of $\{x \in V \setminus W: (x,\si(x)) \in S\}$ is always a difference scheme of finite total dimension; but $S$ is irreducible,
and the closure cannot be represented on $V \times V^\si$.   On the other hand,
the intersection-theoretic and cohomological methods that we 
use do not seem to make sense when one intersects two
$\dim(V)$-dimensional subvarieties of $V \times V^\si \times V^{\si^2}$;
and so we must strain to push the data into $V \times V^\si$ (as in
\ref{dirgen}.)   This is one of a number of places  
 leading one to dream of a broader formalism.

The proofs of the various forms of the main theorem, the applications mentioned
above, and a few others, are gathered in \S \ref{proofs}
and \S \ref{jacobiS}.
\end{subsection}

 \begin{paragraph}{Thanks}   to Dan Abramovich,  Phyllis Cassidy and Ron Livne for
  early discussions of this problem,  and to Gabriel Carlyle,  Zo\'e Chatzidakis, Martin Hils, Emmanuel  Kowalski,   Yves Laszlo, and Thomas Scanlon for their thorough and useful comments. 
Warm thanks to Richard Cohn for information about the Jacobi bound problem.  The paper profited greatly from a seminar
organized by J.-B. Bost, Y. Laszlo and F. Loeser, and a CIRM workshop in Luminy in  March 2013; I am very grateful 
 to all the contributors.

 \end{paragraph}

\end{section}

\begin{part}{Difference schemes }
\begin{section} {Well-mixed rings }

The basic reference for difference algebra is \cite{cohn}.

A {\em difference ring} is a commutative ring $R$ with $1$ and with
a distinguished ring homomorphism $\si: R \to R$.   We write $a^\si = \si(a)$,
and $a^{\sum m_i \si^i} = \Pi _i \si^i(a)^{m_i}$.

A {\em well-mixed ring} is a difference ring satisfying:

$$  a b = 0 \Rightarrow ab^\si = 0 $$

\sk  
Given a difference ring $R$, a {\em difference ideal }  is an ideal such that
$x \in I$ implies $\si(x) \in I$.  If in addition $R/I$ is well-mixed
we say that $I$ is well-mixed.  A {\em transformally prime ideal} is a
prime ideal, such that, moreover, $x \in I \Leftrightarrow \si(x) \in
I$.  Every transformally prime ideal is well-mixed.  

Any difference ring $R$ has a smallest difference ideal $I=rad_{wm}(R)$
 such that $R/I$ is well-mixed.

A difference ring is a
{\em difference domain} if the zero-ideal is a transformally prime ideal. 
\footnote{Ritt and Cohn (\cite{cohn}) use the term {\em prime difference ideal} 
for what we call a transformally prime ideal.  When as we will one considers
difference ideals that are prime, but not transformally prime, this terminology
becomes awkward.   To avoid confusion, we will speak of
{\em algebraically prime} difference  ideal to mean difference ideals that are
prime.  We kept the use of {\em difference domain}, as being less confusing;
though {\em   transformally integral} might be better.  At all events,
difference rings with the weaker property of having no zero-divisors will be referred to as
as {\em algebraically integral}.  Cohn calls an ideal   {\em mixed} if the quotient
is well-mixed.}   %

A difference field is a difference domain, that is a field.
(An endomorphism of a field is automatically injective.)

  \<{lem} \lbl{alg-h}  Let $R$ be a difference ring, $P$ an algebraically prime 
difference ideal, $K$ the fraction field of $R/P$, $h: R \to K$ the natural
map.
\begin{enumerate}
  \item     $P$ is
transformally prime iff there exists a difference field structure on
$K$, extending the natural quotient structure on $R/P$.
  \item Let $R' \leq R$ be a difference subring, $P' = (P \meet R')$.
  Let $K'$ be
the fraction field of $R'/P'$, and suppose $K$ is algebraic over $K'$.
Then $P$ is transformally prime iff $P'$ is.
\end{enumerate}
\>{lem}

\proof (1)  If $P$ is transformally prime, $R/P$ is a difference domain,
and $\si$ extends naturally to the field of fractions: if $b \neq 0$ then $\si(b) \neq 0$, and one can define $\si(a/b)=\si(a)/\si(b)$.  Conversely,
if $L$ is a difference field, and $R \to_h R/P \subset L$
is a homomorphism, if $a \notin P$ then $h(a) \neq 0$ 
so $ h(\si(a)) = \si(h(a))  \neq 0$, and thus $\si(a) \notin P$.

(2)  If $P$ is transformally prime, clearly $P'$ is too.  If $P'$
is transformally prime, then the multiplicative set $h(R' \setminus P') \subseteq K' \setminus (0) \subset K \setminus (0)$
is closed under $\si$, so $\si$ extends to the localization of $h(R)$ by
this set, i.e. to $M = K'(R/P)$.  But since $K' \subseteq M \subseteq K$ and $K$
is algebraic over $K'$, $M$ is a field, hence $M = K$. \qed

\begin{notation} Let $I$ be a difference ideal.  
 $$\srad{I} = \{a \in R:  a^{\sum_{0 \leq i \leq n} m_i \si^i}  \in I, \hbox { some } m_0,m_1,\ldots ,m_n\geq 0 \}$$ 
$$ \sqrt{I} = \{ a \in R: a^n \in I, \hbox { some } n \in \Nn \}$$ \end{notation}

An ideal of $R$ is {\em perfect}  if
$x \si(x) \in I \Rightarrow x, \si(x) \in I$. 

\begin{lem}\lbl{well-mixed1} Let $R$ be a well-mixed ring.  Then $\sqrt{0}$ is 
a radical, well-mixed ideal.   $\srad{0}$ coincides with 
\[ I = \{a \in R: a^{m \si^i} = 0 \hbox{ for some } m \geq 1, i \geq 0\} = \{a \in R: a^{\si^i} \in  \sqrt{0}  \hbox{ for some }  i \geq 0\} \]
 and is a   perfect ideal.

\end{lem}

\proof The statement regarding $\sqrt{0}$ is left to the reader.  If $a^{m \si^i} =
0$, and $m \leq m'$, $i \leq i'$, then $a^{m' \si^{i'}} = 0$.  It
follows that $I$ is closed under $\si$ and under multiplication, and
addition (if $a^{m \si^i} = 0$ and $b^{m \si^i} = 0$ then $(a+b)^{2m
\si^i} = 0$.)  Moreover if $a a^\si \in I$, then $a^{(\si+1)m \si^i} =
0 $.  So $b b^\si = 0$ where $b =a^{m \si^i}$.  As $R$ is well-mixed,
$b^{\si} b^\si = 0$.  Thus $a^{2m \si^{i+1} } = 0$, so $a \in I$. 
\qed

\begin{lem}\lbl{well-mixed1.4} Let $R$ be a difference ring.  Any
intersection of well-mixed ideals is well-mixed.  \end{lem}

\proof:  Clear.

We write $Ann(a/I)$ for $\{b \in R: ab \in I \}$, 
$Ann(a) = \{b \in R: ab=0\}$.

\begin{lem}\lbl{well-mixed1.5}  Let $R$ be a well-mixed ring, $a \in R$, 
  $I(a) $ the smallest well-mixed ideal
containing $a$.  Then $Ann(a)$ is a well-mixed ideal.  If $b \in I(a)
\meet Ann(a)$, then $b^2 =0 $.
\end{lem}

\proof Let $b \in Ann(a)$.  Then $\si(b) \in Ann(a) $ since $R$ is
well-mixed.  Thus $Ann(a)$ is a difference ideal.  Suppose $cb \in
Ann(a)$.  Then $abc = 0$.  So $abc^\si = 0$, as $R$ is well-mixed.  Thus
$bc^\si \in Ann(a)$.  So $Ann(a)$ is well-mixed.

Let $c \in Ann(a)$.  Then $a \in Ann(c)$.  Since $Ann(c)$ is an
well-mixed ideal, $I(a) \subset Ann(c)$.  Thus if at the same time $c
\in I(a)$, we have $c \in Ann(c)$, so $c^2=0$.

\begin{lem}\lbl{well-mixed1.6} Let $R$ be a  well-mixed domain,
$a \in R$,  $J(a) $ the smallest algebraically radical,
well-mixed ideal containing $a$.  Then $Ann(a)$ is an algebraically
radical, well-mixed ideal; and $J(a) \meet Ann(a) =(0)$.
 
\end{lem}

\proof In \ref{well-mixed1.5} it was shown that $Ann(a)$ is a
well-mixed ideal.  If $b^m \in Ann(a)$, then $b^ma = 0$, so
$(ba)^m=0$, so $ba =0$.  Thus $Ann(a)$ is algebraically radical.  Let
$c \in Ann(a)$.  Then $a \in Ann(c)$.  Since $Ann(c)$ is a  well-mixed,
algebraically radical ideal, $J(a) \subset Ann(c)$.  Thus if at the
same time $c \in I(a)$, we have $c \in Ann(c)$, so $c^2=0$, hence $c=0$.

Let $\spe R$ denote the set of prime ideals $P$ of $R$ with $\si^{-1}(P)=P$.  Later, we will
give this set more structure.

\begin{lem}\lbl{well-mixed1.7}  Let $R$ be a well-mixed ring, $0 \neq a \in R$.  For $p \in \spe R$,
let  $a_p$
denote the image of $a$ in the local ring $R_p$. Then $ \{p \in \spe R: a_p \neq 0 \}$  is a nonempty
closed subset of $\spe R$.   
\end{lem}

\proof  By \ref{well-mixed1.6}, $Ann(a)$ is a well-mixed ideal.   
Let $p$ be a maximal well-mixed
ideal containing $Ann(a)$.  Then $p$ is a prime   ideal.  For if $cd \in p$, then 
$p \subset Ann(c/p)$, so $p=Ann(c/p)$ (and then $d \in p$) or $R=Ann(c/p)$ (and then $c \in p$).
Moreover as $p$ is well-mixed, it is a difference ideal, so $p \subset \si^{-1}(p)$; again by
maximality of $p$, $p= \si^{-1}(p)$.  Thus $p$ is a transformally prime ideal.  If $d \notin p$,
then $ad \neq 0$ since $d \notin Ann(a)$.  So $a_p \neq 0$.  \qed

\begin{lem}\lbl{well-mixed1.8} Let $R$ be a well-mixed difference ring,
$S$ a subring.

\begin{enumerate}
  \item Let  $b \in R$, $q=Ann(b) \meet S$.  If $a \in q, \si(a) \in S$ then $\si(a) \in q$.  Same for 
$Ann(b/\sqrt{0}) \meet S$.
  \item If $S$ is Noetherian and $q$ is a minimal prime of $S$, the same conclusion holds.
  \item Call an ideal  $p$  of  $R$ {\em cofinally minimal } if for any finite $F \subset R$
 there exists a Noetherian subring $S$ of $R$ with $F \subset S$ and 
with $p \meet S$   a minimal prime of $S$.  Then any cofinally minimal ideal is a difference
ideal. 
\end{enumerate}
  
 \end{lem}

\proof  \begin{enumerate}
  \item By \ref{well-mixed1}, \ref{well-mixed1.5}, $Ann(b)$ and $Ann(b/\sqrt{0})$ are difference ideals of $R$.  The
statement about the intersection with $S$ is therefore obvious.
\item  Let $S$ be a Noetherian subring of $R$.  Then
$\sqrt{0}_S = q_1 \meet \ldots \meet q_l$ for some minimal primes 
$q_1 \meet \ldots \meet q_l$; and $q=q_i$ for some $i$, say $q=q_1$.  Let $b_i \in q_i$,
$b_i \notin q$ for $i >1$, and let $b = b_2 \cdot \ldots \cdot b_l$.   
Then $q =  Ann(b/\sqrt{0}) \meet S$. 

\item Let $p$ be cofinally minimal, $a \in p$.  Let $F = \{a,\si(a)\}$, and let $S$ be a Noetherian subring of $R$
containing $F$  with $p \meet S$   a minimal prime of $S$.  By (2), $\si(a) \in (p \meet S)$, so $\si(a) \in p$.
\end{enumerate}
 \qed

\smallskip  
\<{remark}\lbl{well-mixed1.8R} Let $k$ be a Noetherian commutative ring, $R$ be a countably
generated $k$-algebra.  Then cofinally minimal ideals exist.  Indeed
if $S$ is a finitely generated $k$-subalgebra of $R$, $p_S$ any
minimal prime ideal of $S$, then $p_S$ extends to a cofinally minimal
$p$ with $p_S = p \meet S$.  \>{remark}

\proof  Find a sequence $S=S_1 \subset S_2 \subset \ldots$ of 
finitely generated $k$-subalgebras of $R$, with $R = \union_n S_n$.  Find inductively
minimal prime ideals $p_n$ of $S_n$ with $p_{n+1} \meet S_n = p_n$, $p_1=p_S$;
then let $p = \union _n p_n$.  Given $p_n$, we must find
a minimal prime $p_{n+1}$ of $S_{n+1}$ with  
$p_{n+1} \meet S_n = p_n$.  Let $r_1,\ldots,r_k$ be the minimal primes of $S_{n+1}$;
then $\meet_i r_i$ is a nil ideal; so $\meet_i (r_i \meet S_n)$ is a nil ideal;
thus $\meet_i (r_i \meet S_n) \subset p_n$.  As $p_n$ is prime,   for some $i$,
$r_i \meet S_n \subset p_n$; as $p_n$ is minimal, $r_i \meet S_n = p_n$. \qed

\<{lem}\lbl{wm1.11} Let $R$ be a difference ring, $I=\sqrt I$ a
well-mixed ideal.  Then $I$ is the intersection of algebraically prime difference ideals. \>{lem}

\proof   We may assume $I=0$.   If $a \neq 0$, we must find
an algebraically prime difference ideal $q$ with $a \notin q$.  By
compactness (see below), we may assume here that $R$ is finitely generated.  Find
 using
\ref{well-mixed1.8R} a cofinally minimal prime $q$ with $a \notin q$.   By
\ref{well-mixed1.8} it is an algebraically prime difference ideal. 
 (Here is the model theoretic compactness argument.   Assume for each finitely generated difference subring $S$ of $R$
we are given an algebraically prime difference ideal $q_S$ with $a \notin q_s$.  One can then  choose an ultrafilter on the set of such subrings,
concentrating for each finite $F$ on $\{S: F \subseteq S\}$.   Then  $R$ embeds into the ultraproduct of the $S$.  But they come
with an algebraically prime difference ideal $q$ = the ultraproduct of the $q_S$, and the pullback is as required.)
 \qed

\<{lem}\lbl{wm1.12} Let $R$ be a well-mixed ring, $a \in R$.  Then $a \in p$
for all $p \in \spe R$ iff $a^n =0$, some $n \in \Nn[\si]$. \>{lem}

\proof  One direction is immediate.  For the other, assume $a^n \neq 0$
for $n \in \Nn[\si]$.  Let $I$ be a maximal well-mixed ideal with $a^n \notin I$,
 $n \in \Nn[\si]$.   By maximality, and using Lemma \ref{well-mixed1},  $\sqrt[\si] I = I$.  So $I = \meet p_j$, $p_j$ transformally prime  
ideals.  Since $a \notin I$, $a \notin p_j$ for some $j$.    \qed

\begin{paragraph}{Difference polynomial rings}

Let $R$ be a difference ring.  A {\em difference monomial} over $R$
is an expression $rX^\nu$, where $\nu = \sum_{i=0}^m m_i \si^i$,
$m_i \in \Nn$.  The {\em order} of the monomial is
the highest $i$ with $m_i \neq 0$.  A {\em difference polynomial } in one variable $X$ (of order $\leq M$)  is 
a formal sum of difference monomials over $R$ (of order $\leq M$.)
The difference polynomials   form 
a difference $R$-algebra $R[X]_\si$.  If $k$ is a field, $k(X)_\si$ denotes the field 
of fractions of $k[X]_\si$.  

\>{paragraph} %

\<{paragraph}{Transformal derivatives} Let $F \in K[X]_\si$ be a difference polynomial in one
variable.  We may write $F(X) = \sum_{\nu} c_{\nu} X^{\nu}$, $\nu \in \Nn[\si]$; where
$\{\nu: c_{\nu} \neq 0 \}$, the {\em support} of $F$, is assumed finite.

Clearly $F(X+U) = \sum_{\nu} F_\nu(X) U^{\nu}$, where $F_\nu$ are certain (uniquely
defined) polynomials.

\<{definition} \lbl{trans-der} The $F_\nu$ will be called 
the {\em transformal derivatives} of $F$, and denoted $\partial_\nu x F = F_\nu$. \>{definition}

Clearly $\partial_\nu x F  = 0$ for all but finitely many $\nu$.  The $\partial_\nu x F $ satisfy rules analogous to those written down by Hasse.  In particular (as can be verified at the level of monomials.)

$\partial_0 (F) = F$

$\partial_\nu (F+G) = \partial_\nu x F + \partial_\nu x G$

$\partial_\nu (FG)  = \sum_{\mu + \mu' = \nu} F_{\mu}(X) F_{\mu'}(Y)$

$\partial _{\si \nu  } \si(F) =  \si( \partial_\nu (F))$

The definitions in  several variables
are analogous.
\>{paragraph}

\end{section}
\begin{section}{Definition of difference schemes}

\<{subsection}{Localization, rings of sections}  \lbl{localization}
  \sk
\<{paragraph} {Localization}
If $R$ is a difference ring, $X \subset R$, $\si(X) \subset X$, $XX \subseteq X$,
and $0 \notin X$, we consider the localization of $R$ by $X$, and write
\[ R[X^{-1}] = \{ \frac{a}{b}: a,b \in R, b \in X \} \]
It admits a natural difference ring structure, with $\si(\frac{a}{b}) = \frac{\si a}{\si b}$.
(This abuses  notation   slightly,  since $R$ need not  inject into $R[X^{-1}]$;   $\frac{a}{b}$
is taken to denote an element of $R[X^{-1}]$, the ratio of the images of $a,b$.)

An element $a \in R$ is said to be $\si$-nilpotent if $0$ lies in the smallest 
$X \subseteq R$ with $a \in X, \si(X) \subset X$, $XX \subseteq X$.  If $a$ is not
$\si$-nilpotent, we denote $R[a \inv]_\si = R[X \inv]$ for this $X$.

If $p$ is a transformally prime ideal, then the localization $R_p$ is
defined to be
\[ R[(R \setminus p)^{-1}] = \{ \frac{a}{b}: a,b \in R, b \not \in p \} \]
\>{paragraph}

 \sk
 The {\em difference spectrum},
$\spe (R)$, is defined to be the set of transformally prime ideals. It is
made into a topological space in the following way:   a
{\em closed} subset of $\spe(R)$
 is the set of elements of $\spe(R)$ extending a given ideal
$I$.

We note immediately that $\spe(R)$ is compact (though rarely Hausdorff.)

\sk
Recall that an ideal of $R$ is {\em perfect}  if
$x \si(x) \in I \Rightarrow x, \si(x) \in I$.  A perfect ideal is an intersection of 
transformally prime ideals.
 There is a bijective correspondence between closed sets in $\spe (R)$
and perfect ideals of $R$. (cf. \cite{cohn}.)

A perfect ideal is well-mixed: modulo a perfect ideal, $ab = 0  
\implies  (ab^\si)(ab^{\si})^{\si} = 0  \implies  ab^\si = 0$.  \sk

 We define
a sheaf of difference rings on $\spe (R)$, called the structure sheaf and denoted
${\cal O}_{\spe R}$, as follows.

If $U$ is an open subset of $\spe R$, a section of ${\cal O}_{\spe R}$, by
definition, is a function $f$ on $U$ such that $f(p) \in R_p$, and such
that for any $p \in U$, for some $b \notin p$ and $a \in R$, $f(q) = \frac{a}{b}$
for any $q \in U$, $b \notin q$.

We let ${\cal O}_{\spe R}(U)$ be the collection of all such functions; it is
a difference ring with the  pointwise operations.  One verifies immediately
that this gives a sheaf  ${\cal O}_{\spe R}$.

The space
$\spe (R)$ together with the sheaf ${\cal O}_{\spe R}$
is called the {\em affine difference scheme} determined by $R$.

If $h: R' \to R$ is a difference ring homomorphism, $X = \spe R, X' = \spe R'$, we have a natural
map $h^*: X \to X'$ of difference-ringed spaces.  At the level of topological
spaces, for $p \in \spe R$, let $H(p) =   h^{-1}(p)$.  This permits us to define  
 a sheaf  on $X'$,  $U \mapsto {\cal O}_{X}(H^{-1}(U))$, and a morphism
 of sheaves, $ {\cal O}_{X'}(U) \to {\cal O}_{X}(H^{-1}(U))$, by $f \mapsto h^*f$,
 $h^*f (p) = h_p  f(q)$, where $q = h^{-1}(p)$, and $h_p: (R')_q \to R_p$ is the natural
 difference ring homomorphism.
 
 Note that if $U$ is an open subset of $X=\spe R$, $p \in U$,  $f \in {\cal O}_X(U)$, then
 for any sufficiently large 
 finitely generated difference subring $R'$ of $R$, there exists an open 
 $U' \subseteq X' = \spe R'$,
 and $f' \in {\cal O}_{X'}(U')$ such that if $i: R' \to R$ is the inclusion, 
 then $p \in (i^*)^{-1}(U') \subseteq U$, and $i^*f' = f |  (i^*)^{-1}(U')$.

\<{lem}  If $f \in {\cal O}_X(X)$, there exists a finitely generated $R'$, $X'=\spe R'$, 
 and $f' \in {\cal O}_{X'}(X')$ with $f = i^*f'$. \>{lem}
 
 \prf  For each $p \in X$, find $b_p \notin p$ and $a_p \in R$   such that
 for $p' \in \spe X$ with $b_p \notin p$, $f(p') = a_p/b_p \in R_{p'}$.  Then no prime
 difference ideal of $R$ can contain all $b_p$.  So the perfect ideal generated
 by some finite set $b_{p_1},\ldots,b_{p_n}$ is improper.  Also given $i,j$
 we have $b_{p_i}a_{p_j} - a_{p_i}b_{p_j} = 0$ in $R_q$, for any
 $q$ such that $b_{p_i},b_{p_j} \notin q$.  So every prime difference ideal
 containing $A(i,j)=Ann(b_{p_i}a_{p_j} - a_{p_i}b_{p_j} )$
 includes 
  $b_{p_i}$ or $b_{p_j}$.  Thus $b_{p_i}b_{p_j}$  
 lies in the perfect ideal generated by $A(i,j)$.  Some finite set $F$
 of elements of $R$ is involved in these various generation processes.  
 Now let $R'$ be 
 a subring with each $b_{p_i}$ and $a_{p_i}$ in $R'$, and containing $F$.  Then any prime
 difference ideal $p'$ of $R'$ avoids some $b_{p_i}$, and we can define
 $f'(p') = a_{p_i}/b_{p_i}$ for this $i$.  The statement is then clear.  \eprf
 
  In other words, each $f$ is determined
 by finite data. 
 
If $Y = \spe R$, 
a ring of the form ${\cal O}_Y(U)$ will be called an {\em affine ring}.  
\sk
\<{paragraph}{Affine rings of global functions}.

  $c \in R$ is a $\si$-unit if $c$ belongs to no transformally prime ideal.
More generally, $c |^\si d$ if for every transformally prime ideal $p$, for some
$b \notin p$, $c | db$.

It is easy to verify that a difference ring
$R$ without zero-divisors is an  affine ring iff for 
all $c,d \in R$, if $c |^\si d$ 
then $c | d$.  

Each of the two properties: well-mixed, affine, imply a property that
one might call ''residually local'': if $a \in R$, and the image of $a$
in every localization $R_p$ at a transformally prime ideal vanishes, then
$a=0$.  (For the well-mixed case,  cf.
\ref{well-mixed1.7}.)

\>{paragraph}

A {\em difference scheme } is a topological space $X$ together
with a sheaf ${\cal O}_X$ of well-mixed rings,
locally modeled on affine difference schemes of well-mixed rings.  In other words
$X$ has an   covering by open sets $U_i$; and there are
 isomorphisms $f_i: \spe (R_i) \to (U_i,{\cal O}_X | U_i)$
for some family of  well-mixed rings $R_i$.

A {\em morphism of difference schemes}  is a morphism of locally ringed spaces,
preserving the difference ring structure on the local rings.

If $Y$ is a difference scheme, a difference scheme {\em over Y} is a
difference scheme $X$ together with a morphism $f: X \to Y$.
If $Y = \spe D$, we also say that $X$ is over $D$.

A {\em pre-} difference scheme is defined in the same way, 
 but with general difference rings in place of well-mixed rings.
 \sk

{\em Gluing} Assume given a space
 $X$ with a    covering by open sets $U_i$ ;
a family of well-mixed rings $R_i$,
and
 homeomorphisms $f_i: \spe (R_i) \to  U_i$;
such that for
 any $i,j$, the map
$$f_j^{-1}f_i: \ f_i^{-1}(U_j) \subset \spe (R_i) \ \to \
f_j^{-1}(U_i) \subset \spe (R_j)$$ is induced by difference ring isomorphisms
between the appropriate localizations of $R_i$, $R_j$.  In this situation there
is a unique difference scheme structure on $X$, such that $f_i: \spe (R_i) \to  U_i$
is an isomorphism of difference schemes for each $i$ (with $U_i$ given the open subscheme
structure.)

If there are finitely many
$U_i$, each of the form $\spe (R_i)$ with $R_i$ a finitely generated difference ring,
(or difference $D$-algebra)
we will say that $X$ is of finite type (over $D$). 

\medskip

\noindent{\bf Remarks and  definitions.}

\st  We will usually consider difference rings $R$ that are finitely generated over
a difference  field, or over $\Zz$, or localizations of such rings. Now by
\cite{cohn}, Chapter 3, Theorem V (p. 89), such rings have no ascending
chains of perfect ideals.     So
their difference spectra are Noetherian as topological spaces.  (Their Cantor-Bendixon rank 
will be $<\omega^2$, not in general $< \omega$ as in the case of algebras.)

\st
A special place is held by difference rings in positive characteristic
$p$
whose distinguished homomorphism is the Frobenius $x \mapsto x^q$,
$q$ a positive power of $p$. Note that for
 such rings, $\si$
is injective iff the ring has no nilpotents, and surjective iff it is perfect.
They are always well-mixed.

\st
 $\spe (R)$ may be empty: e.g. $R= \Qq[X]/(X^2-1)$, $\si(X)=-X$.

However, this does not happen if $R$ is
 well-mixed: by \ref{well-mixed1}, $R$ has a perfect ideal $I \neq
 R$; by \cite{cohn} (or cf.  \ref{well-mixed1.8}) a maximal proper
 perfect difference ideal is prime.

\st 
Let $R$ be a ring, $\si$ a ring endomorphism of $R$.  Then $\si$ acts on
$\spec R$, taking a prime $p$ to $\si^*(p) = \si^{-1}(p)$.   The transformally prime ideals
are precisely the points of $\spec R$ fixed by this map.
(Note that though $\si^*$ is continuous on $\spec R$, the fixed point
set is rarely closed.)

   In some contexts it may be useful to consider not only
points of $\spec R$ fixed by $\si$, but also those with finite orbits, or
even topologically recurrent orbits.   (Cf. \cite{CHP}).
For instance, when $R$ is Noetherian,
the closed subscheme of $\spe R$ corresponding to a difference ideal $J$ may be empty
even if $J \neq R$, while this is avoided if  the wider definition is taken.  In this paper we take
a different approach, replacing arbitrary difference ideals by well-mixed ideals.

\<{sta}  \lbl{wm12}  \>{sta}
Let $R$ be a difference ring,
$\bar R$ be the ring of global sections of $\spe (R)$.
 There is a natural map $i: R \to { \bar R}$ 
  It induces a map
$i^*: \spe {\bar R} \to \spe {  R}$.

There is also a natural map in the opposite direction, $ \spe R \to \spe \bar R    $:

\[  p \mapsto p^*   = \{F \in {\bar R}: F(p) \in pR_p \}   \]

And $ {\bar R}_{p^*} \to R_p $ is defined by :
\[  F/G \mapsto F(p)/G(p) \]

  $p^*$ is the largest prime ideal of $\bar R$ restricting to $p$.
For if $q \meet R = p$, and $F \in q$ and $F(p) \notin p$, say $F(p) = c/d$
with $d \notin p$, $dF \in q$, $dF (p) = c$, $ c \notin  pR_p$,
but $c \in q \meet R$.

The image of the map $^*$ is a closed subset of  $\spe {\bar R}$; it consists of
the transformally prime ideals containing $bF$ whenever $F \in R^*$ vanishes
on every $p \in \spe R$ with $b \notin p$.  I do not know whether the additional primes
are of interest (or if they exist) in general.  In the well-mixed context, they do not: 

\<{lem}\lbl{wm13}  Let $R$ be a well-mixed ring, 
$X = \spe R$,     
$f \in  {\cal O}_X (X)  $, $p,q \in \spe R$.  Suppose $f(p)=0 \in R_p$.
 Then there exists $ b \in R$,
$b \notin p$, with  $bf(q)  = 0 \in R_{q}$.
\>{lem}

\proof Assume first that $X$ is  Noetherian (as a topological space).
Say $f(q) = \frac{c}{d} \in R_q$.  If $ce = 0, e \notin q$,
then $f(q)=0 \in R_q$ and we are done.  If $bc = 0, b \notin p$, then
$bf(q) = 0 \in R_q$ and we are done again.  Thus we may assume
$Ann(c) \subset p \meet q$.  

Since $R$ is well-mixed, so is $ Ann(c)$; by \ref{well-mixed1}, $I=\sqrt[\si] Ann(c)$
is perfect.  Thus 
 $I = p_1 \meet \ldots \meet p_n$ for some   transformally prime ideals $p_i$.

If $p_i \subset p \meet q$, then $f(p_i) = 0 \in R_{p_i}$ (since $p_i \subset p$)
and $f(p_i) = c/d \in R_{p_i}$ (since $p_i \subset q$), hence $c=0 \in R_{p_i}$,
contradicting $Ann(c) \subset p_i$.  Thus for each $i$, there exists 
$a_i \in p_i, a_i \notin q$ or else $b_i \in p_i, b_i \notin p$.  Let $a$ be 
the product of the $a_i$, $b$ the product of the $b_i$.  So $a \notin q,
b \notin p, ab \in I$; thus $(ab)^n \in Ann(c)$, some $n \in \Nn[\si]f$.
  Replacing $a,b$ by $a^n,b^n$, we obtain 
$a \notin q, b \notin p, ab \in Ann(c)$.  So $abc=0$, hence $bc = 0 \in R_q$,
and thus $b f(q) = 0 \in R_q$.  \qed

In general, let $R'$ be a difference subring of $R$ (hence also well-mixed),
$X' = \spe R'$, $f' \in {\cal O}_{X'}(X')$ such that $f$ is the pullback of $f'$
under the natural map $X \to X'$ (See \S \ref{localization},   end of the first paragraph.)
Let $p' = p \meet R', q' = q \meet R'$.  By the Ritt-Raudenbsh theorem and the Noetherian case,
there exists $b  \in R', b  \notin p'$, with $bf'(q')=0 \in R_{q'}$.  It follows
that $b \in R$, $b \notin p$, $bf(q)=0 \in R_{q}$.  \qed

\<{lem} \lbl{wm14}  Let $R$ be a well-mixed ring, 
  $X = \spe R$,   
$f \in  {\cal O}_X (X)  $, $p  \in \spe R$.  Then there exist  $a,b \in R$,
$b \notin p$, with $bf-a = 0$.  (I.e. $bf(a)-a = 0 \in R_q$ for all $q \in \spe R$.)
\>{lem}

\proof  By definition of $ {\cal O}_X (X)$
for some $a,b \in R, b \notin p$, we have:  $f(p) =  \frac{a}{b} \in R_p$. 
  If $bf-a$ satisfies the lemma, then so does $f$.
So we may assume $f(p) = 0 \in R_p$.  By \ref{wm13}, 
for each $q \in \spe R$ there exists $b_q \in R, b_q \notin p$,
with $b_q f (q) = 0 \in R_q$.  But by definition of $ {\cal O}_X (X)$ again,
for some $c,d,d'$ with $d,d' \notin q$ we have 
  $f(q') = \frac{c}{d} \in R_{q'}$ whenever $d' \notin q'$.  The equation
$b_q f(q) = 0 \in R_q$ means that there exists $d'' \notin q$ with
$d''b_q c= 0$.  So $b_q f(q') = 0 \in R_{q'}$ whenever $dd'd'' \notin q'$
Let $U_q = \{q': dd'd'' \notin q' \}$.  By compactness of $X$, finitely many
open sets $U_q$ cover $X$; say for $q_1,\ldots,q_n$.  Let
$b = b_{q_1} \cdot \ldots b_{q_n}$.  Then $bf(q' ) = 0 \in R_{q'}$ whenever 
$q' \in X$.  \qed

  If  $R$ is affine, without 0-divisors, then   $i: R \to \bar{R} $   is an isomorphism;
  so $X \iso \spe \bar{R}$.  This latter statement is true more generally:

\<{prop}\lbl{wm-localization} Let $R$ be a well-mixed ring, $X = \spe R$,
${\bar R} = {\cal O}_X (X)$.   Then $i^*: \spe {\bar R} \to X$ is an
isomorphism of difference schemes.  \>{prop}

\proof  By \ref{well-mixed1.7}, $i: R \to {\bar R} $ is an embedding.    Let $\bar{p} \in \spe \bar{R}$,
$p = \bar{p} \meet R$.  By \ref{wm14}, there 
 exist  $a,b \in R$,
$b \notin p$, with $bf-a = 0 \in \bar{R}$.  So $f \in \bar{p}$ iff $a \in p$.  This shows that $i^*$
is injective, and so induces a bijection of points; similarly \ref{wm14}, \ref{wm12} show  that $i$ induces
an isomorphism $R_p \to {\bar R}_{\bar p }  $.  \qed

\end{subsection}
\<{subsection}{Some functorial constructions}

\st
 Let $R$ be a finitely generated difference ring extension of an existentially closed difference
field $\U$.  The closed  points of $\spe R$ are the kernels of the difference
ring homomorphisms $f: R \to \U$; if $R$ is provided with generators $(a_1,\ldots,a_n)$,
then the closed point $\ker f$ corresponds to the point $(f(a_1),\ldots,f(a_n)) \in \U^n$.
For instance if say $R = \U [ X,\si X, \ldots ]$, then the closed points
of $\spe (R)$ become identified with the points of $\U$ (viewed as the affine line over
$\U$.)   

In general, if $X$,$Y$ are difference schemes, we define a point of $X$ with values
in $Y$  be a morphism $Y \to X$.  The set of $Y$-valued points of $X$ is 
denoted $X(Y)$.

\medskip

The action of $\si$ on $\spe U$ induces an action on the $U$-valued points of $R$
``by composition".  This should not be confused with the action of $\si$ on $\spec R$.
 In particular, $\spe R$ can be viewed as the set of points of $\spec R$ fixed
by $\si$. But this just means we are looking at difference ideals,
and has nothing to do with the fixed field of $\U$. 

This observation yields a construction of difference schemes, starting from an action on
an ordinary (but non-Noetherian) scheme.
I do not know whether or not every difference scheme is obtained in
this way; at least it seems not to be the case via a natural functor adjoint to
 $\fixs$.

\begin{definition} \lbl{fixed}\ 
Let $Y$ be a scheme, not necessarily Noetherian, and let   $\psi: Y \to Y$ be
a morphism of schemes.  Let $\fy$ be the subspace of the topological space
underlying $Y$, whose set of points is $\{p \in Y: \psi(p)=p \}$.  For an open 
$W \subseteq \fy$, let ${\cal F_0}(W) = lim_{W \subseteq U} {\cal O}_Y(U)$ be the 
 direct limit of ${\cal O}_Y(U)$ over
all open subsets $U$ of $Y$ containing $W$.   Let ${\cal F}$ be the sheafification of the presheaf ${\cal F_0}$.  \end{definition}

Notes:

(1)  $Y$ is not assumed to be Noetherian.

(2)  In general, $\fy$ is neither open nor closed in $Y$.  

(3)  For each open $W \subseteq \fy$, ${\cal F_0}(W)$ has a natural difference ring structure:
each ${\cal O}_Y(U)$ is a ring, and the ring homomorphisms
 $$\psi^*_U:  {\cal O}_Y(U) \to {\cal O}_Y(\psi \inv (U))$$
 induce first, by going to the limit on the right,
 maps ${\cal O}_Y(U) \to \lim_{W \subseteq U'} {\cal O}_Y(U') = {\cal F_0}(W)$
 (since $W \subseteq \psi \inv (U)$) and secondly, taking the limit on the left,
 a ring homomorphism $\si_W:  {\cal F_0}(W) \to {\cal F_0}(W)$.

(4) By (3), ${\cal F}$ is a sheaf of difference rings. Given an open $W \subseteq Y$,
   ${\cal F}_W$ denotes the restriction of ${\cal F}$
to $W$.

\lemm{fixed1}  
(1)  $\fy$ has a basis of open sets $W$ such that $(W,{\cal F}_W)$ is an affine difference
scheme.


\>{lem}
\prf

(1)   To verify that $\fy$ is locallly modeled on affine difference scheme, let $p \in Y^\psi$, and let $W_0  = \spec R$ be an open affine neighborhood of $p$ in $Y$.   Then
 $ \psi ^{-1} W_0 \meet W_0$ contains an open affine $W_1=\spec R[f^{-1}]$ of $\spec R$, with
 $p \in W_1$.  The scheme morphism $\psi$ restricts to $\psi: W_1 \to W_0$, corresponding to a ring homomorphism  $\phi_1:R \to R[f^{-1}]$.   Let $R_1= R[f \inv]$, $f_1 = \phi_1(f)$.
 
 Localizing $\phi_1$ yields a homomorphism 
 $\phi_2: R[f^{-1}] \to R[f^{-1},\phi_1(f)^{-1}]$, i.e. $\phi_2: R_1 \to  R_1[f_1 \inv]$.  Let $R_2 = R_1[f_1 \inv], f_2 = \phi_2(f_1)$, etc.  
 Let $R' = R[f^{-1}][f_1 ^{\inv}][f_2 \inv] \cdots$.  The $\phi_i$
 converge to an endomorphism $\si'$ of $R'$.  
 
 On the scheme side, for $n \geq 1$ let $W_{n+1} =W_n \meet \psi \inv(W_n)$.  
 Now  $f$ is a regular function on $W_0$, and $W_1 = \{q \in W_0: f \notin q \}$.  
   If $\theta: X \to W_0$ is any morphism of schemes, it follows that
   $\theta^*(f)$ is a regular function on $X$, and
 $\theta \inv (W_1) =\{q: \theta^*(f) \notin q \}$.  In particular with $X=W_1, \theta=\psi |W_1$ we
 see that  $W_2 = \{q: \theta^*(f) \notin q \} = \spec R_2$.   
 Continuing this way we obtain $W_n = \spec(R_n)$ for $n \geq 1$.  
 So  $\meet_n W_n = \spec (R')$ as subsets of $Y$, and in fact the induced
 topology  on $\meet_n W_n$ as a subspace of $Y$ agrees with the topology of $\spec(R')$
 
   Let $W=W_0 \meet \fy $; 
 then inductively, $W \subseteq W_n$ for each $n$.  So $W$ is an open neighborhood
 of $p$ in $\fy$, and a subspace of $\meet_n W_n = \spec(R')$.  But as a subspace
 of $\spec(R')$, $W$ is identified with the prime difference ideals.  Hence
 topologically $W \cong \spe(R')$.   The isomorphism of sheaves can be seen by looking at stalks.

%


%

%
%
%
\eprf

Most of the difference schemes we
will encounter will admit projective embeddings;
hence they  can all be constructed as $\fixs (Y)$ for some scheme $Y$.

A modification of the above approach does yield all difference schemes.
  Given a difference ring $R$, let
$\spec ' R$ be the set of prime ideals of $R$, with the following topology:
a basic open set has the form  $\meet_n \si^{-n} G$, where $G$ is a Zariski open
set.  Equivalently, a basic open set is the image of $\spec' R'$, where $R' $
is a localization of $R$ by  finitely many elements {\em as a difference ring.} 
Define a structure sheaf, and gluing; obtain a category that one might
call {\em transformation schemes}.  Extend  the functor $\fixs $ to this category,
 and obtain all difference
schemes in the image (but a Noetherian difference scheme need not be 
the image of a Noetherian transformation scheme.)

\smallskip

\begin{paragraph} {Products, pullbacks,   fiber products}
\st  Let ${\cal X},{\cal Y}$ be difference schemes, with underlying spaces
$X,Y$ and sheaves of rings $\cO_X, \cO_Y$ .   Let $\cO_{X,p}$, $\cO_{Y,q}$
denote the stalks at points $p,q$.  Let $R_{p,q} = \cO_{X,p} \tensor \cO_{Y,q}$
be the tensor product, with the natural difference ring structure.  Let
$i_X,i_Y$ denote the maps $(Id \tensor 1): \cO_{X,p} \to R_{p,q}$, resp.
$(1 \tensor Id): \cO_{Y,q} \to R_{p,q}$.

 We define the product ${\cal Z}  = {\cal X} \times {\cal Y}$.
As a set, we let $Z$ be the disjoint union over $(p,q) \in X \times Y$ of
\[ Z_{p,q} = \{ r \in \spe R_{p,q}:  (i_X)^{-1} (r) = p\cO_{X,p}, (i_Y)^{-1} (r) = q\cO_{Y,p} \} \]
If $r \in Z_{p,q}$, we let $pr_X (r) = p$, $pr_Y (r) = q$, $pr(r) = (p,q)$.
 Let $R_r$ be the localization
of $R_{p,q}$ at $r$.

Begin with open $U \subset X$ and $V \subset Y$, and
\[ f \in \cO_X(U)  \tensor \cO_Y(V) \]
$f$ defines a function $F$ on $pr^{-1} (U \times V)$; $F(r) \in R_r$ is the image of
$f$ under the natural homomorphism $\cO_X(U)  \tensor \cO_Y(V) \to R_r$.
A function such as $F$ will be called a basic regular function.
A set such as
\[ W_F = \{r \in pr^{-1}(U \times V): F(r) \notin rR_r  \} \]
will be called a basic open set.
Topologize $Z$ using the sets $W_F$ as a basis for the topology.

 Given an open $W \subset Z$, we let $\cO_Z(W)$
be the set of functions  on $W$ that agree at a neighborhood  $W$ of each point
with  a quotient $F'/F$, where $F,F'$ are basic regular functions and $W \subset W_F$.

It is easy to check that if  ${\cal X} = \spe R$, ${\cal Y} = \spe S$, then ${\cal Z}$
is isomorphic to $\spe (R \tensor S)$.

Similarly, if ${\cal X}, {\cal Y}$ are difference schemes over ${\cal U}$,
we can define ${\cal X} \times_{\cal U} {\cal Y}$.  
Alternatively it can
be defined as a closed difference subscheme of ${\cal X} \times {\cal Y}$,
see (13) below.
If ambiguity can arise as to the map $h: {\cal Y} \to {\cal U}$,
we will write ${\cal X} \times_{{\cal U},h} {\cal Y}$.  

\begin{notation} \lbl{Xy} Let $X$ be a difference scheme over $Y$.  If 
$y$ is a point of $Y$ with values in a difference field $L$, $y: \spe L \to Y$,
we let
$X_y = X \times _{y} \spe L$ 
\end{notation}

\st

 Let $R$ be a difference ring. A
difference-module is a module $A$ together with an additive $\si: A \to A$, such
that $\si(r) (\si(a))$ = $\si(r a)$ for $r \in R$, $a \in A$.
A sheaf of difference modules over a sheaf of difference rings
is defined as in \cite{HA}, II.5, adding the condition that the sheaf
maps respect $\si$.
 If $X$ is a difference scheme,
a   sub(pre)sheaf of ${\cal O}_X$, viewed as a (pre)sheaf of
difference modules over itself, is called a difference ideal (pre)sheaf.

We similarly take over the definition of a quasi-coherent sheaf.

\deff{subschemes}
Given a difference ideal sheaf ${\cal I}$, the stalk ${\cal I}_p$ of ${\cal I}$ at a point
$p \in X$ is a difference ideal of the local difference ring $\cO_{X,p}$.
define the associated closed subscheme
${\cal Y}$
as follows.  The underlying set is
\[ Y = \{p \in X: {\cal I}_p \neq \cO_{X,p} \} \]
with the topology induced from $X$.  Let
$\cO_Y(U)$ be the ring of maps $F$ on $U$ such that $F(p) \in \cO_{X,p} / {\cal I}_p$,
and $U$ admits a covering by open sets $U'$  such that $F | U' $ is represented
 by an element of $\cO_X(U')$ .

Note that ${\cal I}$ can be retrieved from $Y$ as a subscheme of $X$;
${\cal I}(U) $ is the kernel of the map $\cO_X (U) \to \cO_Y(U)$.

\st

Let $f: X \to Y$ be a morphism of difference schemes, and let ${\cal J}$ be a difference ideal
sheaf on $Y$.  Then $f^* {\cal J} $ defined as in \cite{HA} II.5 is a
difference ideal sheaf on $X$.    The closed subscheme associated with $f^* {\cal J} $
is the {\em pullback} of the closed subscheme associated with ${\cal J}$.

\st  Let $f:X \to Z$ and $g: Y \to Z$ be morphisms of difference schemes.  Then
we obtain a map $(f,g) : X \times Y \to (Z \times Z)$.
  Let $Z'$ be the diagonal closed
subscheme of $Z \times Z$ (defined by the obvious equations).
  The pullback of $Z'$ via $(f,g)$ is called the {\em fiber product}
and denoted $X \times_Z Y$.

\st   Let $X$ be a difference scheme, $W$ an open subset of the set of points
of $X$.  Define ${\cal O}_W$ to be the restriction of ${\cal O}_X$ to 
$W$.  Then $(W,{\cal O}_W)$ is a difference scheme.  (The question
is local, so we may assume $X=\spe R$; and also, since $W$ is the union
of open sets of the form $\{p \in X: a \notin p\}$, we may assume $W$ is of
this form.  But then $(W,{\cal O}_W) = \spe R[a^{-1},a^{-\sigma},\ldots]$.)

\end{paragraph}

\

\begin{paragraph}{Difference scheme associated to an algebraic scheme}

\sk   
 For any commutative ring $R$, there
exists a ring homomorphism $h: R \to S$ into a difference ring $S$, with the universal
property for such morphisms: any ring homomorphism on $R$ into a difference
ring factors uniquely as $gh$, with $g$ a difference ring homomorphism.  This universal
ring is denoted by $S = [\si]R$. If $k$ is a difference ring, the same
construction in the category of $k$-algebras is denoted $[\si]_k R$.
 For instance, $[\si]\Qq[X_0] = \U[X_0,X_1,\ldots]$
with $\si(X_i) = X_{i+1}$.

  Let $\cal C$ be the category of
 reduced, irreducible affine schemes over $k$, $\cal C'$
the category of
 reduced, irreducible affine schemes over $k$ with a distinguished endomorphism
compatible with that of $k$.

If $k$ is  an existentially closed difference field (a model of ACFA),
the
map $\spe [\si]_k R \to \spec R$ induces a homeomorphism on the closed points.

 This gives
a functor from $\cal C$ to $\cal C'$, adjoint to the forgetful functor
 ${\cal C'} \to {\cal C}$. Composing with the $\spe$ functor from
$\cal C'$ to difference schemes, we obtain a natural functor from affine varieties
over $k$ to affine difference varieties over $k$.   Call $\spe [\si]_k R$ the
difference scheme associated with $\spec R / \spec k$.

 The functor ${\cal C} \to {\cal C'}$ above
does not naturally extend to a similar functor on irreducible schemes to schemes with
endomorphisms.   (In the affine category, the image of an open subset of a scheme
may not be an open subset of the image; rather a countable intersection of open
subsets.)
Similarly,  while irreducible affine difference schemes are the same as irreducible affine
schemes
with endomorphisms, the notion of gluing is different, reflecting the fact
that the endomorphism is viewed as part of the algebraic structure.

These differences cancel out;   the composed functor, taking an affine variety to the
associated affine difference
scheme, does extend to a functor on schemes.  (By gluing.)   For projective schemes,
we will see   the associated difference scheme can also be described in another way,
defining the difference analog of the functor Proj.

\end{paragraph}

\end{subsection}

\begin{subsection}{Components } \lbl{irred}

For ordinary algebraic schemes, we have two notions, of a {\em reduced}
and an {\em irreducible} scheme.  For difference schemes a third type of
reducibility arises.

A $\si$-nilpotent is an element $a$ of a difference ring $ R$
such that some product of elements
$b_i = \si^{m(i)}(a)$ vanishes.
$R$ is {\em perfectly reduced} if it has no $\si$-nilpotents.
Equivalently, $0$ is a perfect ideal.

$R$ is {\em transformally reduced} if $\si(a)=0$ implies $a=0$ in $R$.
 
By way of contrast, if a difference ring  $R$, viewed simply as a ring,
has no nilpotent elements
(is an integral domain,  has $\spec R$ irreducible)
we will say that $R$ is  {\em algebraically reduced} ( algebraically integral,
algebraically irreducible.)

\begin{lem}\lbl{red-localize}  Let
 $R$ be a difference ring,  $S$ is a subset closed under multiplication
and under $\si$.  Let  $R[S^{-1}]$ be the localization.  If $R$ 
is well-mixed / perfectly reduced / algebraically reduced /
algebraically irreducible, then so is $R[S^{-1}]$.  \end{lem} \proof
Straightforward verification.

\begin{definition}   \
  A difference scheme $X$ is {\em irreducible} if the underlying topological
space is irreducible, i.e. it is not the union of two proper closed subsets.
$X$ is
 {\em
well-mixed / perfectly reduced/ algebraically reduced / algebraically
irreducible / algebraically integral} if for any open $U$, the rings ${\cal O}_X (U)$ have the
property.  $X$ is {\em transformally integral} if $X$ is perfectly
reduced and irreducible.
\end{definition}

Note that transformally integral implies algebraically integral.  Also,  $X$ is well-mixed / perfectly reduced/ algebraically reduced 
iff this property holds for each local ring at a point $p \in X$.   

\begin{notation} If $X$ is a difference pre-scheme, then the ideals
$red_{wm}(R_p)$ (smallest well-mixed ideals of the local rings) generate an ideal sheaf on the
structure sheaf ${\cal O}_X$; it defines a closed subscheme, the
(well-mixed) subscheme $X_{wm}$ of $X$.  We may similarly  
defined the underlying perfectly reduced subscheme, and the somewhat
thicker $X_{wm,red}$ (underlying algebraically reduced, well-mixed subscheme.)
\end{notation}

Every Noetherian   topological space $X$ is a union of finitely many
{\em irreducible components } $X_i$.  The $X_i$ are defined by the following: they are
 closed subsets, no one contained
in another, and $X = \union_i X_i$.

\begin{definition} \ \lbl{component} \end{definition}  Let $X$ be a difference scheme, and let
  $Z$ be an irreducible component of $X$,  defined by a prime ideal sheaf $ {\cal P} \subset
{\cal O}_X$.  Define an ideal sheaf ${\cal I} $ as follows:
$${\cal I} (U) =  \{a \in  {\cal O}_X (U) : ab=0, \hbox{ some } b \notin {\cal P}(U)  \}$$
The {\em sub-difference scheme of $X$ supported on $Z$} is defined to
have underlying space $Z$, and structure sheaf $({\cal O}_X / {\cal I} ) | Z$.

If $Y$ is the sub-difference scheme of $X$ supported on $Z$, then the corresponding
perfectly reduced scheme  will be called the
perfectly reduced subscheme supported on $Z$; and similarly for other reducedness
notions.

\begin{definition} Let $X$ be a  difference scheme over a difference
field $k$ .  We will say that a property holds of $X$ {\em absolutely}
if it holds of $X \times_{\spe k} \spe K$ for any difference field
extension $K$ of $k$.  \end{definition}

We attach here some lemmas on analogs of separability.

\begin{definition} \begin{enumerate}
\item  A difference domain  $D$ is {\em inversive} if $\si_D $ is an automorphism of $D$.
\item  Let $D$ be a difference domain.  Up to isomorphism over $D$, there exists a unique
inversive difference domain $D^{\invv} $ containing $D$, and such that
$D^{inv} =  \union_m D^{\si^{-m}}$.  It is called the inversive hull of $D$.
\item  Let $K \subset L $ be difference fields.  $L$ is {\em transformally separable }
over $K$ if $L$ is linearly disjoint over $K$ from $K^{inv}$.
\item  Recall also the classical definition of Weil:  Let $K \subset L$ be fields.
$L/K$ is a regular extension if
$L$ is linearly disjoint from $K^{alg}$ over $K$.
\end{enumerate}
\end{definition}

\<{lem}\lbl{inv-ext} Let $D$ be a difference domain, $p$ a transformally prime ideal
of $D$.  Then there exists a unique prime ideal $p^{\invv}$ of $D^{\invv}$ such
that $p^{\invv} \meet D = p$.  The natural map $\spe D^{inv} \to \spe D$ is a bijection at
the level of points.\>{lem}

\proof  If $p^{\invv}$ exists, then clearly $a \in p^{\invv}$ iff $\si^n(a) \in p$
for some $n$.  So define $p'$ in this way: $p' = \union_{n \in \Nn} \si^{-n}(p)$.  It
is easy to check that $p'$ is prime, $\si^{-1}(p')=p'$, and $p' \meet D = p$.
 
\begin{lem}  Let $X$ be an algebraically integral
difference scheme over a difference field $k$ \begin{enumerate}
\item  If for some algebraically closed field $L \supset k$,
$X \tensor _k L$ is an integral domain , then $X$ is absolutely
algebraically integral.

\item  Let $K$ be the inversive hull of $k$.  If $X \tensor _k K$ is
transformally reduced, then $X$ is absolutely transformally reduced.

\item  If for some algebraically closed difference field $L \supset k$,
$X \tensor _k L$ is an integral domain, then $X$ is
absolutely transformally integral.
\end{enumerate}
\end{lem}

\proof   The question reduces to the case $X = \spe D$, $D$ a $k$-difference algebra
and a domain.
\noindent{(1)  \ } Let $D' = D \tensor_{k} K$.  Whether or not $D'$ is a domain clearly
depends only on the $k$-algebra structures of $D$ and $K$, and in particular the
choice of difference structure on $K$ is irrelevant.  Moreover for any
field extension $L$ of $k$, we may find $L'$ extending $K$ and containing a copy of $L$;
then by algebra $D \tensor_{k} L'$ is a domain, hence so is $D \tensor _{k} L$.

\noindent{(2)  \ }  Let $D' = D \tensor_{k} K$.  Let $L'$ be a difference field extending $K$.
Let $e = \sum d_i \tensor c_i$ be a nonzero element of
$D' \tensor _K L'$, with $d_i \in D, c_i \in L'$.
We may choose the $c_i$ linearly independent over $K'$.
It follows (using inversivity) that the $\si(c_i)$ are linearly independent over $K'$ .
But then $\si(e) = \sum \si(d_i )\tensor \si(c_i) \neq 0$.  Now if $L$ is a difference
field extension of $k$, the inversive hull $L'$  of $L$ extends $K$.  If
$0 \neq e \in D \tensor _k L$, then the image of $e$ in $D \tensor_k L'$ is nonzero,
since we are dealing with vector spaces.  Thus $\si(e) \neq 0$ as required.

\noindent{(3)  \ } Follows from (1) and (2).

\begin{remark} \ \end{remark}
Note that
 $[\si]\Qq[\sqrt{2}]$ is not a domain.  However let $R$ be a $k$- algebra
and assume $R \tensor _k k^{alg}$ is a domain; then
Then $[\si]_k R$ is also a domain.
 More generally, suppose $D$ is a difference domain, subring of
an algebraically closed difference field $L$. Let $R$ be a
$D$-algebra, and suppose $R \tensor_D L$ is a domain. Then
$ R^* = [\si]_D R$
is a domain.

\begin{definition} \ \end{definition}
Let $k$ be a difference field. An {\it affine difference
variety} over $k$ is a   transformally integral affine difference scheme of finite type
over $k$.

\end{subsection} 

\end{section} 

\begin{section}{Dimensions}
\lbl{dimensions}
\<{subsection}{Total dimension and transformal dimension}

Two types of dimensions are naturally associated with difference equations.
If one thinks of sequences $(a_i)$ with $\si(a_{i})=a_{i+1}$, the
{\em transformal dimension} measures, intuitively,the eventual number of degrees
of freedom in choosing $a_{i+1}$, given the previous elements of the
sequence.  The {\em total dimension} measures the sum of all degrees
of freedom in all stages.   These correspond to transformal transcendence
degree and order in  
\cite{cohn}, but we need to generalize them to difference schemes (i.e. mostly
from difference fields to difference rings.)

The dimensions we consider here and later will take extended natural number values
($0,1,\ldots,\infty$).   We will sometimes define the dimension as the  the maximal integer
with some property; meaning $\infty$ if no maximum exists.

 Since these dimensions are defined in terms of integral domains,
they give the same value to a difference ring  $R$ and to $R/I$, where $I$ is the smallest
well-mixed ideal of $R$.

  If $k$ is a difference field and $K$ is a
difference field extension, a subset $B= \{b_1,\ldots,b_m\}$ of $K$ is {\em transformally
independent over $k$} if $b_1,\si(b_1),\ldots,b_2,\si(b_2),\ldots,b_m,\si(b_m),\ldots$
are algebraically independent over $k$.  The size of a maximal $k$- transformally 
independent subset of $K$ is called the transformal dimension.

Let $k$ be a difference  field, and let 
$R$ be  a difference $k$-algebra.   Consider the set $\Xi$ of triples $(h,D,L)$,
with $D$ an algebraically integral difference $k$-
algebra,    $h: R \to D$ 
 is a surjective homomorphism of difference $k$-algebras, and
$L$ is the field of fractions of $D$.
   Let $\Xi '$ be the set of 
$(h,D,L) \in \Xi$ where  in addition $D$ is a difference 
domain.  In general, $L$ is only a field; but if $(h,D,L) \in  \Xi '$, then $L$ carries
a canonical difference field structure.

The {\em transformal dimension } of $R$ over $k$
is the supremum over all   $(h,D,L) \in \Xi'$ of the transformal dimension of $L$
 over $k$.   
 
The {\em reduced total dimension}  of $R$ over $k$ is the supremum over
all $(h,D,L) \in \Xi ' $ of 
the transcendence degree of $L$ over $k$.

When $R$ is a finitely generated  difference $k$-algebra,
the   {\em total dimension}  of $R$ over $k$ is the supremum over $(h,D,L) \in \Xi$
of the transcendence degree of $L$ over $k$.

\<{lem}\lbl{dim-trans}  Let $k$ be a difference field, $R$ a   difference $k$-algebra.  Then the transformal dimension of $R$ over $k$ is the maximal
$n$ such that the transformal polynomial ring in $n$-variables $k[X_1,\ldots,X_n]_\si$ embeds into $R$ over $k$.  \>{lem}

\proof  If the transformal dimension is $\geq n$, let $(h,D,L) \in \Xi'$ show it; then
$D$ contains a copy of $k[X_1,\ldots,X_n]_\si$.  We have $X_i = h(Y_i)$
for some $Y_i \in R$; and $k[Y_1,\ldots,Y_n]_\si$ must also be a copy
of the transformal polynomial ring.  Conversely, assume
$ k[X_1,\ldots,X_n]_\si \iso _k S \leq R$.  Let $I$ be a maximal well-mixed ideal
with $I \meet S = (0)$ (note the ideal $(0)$ has this property.)  

\claim $I$ is an algebraically prime ideal:   if $ab \in I$, $b \notin I$, then $a \in I$.

\proof    First suppose $a \in S$.  Then $Ann(a/I)$ is a
well-mixed ideal (\ref{well-mixed1.5}), and $Ann(a/I) \meet S = (0)$ since
$S$ is an integral domain.  As $b \notin I$, $Ann(a/I) \neq I$. 
By maximality,   $ 1 \in Ann(a/I)$, i.e. $a \in I$. 

 Now in general:  $Ann(b/I)$ is a well-mixed ideal.  If 
$c \in Ann(b/I) \meet S$, then $bc \in I$, so by the case just covered, as $c \in S$,
 $c \in I$; so $c = 0$.  Thus $Ann(b/I)$ is a well-mixed ideal meeting $S$ trivially,
and $1 \notin Ann(b/I)$; 
so again by maximality, $Ann(b/I) = I$; thus $a \in I$.   

\claim  $I$ is transformally prime.

\proof  Let $J = \si^{-1}(I)$.  Then $J$ is well-mixed (if $\si(ab ) \in I$,
then $\si(a \si(b)) = \si(a) \si^2(b) \in I$.)  Also $J \meet S = (0)$.
 As $I \subset J$, we have $I =J$.

Now $R/I$ shows that  the transformal dimension of $R$ is $\geq n$.  \qed

\begin{lem}\lbl{dim-alg} Let $k$ be a  difference field $k$,
$R$ a difference $k$-algebra.  Then $R$, $R/rad_{wm}(R)$ have the same
total dimension.  If $R$ is well-mixed and finitely generated, then the following numbers are
equal: \begin{enumerate}
 \item The total dimension $t.dim(R)$ of $R$.
\item The maximal 
transcendence degree $t_2$ over $k$ of the fraction field of a 
quotient of $R$ by a prime ideal.   
\item $t_3=$ the maximal $n$ such that the
polynomial ring over $k$ in $n$ variables embeds into $R$.
\item $t_4=$ the maximal Krull dimension of a finitely generated $k$-subalgebra of $R$.
 \end{enumerate} \end{lem}

\proof  Any difference ring homomorphism of $R$ into a
difference ring without zero divisors factors through $R/rad_{wm}(R)$,
so the first point is clear.  Now assume $R$ is well-mixed, of total
dimension $n$.  Clearly $n \leq t_2$,  while $t_2 \leq t_4$ since one
can lift a a transcendence basis from the quotient ring to a set of
elements of $R$, necessarily independent over $k$.    We have $t_4 \leq t_3$ by
standard commutative algebra. Finally, $t_3 \leq
t.dim(R)$: $R$ contains a free polynomial ring $R_1$ of rank $t_3$. 
By \ref{well-mixed1.8R}, there exists a cofinally minimal (prime) ideal
$p$ of $R$ with $p \meet R_1 = (0)$; by \ref{well-mixed1.8}, $p$ is a
difference ideal.  The quotient $R/p$ clearly has field of fractions
with transcendence degree $\geq n$.  \qed

If $R'$ is a finitely generated difference subalgebra of a well-mixed difference $k$-algebra $R$,
by \ref{dim-alg} (3), $t.dim(R') \leq t.dim(R)$.  
 In what follows, the total dimension over $k$ of a well-mixed difference $k$-algebra $R$ is defined to be the supremum over all $R'$ of the total dimension of $R'$, where $R'$ is a finitely generated difference
$k$-subalgebra of $R$.   
Observe that $R$ has the same total dimension as $R/\sqrt{0}$.  If
$R$ is not well-mixed, define the total dimension to be that of $R/J$, with $J=rad_{wm}(J)$ the smallest
well-mixed ideal of $B$.

\<{lem} \lbl{tdim-prep} Let $R$ be a finitely generated difference
$k$-algebra.  Then $t.dim(R) =\sup_{p \in \spe R} t.dim(R_p)$.  
Hence $t.dim(R) \geq t.dim (R[a \inv]_\si)$ when $a \in R$.  \>{lem}

\<{proof}  We may assume $R$ is well-mixed.  If
$t.dim(R_p) \geq n$, let $h': R_p \to D'$ be surjective, $D'$ an
algebraically integral difference $k$-algebra, $L$ the field of fractions fo $D$,
$tr. deg._k(L) =n$.  Let $h$ be the composition $R \to R_p \to D'$, and
let $D = h(R)$.  Any element of $D'$ has the form $h'(a/b)$ for some 
$a,b \in R, b \notin p$.  Since $b$ becomes invertible in $R_p$, $h(b) \neq 0$,
and $h'(a/b) = h(a)/h(b)$.  Thus $L$ is also the field of fractions of $D$, so
$t.dim(R) \geq n$.  Conversely suppose $t.dim(R) \geq n$, and let 
 $(h,D,L) \in \Xi$ show it.  Let $P = \ker (h)$, an algebraically prime difference ideal of $R$.  Let
 $p = \srad{P} = \union_{m \geq 1} (\si^{m}) \inv (P)$.  Then $p \in \spe R$.
 Since $R \to D$ is surjective, and $P \subseteq p$, we have $p = h \inv (q)$
 for some transformally prime ideal $q$ of $D$.  
 If $a \notin p$ then $h(a) \notin q$; so $h$ extends to a surjective homomorphism
 $R_p \to D_q$.  Since $D_q$ has the same field of fractions as $D$, 
 this shows that $t.dim(R_p) \geq n$.  \eprf

  Let $X$ be a  difference scheme over $\spe k$.
The   transformal  (resp. total, reduced total) dimension of $X$ over $k$
is the maximum transformal  (total, reduced total) dimension of $\cO_X(U)$
over $k$, where $U$ is an open affine subset of $X$.  

If $R$ is a finitely generated well-mixed difference $k$-algebra, $X = \spe R$, then
the total dimension of $R$ over $k$ equals the total dimension of $X$.
Indeed let $\bar{R} = {\cO_X}(X)$.  By Proposition \ref{wm-localization}, the sets
of local difference $k$-algebras $\{\bar{R}_q: q \in \spe \bar{R} \}$,
$\{R_p: p \in X \}$ coincide.  So by Lemma \ref{tdim-prep}, $t.dim(\bar{R}) = t.dim(R)$;
whereas any localization (and hence the ring of global sections) has smaller or equal
dimension, so it does not change the supremum.

 If $f: X\to Y$ is a morphism of difference schemes, the transformal (resp. total)
dimension of $X$ over $Y$ is the supremum of the corresponding dimension of
\mbox{$X \times_Y \spe k$}
 over $\spe k$,  where $k$ is a difference field and $\spe k \to Y$ is a $k$-valued point
of $Y$.  (cf.  \ref{total-dim-base-change} for the coherence of these
definitions, when $Y = \spe k$.

 \begin{prop}\lbl{trans-0}  Let $k$ be a difference field, and $X$ a well-mixed 
difference scheme of finite type 
over $k$.  The following conditions are equivalent:
 
\begin{enumerate}
\item $X$ has finite total   dimension
  over $\spe k$
\item $X$ has finite reduced total dimension over $\spe k$
\item $X$ has   transformal dimension zero over $\spe k$ 
\end{enumerate}
\end{prop}

\proof  Evidently  (1) \Implies  (2)  \Implies  (3) .  To prove that     (3) \Implies   (1)  ,
we may assume $X  = \spe R$, $R$ a   finitely generated  
  difference $k$-algebra.   

By \ref{dim-trans}, $R$ does
not contain a copy of the difference-polynomial ring in one variable
$k[t,t^\si,\ldots]$.

Let $r_1,\ldots,r_n$ be generators for $R$, $r_{ij} = \si^j(r_i)$.
Let $\Xi$ be the set of triples $(h,D,L)$ as in the definition of total dimension.
For each $i$, and each $(h,D,L) \in \Xi$, for some $m$, there exists a nontrivial
polynomial relation $F(h(r_i),\ldots,h(r_{i,m})) = 0$, $F \in k[X_0,\ldots,X_m]$.

 By compactness, there exists $m$ and finitely many $F_1,\ldots,F_r \in k[X_0,\ldots,X_m]$ 
such that for all $i$ and all $(h,D,L) \in \Xi$, for some $j \leq r$, 
$F_j(h(r_i),\ldots,h(r_{i,m})) = 0$.   By the claim below, $tr. deg._k(L) \leq mn$,
finishing the proof.  

\claim  Let $k$ be a difference field, $R=k[a,a^\si,\ldots]$ a difference $k$-algebra
with no zero divisors, $L$ the field of fractions of $R$.   Write $a_i = \si^i(a)$,
and assume $F(a_0,\ldots,a_m)=0$, $0 \neq F \in k[X_0,\ldots,X_m]$.  Then
$tr. deg._k (L) \leq (m)$.  

\proof   It suffices to  show that $a_{k } \in k(a_0,\ldots,a_{k-1})^a$ for all $k \geq m$.  We have $F^{\si^{k-m}} (a_{k-m},\ldots,a_k) = 0$.   Let $l$ be least
such  that for some $0 \neq G \in k[X_0,\ldots,X_{l}]$ we have 
$G(a_{k-l},\ldots,a_k) = 0$.  Then $a_{k-l},\ldots, a_{k-1}$ are algebraically
independent over $k$.  (If $H(a_{k-l},\ldots, a_{k-1})=0$ then
$H^\si(a_{k-(l-1)},\ldots,a_k) =0$, lowering the value of $l$.)  So
$a_k \in k(a_{k-l},\ldots, a_{k-1})^a  \subset k(a_0,\ldots,,a_{k-1})^a$.  \qed

\<{cor} \lbl{trans0-r}  If $R$ is a $k$-difference algebra with generators $a_0,
\ldots,a_n$,  and $R$ is an integral domain with field of fractions $L$,
$L_0 = k(\{\si^j(a_i): i \leq n, j \leq d \}) $, 
and $tr.deg._k L_0 \leq d$ , then $tr.deg._k L \leq d$. 
\>{cor}  

\proof  By the Claim of \ref{trans-0}, we have $\si^j(a_i) \in (L_0)^a$ for each $i$ and
all $j$.   \qed

\begin{cor}\lbl{total-krull-0}  Let $k$ be a difference field, 
$R$ be a well-mixed difference $k$-algebra.  Assume $R$ is finitely
generated as a $k$-algebra.  Then the total dimension of $R$ (as a
difference $k$-algebra) equals the Krull dimension of $R$.  If this
dimension is $0$, then $\dimk R < \infty$.
\end{cor}

\proof   
As $R$ is a finitely generated $k$-algebra, the equality follows follows
from Lemma \ref{dim-alg}.  A finitely generated $k$-algebra of Krull
dimension $0$ is   finite dimensional over $k$ (explicitly: let $I$ be the nil
ideal of $R$; $I = \meet_{i=1}^n p_i$, with $p_i$  prime.  
$R/p_i$ is a finite field extension of $k$; so
$\dimk R/p_i < \infty$, hence also $\dimk R/I < \infty$.  
Now $R$ is Noetherian, so
$I$ is finitely generated; thus $I^l/I^{l+1}$ is a finite-dimensional $R/I$-space for each
$k$; so $\dimk I^l/I^{l+1} < \infty$ for each $l \leq r$.  
It follows that $\dimk R < \infty$.) \qed

\begin{lem}\lbl{dim-cl}  Let $X$ be a well-mixed
difference scheme of finite type
over a difference field $k$.  Let $W$ be a   subscheme  of $X$, 
$\bar{W}$   the closure of $W$ in $X$.  Let
$Z = \bar{W} \setminus W$.     Then   
 the transformal (total) dimension of $Z$ and of $W$ is at most that of $X$.
If $W$ has finite total dimension, then $Z$ has smaller total dimension.  
(Hence $W,\bar{W}$ have the same dimension.)
\end{lem}

Here $W$ is a closed subscheme of an open subscheme of $X$.
The closure of $W$ in $X$ is the smallest 
well-mixed subscheme of $X$ containing $W$.

\proof  The statement regarding transformal dimension is obvious;
the one for total dimension follows, say, from \ref{dim-alg} (4).
(A similar statement is true for pro-algebraic varieties.)

\<{example}   One cannot expect a strict inequality for transformal dimension.  \rm 
Consider the subschemes of $\Aa^3$ defined by $x y = \si(x) z$ vs.
$x=0$, or just any $0$-dimensional scheme vs. a point.\>{example}

\begin{remark}\lbl{obviousmonotonicity} \ \end{remark}
If $X$ is a difference scheme over $Y$,
  $\bar{Y} $ is a closed subscheme of $Y$, and ${\bar X}$ is the pullback to $X$,
then the transformal (total) dimension of ${\bar X}$ over ${\bar Y}$ is bounded
by that of $X$ over $Y$. Any closed subscheme of
${  X}$ will also obviously have dimension $\leq m$ over ${Y}$.

\begin{lem} \lbl{total-dim-base-change}  Let $k$ be a difference field, 
$R$ a difference $k$-algebra, of total dimension $l$.  Let 
$K$ be a difference field extension of $k$.  Then $R \tensor_k K$ has
total dimension $\leq l$ over $K$.  \end{lem}

\proof Let $(h,D,L)$ be:  an algebraically integral difference $K$-algebra $D$,
a surjective homomorphism $h: R \tensor_k K \to D$, with $L$ the field
of fractions of $D$.  We must show that $tr. deg._K (L) \leq l$.  Let
$h' = h |R$, $D' = h'(R)$, $L'=$ the field of fractions of $D'$ within $L$.
Then in $L$, $L$ is the field amalgam of $K$,$L'$.  Thus 
$tr. deg._K (L) \leq tr.deg._k (L') \leq l$.  \qed

\begin{remark} \lbl{total-dim-base-change-rem} \ \end{remark}
If $K/k$ is a regular field extension,  
or if $R/k$ is a contained in a regular extension, 
equality holds in \ref{total-dim-base-change}.  

In general it may not, because
of possible incompatibility within $k^{alg}$:  if $k = \Qq$, $a \in R, b \in K$
with $a^2=2,  \si(a)=a, b^2=2, \si(b)=-b$, then regardless of the total dimension
of $R$, the total dimension of $R \tensor _{\Qq} K$ equals zero.

If one counts total dimension only with respect to difference domains lying
within a given universal domain (model of ACFA), and this universal domain
contains $K$, then again this dimension is  base-change invariant.

\begin{lem}  \lbl{inv-base-change} Let $k$ be a difference field, 
$R$ a difference $k$-algebra  of total dimension $l$.  Then:

(1)  There exists a difference field $\bar{k}$ extending $k$, with underlying field
structure equal to $k^{alg}$, 
such that $R \tensor _k \bar{k}$ has total dimension $l$ over $\bar{k}$.

(2)    $R \tensor _k k^{inv}$
has total dimension $l$ over $k^{inv}$. \end{lem}

\prf  In both cases, the total dimension is at most $l$ by  Lemma \ref{total-dim-base-change}.
  Thus we must show that it equals at least $l$.  We may pass to 
to a quotient of $R$ demonstrating that $R$ has total dimension $\geq l$;
i.e. we may assume $R$ is an integral domain, whose field of fractions $K$ has
$k$-transcendence degree $\geq l$.  

To prove (1), let $R'$ be the ring generated by $R$ and $k^{alg}$ within $K^{alg}$.  
Then $R'$  is an integral extension of $R$.
By the lemma on extending homomorphisms, the homomorphism $\si: R \to R \subset K$
admits an extension $\bar{\si}: R' \to K^{alg}$.  But
$\bar{\si} (R) \subseteq R$, $\bar{\si}(k^{alg}) \subseteq k^{alg}$, so  
$\bar{\si}(R') \subseteq R'$.   Let $\bar{k} = (k^{alg}, \bar{\si})$.  Then $R'$ is a difference ring of total dimension
$l$ over $\bar{k}$.  Since $R'$ is a homomorphic image of $R \tensor_k \bar{k}$, the same
is true of the latter.  

(2)  If $K/k$ is a regular field extension, then $R \tensor _k k^{inv} \leq K \tensor_k k^{inv}$ is an integral domain, and the assertion is clear.  In particular, (2) is true when $k$ is algebraically
closed.   Thus 
$t.dim (R \tensor_k \bar{k}) \tensor_{\bar{k}} \bar{k}^{inv} = t.dim (R \tensor_k \bar{k}) = l$.
But $(R \tensor_k \bar{k}) \tensor_{\bar{k}} \bar{k}^{inv} =  R \tensor_k \bar{k}^{inv} 
= (R \tensor_k k^{inv} ) \tensor_{k^{inv}} {\bar{k}^{inv}}$.  Hence by 
Lemma \ref{total-dim-base-change}, $t.dim (R \tensor_k k^{inv} ) \geq l$.  \eprf

Note that  the reduced total dimension can certainly go down upon base change
to $K^{inv}$, even if
$k=k^{alg}$; e.g.  $R = k(t)$, $K=k(\si(t)) \subset R$.

\begin{paragraph}{Embedding in the transformal line}

A difference domain $R$ is {\em twisted-periodic} if for some $n$, $p$,$m$, $R \models (\forall x) \si^n(x)=x^{p^m}$.  If a subring $R$ of a difference field $K$ is 
not   twisted-periodic, then   $R^n$ is Ritt-Cohn dense in 
$K^n$, i.e. every $K$-difference  variety containing $R^n$ equals $K^n$ (cf. \cite{cohn}.  Note that if every element of $R$ satisfies a difference equation of
total dimension $n$, then every element of $K$ satisfies one of total dimension $\leq 2n$.
In one dimension, polarize.)

Below, a morphism $f:X \to Y$ of difference schemes over a field $K$ will be
said to be
{\em point- injective } if for every difference field extension $K'$ of
$K$, $f$ induces an injective map $X(K') \to Y(K')$.

\<{lem}\lbl{tpl}
  Let $R$ be a   difference domain, with field of fractions $K$; let $r: R \to k$
be a surjective homomorphism into a difference field $k$, with $k$ not twisted-periodic. Assume $R$ is a valuation ring.
Let $X \subset \Aa^n$ be an   difference scheme over $R$, of finite total dimension.
 Then there exist $c_1,\ldots,c_n \in R$, such that if $H(x) = \sum c_i x_i$, $h(x) = \sum {\bar c}_i x_i$,
then $H$ and $h$ are point- injective.
   \>{lem}

\proof  If $\sum c_i x_i = \sum c_i y_i$, $x,y \in X $ distinct, then $c = (c_1,\ldots,c_n)$ has
transformal dimension $<n$ over $K(x,y)$.  Thus $\{c: (\exists x \neq y \in X) c \cdot x = c \cdot y  \}$
has transformal dimension $\leq (n-1) + 2 trans. dim (X) = n-1$.  So for any $c$, outside
a proper difference subscheme of $\Aa^n$, $H(x) = \sum c_i x_i$ is injective.   \qed

\<{cor}\lbl{fin-dim-1}  Let $X$ be a projective difference scheme of finite total dimension over a field $k$,
with $k$ not twisted-periodic.
Then there exists a linear projection to $\Pp^1$, defined over $k$,    
point- injective on $X$.  \>{cor}

\proof  Say $X \subset \Pp^n$.  The set of hyperplanes passing through some point of $X$ is contained
in a difference scheme of transformal dimension $n-1$.  So there exists a hyperplane defined over $k$
and avoiding $X$.  Thus we may assume $X \subset \Aa^n$.  Now \ref{tpl} applies.  

(Improve this to:  $X$ is isomorphic to a difference subscheme of $\Pp^1$.)

\end{paragraph}

\end{subsection} 

\begin{subsection}{Dimensions of generic specializations}

\begin{definition}  Let $X$ be a difference scheme over a difference domain
$D$ with field of fractions $K$.  A property   will be said to hold {\em
generically } of $X$ if it is true of $X \times_{\spe D} {\spe K}$.

   For example, if
$R$ is a  difference  $D$-algebra, the {\em generic} transformal (total)
dimension of $R$ over $D$ is the transformal (total) dimension of
$R \tensor_D K$ over $K$.
\end{definition}

\begin{lem}\lbl{5.11} Let $D$ be a difference domain, $R$ a finitely generated
difference $D$-algebra, of generic transformal (total) dimension $d$ over $D$.
Then for some dense open $U \subset \spe(D)$, for all $p \in U$,
if $h: D \to K$ is a map into a difference field with kernel $p$,
then $R \tensor_D K$ has transformal (total) dimension $\leq d$
 over $K$. 
 \end{lem}
\proof   It is easy to give an effective argument (cf. \ref{si-dim}.)  We give a qualitative one here.
Suppose the lemma is false; then there exist difference ring homomorphisms
$h_i: R \to D_i \subset L_i$ such that $h_i(R)=D_i$ is a domain with field of fractions $L_i$,
$h_i(D)$ has field of fractions $K_i$, and  $L_i$ is a difference field of  transformal dimension $>d$ over $K_i$,
(respectively, $tr. deg._{K_i}L_i > d$); and such that $p_i = ker(h_i | D)$ approaches the generic point of $D$.
Let $F$ be a finite set of generators of $R$ as a $D$-algebra.  Then for each $i$, 
for some $r_0,\ldots,r_d \in F$ (respectively:  $F \union \ldots \union \si^d(F)$)
the images $h_i(r_0),\ldots,h_i(r_d)$ are transformally independent over $K_i$.
(In the total dimension case, use \ref{trans0-r}.)
We may assume it is always the same sequence $r_0,\ldots,r_d$ (by refining the limit.)
Taking an ultraproduct we obtain $h_*: R \to D_* \subset L_*$, injective on $D$, with
the analogs of the above properties; in particular $h_*(r_0),\ldots,h_*(r_d)$ are 
 transformally (resp. algebraically) independent over $K_*$, contradicting the assumption on generic dimension. \qed

\begin{rem} \ \end{rem}
(1) \  One can find a difference domain $D'$,
$D \subset D'$,    $D'$ a finite integral extension
of $D$, and a dense open $U \subset \spe(D)$, such that
 for all $p \in U$,
 if $h: D  \to L$ is a map into a difference field with kernel $p$,
and $h$ extends to a difference ring homomorphism $h': D' \to L$,
 then $R \tensor_D L$ has transformal dimension precisely $d$
  over $L$. The proof of this fact, that we will not require, uses
  \ref{direct2} and \ref{direct-ritt-0} below.

\smallskip \noindent
(2) The base extension $D'$ in (1) cannot be avoided:
   
 for instance let $D = \Zz$, $R= D[X,\si(X),\ldots]/(X^2+\si(X)^2, \si^2(X)-X)$.
  Then $R$ has transformal dimension $1$ over the field of fractions of $D$.
  But if $p$ is any rational prime with $p=1 \hbox{\rm \ mod\ }  4$, and $h: R \to L$
 is a map into a difference field of characteristic $p$, then the fixed
 field of $(L,\si)$ has an element $i$ with $i^2=-1$, so $h(X)=\pm ih(\si X)$.
 Applying $\si$, $h(\si X) = \pm i h(\si^2 X)$ so $h\si^2(X) = - hX$,
 yielding $X=0$. It follows that $R \tensor_D (Z/pZ)$
 has transformal dimension $0$.

\<{lem}\lbl{wta}  Let $D$ be a difference domain, with
field of fractions of $D$
of transcendence degree $\geq k$ over a difference field $k$.
Let $R$ be a $D$-difference algebra, of generic total dimension $m$.  Then $R$ has total dimension $\geq k+m$
over $k$. \>{lem}
\proof There exists a surjective homomorphism   of transformal
$K$-algebras $h: R \tensor_D K \to D' $, $D'$ has field of fractions $K'$, $tr. deg._K K'  = m$.  So $tr. deg_k K' \geq k+m$,
and $h(R)$ generates $K'$ as a field.  \qed

 \begin{subsection}{Difference schemes and pro-algebraic varieties}
\lbl{pro-alg-sec}

Let $X$ be a  difference subscheme of
 an algebraic variety $V$ over a difference field 
$k$; more precisely, of the difference scheme
$[\si]_{k}V $.  Let  $ V_n= V \times \ldots \times \ldots \times
V^{\si^{n}} $.    (So   $V=V_0$.)

We associate with $X$   a sequence of   algebraic
subschemes $X[n]$ of $V_n$; $X[n]$ will be called 
the $n$'th -order {\em weak Zariski closure of $X$}.   We describe these locally
on $V$.  Let $U$ be an open affine subset of $V$.  So $U$ can be
identified with $\spec A$, with $A = {\cal O}_{V}(U)$, a $k$-algebra. 
The inclusion $X \to [\si]_{k }V$ gives a $k$-algebra homomorphism
$[\si]_k A \to {\cal X}(U)$, with kernel $I$ (so $I$ is a difference
ideal of $\spe [\si]_k A$, and $X \meet U \iso \spe ([\si]_k A) / I$.) 
Let $A_n$ be the sub $k$-algebra of $A$ generated by $A \union \ldots
\union \si^{n}(A)$; so $\spec A_n = U \times \ldots \times \ldots
\times U^{\si^{n}}$.  Finally define $(X  \meet U)_n = \spec A_n / ( I
\meet A_n)$.

We also let $X_\omega$ be the projective limit of the $X[n]$; it can be viewed
as a scheme, or as a pro-(scheme of finite type).  
In particular, $[\si]_kV_\omega = \Pi_{n\geq 0} V^{\si^n}$.  Note
that $[\si]_kV_\omega$ is isomorphic to the
   scheme $\Pi_{n\geq 1}  V^{\si^n}$ (by an isomorphism intertwining with
$\si: k \to \si(k)$).  At the same time we have the projection 
  $r: \Pi_{n\geq 0} V^{\si^n} \to  \Pi_{n\geq 1}  V^{\si^n}$.   

\<{remark}\lbl{pro-alg}  A subscheme $Y$ of $\Pi_{n\geq 0} V^{\si^n}$ has the form $X_\omega$
for some difference subscheme $X$ of $[\si]_kV$ iff $Y$ contains
${r Y}^{\si^{-1}}$ as a  scheme \>{remark}

\proof Locally, $Y$ is defined by an ideal $I$.  $I$ is a difference ideal
iff $I \subset \si^{-1}(I \meet R^{\sigma})$.  

\<{definition} $X$ is {\em weakly Zariski dense in $V$ } if
the 0-th-order weak Zariski closure of $X$ equals $V$.  \>{definition}

We will also say that $X$ is {\em Zariski dense} in $V$ if for some difference field
extension $K$ of $k$, $X(K)$ is Zariski dense in $V$.   (It is safer to use this only
for $k=k^{alg}$ and $X/k$ irreducible.)

It is quite possible that $X$ be weakly Zariski dense in $V$, while the set of
points of $X$ is not Zariski dense in $V$.  For instance, when $V=\Aa^{1}$,
the subscheme $X$ cut out by $\si(x)=0$ has this property.

If $X$ is algebraically integral, each $X[n]$ is an irreducible  algebraic
variety over $k$, of dimension $\leq n \dim V$.   The natural map $X[{n+1}] 
\to X[n]$ is dominant.

\end{subsection}
 \begin{subsection}{Transformal degree}

Let notation be as above:  $k$ is a difference field, $X[n]$ is the
$n$'th-order weak Zariski closure of $X$.

\begin{lemdef} \lbl{dimgrowth-int} Assume $X$ is an algebraically integral difference
subscheme of an algebraic variety $V$ over $k$.  There exist integers
$a,b$ such that for all sufficiently large $n \in \Nn$,
$$ \dim(X[n]) = a(n+1) +b$$

$a$ is the transformal dimension of $X$, while $b$ is called
the {\em  dimension growth degree} or {\em transformal degree}.  If $a=0$ then
$b$ is the total dimension of $X$.

\end{lemdef}

\proof %
 One proof uses intersections with generic hyperplanes (a Bertini argument.)  
We can assume $V$ is projective.  We define the notion of a generic
hyperplane section $X'$; it is the intersection of $X$ with a linear equation,
whose coefficients are generic in the transformal sense.  Show
that if $X$ has positive transformal dimension, then $\dim X'[n] = \dim X[n]-1$
for large $n$.  In this way reduce to the case of transformal dimension $0$, where
we must show that $\dim(X_n)$ is bounded (hence 
eventually constant), with eventual value equal to the total dimension.  This is clear
by Proposition \ref{trans-0} and Lemma \ref{dim-alg}.

Here is a more direct proof.  Let $(c_0,c_1,\ldots)$ be such that
$(c_0,\ldots,c_n)$ is a generic point over $k$ of $X[n]$, in some field
extension of $k$.  (This makes sense since $X$ is algebraically integral;
so $X[n]$ is an irreducible variety over $k$; and $X[n+1]$ projects
dominantly to $X[n]$.)  So $(c_0,\ldots,c_n) \to (c_l,\ldots,c_{n+l})$ is
a (Weil) specialization over $k$, and thus (with $k[l,n] = k(c_l,\ldots, c_{n+l}) $)

$a(l,n) = tr. deg._k k[l,n]$
is non-increasing with $l$.  So for each $n$, for some $\beta(n)$, for $l \geq \beta(n)$,
$a(l,n)=a(l+1,n)=\ldots =_{def} e(n)$.  Take the least possible value of  $\beta(n+1)$,
subject to $\beta(n+1) \geq \beta(n)$. 
 
Now $\epsilon(n) = e(n+1)-e(n)$ is non-increasing:  for large enough $l$,
$e(n+1)-e(n) = tr. deg. _{k[l,n]}  k[l,n+1] $;  so 
\[ e(n+2)-e(n+1) = tr. deg. _{k[l,n+1]}  k[l,n+2]  \leq  tr. deg. _{k[l+1,n]}  k[l+1,n+1] =  e(n+1)-e(n) \] 

Thus $\epsilon(n)$ is eventually constant, with value $a$.  So for all
$n \geq \gamma$, for all $l \geq \beta(n)$, $a(l,n)= \alpha(n) + an$
for some $\alpha(n)$.  

Note that $\beta(n)=\beta(n+1)$ unless $e(n) > e(n+1)$.  Thus $\beta(n)$ also stabilizes at
some maximal $\beta$.

It follows that for $l \geq \beta, n \geq \gamma$, $a(l,n) = \alpha + an$. 

Finally,  $tr. deg. _{k(c_{\beta},c_{{\beta}+1},\ldots)} k(c_0,c_1,\ldots) = \delta=
tr. deg. _ k(c_{\beta},c_{{\beta}+1},\ldots,c_L) k(c_0,\ldots,c_L)$, for some sufficiently 
large $L \geq \beta + \gamma $.  

So for $n \geq L$, $tr. deg. _k k(c_0,\ldots,c_n) = \alpha + \delta + an$.

For the identification of $a$ with the transformal dimension, we argue as follows.
 
 Let $R$ be an integral domain, and a finitely generated difference $k$-algebra.  If $B$
 is a finite set of generators  of $R$, let $R[n;B]$ be the $k$-algebra generated by $B \union \ldots \union \si^n(B)$, and let $a(n, R,B)=\dimk R$, the Krull dimension of $R[n;B]$.   We have seen that
 $\lim_{n \to \infty} (1/n) a(n,R,B)$ exists; denote this $a(R,B)$.    
  
 We will use two statements about the dimension of  finitely generated $k$-algebras $S,S'$ without zero-divisors, with fields of fractions $F,F'$:
 
(*)  If $S \leq S'$ then $\dimk(S) \leq \dimk(S')$.

(**) If in addition $F'$ is algebraic over $F$, then $\dimk(S) = \dimk(S')$.  

(In terms of varieties, if $V' \to V$ is dominant then $\dim(V) \leq \dim(V')$, and 
if in addition $V' \to V$ is generically finite, then equality holds.)
 
 Now:
 
1)  If $R \subseteq R'$, and $B \subseteq B'$, then $R[n;B] \subseteq R'[n;B']$. 
By (*),  $a(n,R,B) \leq a(n,R',B')$ and in particular $a(R,B) \leq a(R',B')$.

2)  If $B \subseteq B' \subseteq R$ ($B'$ finite), then $B' \subseteq R[k;B]$
for some $k$, $\si^l(B') \subseteq R[k+l,B]$, 
and so $R[n;B] \subseteq R[n;B'] \subseteq R[n+k; B]$.  Thus 
$a(n,R,B) \leq a(n,R,B') \leq a(n+k,R,B)$.  So 
$a(R,B)=a(R,B')$, and $a(R,B)$ does not depend on the choice of finite generating set; we can write $a(R)$ for this.

3)  Let $D$ be a difference domain, and a finitely generated difference $k$-algebra.  Let $R$
be a finitely generated difference $D$-algebra, not embedding $D[X]_\si$ (the difference polynomial 
ring in one variable over $D$.)  Then $a(D)=a(R)$.  Indeed, $a(D) \leq a(R)$ by (1).  For
the other inequality let $F$ be the field of fractions of $D$; it is a difference field.  
Let  
$B_0$ be a finite set of generators for $R$ as a difference $D$-algebra.  For $b \in B_0$,
for some $n=n(b)$, $\si^n(b) \in F(b,\ldots,\si^{n-1}(b))^{alg}$ (otherwise
the homomorphism $D[X]_\si \to R$, $X \mapsto b$, is an embedding.)  Let
$B_1 = \{\si^i(b): b \in B_0, i < n(b)\}$.  Then $B_1$ is a finite set of generators for $R$
as $D$-algebra, and $\si(B_1) \subseteq F(B_1)^{alg}$.  Let $B'$ be a finite
subset of $D$ contains the coefficients of the
relevant polynomials, so that  if $F'$ is the field of fractions of the $k$-algebra generated by $B'$,
then $\si(B_1) \subseteq F'(B_1)^{alg})$.   Finally let $B = B_1 \union B'$.  Then the field of
fractions of $R[n,B]$ is algebraic over the field of fractions of $D[n,B'][B_1]$.
(For $n=1$ by definition; for $n+1$, since $\si^{n+1}(B) = \si^{n+1}(B_1) \union \si^{n+1}(B')$,
and $\si^{n+1}(B_1) \subseteq \si^n(F'(B_1)^{alg})$, and induction.)  Thus by (**), 
$a(n,R,B) = \dimk D[n,B'][B_1] \leq a(n,D,B') + e$ where $e$ is the number of points in $B_1$.
Since $e$ is constant, $e/n \to_{n \to \infty} 0$, and  $a(R) = a(R,B) = a(D,B') = a(D)$.  

Now assume $R$ is well-mixed, and use Lemma \ref{dim-alg}; 
Applying (3) to $D= k[X_1,\ldots,X_t]_\si$ with $t$ the transformal dimension of $R$, we
see that $a(R)=s$.   
  \qed

\<{example}\lbl{growth-exmpl} The eventual dimension growth formula of \ref{dimgrowth-int} need not
hold for all $n$;
 the differences $\dim(X[n+1])-\dim(X[n])$ need not be monotone. \>{example}
Let $D=\Qq[z]_\si$, $R= D[x]$ with $\si(x)=0$, and let $y=xz^\si$.  Let $R[0]=\Qq[x,y,z] \leq R$,
$V= \Aa^3 = \spec R[0]$, $X = \spe R$ embedded into $V$ dually to $R[0] \to R$.  
Then $X[n] = \spec  \Qq[x,y,z,z^\si,\ldots,z^{\si^n}]$; letting $a_n = \dim X[n]$ we have $(a_0,a_1,a_2,a_3,\ldots) = (3,3,4,5,\ldots)$.    Taking the product of two copies of this
with one (such as $\si^3(t)=t$)  of growth $(1,2,3,3,\ldots)$, we obtain the dimension growth sequence  
$(7,8,11,13,15,\ldots)$ with differences $(1,3,2,2,\ldots)$.

\begin{rem}\lbl{dimdeg}    Assume $X \subset Y$ are two algebraically integral
difference subschemes of the algebraic variety $V$.  If 
 $X,Y$ have
the same transformal dimension and dimension growth degree, then $X=Y$. \end{rem}     

\proof  
In fact, for every sufficiently large $n$, $X[n] \subset Y[n]$,
and $\dim(X[n]) = \dim(Y[n])$; so $X[n]=Y[n]$; and so $X=Y$. \qed

 \begin{cor}\lbl{rm1} Let $k$ be a difference field, $D$ an
algebraically integral finitely generated $k$-difference algebra. 
Then for some $n $, $\si^n(D)$ is a difference domain.

In fact if $R$ is a finitely generated $k$-subalgebra of $D$, generating
$D$ as a difference ring, and $b$ is the dimension growth degree of $\spe (D)$
as a difference subscheme of $\spec (R)$, then
one can take $n \leq b$.
\end{cor}
  
\proof     Let $a$
be the  transformal dimension of   $\spe D$ over $k$; this equals
the transformal dimension of $\spe \si^n(D)$ over $\si^n(k)$, 
for any $n$. Let $b_n$ be the dimension growth degree of $\spe \si^n(D)$
as a difference subscheme of $\spec \si^n(R)$, a variety over
$\si^n(k)$.   In other words, let $R_{n,m}$ be the subring of $D$
generated by $\si^n(R) \union \ldots \union \si^{n+m-1}(R)$, $K_{n,m}$
the field of fractions of $R_{n,m}$, $c(n,m)$ the transcendence
degree of $K_{n,m}$ over $\si^n(k)$.  Then for any $n$, for large
enough $m$, $c(n,m) = am+b_n$.   We have $b=b_0 \geq b_1 \geq \ldots$,
so $b_i=b_{i+1}$ for some $i \leq b$.  The lemma now follows by
applying the following claim to $\si^i(D)$.

\claim    If $b=b_1$,  then $D$ is a difference domain.

\proof  Pick any $c \in D, c \neq 0$, in order to show that $\si(c) \neq 0$.
Let $m$ be large enough so that $c \in R_{0,m}$, that 
$c(0,m)=am+b$ and $c(1,m)=am+b_1$.  As $b=b_1$, 
$c(0,m)=c(1,m)$, i.e. $tr. deg._k K_{0,m} = tr. deg. _{\si(k) }K_{1,m}$.
Thus by Krull's theorem, the Krull dimensions of $R_{0,m}$ and
of $R_{1,m}$ are equal.   Now the surjective homomorphism
$\si: R_{0,m} \to R_{1,m}$ has a prime ideal $p$ for kernel (as $R_{1,m}$
is an integral domain.)   .
Thus $R_{0,m}/p \cong R_{1,m}$ has the same Krull dimension as
$R_{0,m}$.  As $0$ is also a a prime ideal of $R_{0,m}$, this forces $0=p$.
So $c \notin p$, and thus $\si(c) \neq 0$. \qed

\end{subsection} 

\begin{subsection}{Finite presentation}

We will need a certain uniform version of Hilbert's basis theorem;  for lack of reference
we give the proof.    Let $\cS$ be a  family of subrings of some commutative ring $S$ (with $S \in  {\cS}$).
 Call  ${\cS}$ {\em uniformly Noetherian}  if for any ideal $I$ of $S$,
for some finite $F \subseteq S$,   for any   $R \in \cS$ containing $F$,
$F \meet I$ generates $I \meet R$ as an ideal of $R$.  

\<{lem} \lbl{hilbert}  Let $K$ be a field, or $K=\Zz$, $n \in \Nn$, and let $\cS(K,n)$ be the family of 
polynomial rings $R[X_1,\ldots,X_n]$, with $R$ a subring of $K$.  Then $\cS(K,n)$
is uniformly Noetherian.  \>{lem}

\<{proof}
Note that the lemma is true for $n=0$ ( for a  field $K$, $\cS(K,0)$ is the family of subrings of $K$; one can  always take $F=\{0,1\}$.)  Thus it suffices to show that 
if $\cS$ is a uniformly Noetherian family of subrings of a ring $S$, then 
$\cS[X] := \{R[X]: R \in \cS\}$ is a uniformly Noetherian family of subrings of $S[X]$.  
The usual
proof of the Hilbert basis theorem shows this:   Let $I$ be an ideal of $S[X]$, and let 
$J_k = \{a \in S: (\exists f)(\deg(f)< k, \ aX^k + f \in I)\}$, an ideal of $S$. 
  We have
$J_0 \subseteq J_1 \subseteq \ldots J_m = J_{m+1} = \ldots$ for some $m$.
Let $F_k$ be a uniform set of generators for $J_k$, i.e. 
$F_k \meet J_k$ generates $J_k \meet R$ for any  $R \in \cS$ containing $F_k$.
   For each $a \in F_k$,
find $f=f_a^k \in I$ of the form $aX^k + g$, $\deg(g) < k$.   
Let $F$ be the union of all $f_a$, with   $a \in F_k$ , for $k \leq m$.  So $F$ is finite.   Let $R \in \cS$ with $F \subseteq R[X]$.
 Then all coefficients of all $f \in F$ are in $R$; in particular each $F_k \subseteq R$.  Hence $F_k \meet J_k$ 
 generates $J_k \meet R$ for $k \leq m$.   Let $I'$ be the ideal of $R[X]$
 generated by $F \meet I$.  To show that $I' = I \meet R[X]$,    let $h \in I \meet R[X]$, $\deg(h)=l$; use induction on $l$.  Write $h=aX^{l} +g$, $\deg(g)<l$.  Then $a \in J_l$,
 so $a \in J_k$ with $k= {\min(l,m)}$.   By assumption, $a=\sum r_i a_i $ for some $r_i \in R$
 and $a_i \in F_k$.  So $h' = h - X^{l-k} (\sum r_i f_{a_i}^k)$ has degree $<l$.  By induction
 $h' \in I'$; so $h \in I'$.  \eprf

\<{cor}  \lbl{hilbert-c} Let $D$ be an integral domain, $R$ a finitely generated $D$-algebra.
Then for some $a \in D, a \neq 0$, $R[a \inv]$ is a finitely presented $D[a \inv]$-algebra.

In fact, the presentation is uniform over   all sufficiently large localizations of subrings of $D$,
in the following sense.  
Let  $K$ be the field of fractions of $D$, and let
$h: D[X] \to R$ be a surjective homomorphism of $D$-algebras, 
$D[X]$ a polynomial ring; extend naturally to $h_K: K[X] \to R_K =  R \tensor_D K$.
Then for some  finitely generated $R_0 \subseteq D$ and $0 \neq a \in R_0$, for 
some finite $F_0 \subseteq R_0[a \inv][X]$, for 
any subring $R'$
of $K$ with $R_0[a \inv] \subseteq R'$, the restriction of $h_K$ to $R'[X]$ has kernel generated
by $F_0$.  
\>{cor}

\<{proof}   
Let $F$ be a uniform set of 
generators for $\ker h_K$ in the sense of Lemma \ref{hilbert};
 so that $F \meet \ker h_K$ generates $R'[X] \meet \ker h_K$ 
whenever $F \subseteq R' \subseteq K$.     Find $R_0,a$ such that $F \subseteq R_0[a \inv]$;
let $F_0 = F \meet \ker h_K \subseteq R_0[a \inv]$.  
If $R'$ contains $R_0[a \inv]$, let $h' = h_K | R'$.  Then $h'$ vanishes on $F_0$.
Conversely, if $f \in \ker h'$
then $f \in \ker h_K \meet R'[X]$, so    $f \in R'[X]F_0$.   So $\ker h'$ is generated by $F_0$.
 \eprf 

%
%

\<{cor}\lbl{hilbert-c-c}  Let $D,R$ be as in \ref{hilbert-c}, and let $I$ be an ideal of $R$.  Then  for some $d \in D, d \neq 0$, and some finitely generated $I_0 \subseteq I$, 
   for any $r \in I$, for some $N \in \Nn$, $rd^N \in I_0$.  \>{cor}
    
 \proof 
Applying \ref{hilbert-c} in $R/I$, we   find $d,I_0$ such that $I R[d \inv] = I_0 R[d \inv]$.  
But $I_0 R[d \inv] = \{a / d^n: a \in I_0\}$.  So for any $a \in I$, for some $b \in I_0$,
 $a=b/d^n$; $ad^n -b =0$; $ad^{n+m} -bd^m =0$; $ad^N \in I_0$.  \qed

\<{cor} \lbl{noether-1}  Let  $k$ be a difference field, or $\Zz$.  $D$ be a difference domain, finitely generated as a difference $k$-algebra.  Then there exists $0 \neq d \in D$ such that $D[d \inv]_\si$
is finitely presented as a difference $k$-algebra.  \>{cor}

\proof  Let $K$ be the field of fractions of $D$, $K ^{inv}$ the inversive hull of $K$.  Let $a_1,\ldots,a_r$ be     generators of $D$, $a = (a_1,\ldots,a_r)$.   In this proof, 
for $Y \subseteq K^{inv}$, write $k[Y]$ for 
the $k$-algebra of $K ^{inv}$ generated by $Y$, $k[Y]_\si = k[Y,\si(Y),\ldots]$.
Let $D^- = k[\{\si^n(a): n < 0 \}]$.  By Corollary \ref{hilbert-c},  there exists $0 \neq e \in D^-$ 
and a finite $f \subset D^-[e \inv][X]$ 
such
that the natural map induces an isomorphism

 (*) $D^-[e \inv, a] \cong D^-[e \inv][X]/(f)$,
 
  and indeed 
  $R'[Y \inv][a] \cong R'[Y \inv][X]/(f)$ whenever $R \subseteq R' \subseteq D$ and $e \in Y \subseteq R' \setminus (0)$.  
  
The element $e$ as well as all coefficients of the polynomials in $f$   lie in 
$ k[\si \inv (a),\ldots, \si ^{-l}(a)]$ for some $l$.  Let $d=\si^l(e)$.  For $n \geq 0$, let 
$R_{-n} = k[a^{\si^{-1}},\ldots, a^{\si ^{-l-n}},e \inv, \ldots, \si^{-n}(e \inv)]$.  By (*), 
   $$R_{-n}[a] \cong R_{-n}[X]/ (f)$$ 

Note that $\si (e) \in R_{-l}[a]$; say $\si (e) = h(a)$ with $h \in R_{-l}[T]$.
  Let $c = (\si  (e))^{-1}$, and let $g =  (f,Yh(X)-1)$.    Then

$$(**)  \ \ \ \ \ \ \ \ \ \      R_{-n}[a][c] \cong R_{-n}[X,Y] / (g)$$

 Let $$R_n = \si^{l+n}(R_{-n})   
 =k[a,\ldots,a^{\si^{l+n-1}},d \inv, \ldots,\si^n(d) \inv] $$
Applying $\si^{n+l}$ to (**), we see that $R_{n+1}   \cong R_n[X,Y]/(\si^{n+l}(g))$ naturally.
Let $H:k[X,Y]_\si \to D[d \inv]_\si$ be the  difference ring homomorphism with $X \mapsto a, Y \mapsto d \inv$.  Then the kernel of $H$ is generated by the kernel of $H | k[X,Y,\ldots,X^{\si^l},Y^{\si^l}]$ together with $\si^l(g),\si (\si^l(g)), \ldots$.  As a difference ideal, this is
finitely generated.    \qed

 \def\bR{{\bar{R}}}
 \<{cor}\lbl{hilbert-c-1}  Let $R$ be a finitely generated   difference $k$-algebra, with $k$ a difference field or $\Zz$.   Let $P$ be a transformally prime ideal of $R$.   Then there exists
 $a \notin P$ such that with $\bR = R[a \inv]_\si$, $P \bR $ is finitely generated as a difference ideal
in $\bR$.  \>{cor}
 
\proof  This  follows from \ref{noether-1}, applied to $R/P$.  \qed
 
\begin{definition}\lbl{radicial}  Let $f:R \to S$ be a morphism 
of difference rings.  $S$ is
{\em transformally radicial} over $R$ if for any $s \in S$, for some
$m \geq 0$, $s^{\si^m} \in f(R)$. \>{definition}

\begin{lem} \lbl{rad-fg}   If $S$ is a finitely generated
$R$-difference algebra and is transformally radicial over $R$, then $S$
is a finitely generated $R$-algebra. 

 If in addition  $R$ is a difference domain, then for some $0 \neq a \in R$,
$S[a \inv]_\si$ is a finitely presented $R[a \inv]_\si$-algebra, hence
a  finitely presented $R[a \inv]_\si$- difference algebra.  \end{lem}  

\proof  Let $H$ be a finite set of generators of $R$ as a difference $k$-algebra.
For $a \in H$, $\si^m(a) \in k \cdot 1$ for some $m$; let 
$T(a) = \{a,\si(a),\ldots,\si^{m-1}(a)\}$.   
Then $H' = \union_{a \in H} T(a)$ is a finite set of generators of $R$ over
  as a $k$- algebra.

Now assume  $R$ is a difference domain.
By Corollary \ref{hilbert-c}, for some $0 \neq a \in R$,  $S[a \inv]_\si$ is a finitely presented $R[a \inv]_\si$-algebra.

Now in general, if $S''$ is a difference $R''$-algebra, and a finitely presented 
$R''$-algebra, then it is a finitely presented $R''$-difference algebra:  if $S''=R''[b]$, $b=(b_1,\ldots,b_n)$, 
then $\si(b) = f(b)$ for some $f \in R[X]$, and it suffices to add to the $R''$-algebra
presentation the relation $\si(b) = f(b)$.     \qed

 \<{lem} \lbl{noether-2}      Let $R$ be a finitely generated algebraically reduced, well-mixed difference $k$-algebra, 
with $k$ a difference field or $\Zz$.   Assume $R$ has a unique minimal transformally prime ideal $P$.  Then there exists $a \notin P$ such that $\bR = R[a \inv]_\si$ is   a finitely presented difference $k$-algebra.   In fact for some $n$, 
 $D=\si^n(\bR)$ is a difference domain, and $\bR_\si$   is   a finitely presented difference $D$-algebra.
  
\>{lem}

\<{proof}   

By Corollary \ref{hilbert-c-1}, for some $a \notin P$, $P R[a \inv]_\si$ is a finitely generated difference
ideal.  Passing from $R$ to $ R[a \inv]_\si$, we may assume $P$ is finitely generated.   
By Lemma
\ref{well-mixed1}  $\srad{0}$ is a perfect ideal by , hence an intersection of transformally prime ideals; by assumption on $P$, 
$\srad{0} = P$.  Since $P$ is finitely generated, and $\srad(0) =  \union_n \si^{-n} (0)$, we have
 $\si ^{-n} (0) = P$ for some $n$.   Thus  if $a \in P$ then $\si^n(a)=0$.  

\claim   $D=\si^n(R)$ is a difference domain.

Indeed if $a,b \in \si^n(R)$ and $ab=0$, write $a = \si^n(c), b=\si^n(d)$; then $\si^n(cd)=0$,
so $c \in P$ or $d \in P$; thus $a =0 $ or $b=0$.  Similarly if $\si(a)=0$ then $a=0$.  Thus 
$\si^n(R)$ is a difference domain.

 Now $R$ is a finitely generated radicial $D$-difference algebra.  
  By Lemma \ref{noether-1},
 for some $0 \neq a \in D$, $D[a \inv]_\si$ is
 finitely presented as a difference $k$-algebra. 
  By  Lemma \ref{rad-fg}, for some $0 \neq b \in D$ (with $a | b$),
 $R[b \inv]_\si$ is a finitely presented difference  $D[b \inv]_\si$-algebra.  
Thus $R[a \inv]_\si$ is finitely presented as a difference $k$-algebra.   
  \eprf

%

Lemma \ref{noether-2} applies in particular to  algebraically integral 
finitely generated  $k$-algebras; since in this case $\srad{0}$ is the unique
minimal transformally prime ideal.  Hence \ref{hilbert-c-1} extends to:

\<{cor} \lbl{hilbert-c-2}  
 Let $R$ be a finitely generated   difference $k$-algebra, with $k$ a difference field or $\Zz$.   Let $P$ be an algebraically prime difference ideal of $R$.   Then there exists
 $a \notin \srad{P}$ such that with $\bR = R[a \inv]_\si$, $P \bR $ is finitely generated as a difference ideal
in $\bR$.  \qed \>{cor}

\end{subsection}

\begin{subsection}{A chain condition}

We present two variants 
of the Ritt-Raudenbush finite basis  theorem for 
perfect difference ideals; one (\ref{prime-noether})  for algebraically integral well-mixed ideals,
the other (\ref{rad-noether-0}) for algebraically radical ones in finite total dimension.    Yves
Laszlo found a mistake in an earlier proof of \ref{rad-noether-0}, stated without 
the finite total dimension assumption; I do not know if the statement is valid without it.

 \begin{lem}\lbl{prime-noether} There is no strictly descending infinite chain
of algebraically integral difference subschemes of an algebraic variety
$V$ \end{lem}

\proof  In such a descending chain, the transformal dimension
must be non-increasing, and eventually stabilize; after that point,
the  dimension growth degree cannot increase, and
eventually stabilizes too.   By \ref{dimdeg}, if  
$X \subset Y \subset V$ are algebraically integral, and have
the same transformal dimension and dimension growth degree, then $X=Y$.  \qed

Let $\rwm(R)$ be the set of  an algebraically radical, well-mixed ideals of a difference ring $R$.   When  $X \subseteq R$, let    $<X>= <X>_R$ 
be denote the smallest $I \in  \rwm(R)$ containing $X$.   $I \in \rwm(R)$ is {\em finitely generated}   if $I = <X>$ for some finite $X$.
Let $k$ be a difference field, or $\Zz$, and let $R$ be any algebraically reduced, well-mixed $k$- algebra.   For the purposes of the present proof, we call 
   $R$   {\em finitely presented}
    if $R \iso k[X]_\si / I$ for some finitely
generated $I \in \rwm(k[X]_\si)$, where $k[X]_\si$ is the 
difference polynomial algebra over $k$ in finitely many variables.  
  $\Nsi$ denotes the ring of expressions $\sum_{i=0}^m m_i \si^i$, with $m_i \in \Zz, m_i \geq 0$.
 
\<{lem}  \lbl{rad-noether-0}  Let $R$ be a finitely generated well-mixed $k$-algebra.  Assume $R$ has  finite total dimension, or the weaker assumption:

($\alpha$)  For any well-mixed quotient $\bar{R}$ of $R$, there is a bound $m$
on the transformal degree of ${\bar{R}}/p$, as $p$ ranges over the minimal algebraically prime difference ideals of $\bar{R}$.
 
Then   every element of  $\rwm(R)$ is finitely generated. \>{lem}

The first part of the proof (through Claim 2) just recovers the Ritt-Raudenbush basis theorem, and does not 
make use of finite total dimension.  We may factor out $\sqrt{0}$, so assume $R$ is algebraically reduced.

If the union of a increasing chain of $\rwm(R)$ is finitely generated, then cofinally many
elements of the chain are already finitely generated.     Thus if not every element of $\rwm(R)$ is finitely generated, by Zorn's lemma,
there exists a maximal non-finitely-generated  $I \in \rwm(R)$.  We will reach a contradiction from this assumption.  Note that  every element of $\rwm(R/I)$ is finitely generated
(either it is $0$, or it lifts to a finitely generated element of $\rwm(R)$.)  Thus every increasing
chain of elements of $\rwm(R/I)$ is eventually constant.

\<{claim}{1}    $I$ is algebraically prime.  Hence (\ref{well-mixed1}) $P= \srad {I}$ is transformally prime.   

\proof    If $ab \in I$, let $I_1
= Ann(a/I)$; then $I_1$ is well-mixed and radical.  Let $I_2 =
Ann(I_1/I)$; then $I_2$ is well-mixed and radical, and $I^2 \subseteq
I_1 I_2 \subseteq I$.  If $I\neq I_1$ and $I \neq I_2$ then $I_1,I_2$
are finitely generated as well-mixed radical ideals, say by $X_1,X_2$. 
Any algebraically prime difference ideal containing $X_1X_2$ must
contain $I_1$ or $I_2$; hence also $I$.  Thus any radical well-mixed
difference ideal containing $X_1X_2$ must contain $I$.  So $I$ is
finitely generated, a contradiction.  Thus $I=I_1$ or $I=I_2$; so $b
\in I$ or $a \in I$.  Thus $I$ is prime.  \qed \end{claim} 

Since $\si^{-n}(I)$ stabilizes,  $\si ^{-n} (I) = P$ for some $n$. 

By Lemma \ref{noether-1}, there exists $d \in R \setminus P$ and a finite $X_0 \subseteq R$ 
such that  $I R[d \inv]_\si  \in \rwm(R[d \inv]_\si)$ is generated by $X_0$.  But
$<X_0>_{R[d \inv]_\si} = \{a/d^n: a \in <X_0>, n \in \Nsi \}$.
  It follows that 
for any $a \in I$, for some $n \in \Nsi$, $d^na \in <X_0>$.  So
$d \in \srad{Ann(a/<X_0>)}$.  

But
by the maximality of $I$,   $<I,d>$ is finitely generated; so  
$<I,d> = <X_1 \union \{d\}>$ for some finite $X_1 \subseteq I$.   

 Replacing 
$R$ by $R/ <X_0 \union X_1>$,  we may assume: 

(*) For any $a \in I$, $d \in  \srad{Ann(a)}$  

(**) $<I,d> = <d>$. 
 
By (*) and (**), for any $a \in I$, $<I,d> \subseteq  \srad{Ann(a)}$; in particular
$a \in \srad{Ann(a)}$; so $a^{\si^n} \in Ann(a)$ for some $n$, i.e. $a^{\si^n}a=0$, 
  as $R$ is well-mixed and reduced,  $a^{\si^m}=0$ for some $m$.  

We thus have:

\claim{2}   For any $a \in I$, for some $n \in \Nn$,  $a^{\si^n}=0$.   \qed

\claim{3}  Say  $R/I$ has  transformal dimension $n$.  Then there  exists a finitely generated  
$I_0 \in \rwm{R}$,  $I_0 \subseteq I$ such that $R/I_0$ has transformal dimension $n$.    
  
\proof   Let $a_1,\ldots, a_N$ be   generators of $R$; for each subset $s$ of $\{X_1,\ldots,X_N \}$ of size $n+1$, there exists a nonzero $f_s \in k[X_1,\ldots,X_N]$ using the variables
in $s$ only, with $f_s(a_1,\ldots,a_N) \in I$.  let $I_0 = <\{f_s(a_1,\ldots,a_N): s\}>$.
Then in $R/I_0$ the image of $a_1,\ldots, a_N$  is a set of generators, with no transformally independent subset of size $n+1$; so the transformal dimension is at most $n$, hence $n$. \qed

Again we may replace $R$ by $R/I_0$, so we may assume $R$ has transformal dimension $n$.   
Let $Q$ be the set of cofinally  minimal prime ideal of $R$.  
We now use the assumption ($\alpha$);  so for each $p \in Q$,  $R/p$ has  transformal degree $\leq m$.   
By Lemma \ref{rm1},   $\si^{-m}(p)$ is transformally
prime.  By Claim 2,  $P \subseteq \si^{-m}(p) $.

Now using \ref{wm1.11},
 $\si^{-m}{0} = \si^{-m}( \meet_{p \in Q}  p )= \meet_{p \in Q} \si^{-m}(p) \supseteq P$.
 
So $P = \si^{-m}(0)$.  Let $D = \si^m(R)$. $D \iso R/\ker(\si^m) = R/P$, so $D$
is a difference domain. 

By Lemmas \ref{rad-fg} and \ref{hilbert-c-c}, for some finitely generated $I_1 \in \rwm{R}$, and some $d \notin P$,
for any $a \in I$ we have $ad^n \in I_1$ for some $n \in \Nn$. 
But $I_1 = \sqrt{I_1}$, so $ad^n \in I_1$ implies $(ad)^n \in I_1$ and $ad \in I_1$.  Thus, 
 replacing $R$ by $R/I_1$, 
we have:  

\claim{3}  For any $a \in I$,  $ad = 0$.

By Claim 3, for any $a \in I$, $d  \in Ann(a)$.  On the other hand by (**), $a \in <d>$ .  So  $a \in Ann(a)$.  $ \sqrt Ann(a)$ is a radical well-mixed ideal containing $d$, so $a \in <d> \subseteq \sqrt Ann(a)$, so $a^n a = 0 $ for some $n \geq 1$;  hence $a=0$.   So $I=0$.  This contradicts the assumption that $I$ is not finitely generated. \qed

 \begin{cor}\lbl{direct-ritt-0}   Let $X$ be an affine or projective
difference scheme over $\spe k$, $k=\Zz$ or a difference field,
of finite total dimension.
Then $X$ may be embedded as a closed subscheme of a directly presented (over $k$)
difference scheme $\tilde{X}$.  Moreover, $X, \tilde{X}$ have the same
topological space and the same underlying algebraically reduced
difference scheme $X_{red} = \tilde{X}_{red}$.  
\end{cor}

\proof (Affine case.)   We may write $X = \spe ([\si]_D D[X]/ J)$, where
$D[X]$ is a polynomial ring over $D$ in finitely many variables, and
$J$ a prime ideal of $[\si]_D D[X]$.  By Proposition \ref{rad-noether-0},
$J$ is
finitely generated as a well-mixed, algebraically reduced ideal.  Adding
variables, we may assume $J$ is so generated by $J_0 = J \meet
D[X,\si(X)]$.  Let $\tilde{X} = \spe ([\si]_D D[X]/ J_0)$.  \qed.

\begin{cor}   \lbl{rad-prime-0}
Let $k$ be a  difference field,
$R$ a finitely generated difference $k$-algebra of finite total dimension.   If
$p$ is a well-mixed ideal and is radical, then there exists a finite
number of algebraically prime difference ideals whose intersection is
$p$.    \end{cor}
 
\proof Using Lemma \ref{rad-noether-0}, we can prove the statement
by Noetherian induction on well-mixed, radical ideals.  Let $p$ be
such an ideal.  If $p$ is algebraically prime, we are done; otherwise,
say $ab \in p$, $a,b \notin p$; let $p_1 = Ann(a/p)$, $p_2 = Ann(p_1/p)$; then 
$p_1,p_2$ are  well-mixed,
radical ideals bigger than $p$, with $p_1 p_2 \subseteq p \subseteq p_{1} \meet
p_{2}$.  Let $q = p_{1} \meet p_{2}$; then $q^{2 } \subseteq p_1 p_2
\subseteq p$, so as $p$ is radical, $q \subseteq p$; and thus $p =
p_{1} \meet p_{2}$.  Now by induction, each $p_{i}$ is the
intersection of finitely many algebraically prime difference ideals;
hence so is $p$.  \qed

The corollary can be restated as follows:  $S=R/p$ has a finite number of minimal
algebraically prime difference ideals; their intersection equals $0$ in this ring.
(Indeed if $p_1,\ldots,p_m$ are algebraically prime
difference ideal of the difference ring $S$, no one contained in another, and
with $ p_1 \meet \ldots \meet p_m = (0)$ then each
is minimal.)   

Let $P_i  = \{a: \si^m(a) \in p_i , \hbox{ some } m \}$.
By Corollary \ref{rm1} applied to $R/p_i$,  for some fixed $m$, $P_i = \{a: \si^m(a)  \in p_i \}$.
Thus $P_i$ is a   transformally prime ideal.  Every transformally prime ideal contains
some $p_i$ and hence some $P_i$.

\<{rem} \lbl{t-emb} \ \>{rem} Every minimal transformally prime
ideal must therefore equal some $P_i$; but
  not every $P_i$ must be minimal.  (Consider
$\Qq(x,y)_\si / (xx^\si,xy^\si)$.)  
 
The prime ideals of $R[a^{-1} ]_\si$ correspond to primes $p$ of $R$ with $a^{-\si^n} \notin p$ for all $n$.
This correspondence preserves inclusion, and also preserves the set of difference ideals. The minimal 
algebraically prime difference ideals of $R[a^{-1}]_\si$ are thus precisely the proper $p_i' = p_i R[a^{-1}]_\si$,
and $p_i' \meet R = p_i$.  Thus for any irreducible difference scheme $X$ of finite total dimension,
one obtains canonically a finite set of algebraically prime difference ideal sheaves;
the corresponding closed subschemes are called the {\em algebraic components} of $X$.  

\begin{definition} \lbl{wct} Let $X$ be a  difference scheme.  An {\em algebraic component} of $X$ is a maximal algebraically integral 
subscheme $Y$ of $X$.
\end{definition}

From the discussion above, we obtain:  

\<{lem} \lbl{power}  Let $X$ be a difference scheme of finite type, and of finite total dimension.
Then every algebraically integral subscheme of $X$ is contained in some algebraic component.
If $C$ is an algebraic component of $X$, it has a unique maximal 
 transformally integral difference subscheme 
 $\bar{C}$.  For some $n$,  for any algebraically integral difference ring $D$, if $a \in C(D)$
 then $\si^n(a) \in \bar{C}(D)$.   
\qed \>{lem}

(In the language of the next section, \ref{redseq}, we can say that $B_n(C)$ is a component of $B_n(X)$.)




\end{subsection} 

\begin{subsection}{Directly presented  difference schemes}
\lbl{dirpres}

\begin{definition} \ \end{definition} Let $D$ be a difference
(well-mixed) ring.  Let $F =  D [X_1,\ldots,X_n]_\si$ be the
difference polynomial ring over $D$ in  variables $X_1,\ldots,X_n$.
  Let $I$ be a difference (well-mixed)
ideal of $F$.  $I$ is {\em directly generated} if $I$ is generated as
a difference (resp.  well-mixed) ideal by $I \meet F_1$.  A surjective
difference homomorphism $h: F \to R$ is said to be a {\em direct
presentation of $R$ as a difference (resp.  well-mixed) $D$-algebra} via
$R_1 = h(F_1)$ if $Ker(h)$ is directly generated.

\smallskip

In other words, $R$ is generated by $R_1$, and the relations
between the generators can all be deduced from relations between
$R_1$ and $\si(R_1)$.  Every finitely presented difference $D$-algebra
admits a direct presentation.  We will say that the {\em direct
total dimension} of $R$ (for this presentation) is the Krull dimension
over $D$ of the ring $R_2$ generated by $R_1 \union \si(R_1)$, and
that $R$ is {\em directly reduced (irreducible, absolutely
irreducible)} if $f=0$ whenever
$f \in R_1$ and $f^n \si(f)^m=0, \ n,m >0$ (respectively, if $R_2$ has no zero-divisors,
$R_2 \tensor_D L$ has no zero-divisors for some algebraically closed
field $L$ containing $D$.)

A similar definition can be made for difference schemes.  Let $V$ be an
algebraic scheme over a difference ring $D$, $S$ a subscheme of
$Y = (V \times V^{\si})$.  We let $\Sigma$ denote the graph of 
$\si$ on $[\si]_D V \times_{\spe D} [\si]_D V^\si$; more precisely, 
in any affine open neighborhood where $V = \spe R$,
$R$ a difference $D$-algebra, we let $I_\si$ be the   ideal
in $R \tensor_D R$ generated by the elements $\si(r) \tensor 1 - 1 \tensor r$;
$\Sigma$ is the corresponding closed difference  subscheme.  We also let
$S \star \Sigma$ be the projection to $[\si]_D V $ of the difference scheme
$[\si]_D S \meet \Sigma$.  It is isomorphic to $[\si]_D S \meet \Sigma$.

A direct presentation of $X$ is an embedding of $X$ into $[\si]_D V$,
$V$ an algebraic scheme of finite type over $D$, such that the image
of $X$ has the form $S \star \Sigma$ as above.

\

 
\begin{definition} (The limit degree, cf. \cite{cohn}).   Let $Z$ be an irreducible
difference variety over a difference field $K$, of transformal
dimension $0$.  Let $a$ be a generic point of $Z$ over $K$.  Then $K_n
= K(a,\si(a),\ldots,\si^n(a))$ is a field; for large $n$, $K_{n+1}$ is
a finite algebraic extension of $K_n$; and the degree $[K_{n+1}:K_n]$
is non-increasing with $n$.  Let $\deg_{lim}(Z)$ be the eventual value
of this degree.   \end{definition}

\smallskip

\begin{lem}\lbl{direct2}   Let $V$ be a smooth algebraic variety over a
 difference field $K$, $\dim(V) = d $.
   Let $S \subset V \times V^\si$ be an absolutely irreducible  subvariety,
and assume
either the projection $S \to V$ or $S \to V^\si$ is smooth, of relative dimension $e$.   Let
$Z =  S \star \Sigma \subset [\si]_K V$ be the
difference subscheme of
$[\si]_KV$ cut out by: $(x,x^\si) \in S$.
  Then 
\begin{enumerate}
  \item $Z$ is a disjoint union  of components.
  \item  Each component $W$ is actually a component of the pro-algebraic
variety $$S[\infty] = S \times_{V^{\si}} S^\si \times_{{V^{\si^2}}} \times \ldots$$
 in the strong sense that $W[n]$ is a component of 
$$S[n] = S \times_{V^{\si}} S^\si \times_{V^{\si^2}} \times \ldots S^{\si^{n-1}}$$
  \item  Each component $W$ is Zariski dense in $V$,  and has    $\dim(W[n]) = d+ne$ for all $n \geq 1$.  
  \item  If $e=0$,   
\[ \deg_{cor} S = \sum_{W} \deg_{lim} W \] 
where  the sum is taken
over all components $W$ of $Z$.
\item $Z$ is perfectly reduced.
\end{enumerate} 
\end{lem}

\proof    We take the case: $S \to V$ is smooth (the other being similar.)

  Since
$S^{\si^n} \to {V^{\si^n}}$ is smooth, and smooth morphisms are preserved under base change, $S[n+1] \to S[n]$ is also smooth.  Compositions of smooth morphisms are also smooth,
so for each $n$, $S[n] \to V$ is a smooth morphism, and indeed $S[n]$ is smooth over $K$.

It follows that each $S[n]$ is a disjoint union of irreducible components.  By the dimension theorem (in the smooth variety $V \times   V^\si \times \ldots \times V^{\si^n}$), 
each component has dimension $\geq d+ne$.  On the other hand 
by adding up relative dimensions we see that $\dim S[n] \leq d+ne$, so equality holds.

Since the components of the $S[n]$ are disjoint, letting $S[\infty]$ be the projective
limit of the $S[n]$, we see that the   pro-algebraic variety $S[\infty]$ is a disjoint union
of irreducible pro-algebraic varieties $C$; $C$ is an inverse limit of components
$C[n]$ of $S[n]$.  Now by definition, $C[n+1]$ projects to $C[n]$; but
one can also project to $V^\si \times \ldots \times V^{\si^n}$; the image is a component
which is either contained in, or disjoint from, $C[n]^\si$.  In the former case, 
for any difference subvariety $W$ of $Z$, $W[n+1] \meet C[n+1] = \emptyset$;
call such a pro-component {\em empty of difference points}.    If the latter case holds for
each $n$, then $C$ defines  a perfect difference scheme, and so a component $W(C)$ of $Z$.    Thus $Z$ is the disjoint union of the various (nonempty) $W(C)$.  
 It follows
that $Z$ is perfectly reduced.

Finally,   assume $e=0$.  The additivity of degree
may be proved by induction on $\deg_{cor} (S)$.
If for all $m$, the $S[m+1]$ is absolutely irreducible,
 then $Z$ is transformally integral, and
$\deg_{lim}(Z) = \deg_{cor}(S)$.  Otherwise, consider the minimal $m$
such that $S[m+1]$ is reducible.  Let $Y = S[m]$, and let
 \[T =  ((a_0,\ldots,a_m),(a_1,\ldots,a_{m+1})) \in (Y \times Y^\si):
(a_0,\ldots,a_{m+1}) \in S[m+1] \} \] 
Let $\{T_j\}$ be the components of
$T$ projecting dominantly to $Y$.  They map finite-to-one to $Y^\si$,
hence dominantly.  Then \[ \sum_j \deg_{cor} (T_j) = \deg_{cor}(T) \]
Since $\deg_{cor} (T_j) < \deg_{lim} (Z)$ for each $j$,
 \[\deg_{cor}(T_j) = \sum \deg_{lim}(Z') \] where the $Z'$ here range over
the components of $Z$ , Zariski dense in $X$, such that for $a \in Z'$
\[ (a,\si(a),\ldots,\si^{m+1}(a) ) \in T_j \] Since each Zariski dense
$Z'$ falls into a unique $T_j$ in this sense, the required equation
follows.

\qed


Let $f[1](S) = \{a:  \dim(f^{-1}(a) \meet S) \geq 1 \}$.  
Note also, in a converse direction to \ref{direct2}:

\<{lem} \lbl{zd} Let $S \subset V \times V^\si$ be absolutely
irreducible $k$-varieties of dimension $d$, $Z = [\si]_k pr_0[1](S)$.
Let $X$ be a difference subscheme of  $ S \star \Sigma$; $[X]$ is the set of points of $X$. 

(1)  If $[X] \meet Z = \emptyset$ or $\si [X] \meet  pr_1[1](S)= \emptyset$, then $X$ has total dimension $\leq d$. 

  Assume $X$ is algebraically integral.  
  
 (2)  If $X$ is weakly Zariski dense in $V$, then $X$ is a component of $S \star \Sigma$,
and $X$ has total dimension $d$.   

(3)    If $t.dim(X)=d$, and $[X] \not \subseteq [Z]$,  then $[X]$  
 is  Zariski dense in $V$.

\>{lem} 

 \prf  (1) If $b \in X(D)$, $D=k[b]_\si$  an algebraically integral difference $k$-
algebra, with fraction field $K$, then $(b,b^\si) \in S$, $b \notin pr_0[1](S) $,
 so $b^\si \in  k(b)^{alg}$.  Inductively, $k[b]_\si \subseteq k(b)^{alg}$.  Thus $tr.deg._k K = tr.deg._k k(b) \leq \dim(V) = d$.  This shows that $X$ has total dimension $\leq d$.  
 
 If $b^\si \notin pr_1[1](S)$, then $b \in k(b^si)^alg$.  But $b^\si$ is a specalization of $b$.
 So $tr.deg._k k(b) \geq tr. deg._k k(b^\si)$.  It follows that $b^\si \in k(b)^{alg}$, and we continue as before.

For (2,3) assume  $X$ is algebraically integral, and let $(a_0,a_1,\ldots)$ be a generic point of the pro-algebraic variety corresponding to  $X$, i.e.  $(a_0,\ldots,a_n)$ is a generic
point of  $X[n]$ for
each $n$.

(2)  Assume$X$ is weakly Zariski dense
in $V$, i.e. $X[0]=V$.  
 Then $X[1] \to X[0]$ is dominant, so $d= \dim (S) \geq \dim X[1] \geq X[0] = d$.
So $(a_0,a_1)$ is a generic point
of $S$, and thus $a_0 \in k(a_1)^{alg}$.  So  the specialization
$$(a_0,\ldots,a_n) \to (a_1,\ldots,a_{n+1})$$
 does not lower transcendence degree, hence it
 is an isomorphism.    Thus
$k(a_0,\ldots)$ is a difference field, and hence the prime ideal corresponding
to $X$ is
a transformally prime ideal, and $X$ is a component.  Also, 
since $k(a_0),k(a_1)$ have the same transcendence degree and $a_0 \in k(a_1)^{alg}$,
we have $a_1 \in k(a_0)^{alg}$, and since the  displayed specialization is an isomorphism,
$a_{n+1} \in k(a_n)^{alg}$; so $tr.deg._k k(a_0,a_1,\ldots) = d$.  As this is a generic point,  the
algebraically integral difference scheme $X$ has total dimension $d$.   

(3) The point  $a_n$ cannot be in the closed set $\si^n(pr_0[1](S))$: if this were the case for some
$n$, it would be true for all larger values of $n$.  But by \ref{power}, 
 for sufficiently large $n$, $a_n \in \si^n([X])$; and $\si^n([X]) \meet
\si^n(pr_0[1](S)) = \emptyset$.)  
Now $(a_n,a_{n+1}) \in S^{\si^n}$, so 
 $a_{n+1} \in k(a_n)^{alg}$ for each $n$.  Thus   $tr.
deg._k k(a_0) =   d $.  In other words, (2) holds; and (3) follows.       \eprf

\<{lem}\lbl{dim-const}  Let $Y$ be a   difference scheme 
with no strictly decreasing chain of  perfectly reduced subschemes, 
and let $X$ be a difference scheme of finite type over $Y$ (via $f: X \to Y$.)  Let $m$
be an integer.  Then:

\begin{enumerate}
  \item  there exists a  perfectly reduced difference subscheme $D_m(f)=D_m(X/Y)$ of $Y$, such that $X_a$ has total dimension
$\geq m$ if $a$ is a generic point of (any component of) $D_m(X/Y)$, and has total dimension $< m$ (or is empty) if $a \notin D_m(X/Y)$.    
  \item The above properties characterize $D_m(X/Y)$ uniquely.
  \item  If $X$ has finite
total dimension $d$, then $D_m(X/Y)$ has total dimension $\leq d-m$.
\item If equality holds in (3), then $f^{-1}D_m(X/Y)$ 
contains an algebraic component of $X$. (cf. \ref{wct}).
\end{enumerate}     \>{lem}

\proof  The uniqueness (2) is clear, since if $Z,Z'$ are two candidates, each generic point of $Z$ must lie on $Z'$, and vice versa.  To show existence (1),
and (3), we use  Noetherian induction on $Y$.  We may assume $Y$ is perfectly reduced.  We may also assume
$Y$ is irreducible, since if it has a number of components $Y_j$, then we can let $D_m(X/Y) = \union_j D_m(f^{-1}Y_j/Y_j)$.
If $X/Y$ has generic total dimension $\geq m$, let $D_m(X/Y) = Y$.  By \ref{wta}, $Y$ has total dimension
$\leq d-m$ in this case.  Otherwise, $Y$ has total dimension $\leq m-1$.  By \ref{5.11}, there exists a proper closed difference
subscheme $Y'$ of $Y$ such that for $a \in Y \setminus Y'$, $X_a$ has total dimension $\leq m-1$.  Now $D_m(f^{-1}Y'/Y')$
exists by Noetherian induction, and we can let $D_m(X/Y)=D_m(f^{-1}Y'/Y')$. 

Finally, (4) is clear since $f^{-1}D_m(X/Y)$ is a difference subscheme of $X$
of the same total dimension $d$ (by the first part of (1).)
 \qed

{\bf Caution:}  This is not preserved under base change $Y' \to Y$.
$D_m(X'/Y')= \emptyset$, $D_m(X/Y)=Y$ is possible.

From a logical point of view, ACFA does not eliminate quantifiers.  
$D_0(X/Y)$ for instance is therefore not the projection,
but rather a difference-scheme theoretic closure of the projection.

\end{subsection}

 \end{subsection}

\end{section} 

\begin{section}{Transformal multiplicity}
\lbl{trans-mult}

Recall Definition \ref{radicial}:  if $f:R \to S$ is a morphism of  difference rings,  $S$ is
 transformally radicial  over $R$ if for any $s \in S$, for some
$m \geq 0$, $s^{\si^m} \in f(R)$.  We will extend this to difference schemes, and 
analyze transformally radicial maps $X \to Y$ by dividing $Y$ according to the
relative total dimension of $X$.

\begin{remark}\lbl{red-rad} Let $k$ be an inversive difference
field, let $R \subset S$ be difference
$k$-algebras, and assume
$S$ is  transformally radicial over $R$.   Then $R$,$S$ have
the same reduced total dimension over $k$.  \end{remark}
\proof Let $I = \{r \in R: (\exists m \geq 1) \si^m(r)=0\}$,
$J = \{s \in S: (\exists m \geq 1) \si^m(s)=0\}$.  Any difference ring
homomorphism on $R$ into a difference domain must factor 
through $R/I$; thus $R$,$R/I$ have the same reduced total
dimension, and similarly so do $S,S/J$.  But $I = J \meet R$,
and $R/I \iso_k S/J$.  \qed

\<{definition}
 A difference scheme $X$ over
a difference scheme $Y$ is {\em transformally radicial} if
there exists an open affine covering $\{U_j\}$ of $Y$ and $\{V_{ij}\}$
of $X$, such that ${\cal O}_X(V_{ij})$ is transformally radicial
over ${\cal O}_Y(U_j)$  \end{definition}

If a difference
scheme $X$ of finite type over $Y$ is transformally radicial over $Y$, then $X$
has finite total dimension over $Y$; this follows from the first statement of  Lemma \ref{rad-fg}.

We begin by a look at the generic behavior; for
this purpose we may pass to the fields of fractions.  

\begin{subsection}{Transformally separable extensions}
 
We define an invariant $\iota$
dual to the limit degree, and a purely inseparable variant $\iota'$.  
The invariant  $\iota'$  
  of \ref{towers} below   figure in the asymptotics in $q$ of
the multiplicity of the points of $M_q(X)$, the reduction of $X$
to a difference scheme with structure endomorphism the $q$-Frobenius.
 Most likely these lemmas appear in \cite{cohn}.

\<{lem} \lbl{inv-fin} Let $F$ be an inversive difference field, $K$ a difference field extension of $F$
with $tr. deg._F (K) < \infty$.  Then $K^a$ is inversive. \>{lem}

\proof $F \subset \si(K) \subset K$, and $tr. deg._F (\si(K)) = tr.deg._F (K)$. \qed

\begin{lem} \lbl{towers} Let $K \subset L \subset M$ be difference fields, with
$M$ of finite transcendence degree over $K$.
Consider  $L,K,M$  as subfields of $M^{inv}$.

  \begin{enumerate}
\item   $L^{inv}$ is an algebraic extension of $K^{inv}L$.
Thus if $L$ is finitely generated over $K$ as a difference field, then
$LK^{inv}$ is a finite extension of $\si(L)K^{inv}$.

Write $\iota(L/K)=[LK^{inv}: \si(L)K^{inv}]$,
$\iota'(L/K)=[LK^{inv}: \si(L)K^{inv}]_{insep}$
\item  Assume $M/L$ is a regular field
extension.  Then $M/K$ is transformally separable
iff $M/L$ and $L/K$ are transformally separable.
\item  Assume $M$ is finitely generated and transformally separable over $K$.  Then
 $\iota(M/K) = \iota(M/L) \iota(L/K)$, and similarly for $\iota'$.

 \end{enumerate}
\end{lem}
\proof

By Lemma \ref{inv-fin}, we have:

\claim  Let $K$ be an inversive difference field, $ L$ a difference field extension of $K$
of finite transcendence degree.  Then $L^{inv}$ is an algebraic extension of $L$.

(1) follows:
since the   transcendence degree of $L$ over $K$ is finite, so is that of $LK^{inv}$
over $K^{inv}$, and the claim applies.

\smallskip

(2)  The "if" direction follows from the transitivity properties
of linear disjointness. Assume $M/K$ is transformally separable.  Then
$L/K$ is a fortiori transformally separable.  We need to prove that $M/L$
is transformally separable.
 $M$ is linearly disjoint from the compositum $K^{inv}L$ over $L$.
We must show that $K^{inv}M$ is linearly disjoint from $L^{inv}$ over $K^{inv}L$.
Since $M/L$ is a regular field extension , and
$M$ is linearly disjoint from $LK^{inv}$ over $L$,
$MK^{inv}/LK^{inv}$ is a regular extension.   In other words $MK^{inv}$ is linearly disjoint
over $LK^{inv}$ from the algebraic closure of $LK^{inv}$.  By (1), $L^{inv} \subset (LK^{inv})^{alg}$.

\smallskip
(3) We may assume here that $K$ is inversive.  We have
$$ \iota(M/K) = [M:M^\si] = [M: M^\si L] [M^\si L : M^\si] $$
Yet $M$ is linearly disjoint from $L^{\si^{-1}}$ over $L$, so
$M^\si$ is linearly disjoint over $L^\si$ from $L$, and
$$[M^\si L : M^\si] = [L : L^\si] = \iota(L/K)$$
Similarly $$[M:M^\si L] = [M L^{inv} :M^\si L^{inv}] = \iota(M/L)$$
\qed 
 
\end{subsection}

 \begin{subsection}{The reduction sequence of a difference scheme}
\lbl{redseq}
Let $X$ be a difference  scheme.    We define a sequence of functors
$B_n$, and maps giving a sequence $X \to B_1X \to B_2 X \mapsto \ldots $.
Denote $B_n X$ by $X_n$ in this section.  We will also define functorially maps
$r_n: X \to X_n$ as well as $i_n: X_n \to X$.  The maps $i_n$ will allow us to identify
$X_n$ with a subscheme of $X$, but it is the sequence of maps $r_n$ that will really
interest us.

When $X = \spe R$ is affine, we let $X_n = \spe \si^n(R)$.   Let $r_n$ be induced
by the inclusion $\si^n(R) \subset R$.  Let $i_n$ be induced by the surjective homomorphism
$\si^n : R \to \si^n (R)$.

If $f: R \to S$ is a ring homomorphism,, we let
 $B_n(f) : \si^n(R) \to \si^n(S)$
be the restriction of $f$; and  with $f^*: \spe S \to \spe R$ the corresponding
map of difference schemes, $B_n(f^*) = B_n(f)^*$.

If $X$ is a multiplicatively and transformally closed subset of $R$, $0 \notin X$,
let $S = R[X^{-1}]$.
We compare $\si^n ( S)$ to the localization $R'[\si^n(X)^{-1}]$, where
$R' = \si^n(R)$.   The natural map $R'[\si^n(X)^{-1}] \to \si^n(S)$ is clearly surjective,
and since $R$ is assumed to be well-mixed, injective too.  Thus \[
B_n(j): B_n(\spe S) \to B_n( \spe R) \] is compatible with
localizations.  By gluing we obtain a functor $B_n$ on well-mixed
schemes.

Observe that $R[X^{-1}]$ and $R[\si^n(X)^{-1}]$ may not be the same ring, but they
have the same affinization, i.e. the same difference spectrum.  Indeed an element $a \in X$ may not be invertible in
$R[\si^n(X)^{-1}]$, but $\si^n(a)$ will be invertible, and therefore $a$ will be a $\si$-unit
in this ring.

The maps $ i_n: \spe R \to B_n ( \spe R)$ and $r_n: B_n (\spe R) \to \spe R$,
also globalize.

\begin{remark}  \ \end{remark} If $k$ is inversive, $B_n(k) = k$, and $B_n$ induces a
functor on difference
schemes over $k$.  The maps $r_n$ are maps of $k$-difference schemes, but  the
maps $i_n$ are not (so that the identification of $X_n$ with a subscheme of $X$ involves
a twisting vis a vis $k$.)

\medskip

 Let $X$ be a difference scheme of finite type (or, of finite type over
 an inversive difference field $K$.)  It follows from Lemma
 \ref{rad-fg} that each ring ${\cal O}_X(U)$ is finitely generated
 over ${\cal O}_{B_n(X)} (U)$.  In particular it has finite relative
 total dimension $\tau_n$.  If $X$ has finite total dimension, then all
 the $\tau _n$ are bounded by this dimension.

\begin{definition}  The {\em transformal multiplicity} of a difference
scheme $X$ is the supremum of
the total dimensions of the morphisms $X \to B_n(X)$.   

If $X$ is a difference scheme of finite type over a difference field
$k$, we define the {\em transformal multiplicity of $X$ over $k$} to
be the transformal multiplicity of the difference scheme $X \tensor_k
K$ , where $K$ is the inversive closure of $k$.

 If $X$ is a difference scheme of finite type over a difference scheme
 $Y$, define the {\em relative transformal multiplicity of $X$ over
 $Y$} to be the supremum of the transformal multiplicity of $X_y$ over
 $L$, where $L$ is a difference field and $y$ is an $L$-valued point of
 $Y$.
\end{definition} 

Observe that as a map of points, $\spe X \to \spe B_n(X)$ is bijective.  (Every
transformally prime ideal $p$ of $\si^n(R)$ extends uniquely to a transformally prime ideal
of $R$, namely $\si^{-n}(p)$.)  While $B_n(X)$ may not be Noetherian as a scheme, $X \to B_n(X)$ is
a map of finite type of schemes, and by Lemma \ref{total-krull-0}, the transformal 
multiplicity is  the maximum dimension of a fiber of this map (an algebraic scheme)
over a point $p \in \spe B_n(X)$.

\begin{example} \ \end{example} Consider the   subscheme $X$ of $\Aa^1$ defined
over an algebraically closed difference field $K$ by 
$f(X^\si,\ldots,X^{\si^{m+1}}) = 0$, $f$ an irreducible polynomial.
If $K$ is inversive, then we may write $f = g^\si$, and the
corresponding reduced scheme will be given by
$g(X,\ldots,X^{si^m}) = 0$; the transformal multiplicity is $1$.
  In general $X$ may be transformally reduced; but
it becomes reducible after base change, and the transformal multiplicity is at all events
equal to $1$.

\smallskip

\begin{prop}\lbl{si-dim3.2}  Let $K$ be an inversive difference field,
 $X/K$ a 
difference scheme of finite type and of total
 dimension $d$.  Let $k \geq 1$.  There exist canonically defined
closed subschemes $Mlt_k
 X$ of $X$ such that:
\begin{enumerate}
\item  $X \setminus Mlt_kX$ has transformal multiplicity $< k$
 over $K$.
\item $Mlt_k X  $ has reduced total dimension
 $\leq d-k$ over $K$.
\item If $Mlt_k X$ has  reduced total dimension $  d-k$, then it contains
an algebraic component of $X$.     In this case, $X$ has  transformal multiplicity $\geq k$.
\end{enumerate}
    \end{prop}

\proof    

Let $Y= B_{k+1}X$, and
let $r=r_{k+1}: X \to Y$  be the reduction sequence map. 

Let $Y_k = D_k(r)$; cf. \ref{dim-const}.
    
Let  $Mlt_kX = r^{-1}Y_k$.  (Or the underlying reduced scheme).

(1)  By Lemma \ref{trans-mult-0} below, we have to show that 
 if $L'$ is an inversive difference ring and 
$y$ is an $L'$-valued point of
$B_{k+1}(X \setminus Mlt_kX) $, then $X_y$ has total dimension 
$< k$ over $L'$.  Now $y \notin Y_k$; by definition of $Mlt_k(Y)$, and \ref{dim-const},
 the total dimension of $X_y$ is $<k$.

(2),(3) come from     \ref{dim-const} (3), (4).  
  
\qed

\<{notation} \lbl{z-mlt-not}  $Z_0X = X \setminus Mlt_1 X$.  For $k \geq $, 
$Z_kX = Mlt_k X \setminus Mlt_{k+1}X$. \>{notation}

  \begin{rem}  \lbl{si-dim3.2R}  The assumption in \ref{si-dim3.2}
that $K$ is inversive is not necessary.
    \end{rem}

\proof  Let $X' = X \times_K K'$, and let $X'_k = Mlt_k(X')$,
satisfying the conclusion of Proposition \ref{si-dim3.2}.  Let
$j: X' \to X$ be the natural map of difference schemes.  Then
the radicial map 
$j$ induces a bijection between the points of $X'$ and of $X$,
or between the perfectly reduced subschemes of $X'$ and of $X$,
preserving reduced total  dimension.
Let
$Mlt_k X  $ be the perfectly reduced subscheme of $X$ 
corresponding to $Mlt_k X'$.  Then (1),(2),(3) are clear. 
(cf.   \ref{inv-base-change}).   \qed

The next lemma, \ref{trans-mult-n}, falls a little 
short of concluding, when $X$ itself has transformal multiplicity $\leq n$,
 that $B_n(X)$
 must have transformal multiplicity $0$.  
The obstacle  can be explained in terms of sheaves of difference
algebras;  $(B_n)_*({\cal O}_X)$ need not coincide with 
${\cal O}_{B_n(X)}$.

\begin{lem}\lbl{trans-mult-n} Let $X=\spe R$ be a well-mixed
difference scheme, and assume $X \to B_{n'}(X)$ has total dimension
$\leq n < n'$.  Then for any difference field $L$ and any $L$-valued
point $y_{n'}$ of $B_{n'}(X)$, and $y \in X(L)$ lifting $y_{n'}$,
$\si^n( {\cal O}_{X,y} \times_{B_{n'}(X),y_{n'}} L)$ has total
dimension $0$.  \end{lem}

\proof   
By \ref{inv-base-change},
we may take $L$ to be inversive. 

We may assume $X = \spe R$, $R$ a well-mixed difference ring.  Let
$R_n = \si^n(R)$.  Let $y$ be the unique extension of $y_{n'}$ to a
difference ring morphism $y: R \to L$; let $y_m = y | R_m$.  Let
$\tensor'$ denote the tensor product in the well-mixed category; $A
\tensor '_B C$ is the quotient of $A \tensor _B C$ by the smallest
well-mixed ideal of that ring.  Let $S = R \tensor ' _{R_{n'},y_{n'}}
L$, $z: S \to L$ the induced map, $S_n = \si^n(S)$, $z_n = z | S_n$.

We have to show that $S_n$ has total dimension $0$.  Now $S$ is a
well-mixed $L$-algebra, and $\si^{n'}(S) \subset L$.  By assumption,
$S$ has total dimension $\leq n$.  We are reduced to showing:

{\bf Claim}  Let $L$ be an inversive difference field,
 $S$ a finitely generated well-mixed difference $L$-algebra, with
 $\si^{n'}(S) \subset L$.  Assume $S$ has total dimension $ n < n'$. 
 Then $\si^n(S) $ has total dimension $0$.

   There is no harm factoring out the nil ideal of $S$ (as a $k$-algebra), as this
is a difference ideal, and the total dimension of $S_n$
will not be effected.  As $S$ is well-mixed, and Noetherian, $0$ is
the intersection of finitely many prime ideals $p_1,\ldots,p_r$, and
they are difference ideals.  We have $\meet_{i=1}^r (p_i \meet S_n) =
0$, so it suffices to show that $S_n/(p_i \meet S_n)$ has Krull
dimension $0$ for each $i$; for this we may work with $S/ p_i$.  Thus
we may assume $S$ is an integral domain.  At this point we will show
that $S=L$.  Otherwise let $a \in S \setminus L$.  As $L$ is
inversive, $\si^{n'}(a)=\si^{n'}(b)$ for some $b \in L$, so $a-b
\notin L$, and $\si^{n'}(a-b)=0$.  Thus it suffices to show that for
$n \leq m < n'$, if $a \in S$ and $\si^{m+1}(a)=0$ then $\si^m(a)=0$. 
$S$ is a f.g. domain of Krull dimension $\leq n$.  If $a \in S$, then
$a,\si(a),\ldots,\si^n(a)$ are algebraically dependent over $L$ in the
field of fractions of $S$; so there exists an $F \neq 0 \in
L[X,X_1,\ldots,X_n]$ with $F(a,\ldots,\si^n(a))=0$.  Take $0 \neq F
\in L[X,X_1,\ldots,X_m]$ of smallest number of monomials, and then
least degree, such that $F(a,\ldots,\si^m(a))=0$.  Let $i_0$ be least
such that $F$ has a monomial from $L[X,X_1,\ldots,X_{i_0}]$.  Apply
$\si^{m-i_0}$ to $F$.  The monomial involving $X_{i}$ for some $i >
i_0$ disappear when applied to $a$ (as $\si^{m+1}(a)=0$); deleting
them, we obtain a shorter polynomial (though in higher -indexed
variables ) vanishing on $(a,\ldots,\si^m(a))$.  This is impossible;
so no monomials disappear; i.e. all monomials are from
$L[X,X_1,\ldots,X_{i_0}]$; but none are from
$L[X,X_1,\ldots,X_{i_0-1}]$.  So $F$ has the form $X_{i_0}F'$; and
being irreducible, it is just a multiple of $X_{i_0}$.  So
$\si^{i_0}(a)=0$.  \qed

\begin{cor}\lbl{trans-mult-0} A difference scheme $X$ has transformal multiplicity $\leq n$ iff
  $X \to B_{n+1}X $ has total dimension $\leq n$. \end{cor}

\proof  It suffices to show that if $n<n'$ and
$X \to B_{n'}(X)$ has total dimension $\leq n$, then so does 
$X \to B_{n'+1}(X)$.   Let $k$ be a difference field, $y$ an $k$-valued point
of $B_{n'+1}$.  In terms of local rings, we have 
$R \supset \si^{n'}(R) \supset \si^{n'+1}(R)$, and $z: \si^{n'+1}(R) \to k$.
We have to show that $S = R \tensor _{\si^{n'+1}(R),z}  k$ has total dimension $\leq n$.
Let $(h,D,L)$ be as in the definition of total dimension.  
Let $S_m$ denote the image of $\si^m(R) \tensor _{\si^{n'+1}(R),z} k$ in $S$.

Observe in general that if $ X \to_f Y$ has relative total dimension $\leq n$, then 
so does the induced map $BX \to BY$.  Indeed, by Remark \ref{obviousmonotonicity},
the pullback $f^{-1}BY \to BY$ has total dimension $\leq n$, and $BX$ is a 
subscheme of $f^{-1}BY$, so \ref{obviousmonotonicity} applies again.  In 
particular, as  $X \to B_{n+1}X$ has total dimension $\leq n$, so does 
$B_{n'-n}X \to B_{n'+1}X$ .

By  lemma \ref{trans-mult-n},  $\si^n(\si^{n'-n}R \tensor _{\si^{n'+1}(R),z}  k)$ 
has total dimension $0$.  Thus $h(S_{n'})$ is finite dimensional over
$k$, so it is a difference subfield of the domain $D$, call it $k'$.   So
$h(\si^{n'}(R)) \subset k'$.  As $X \to B_{n'}(X)$ has total dimension $\leq n$,
$h(R)$ is contained in a field $F$ of $k'$-transcendence degree $\leq n$.
But $[k':k]< \infty$, so  $tr. deg._{k}(F) = tr. deg_{k'} (F) \leq n$.  \qed

\begin{paragraph} {The reduction sequence: two side remarks}
 
\begin{lem} \lbl{rm2} If $X$ is algebraically integral,of finite type over
a field $K$, then
the   $X_n$ stabilize as subschemes of $X$:  for some $n$, $X_n = X_{n+1} = \ldots$.
\end{lem}
\proof This reduces to the affine case by taking an open affine cover; there it follows
from Corollary \ref{rm1}.  \qed

( What can be said without the algebraic integrality assumption?
 
 If $X$ is transformally integral, the schemes $X_n$ are all isomorphic to $X$.
The maps $r_n: X \to X_n$ need not be isomorphisms however.

   If $X$ is a
scheme over an inversive difference field $K$, then $X_n$ is not necessarily isomorphic
to $X$ over $K$, but is the transform of $X$ under $\si^n$.
In this case we can also define $X_n$ for negative values of $n$, obtaining a sequence
$ \ldots X_{-2} \to X_{-1} \to X \to X_1 \to \ldots $.)

\end{paragraph}
\end{subsection}

\<{subsection}{$\si$-primary ideals}
 \<{definition}  Let $R$ be a difference ring.  A well-mixed ideal $I$ is {\em $\si$-primary} if
 whenever $a,b \in R$ and $ab \in I$, either $a \in I$ or $b^{\si^n} \in I$ for some $n \in \Nn[\si]$. \>{definition}
 
   We have:
 if $ab \in I$, $b \notin P$ then $a \in I$.  If $\srad{I} = P$ and this condition holds, we say
 that $I$ is $P$-$\si$-primary.  Any finite intersection of $P$-$\si$-primary ideals is again $P$-$\si$-primary.
 
 If $P$ is a transformal prime, we say that a well-mixed ideal $I$ is $P$-$\si$-primary
 if $P =  \srad{I}$ and whenever $ab \in I$, $b \notin P$ we have $a \in I$.
 
\<{lem} \lbl{primary}  Let $I$ be a well-mixed ideal, and let $P = \srad{I}$, and 
$R_P$   the localization of $R$ outside $P$.  Let $\iota: R \to R_P$ be the natural map.
The following conditions are equivalent:
\begin{enumerate}
  \item $I$ is $\si$-primary
  \item $P$ is a transformal prime, and $I$ is   $P$-$\si$-primary
  \item  $P$ is a transformal prime, and  $I =\iota \inv (I')$  for some well-mixed ideal $I'$ of $R_P$.
\end{enumerate} \end{lem}

 \proof   Since $I$ is a well-mixed ideal, $\srad{I}$ is a perfect ideal (Lemma \ref{well-mixed1}.)
 If $I$ is $\si$-primary, then $P$ is a prime ideal, hence a transformal prime.  Thus
 (1) implies (2).  Assume (2).  Let $I' = \{a/b: a \in I, b \notin P\}$.  This is a well-mixed
 ideal of $R'$:  if $(a/b)(a'/b') \in I'$, then $(a/b)(a'/b') = a''/b''$ with $a'' \in I, b'' \notin P$.
 So $aa'b''=a''bb'$ in $R'$; upon multipliying denominators and numerators by an 
 appropriate element outside $P$, we may assume $aa'b''=a''bb' $ in $R$.  As $I$ is $P$-primary, and $a'' \in I$, we have $aa' \in I$; as $I$ is well-mixed,  $a^\si a' \in I$. 
 Now $(a/b)^\si (a'/b') = (a^\si a') / (b^\si b') \in I'$.  So $I'$ is well-mixed.  We have $I \subseteq \iota \inv (I')$, and if $a/1 \in I'$, then $a = a'/b$ for some $a' \in I$, $b \notin P$,
 so $(ab-a')b'=0$ for some $b' \notin P$, thus $ab-a' \in I$, so $ab \in I$, so $a \in I$.  This proves (3).  Now assume (3).  If $ab \in I$, and $b \notin P$, 
 then $\iota(b)$ is invertible; so $a/1 \in I'$, thus $a \in I$.    \qed 
 
 \<{lem}  Let $R$ be a finitely generated difference $k$-algebra, $k$ a difference field.
  Then there is no strictly increasing chain
 of    algebraically radical
 $\si$-primary well-mixed ideals.  \>{lem}

 \proof  Suppose $I_1 \subset I_2 \subset \ldots \subset P$ is such a chain.  Let $P_k = \srad{I_k}$.  Then $P_1 \subseteq P_2 \subseteq \ldots$, so by Ritt-Raudenbush 
  the $P_k$
 are eventually equal, and we may assume $\srad{I_k}=P$.  Let $I = \union_k I_k$.
 If $b \in R \setminus P$, $R' =  R[b \inv]_\si$, $I_k' = I_k R[b \inv]_\si$, then the $I_k'$
 form a strictly ascending chain, by Lemma \ref{primary} (3).  Thus, replacing $R$ by 
 an appropriate $R'$ using Lemma \ref{noether-2}, we may assume $I = \union_k I_k$ is finitely  generated;
but it follows that $I=I_k$ for large enough $k$, a contradiction.   \qed

\>{subsection}

\end{section} 

\<{section}{Transformal valuation rings in transformal dimension one}  \lbl{td1}

\<{subsection}{Definitions}

 \<{notation}  When $K$  is a valued field, $\OK,\M_K, {\val}(K),K_{\res}$
will denote
 the valuation ring, maximal ideal;  
 value group, residue field  of $K$; but we will often denote 
$R=\OK, M = \M_K, \bar{K} = K_{\res},  \Gamma=  {\val}(K)$. \>{notation}

\<{definition} \lbl{tvr-def}
 \begin{enumerate}
 \item A transformal valuation ring  (domain) is a valuation ring $R$ that is also a difference ring,
  such that $\si(M) \subset M$ (and $\si$ injective.)

\item The valuation $v$
 will be said to be  $m$-increasing
 (resp. strictly increasing)
 if  for all $a \in R$ with $v(a)>0 $, $ v(\si(a)) \geq m \cdot v(a)$
  (resp. $v(\si(a))>v(a)$).  It is $\omega$-increasing if $m$-increasing for all $m$.

 \item A weakly
transformal valued field is an octuple $(K,R,M,\si,\Gamma,\bar{K},\val,\res)$ such that $(R,M,\si)$ is a transformal valuation
ring, $K$ is the field of fractions of $R$, $res:  K \to \bar{K}$ is the residue map, $val:  K^* \to \Gamma = G_m(K)/G_m(R)$
is the valuation map.  
\footnote{$G_m(R)$ is the group of invertible elements of $R$.}
We will also apply the term to parts of the data.  If $\si$ extends to an endomorphism of $K$, we
denote it too by $\si$, and we say that $K$ is a transformal valued field.

\item  Assume given a distinguished  $t \in R $,
or at least a distinguished $\tau \in \val(R)$ ($\tau = \val(t)$). Let 
$$\Delta =    \{v \in \val(K):  -\tau < nv < \tau, \ n=1,2,\ldots \} $$

Write $a << b$ if $ma <b$ for all $m \in \Nn$.
    $K$  has {\em   transformal ramification dimension $rk_{ram}(K) = r$} (over
$\tau$)
if $r$ is the length of a maximal chain $v_0,\ldots,v_r \in \val(R)$ with
$0 = v_0 << v_1 << \ldots << v_r << \tau $.  When $K$ is $\omega$-
increasing, 
$rk_{ram}(K) \leq \dim_{\Qq} (\Qq \tensor \Delta)$.  
(In   transformal dimension one, equality holds; cf \ref{vt1c}.)  When $K$ extends $k(t)_\si$,
taking $\val(t)$ distinguished, we define 
the {\em valuative rank} $\vrk(K/k(t)_\si)$ of $K$ over $k(t)_\si$ to be $rk_{ram}(K) + tr. deg._k \bar{K}$.  

\end{enumerate}

\end{definition}

Weakly transformal valuation rings will not be used in the proof of the main results.

\<{notation} \ \lbl{not-moving} \>{notation}
A $\Zz$-polynomial in one variable $F(X)$ is said to be {\it positive at $\infty$} if $F(t)>0$ for sufficiently large real
$t$.

Given a difference ring $k$, let $ k[t,t^{-1}]_\si$ be the transformal localization by $t$ of the transformal polynomial ring
$k[t]_\si$, and let $\kt'$ be the sub-difference ring $$\kt'= k[t^{F(\si)}:  F \in \Zz[X], F(\infty)>0 ] \leq
k[t,t^{-1}]_\si$$ Write $t_n= \si^n(t)$.  Note the homomorphism $\kt \to k$, with kernel generated by the $(t_n)$.

Alternative description:  Let $k$ be a difference field.  Let $k(t)_\si = k(t_0,t_1,\ldots)$, with $\si(t_n)=t_{n+1}$,
$t=t_0$.  Then $k(t)_\si$ admits unique a $k$-valuation, with $0 < val({t_i} ) << val(t_j) $ whenever $i<j \in \Nn$.  (Here
$\alpha << \beta$ means:  $m \alpha < \beta$ for all $m \in \Zz$.)  $\kt$ is the associated valuation ring, and $\Ak = \spe
\kt$.  Note that $\kt$ is the localization of $\kt'$ at the prime $t=0$.

\<{lem} \lbl{tvrl} Let $(K,R,M,\si,\Gamma,\bar{K},val)$ be a weakly transformal valued field.  \begin{enumerate}

\item $\si^{-1}(M) = M$.  \ $\si(M) = M \meet \si(R)$.  If $K$ is $1$-increasing,
then every ideal of $R$ is a well-mixed difference ideal. 

\item There exists a unique difference field structure on $\bar{K}$, such that the residue map $R \to \bar{K}$ is a morphism
 of difference rings.  

\item There exits a convex subgroup $\Gamma'$ of $\Gamma$ and a homomorphism $\si_{\Gamma}:  \Gamma' \to \Gamma$ of ordered
Abelian groups, such $val(a) \in \Gamma'$ and $\si_{\Gamma}(val (a)) = val(\si (a))$ when $a \in R$, $ \si(a) \neq 0$.

If the valuation is $k$-increasing, then $k v \leq \si_{\Gamma}(v)$ for $v \in \Gamma, v \geq 0$.  

  \item Let $L$ be a subfield of $K$.  Then $\si(L^a \meet R) \subset (\si(L \meet R))^a$.  (where $S^a$ denotes the algebraic
closure of the field of fractions of $S$.)  

\item If $L$ is a subfield of $K$, and $\si( L \meet R) \subset L^a$, then $L^a \meet R$ is a difference subring of $R$.  

\item Let $\Delta$ be a convex subgroup of $\Gamma$ with ${\si_\Gamma}^{-1}( \Delta) \subset \Delta$.  (An automatic condition
  for    $1$-increasing valuations.)  Let $\hat{v}$ be the valuation $ val(a) + \Delta$; with value group $\Gamma/\Delta$, residue field
  $\hat{K}$.  Then a weakly transformal valued field structure is induced on $\hat{K}$, with residue field $\bar{K}$, value
  group $\Delta$.  

\end{enumerate} \>{lem} \proof \begin{enumerate}

\item $\si(M) \subset M$ by assumption, while $\si^{-1}(M) \subset M$ since $M$ is the unique maximal ideal of $R$.  Thus
$\si^{-1}(M) = M$.  If $\si(r) \notin \si(M)$, then $r \notin M$, so $r$ is a unit of $R$, $rr'=1$, so $\si(r)\si(r')=1$, and
thus $\si(r) \notin M$.  

\item $\si$ induces an endomorphism $\si:  {\bar K} = R/M \to \si(R) / \si(M) = \si(R) / (M \meet \si(R)) \subset {\bar K}$ of
${\bar K}$.  

\item Let $R' = \{r \in R:  \si(r) \neq 0 \}$, $\Gamma' = \{ \pm val(a):  a \in R, a, \si(a) \neq 0 \}$.  If $a,b \in R$ and
$ab \in R'$ then $a \in R'$, so $\Gamma'$ is convex.  If $a,b \in R'$ and $val(a)=val(b)$, then $a=cb,b=da$ for some $b,d \in
R$, so $\si(a)=\si(c)\si(b),\si(b)=\si(d)\si(a)$, and thus $val \si(a) = val \si(b)$.  Define $\si_{\Gamma}( val(a) ) = val
(\si(a))$, and observe that $\si:  val(R') \to \Gamma$ is a homomorphism of ordered semi-groups.  Extend it to an ordered
group homomorphism $\Gamma' \to \Gamma$.

  \item Let $b \in L^a \meet R$.  There is a nonzero polynomial $P \in L[Y]$ such that $P(b)=0$.  Dividing by the coefficient
of lowest value, we can assume the coefficients of $P$ are in $R$, and at least one of them has value $0$.  By (1), it follows
that $P^\si \neq 0 \in \si((L \meet R))[X]$, and $P^\si(b^\si)=0$.  

\item From (4):  $\si(L^a \meet R) \subset (\si(L \meet R))^a \subset L^a$, so $\si(L^a \meet R) \subset (L^a \meet R)$.

\item Let $R_\Delta = \{a \in K:  (\exists \delta \in \Delta) (\delta \leq val(a) )\}$, $M_{>\Delta} = \{r \in K:  (\forall
\delta \in \Delta) \, \delta < val(r) \}$, $\pi:  R_\Delta \to \hat{K}$ the natural map.  Let $\hat{R} = \pi(R)$; it is a
valuation ring of $\hat{K}$, with maximal ideal $\hat{M}=\pi(M)$.

If $r \in M_{>\Delta}$, then $\si(r)=0$, or $val (\si(r)) = \si_{\Gamma}(val(r)) \notin \Delta$ (otherwise $val(r) \in
\Delta$.)  So $\si(M_{>\Delta}) \subset M_{>\Delta}$, and $\si$ induces an endomorphism $\si:  \hat{R} \to \hat{R}$.  $\pi$ is
a morphism of difference rings, and $\si (\hat{M}) = \si(\pi(M)) = \pi(\si(M)) \subset \pi(M) = \hat{M}$.

 We have $val(\si(r)) \geq val(r)$, implying in particular when $val(r) \in \Delta$ that $val(\si \hat{r}) \geq
val_{\hat{K}}(\hat{r})$, so that $(\hat{R},\si)$ is  $1$-increasing if $(R,\si)$ is.  \end{enumerate}

\>{subsection}

\<{subsection}{Transformal discrete valuation rings}

We will be interested in $\omega$-increasing transformal valued fields $L$, finitely generated and of transformal dimension one over trivially valued difference subfield $F$.   We will call these 
 {\em transformal discrete valuation rings}.     Picking any $t \in L$, $\val(t)>0$, we can view 
$L$ as an extension of $F(t)_\si$ of finite transcendence degree.

Let $\Z-si = \Zz[\si]$,  $\Q-si = \Qq[\si]$, viewed as   ordered $\Zz[\si]$-modules,
(with $\Qq < \Qq \si < \ldots$.  ).  $\Z-si$,$\Q-si$ are   the value groups of $F(t)_\si$,
$(F(t)_\si)^a$.

{\em  In this section, all transformal valued fields are assumed to be
$\omega$-increasing, with value group contained in $\Q-si$.      }

Let $L$ be a transformal valued field.   . 

\<{lem} \lbl{hensel-alg} Let $L^h$ be the Henselization of $L$ as a valued field.  The endomorphism $\si$ of $L$ 
lifts uniquely to a valued field endomorphism of  $L^h$.  If $L$ is $\omega$-increasing, so is $L^h$.   \>{lem}

\proof $\si: L \to \si(L)$ is an isomorphism of valued fields; by the universal property of the
Henselization, for any Henselian valued field $M$ containing $\si(L)$, $\si$ extends uniquely to 
an embedding $L^h \to M$; in particular, with $M=L^h$, $\si$ extends to $\si: L^h \to L^h$.  The property of being $\omega$-increasing 
depends on the value group, which does not change. \qed

  One can also  canonically define a  {\em transformal Henselization} of a transformal valued field, 
 using   difference polynomials and their derivatives; cf. the remarks
 following \ref{hensel-c}, as well as \ref{res-surj}. 
Since we will only work with the transformal analogue of discrete valuation rings, we will be able to use
the somewhat softer notion of topological closure.

\<{paragraph}{Completion and closure}

Let $L$ be a transformal  valued  field, with value group contained in $\Qq_\si$.
We define a topology on $L^h$, a basic open set being a ball of nonzero radius.   This topology
is in general incompatible with valued field extensions.   However, according to lemma \ref{vt1},
all the valued field extensions we will consider here have cofinal value groups; so the inclusions
of valued fields are continuous, and the induced topology on a subfield coincides with the intrinsic one.  

We can construct the completion $\hat{L}$ of $L^h$ for this topology; an element of $\hat{L}$ is
represented by a sequence $a_n$ of elements of $L^h$, with $\val(a_{n+1}-a_n) \geq \si^n(v)$
for some $v \in \val$, $v>0$.  This completion carries a natural $\omega$-increasing transformal valued field structure.  

  If $K \leq L$ are transformal discrete valuation rings over $F$,  
the topological closure of $K^h$ within $\hat{L}$ can be identified with $\widehat{K}$.

The Henselization and completion processes do  not change the value group or residue field.  By Lemmas  \ref{vt1}, \ref{vt2},
the residue field is an extension of $F$ of finite  transcendence degree, and the value group of the algebraic clsoure is $\Qq[\si]$.  
 We will obtain a theory closely analogous to the theory of discrete valuation rings in algebraic geometry.  
 
\>{paragraph}

\>{subsection}
\<{subsection}{The value group}

\<{lem} \lbl{vt1}   Let  $L$ be 
an extension of $F(t)_\si$ of finite transcendence degree. 
Then 
the value group $\val (L^a)$  of $L^a$ is isomorphic to $\Q-si$, as an ordered
  $\Zz[\si]$-module.  Let $K \leq L$ be a difference subfield, nontrivially valued. Then
 $\val(K)$ is cofinal in $\val(L)$.  
   $\val(L^a)/\val(K^a)$ is a principal, torsion
$\Qq[\si]$-module.   \>{lem}

\proof   Pick $t \in K$, $\val(t)>0$.  Since $tr. deg. _{F(t)_\si} L < \infty$, $\val(L^a)/ \val({F(t)_\si}^a)$ is a finite-dimensional $\Qq$-space.  Thus $\val(L^a)$ is a   finitely generated   $\Qq[\si]$-module.   Since $\Qq[\si]$ is a principal ideal domain,
and $\val(L^a)$ is torsion free,  
$\val(L^a)$ is a free  $\Qq[\si]$-module.  Now the quotient $\val(L^a)/ \val({F(t)_\si}^a)$ is finite-dimensional,
so the rank of $\val(L^a)$ must equal one.  Let $h:  \Q-si \to \val (L^a) $
be an isomorphism of  $\Zz[\si]$-modules.    Replacing $h$ by $-h$ if necessary,
 we may assume $h(1)>0$.  It follows that $h$ is order preserving.   
Any nonzero $\Zz[\si]$ - submodule
of $\Q-si$  is cofinal, and co-torsion.  
 \qed

In particular, since $\val(L) \leq \Q-si$ as an ordered Abelian group,
there are no "irrational" values in $\val(L)$; for any $0 < u < v \in \val(L)$,  
for some (unique)  $0 \leq q \in \Qq$, $|qv-u| << v$.   This follows from the fact that the
 divisible ordered Abelian group $\Q-si$,   is a lexicographic direct sum of copies of $\Qq$.
 

 \<{cor}\lbl{vt1c}   Let  $L$ be 
a transformal valued field extension of $F(t)_\si$ of finite transcendence degree. 
Then $rk_{ram}(L / F(t)_\si) = rk_{\Qq} (\val(L)/ \val F(t)_\si)$.   
 \>{cor}
 
\<{proof}  This is in fact a statement about the ordered Abelian group $A=\Q-si$.  Let 
$0 < \tau \in A$, so that $B = \Qq[\si] \tau$ is a copy of $\Q-si$ embedded
in $A$; and let $0<<a_1 << \ldots << a_n << \tau$ with $n$ maximal possible. 
We have to show that $rk_{\Qq} (B/A)=n$.  Clearly  $a_1,\ldots,a_n$
are $\Qq$-linearly independent over $B$, so $n \leq  rk_{\Qq} (B/A)$.  Factoring out
the (convex) $\Qq$-space generated by $a_1,\ldots,a_n$, we may assume $n=0$.
In this case we have to show that $A=B$.  Indeed otherwise let $a \in A, a \notin B, 0<a$.
As in \lemref{vt1}, $B$ is cofinal in $A$.  Let $m$ be least such that for some 
 $a \notin B, 0<a < \si^m(\tau)$.  Then for some rational $q \geq 0$, $a' =a- q \si^m(\tau) <<
 \si^m(\tau)$; but $a' \notin B$, contradicting the minimality of $m$.  \eprf 

\>{subsection} 


  \<{subsection}{The residue field}

\<{lem}\lbl{vt2}  Let  $L$ be an $\omega$-increasing transformal valued field of 
finite transcendence degree over $F(t)_\si$.    Then the residue field $\res (L)$ of $L$ has finite transcendence degree over $F$. \>{lem}

\proof   The residue field of $F(t)_\si$ is $F$ itself, and  $L$ has  finite transcendence degree over
$ F(t)_\si$.  \qed

This raises the question, that I did not look into:   If in \ref{vt2} $L/F(t)_\si$   is 
finitely generated, is $\res L$   finitely generated
over $\res(K)$ as a difference field
 (up to purely inseparable extensions )?  It seems likely that this can be proved by analogy with \lemref{vt1fg}.  
\end{subsection} 

\<{subsection}{Valued field lemmas}

We will need an observation from the theory of valued fields.
(Compare \cite{HHM1}  Part I, \S 2.5 ( ''Independence and orthogonality for unary types".)   )

Let $K \leq L$ be an inclusion of valued fields, $c \in L$.  We let 
$$T(c/K) =  \{val(c-b): b \in K \}$$
Also let $E(c/K)$ be the stabilizer of $T=T(c/K)$ in $\val(K)$, i.e.
$$E = \{e \in \val K: e+T=T \}$$
If $T(c/K)$ has no greatest element, then
$\val K(c) = \val K$  (see Lemma \ref{val-lem}).  In this case $T(c/K) $ is a downwards-closed
subset of $\val K$.  Hence $E$ is a convex subgroup of $\val K$.   

\<{lem}\lbl{val-lem} Let 
 $K$ be an algebraically closed valued field, $L$ an extension of 
transcendence degree $1$.   Let $c \in L \setminus K$.  Then the following conditions are equivalent:
\begin{enumerate}
  \item  $\res(K) \neq \res(L)$ or  $\val(K) \neq \val(L)$
  \item  $\res(K(c)) \neq  \res(K)$ or  $\val(K(c)) \neq \val(K)$
  \item  $T(c/K) = \{val(c-b): b \in K \}$ has a maximal element. 
\end{enumerate} 
In particular, the third condition depends on $K,L$ alone and
not on the choice of $c$. 

$T(c/K)= \val K$ iff $L$ embeds into the completion $\widehat{K}$ as
a valued field.
  \>{lem}
\proof:    Let $\Oo_L$ be the valuation ring of $L$.
(1) implies (2) since if $\res(K(c))=\res(K)$,
and $\val (K(c)) = \val(K)$, 
then $\res(K(c))$ is algebraically closed and $\val(K(c))$ is divisible;
so they cannot change under algebraic field extensions.  Now assume 
$\val(K(c)) \neq \val(K)$.  Then $\val f(c) \notin \val(K)$ 
 for some $f \in K[X]$.    Splitting $f$ into linear factors, wee see that 
we can take it to be linear.  So $\val(c-b) \notin \val(K)$ for some $b \in K$.
It follows that $\val(c-b') \leq \val(c-b)$ for all $b' \in K$.  (Otherwise
$\val(c-b)=\val(b-b')$.)  Next suppose $\val(K(c)) = \val(K)$, but
 $\res f(c) \notin \res (K)$ for some $f \in K[X]$ with $f(c) \in \Oo_L$.  
Taking into account $\val (K(c)) = \val(K)$, we can split $f$ into linear
factors $f_i$ such that each $f_i(c) \in \Oo_L$.  Thus $ac-b \in \Oo_L$,
$\res(ac-b) \notin \res(K) $ for some $a,b \in K$.  It follows that
$\val(ac-b') \leq 0$ for all $b' \in K$.  So $\val(c-b/a)   \geq \val(c-b'')$
for all $b'' \in K$.  This proves that (2) implies (3).   

Now assume (3): $\val(c-b) \geq \val(c-b')$
for all $b' \in K$.  If $\val(c-b) \notin \val(K)$, (2) holds.  If $\val(c-b) = \val(d)$,
$d \in K$, then $(c-b)/d \in \Oo_L$, and $\val((c-b)/d-d') \leq 0$ for all $d'$;
so $\res((c-b)/d ) \notin \res(K)$.  Thus (1).  \qed

\<{example} \lbl{vt4-e}  Let $F$ be a   field, $K=F(t_0,t_1,\ldots)^a$, valued over $F$ with $0<\val(t_0)<<\val(t_1)<< \cdots$.  
Let $K' = F(t_e,t_{e+1},\ldots)^a$, and let $c \in \widehat{K}  \setminus \widehat{K'}$.  Then 
$T(c/K') = \{ \val(c-b): b \in K' \}$ has a maximal element.  \>{example}

\proof  If $T(c/K')$ is unbounded in $\val(K')$, then
$c \in \widehat{K'}$.  Otherwise, for all $b \in K'$, 
$\val(c-b) \leq \alpha$.  Find $c' \in K$ with $\val(c-c') > \alpha$.
Then $T(c'/K')  = T(c/K')$.  Now $K^a/(K')^a$ is generated (as an algebraically
closed field) by $e$ elements,
whose values are $\Qq$-linearly independent over $\val(K')$.  Thus
  $\val K'(c') \neq \val K'$; so  \ref{val-lem} applies,
and shows $T(c'/K')$ has a maximal element.  \qed

\<{lem}\lbl{val-ort}  Let 
 $K$ be an algebraically closed valued field, $L$ an extension of 
transcendence degree $1$, $c,d \in L \setminus K$.  Then  
$E(c/K)=E(d/K)$.  \>{lem}

\proof  
 
If one of $T(c/K),T(d/K)$ has a last element, then by \ref{val-lem}
so does the other, and $E(c/K)=(0)=E(d/K)$.  Thus we may 
assume $T(c/K),T(d/K)$ have no last element, so that $\val K = \val L =: \Gamma$.
Further, since $E(c,K),E(d/K)$ are convex subgroups, one is contained in the other, so we may assume   that $E(c/K) \subseteq E :=E(d/K)$.

In case  $E  = \Gamma$,  , by the last statement of \ref{val-lem},  $L$ embeds   into  
$\widehat{K}$, and hence $E(c/K)=\Gamma$ too.

In general, let $\Gamma' = \Gamma / E $,
$r: \Gamma \to \Gamma'$ the quotient map,  and let $\val': L \to \Gamma'$
be the   valuation obtained by composition, $\val' = r \circ \val$.   Note that $\val'(L)=\val'(K) = \Gamma'$.
Call the  residue fields $K',L'$.  Let
$T'(y/K) = \{val'(y-b): b \in K \}$. 

  If $E(c/K) \neq E(d/K)$, then
 $T'(c/K)$ has a last element.  (Indeed let $\gamma \in E(d/K) $
 with $\gamma > 0$ and $\gamma \notin E(c/K)$.  Then there
 exists $\alpha \in T(c/K)$ with $\alpha + \gamma \notin T(c/K)$.
 Let ${\alpha}'$ be the common image of $\alpha,\alpha+\gamma$
 in $\Gamma'$.  Then clearly $\alpha'$ is the greatest element of
 $T'(c/K)$.)  By \ref{val-lem}, $T'(d/K)$ has a greatest element too.
Effecting additive and multiplicative translations, we may assume   these greatest
elements are $\val'(c)=0,\val'(d)=0$.  Let $c',d'$ be the residues
of $c,d$.

Since  $\val'(L)=\val'(K) = \Gamma'$, by \ref{val-lem} we must
have $K' \neq L'$.   We have an induced valuation  $\val'': L' \to E $,
$\val'' (\res' (a) ) = \val(a)$ for $a$ with $\val(a) \in E $.  
Clearly $T''(c'/K') = \{\val''(c'-b': b' \in K') \} = E \meet  T(c/K)$, and
similarly for $d$; thus $E''(c'/K') = E(c/K), E''(d'/K') = E(d/K)$.  But
now $T''(d/K)$ is the entire value group $E(d/K)$, so by the special case of $E=\Gamma$ considered above we have 
$E''(c'/K') = E''(d'/K')$.  \qed

\<{notation}  \lbl{calHnotation}  For a linear ordering $L$ and $X \subseteq L$, let ${\cal H}_L(X) = \{y: (\exists x \in X) (y \leq x)\}$ \>{notation}

\<{rem}\lbl{vt1c2}  Let $M \subset \Rr[\si]$ be a $\Zz[\si]$-module.     The nonzero convex subgroups form a single $\si$-orbit $\{E_n: n=1,2,\ldots\}$, in the sense that
$E_{n+1}$ is the convex hull of $\si(E_n)$; $E_n = \si^{-1}(E_{n+1}).$  \>{rem}   
 
\proof   Let $(0)=C_0 \subset C_1 \subset \ldots$
be the proper convex subgroups of $\Rr[\si]$ (with $\si^{-1}(C_{n+1})=C_n $).  Then for some $m$, 
 $M \meet C_{m+1}  \neq (0)$,   $M \meet C_{m}  = (0)$. 
So for $n \geq 0$,  $E_n = C_{n+m} \meet M$ are convex subgroups of $M$.
If $E$ is any proper convex subgroup of $M$, then the convex hull of $E$ in $\Z-si$
must be some $C_{n+m}$, and it follows that $E=E_n$.    \qed

\<{cor}\lbl{vt1r1} Let $K'=K(c)^a, K''=K(d)^a$ be two 
valued field extensions of an algebraically closed valued field $K$,  with $E(c/K) \neq E(d/K)$.  Then $K',K''$ have a unique valued field
amalgam. 

We have 
$T(c/K'') = {\cal H} T(c/K)$ (where  ${\cal H}(X) = \{y: (\exists x \in X) (y \leq x)\}$ ).

In particular if $\val K'' = \val K$ then $T(c/K'')=T(c/K)$,
and $E(c/K'') = E(c/K)$.
 \>{cor}

\proof  Assume $K',K''$ are embedded in $L = K'K''$.  Let
 $a \in K', b \in K''$.  Then $E(a/K) = E(c/K) \neq E(d/K) = E(b/K)$, so $T(a/K) \neq T(b/K)$.  One of $T(a/K),T(b/K)$ is bigger, say
$T(a/K) \subset T(b/K)$.  Find $e \in K$ with $\val(b-e) > T(a/K)$.
Then $\val(a-b) = \val(a-e)$.  This shows that the values of elements
$a-b$ are all determined.  In particular the values $\val(c-b)$
for $b \in K(d)$ are determined; this determines $tp(c/K'')$ and 
hence $tp(K'/K'')$.   

Moreover we see that
every element $\val(c-a)$ of $T(c/K'')$ either equals $\val(c-e)$
for some $e \in K$, so that it lies in $T(c/K)$, or else 
equals $\val(a-e) $ where $\val(c-e) > \val(a-e)$; so
in either case it lies in ${\cal H} T(c/K)$.  
\qed

This proof uses the stationarity theorem of \cite{HHM2} (cf. Theorem 2.11); 
alternatively one can prove the corollary directly, along the lines
of the proof of \ref{val-ort}.

\<{prop}\lbl{vt1r3}  Let $L$ be a valued field extension of an algebraically closed valued field $K$.  
 Let $a_i \in L \setminus K$,   
$$T_i = T(a_i/K) = \{val(a_i -b): b \in K \} $$
Assume  that the convex subgroups 
 $$E_i = E(a_i/K)=  \{e \in \val K: e+T_i=T_i\}$$
 are distinct.  
   Then $a_1,\ldots,a_n$ are algebraically independent over $K$.
The valued field structure of $K(a_1,\ldots,a_n)$ is uniquely determined
by the $T_i$.
 \>{prop}

\proof  By assumption, $a_i \notin K^{alg}$. Use induction
on $n$.  Some
$E_i$ is nonzero; say  $E_1 \neq (0)$.  Let $K'  = K(a_1)^a$.  By \ref{val-lem},
$\val K' = \val K$.  By \ref{val-ort}, $a_i,a_1$ are algebraically
independent over $K$ for $i \neq 1$.  By  \ref{vt1r1}, the type
of $(a_i,a_1)/K$ is determined, and 
$E(a_i/K') = E(a_i/K)$.  By induction, $a_2,\ldots,a_n$ are
independent over $K'$, and their type over $K'$ is determined.
The conclusion follows.  \qed

 Let $L$ be a model, and let
 Let $F \leq K_{1},\ldots,K_{n} \leq L$.  We say that $K_{1},\ldots,K_{n}$ are almost orthogonal over
 $F$ if for any $a=(a_{1},\ldots,a_{n})$, $a_{i}$ a tuple of elements of  $K_{i}$, $\union_i tp(a_i/F)$
 implies $tp(a/F)$.  

In the case of  an algebraically closed valued field $F$,
 $K_{1},K_{2}$ are almost orthogonal over
 $F$ if whenever $f_{i}:K_{i} \to M$ are embeddings 
 of valued fields, with $f_{1}| F =  F_{2} | F$, there exists
 a valued field embedding $f: K \to M$, with $f | K_{i}=f_{i}$. 

In this language, \ref{vt1r3} asserts the almost-orthogonality of certain
unary types.

   Let $F,K,L,K',L$ denote algebraically closed valued fields.  $\widehat{F}$
 denotes the completion of $F$.  

  When $K \leq L$ is an extension of valued fields, $K=K^a$, 
and $tr. deg. _K L = tr. deg. _{K_{\res}} L_{\res}$, we will say
that $L/K$ is {\em purely inertial}.   If instead 
$tr. deg. _K L  = rk_{\Qq} (\val(L)/\val(K))$, we will say that
$L/K$ is {\em purely ramified.}

The following lemma will use Theorem 2.11 of \cite{HHM2}, or rather the following corollary:
if $F(a),K$ are almost orthogonal over $F=F^{a}$, then so are $F(a)^{a},K$.
(Note that this is immediate for   perfect Henselian closure in place of
algebraic closure.)

 \<{lem}\lbl{acf-acvf-comp}  
  
 Let $F \leq K_{1},K_{2}$ be algebraically  
  closed valued fields.  Assume $K_2$ can be embedded into $\hat{F}$ over $F$;
  and that in any embedding of $K_1,K_2$ into a valued field extension $L$ of $F$, the images
  are linearly disjoint over $F$.  
   Then $K_{1},K_{2}$ are almost-orthogonal over
 $F$.
 
 In particular, if  $F \leq K_{1},K_{2} \leq \hat{F}, $ then $K_{1},K_{2}$ are almost-orthogonal over
 $F$ if and only if  they are  linearly disjoint over $F$.
 
  \>{lem}

 \proof  We may assume here that $tr. deg. _{F }K_{1} < \infty$. 
 If $F \leq F' \leq K_1$ then the hypotheses hold for $F,F',K_2$ and then also
 for $F',K_1,F'K_2$; so using transitivity of almost orthogonality 
    the lemma reduces to the case 
    $tr. deg. _{F }K_{1} =1$.  Fix some embedding of $K_1,K_2$ into a valued field $L$.
    Let $a \in K_{1} \setminus F$.  
    Then $a \notin K_{2}$, by the linear disjointness assumption.
      Let $B$ be a $K_{2}$-definable ball
    with $a \in B$.  So $B = B_{\rho}(b)= \{y: \val(y-b) \geq \rho \}$,
    $\rho  \in \val(K_{2}), b \in K_{2}$.  
As $a \notin K_{2}$, $\rho < \infty$.
Since $b  \in \widehat{F}$
    and $\rho \in \val(\widehat{F}) = \val(F)$, there exists $b' \in F$
    with $\val(b-b')> \rho$.  So $B$ is $F$-definable.  Thus
    $tp(a/F)$ implies $a \in B$.  Since $B$ was arbitrary, and $K_{2}$ is 
    algebraically closed,  the formulas $x \in B$
    of this kind generate $tp(a/K_{2})$.  Thus
    $tp(a/F)$ implies $tp(a/K_{2})$.  So $F(a),K_{2}$ are almost 
    orthogonal over $F$. By Theorem 2.11 of \cite{HHM2}, $K_{1}
    ,K_{2}$ are almost-orthogonal over $F$.    
   \qed

\lemm{complete-max1}  Let $F$ be a complete valued field, $K$ an immediate extension of
$F$, and $a_1,\ldots,a_n \in K$ linearly independent over $F$.  Then there exists $\lambda=\lambda(a_1,\ldots,a_n)  \in \G(F)$
such that if $b_1,\ldots,b_n \in F$ and $\min_i \val(b_i) =0$ then $\val(\sum b_i a_i) < \lambda$. \>{lem}

\prf  By induction on $n$.  The case $n=1$ is evident, with $\lambda(a_1)=\val(a_1)$.
Suppose for contradiction that for each $\lm \in \G$ there exist $b_1^{\lm}, \ldots, b_n ^{\lm} \in F$ with  $\min_i \val(b_i^{\lm}) =0$ and $\val( \sum b_i^{\lm} a_i ) > \lambda$.  For some $i \leq n$, for cofinally many
$\lm$ we have  $\val(b_i^{\lm})=0$; so we may assume $\val(b_n^{\lm}) = 0$;
and so in fact we may assume $b_n=-1$.   Now for $\lm' > \lm$ we have
 $\val(\sum_{i=1}^{n-1} (b_i^{\lm '} - b_i ^{\lm})a_i) > \lm$.   
 Write $  (b_i^{\lm '} - b_i ^{\lm}) = c_i(\lm,\lm') e(\lm,\lm')$ with $\min_i c_i(\lm,\lm') = 0$.
 Then $\lm(a_1,\ldots,a_{n-1}) > \sum c_i(\lm,\lm') a_i  > \lm - \val(e(\lm,\lm'))$
 so $$\val(b_i^{\lm '} - b_i ^{\lm}) \geq \val(e(\lm,\lm')) > \lm - \lm(a_1,\ldots,a_{n-1}) $$
 Hence $b_i^\lm $ convereges; let $b_i^{\infty} = \lim_{\lm \to \infty} b_i^{\lm}$.  
 Clearly $\sum_{i=1}^n b_i^{\infty} a_i = 0$.  But this contradicts the $F$ linear independence of 
 $a_1,\ldots,a_n$.
\eprf

 \lemm{complete-max}    Let $F =  \widehat{F}$ be an 
  algebraically     closed valued field.  Let $K_1$ be an immediate extension of $F$,
  and let $K_2$ be another extension of $F$, almost orthogonal to $K_1$ over $F$.
  Then $K_1,\hat{K_2}$ are almost orthogonal over $F$.  Moreover the compositum
  $K_1 \hat{K_2}$ within any common valued field extension is an immediate extension of $\hat{K_2}$     \>{lem}
 
\prf  Using transitivity of almost orthogonality as in  \lemref{acf-acvf-comp},
and quoting   \lemref{acf-acvf-comp} for $(K_1K_2,K_2,\hat{K_2})$ in place of $(K_1,F,K_2)$, 
 it suffices to show that 
for any valued field $L \supset F$ and valued field embeddings $K_1 \to L$, $ \hat{K_2} \to L$ over $F$,
the images of $K_1,\hat{K_2}$   are linearly disjoint over $F$.    We may take these embeddings to be the inclusions.   

Let $a_1,\ldots,a_n \in K_1$ be linearly independent over $F$.  Let $\lm = \lm(a_1,\ldots,a_n)$
from Lemma \ref{complete-max1}.   Suppose, for contradiction, that $b_1'',\ldots,b_n'' \in \hat{K_2}$
and $\sum b_i''a_i =0$.  Then we can approximate $b_i''$ by $b_i' \in K_2$ to within valuation $>\lm -\val(a_i)$, so that $\val(\sum b_i' a_i) > \lm$.  Let $b' = (b_1',\ldots, b_n')$.  By the almost orthogonality, $tp(b'/F)$  implies $\val( \sum x_i a_i) > \lm$.  So some formula
$\phi(x) \in tp(b'/F)$ implies the same.  But $F \models ACVF$, so there exists
$b=(b_1,\ldots,b_n) \in F^n$ with $F \models \phi(b)$.  Thus $\val (\sum b_i a_i) > \lm$.  But this
contradicts the definition of $\lm$.   

For the last point, note that both $K_1$ and $\hat{K_2}$ embed into the 
maximal immediate extension of $\hat{K_2}$, over $F$.  By almost orthogonality
the compositum within this extension is isomorphic over $F$ to the compositum within any
other. \eprf

\<{lem}\lbl{domin-val} Let $K_0 \subset K_1 \subset  \cdots \subset K_n$ be algebraically closed valued fields, $K=\hat{K_0}, L=\hat{K_n}$.    Assume: 
$\val(K_0)$ is cofinal in $\val(K_n)$; and for each $i$, $K_{i+1}/K_i$ is 
purely ramified, or purely inertial, or $K_{i+1} \subset \widehat{K_i}$.  Let  
$\bK$ be an immediate extension of $K$.    Then:
   $\bar{K},L$ are almost orthogonal over $K$, and $M=L \bK$ is  an immediate extension of $L$.  
 \>{lem}

\proof  
\def\L-{{L^{-}}} \def\M-{{M^{-}}}
(1) Let $\L- =\hat{K_{n-1}} ,  \,  \M-   = \L-\bK$.  We may 
assume inductively that $ \bar{K},\L- $ are almost orthogonal over $K$, and $\M-$  is  an immediate extension of $\L-$, hence of ${K_{n-1}}$.  We will show that $\M-,L$ are
are almost orthogonal over $\L-$, and $M$  is  an immediate extension of $L$.  
From this (1) follows by transitivity.  

If $K_n \subseteq \hat{K_{n-1}}$, then $\L-=L$ and there is nothing to prove.  
Otherewise $K_n/K_{n-1}$ is  purely ramified, or purely inertial.  Since $\M-$
  is  an immediate extension of   ${K_{n-1}}$, it follows that
$\L- K_n / \L-$ is purely ramified, or purely inertial; for instance in the purely inertial case,
using $(K_{n-1})_{\res} = (\L-)_{\res}$ we have
$$tr. deg._{\L-} \L- K_n   \leq tr. deg._{K_{n-1}} K_n  =  tr.deg._{(K_{n-1})_{\res}} (K_n)_{\res} \leq 
tr.deg._{(\L-)_{\res}} (\L- K_n)_{\res}  $$
so the numbers are equal. 
This also shows that   $\L-K_n  $ is
almost orthogonal to $\M-$ over $\L-$.   By \lemref{complete-max}, $L$ is 
almost orthogonal to $\M-$ over $\L-$, and
$\M- L$ is an immediate extension of $L$.  But $M=\M-L$.                 \qed

{
{\noindent \bf Remark}  Taking $\bar{K}$ to be a maximal immediate extension of $\hat{K}$, it follows 
from    \cite{HHM2}, Theorem 6.13 that $L$ is dominated by $L_\res \union \{a_i\}$  over $K \cup \val(L)$;
where $a_i$ are elements of $L$ such that $\val(a_i)$ form a $\Qq$-basis for $\val(L)$ over $\val(K)$.
}

\<{lem} \lbl{gauss} Let $k \leq K$ be   algebraically closed fields.  Let $v$ be a valuation of
$K/k$, $R$ the valuation ring.  Let $X$ be a quasi-projective scheme over $R$, with $X(K)$ and $X(k)$ finite.
Consider the residue homomorphism $R \to k$, and the induced map $r: X(R) \to X(k)$.
If $p \in X(k)$ and $r^{-1}(p)$ has $n$ distinct points, then $p$ has geometric multiplicity 
$\geq n$ on the scheme $X_0 = X \tensor_R k$.  More generally, this holds true if $r^{-1}(p)$
has $n$ points counted with multiplicities on $(X \tensor_R K) $.  If $X$ is flat over $R$ (or, is reduced and has no component lying over the special fiber of $\spec R$) equality holds.      \>{lem}

\proof  Find an affine open subscheme containing the said $n$ points $q_1,\ldots,q_n$, and thus   reduce   to $X = \spe A$, $A$ a finitely 
generated  $R$-algebra.   Let $B$ be the image  of $A$ in $A \tensor_R K$.     As $X$ is $0$-dimensional, $A \tensor_R K$ is a finite dimensional $K$-space, and so the finitely generated $R$-submodule $B$ is a free $R$-module of finite rank.   
  It follows that  $\dim_k (A \tensor_R k) \geq 
\dim_k  ( B \tensor _R k) = \dim_K (B \tensor _R K)  
= \dim_K(A \tensor_R K)$.

If $X$ is flat over $R$,  then $A$ is a flat $R$ module, so $A \to A \tensor_R K$ is an injective map of modules, and hence $B \cong A$.

 \qed

 \end{subsection}

\<{subsection}{Structure theory}

\<{example}  Let $K$ be a transformal discrete valuation ring, with value group 
 $\Zz[\si]$, and   residue field   $F(a)_\si$, where
  $g(a)=0$ for some difference polynomial $g$ over $F$ while $g'(a) \neq 0$.   Then  there exists a $1$-generated
  transformal discrete valuation ring over $F$, with the same value group and residue field.  \>{example}
  
\<{proof}  We may assume $K$ is Henselian.  Pick $t$ such that $\val(t)$
generates the value group as a $\Zz[\si]$-module.  Lift  $g$ to
$\OK(x)$, and lift $a$ to $\OK$.  Then $\val(g(a)) > 0$, so $\val ( (g-t)(a) > 0)$.  
By one step of Hensel's lemma, we can perturb $a$ so that $\val( (g-t)(a)) > \val(t)$.  Then 
$\val g(a) = \val (t)$.  So $F(a)_\si$ has the same residue field and value group as $K$.
\eprf

\<{example} \ \>{example}\begin{enumerate}
  \item   Let $F$ be a valued field.  Form $F(t)_\si$, and let $L$ be the completion.
Let $a = \sum_{n=0}^\infty  t^{\si^n} \in L$, $K = F(t,a)_\si$.  We have $\si(a)-a-t = 0$.  
 
  \item   Let  
  $y^\si-y = t^{-1} $.  Then $(F(t)_\si)^{inv}(y)$ is an immediate extension of $(F(t)_\si)^{inv}$:
$$ y = t^{-1/\si} + t^{-1/\si^2} + \ldots $$
  $F(t,y)_\si$ is a "ramified" extension of $F(t)_\si$, in the sense
that $dim_{\Qq}(\Qq \tensor ( \val (F(t,y)_\si) / \val( F(t)_\si)) )=1 $.
\end{enumerate}

\<{prop}\lbl{vt4}  Let $K=K^a$ be a transformal valued field with value group 
$\val(K) \subset \Q-si$.  
Assume:  (*) for any $c \in K$,
$\{\val(c-b): b \in \si(K) \}$ has a greatest element.

Let $L$ be a valued difference field extension of $K$, with 
$  tr. deg._K(L) < \infty$, $\val(L)=\val(K)$,
and $L_{\res} = K_{\res}$.  Then $\widehat{L} =  \widehat{K}$.      \>{prop} 

\proof  Let $a \in L$; we will show that $a \in \widehat{K}$.  We may assume $\val(a)>0$.  Let
$a_n = \si^n(a)$, $C_n =  \{ \max(0,\val(a_n-b)): b \in K  \}$, $C=C_0$.
 Then $C_n$ is closed downwards in the non-negative elements $\val(K) $.
$C_n$ cannot have a maximal element, since
$L$ is an immediate extension of $K$.    
  
  If $C$ is unbounded,
then $a \in \widehat{K}$.  Suppose for contradiction that $C$ is
bounded.  We have $\Q-si = \Qq v  \oplus \Qq v^\si + \ldots$; let 
$E_m = \Qq v \oplus \ldots \oplus \Qq v^{\si^m}$. 
Let $m$ be least such that $C \leq E_m + c_0$ for some $c_0 \in C$ (here we write $X \leq Y$
if $(\forall x \in X)(\exists y \in Y) (x \leq y)$.)  
Since $\val(L)=\val(K)$, dividing $a$ by an element of $K$ with value
$c_0$, we may assume 
$C \leq E_m$; as $E_m$ is convex, we have $C \subset E_m$.  
Then $\si(C) \subset E_{m+1}$.  Let $H_1$ be the convex hull of $\si(C)$.

\claim $E_m + H_1 = H_1$.

\proof Let $0<e \in E_m$, $0<c \in C$.  Since $C \not \leq c+E_{m-1}$
there exists $c' \in C$ , $l \in \Nn$ with $c' -c > e  / l$.  Then $\si(c'-c) > l (c'-c) > e$.
So $e+ \si(c) \leq \si(c')$.

\claim $H_1=C_{1}$

\proof  Let $c_1<c_2<\ldots$ be cofinal in $C$.  Let $b_n \in K$,  $\val(a-b_n) = c_n$; 
let $B_n = \{x: \val(x-b_n) \geq c_n \}$; $P = \meet _{n \in \Nn} B_n$;
$ {B_n}^\si = \{x: \val(x-\si(b_n)) \geq \si(c_n) \}$;
$P^\si = \meet_{n \in \Nn} {B_n}^\si$.
 
Then $P(K) = \emptyset$.  So $P^\si(K^\si) = \emptyset$.   
By the assumption (*) on $K/K^\si$, $P^\si(K) = \emptyset$.  Since
$a_1 \in P^\si$, it follows that there can be no $b \in K$
with $\val(b-a_1) > \si(C)$.  Thus $C_1 = H_1$.

Now $\si^n$ induces an isomorphism $(K,\si(K)) \to (\si^n(K),\si^{n+1}(K))$.
So the convex hull of 
$\si^n(C)$ is  $C_n$ for all $n$; $E_{m+n} + C_{n+1} = C_{n+1}$;
and $C_n \subset {\cal H}(E_{m+n})$ (\ref{calHnotation}).
   Thus the hypotheses of Lemma \ref{vt1r3}
are valid of the $C_m$.  

But now by \ref{vt1r3}, the elements $a_n$ must be
algebraically independent over $K$.  
This contradicts the assumption of finite transcendence degree. 
  \qed

The condition (*) is easily seen to apply in the fundamental case: $K=F(t)_\si^a$.
(See \ref{vt4-e}.) 

\<{example} The composition of difference  polynomials extends to $F(t)_\si$, and by
continuity to the completion $L$ of $F(t)_\si$. 
 By \ref{vt4}, difference polynomials with nonzero linear term have a left compositional  inverse in $L$.  \>{example} 
For instance, the compositional inverse of $t+t^{\si+1}$ is the repeated fraction:
\def\ff2{{1+\frac{t^{\si^2} }{1+\cdots}    }}
$$\frac{t}  {1+\frac{t^\si }{\ff2}    }$$

\<{paragraph} {Transformally Henselian fields}

\<{lem}\lbl{hensel-c}  Let $K$ be a  complete, algebraically Henselian transformal valued field over $F$, with value group 
$\val(K) \subset \Q-si$.  Let ${f} \in \OK[X]_\si$ be a difference polynomial,
 ${f_\nu}  = \partial_{\nu} {f}$.  
Suppose  $a \in \OK$, $v=\val ({f}(a)) > 2 \val ( f_1(a)) =2v'$.   Then 
there exists $b \in K$, $f(b)=0$, $\val(a-b) \geq v-v'$. \>{lem}

\proof   Let $e_0 = f(a)$.
We will find
$a_1 \in \OK$, $\val(a-a_1) \geq v-v'$, $\val f(a_1) \geq  \si(v)$.  Since
$\val ( f_1(a) - f_1(a_1)) \geq v-v' > v' = \val f_1(a)$, it follows that $\val f_1(a_1) = v'$.
 Iterating this we get a sequence
$a_n$, with $v_n =_{def} \val f(a_n) \geq \si(v_{n-1})$, and $\val (a_{n-1} - a_n) \geq v_n - v'$. 
As $K$ is complete, there exists a unique $b \in K $ such that the sequence $a_n$ converges to $b$;
by continuity, $f(b)=0$.  

To obtain $a_1$, write $a_1= a + re$, with $e=e_0 / f_1(a)$, and $r \in \OK$ to be found.  
    Then 

$$f(a_1)=f(a+re) = e_0 + \sum_{i=1}^m f_i(a) (re)^i + \sum_{\nu} f_\nu(a) (re)^\nu$$

Here the $f_i$ are transformal derivatives, cf. \reff{trans-der}
 $\nu$ runs over indices $\geq \si$, while $i$ ranges over the nonzero finite indices.  Note
that $\val    f_\nu(a) (re)^\nu \geq \val(\si(e)) = \si(v)$.  Thus it suffices to find $r$ with 
$e_0+\sum_{i=1}^m f_i(a) (re)^i =0$, or with $r$ a root of 
$$g(x)= 1+ x + f_2(a) (e^2/e_0) x^2 + \ldots + f_m(a) (e^{m}/e_0) x^m $$
(Divide through by $e_0$, noting  $f_1(a)(e/e_0) =1$.)
This (ordinary) polynomial has derivative 
$$g'(x) = 1 + 2f_2(a)(e^2/e_0)x + \ldots $$
Now   for $k \geq 2$, $\val (e^k / e_0) = k (v-v') -v  
  \geq 2(v-v') - v = v-2v' > 0$.  Thus $\res g'$ is a nonzero constant polynomial.  By the
ordinary Hensel lemma, there exists $r \in \OK$ with $g(r)=0$, near any $c \in \OK$ with $\val g(c) > 0$.
Letting $c = -1$ 
 we see that there  exists $r \in \OK$ with $g(r)=0$.  \qed

We will tentatively call a transformal valued field 
satisfying the conclusion  of \ref{hensel-c}
of characteristic $0$
{\em transformally henselian}.   In characteristic $p>0$, we demand
also that the field $L$ be perfect, and that all the Frobenius
twists $(L,\si \circ (\phi_p)^n)$ also satisfy the conclusion  of \ref{hensel-c}. 
( $\phi_p(x)=x^p$).

\>{paragraph}

 \<{paragraph}{Lifting the residue field}

 \<{cor} \lbl{res-lift-lem-0}  Let $K =K^a  \leq L $ be transformal valued
fields, with $L$   transformally henselian.  Let 
 ${f} \in K_{\res}[X]_\si$ be a difference polynomial
of order $r$, ${f}' = \partial_1 {f}$. 
  Let $\ba \in L_{\res}$, ${f}(\ba)=0$, ${f}'(\ba) \neq 0$,   
$tr.deg._{K_{\res}} K_{\res}(\ba) = r$.
Then   there exists a purely inertial extension $K'=K(b)_\si$ of $K$,
$K' \leq L$, $\ba \in \res(K')$.    \>{cor}

\proof  Pick any $a \in {\cal O}_L$ with $\res(a)=\ba$, and also lift $f$ to $\OK[X]_\si$.  Then $\val f(a) > 0$,
$\val f'(a) = 0$, so by the transformal Hensel property (cf. \ref{hensel-c})  there exists $b \in \OK$ with $f(b)=0$, $\res(b)= \ba$.  Clearly $K'=K(b)_\si$
is purely inertial over $K$.  \qed

Let $K$ be an $\omega$-increasing transformal valued field over $F$.   The set of elements
of $K$ satisfying nontrivial difference polynomials over $F$ forms a difference subfield $F'$ of $K$; it is the union
of all subfields of   transformal transcendence degree $0$ over $F$.  The residue map $\res$ is injective
on $F'$:  to see this, suppose  $\res(a)=0$ and $a \neq 0$.  Then $\val(a) > 0$, then $0<val(a) << \val(\si(a)) << \ldots$, so $a,\si(a),\ldots$
are algebraically independent over $F$, so $a \notin F'$.

\<{cor} \lbl{res-lift-lem}  Let $K$ be an $\omega$-increasing   transformally henselian  valued field over $F$.
  Let $\ba \in K_{\res}$, ${f}(\ba)=0$, ${f}'(\ba) \neq 0$, ${f} \in F[X]_\si$ a difference polynomial, ${f}' = \partial_1 {f}$.  
Then $F(\ba)_\si$ lifts to $K$:  there exists a subfield $F(b)_\si$ of $K$ such that the residue map 
restricts to an isomorphism $F(b)_\si \to F(\ba)_\si$ (of difference $F$-algebras.)  \>{cor}

\proof  Pick any $a \in \OK$ with $\res(a)=\ba$, and also lift $f$ to $\OK[X]_\si$.  Then $\val f(a) > 0$,
$\val f'(a) = 0$, so by \ref{hensel-c}, there exists $b \in \OK$ with $f(b)=0$, $\res(b)= \ba$.  The difference
field $ F(b)_\si$ has finite transcendence degree $m$ over $F$, and hence contained in $F'$; 
thus it is trivially valued, so the residue map is injective on $F(b)_\si$.  since $\res(b) = \ba$,
it   carries $ F(b)_\si$ into and onto $F(\ba)_\si$.  \qed

\<{prop} \lbl{res-lift}  Let $K$ be a perfect $\omega$-increasing   transformally henselian  valued field
over a  difference field $F$.
 
Assume $F^{inv} \subset F^a$, and $K_{\res}$ has transformal transcendence degree $0$ over $F$.  Then
  $\res: F' \to K_{res}$ is surjective.  
\>{prop}

\proof   By Lemma \ref{inv-fin}, $(F')^{inv} \subset (F')^a$; and clearly $F'$ is perfect.
We may   assume $F'=F$, and show that $\res(K)=F$.       Suppose
$K_{\res} \neq F$.  In characteristic $0$, choose a difference polynomial $f \in F[X]_\si$ of smallest possible
transformal order and degree, and  $\ba \in K_{\res}$
with $f(\ba)=0$.  If $f \in F[X^\si]_\si$, then replacing $\ba$ by $\si(\ba)$ we may
lower the degree of $f$; unless $\si(\ba) \in F$.  But in this case, $\ba \in F^a$,
so the order of $f$ must be $0$, so $f \in F[X^\si]$ is impossible.  Thus $f \notin F[X^\si]$.
It follows that $f' \neq 0$ and $f'$ has smaller order or degree.  So the hypothesis
of \ref{res-lift-lem} is met, and we can increase $F'$.

In positive characteristic, this may not be the case, because of polynomials such as $x^\si - x^p$.  
However, after replacing $\si$ by $\tau = \si(x^{1/p^m})$ for large enough $m$, it is possible to 
find such an $\ba$ for $\tau$, lowering the degree of $f$ further.  Assuming the residue field of finite total dimension, this permits lifting the residue field; we may then   return from $\tau$ to $\si$.  \qed
 
\<{rem}\ \>{rem}  In positive characteristic, when $K$ is a 
transformal discrete valuation ring
over $F$, 
the Proposition applies to the perfect closure of $K$.  But  we can   lift $K_{res}$
to a subfield of  the completion of  a smaller extension $L$ of $\hat{K}$ within the perfect closure,
such that $L$ is still the completion of a  transformal discrete valuation ring, and in particular
$\val(L) \subset \Z-si$;  the details are left to the reader.  
 
\<{example}  \ \lbl{res-lift-ex} \>{example}  The assumption in \ref{res-lift} that $F^a$ is inversive is
necessary:  take $F = \Qq(b,b^\si,\ldots), K=F(t,c), \si(c)=b+t$.
\>{paragraph} 

\<{paragraph}{Characterizing the completion}

Here is one version of a  theorem characterizing completions of 
transformal discrete valuation rings.  

\<{lem}\lbl{char-comp} Let $L=L^a$ be a transformal valued field over $F=F^{inv}$, 
with value group $\subset \Qq_\si$.    Assume $L$ is complete.  Then the following conditions 
are equivalent:  \begin{enumerate}
  \item $L=\widehat{L_0}$ for some transformal valued field $L_0$ of transformal dimension $1$ over $F$.
  \item Whenever $L' = (L')^a  \leq L$ is complete with the same value group and
residue field,   $L'=L$.
  \item $L \iso \widehat{K}$, where $K = (L_{\res}(t)_\si)^a$.  
\end{enumerate}
 
\>{lem}

\proof  Assume(1), and let $L'$ be as in (2).  By \ref{vt2},  $L$ has value group $\Qq_\si$, and   $L_{\res}$ has
transformal transcendence degree $0$ over $F$.  So $L_{\res}$ is inversive.  

Let $t \in L'$ be such that $\val(t)$ generates $\val(L)$ as a $\Qq[\si]$-module.
By \ref{res-lift}, $L_\res$ lifts to a trivially valued difference subfield $F'$ of $L'$.
Let $K = (F'(t)_\si)^a$.  Then $K_{\res} =   L_{\res}$, and $\val(K) = \val(L)$.  

We have $\val(K) \neq \val(K^\si)$, while $tr. deg._{\si(K)} K = 1$.  By
\ref{val-lem}, the  hypothesis (*) of \ref{vt4} holds for $K$.  Thus by 
\ref{vt4}, $L =  \widehat{K}$.  Since $\widehat{K} \leq L' \leq L$, we have also $L'=L$.  
This proves (2).  

The same proof   shows that $K$ embeds into $L$; so assuming (2), we conclude (3),
using the assumption of (2) in place of \ref{vt4}.

(3) implies (1) trivially.
  \qed

 \>{paragraph}
 \<{paragraph}{Lifting the value group}
Let $K=K^a$ be a valued field. 
If $K'/K$, $K''/K'$ are purely ramified,
then clearly $K''/K'$ is purely ramified.  Thus $K$ has a maximal purely
ramified extension $K'$ within a given extension $L$ (not necessarily unique.) 

 
Let $L$ be a transformally valued field.  Let 
$$H_e =H_{e}(L) =  \{v \in \val(L): (\forall u > 0) (|v| < \si^{e} (u)) \}$$
$H_e$ is a convex subgroup of $\val(L)$;  in case $\val(L) = \Qq_\si$,   $H_e$ is the $e$'th nonzero convex 
subgroup.  When $v > c$ for every $c \in H_{e}$, we write: $v>H_{e}$.
Let $h_e: \val(L) \to \val(L)/H_e$ be the quotient map, and let 
$$VAL_e = h_e \circ \val :  \ L^* \to \val(L)/H_e $$

We obtain auxiliary valuations, with residue map denoted $\RES_e$.   Moreover we obain a valuation 
$VAL^*_e: \RES_e \to H_e$, with residue field $L_{\RES}$.
Neither of these is  a transformal valuation.   


When  $L =  {F(t)_\si}^a$, $\RES_e$ induces an isomorphism  $  F(t,t^\si,\ldots,t^{\si^{e-1}})^a \to \RES_e(L)$ naturally.
When $L=\widehat{F(t)_\si}^a$ the same is true, since any element of $L$ is close to an element
of $F(t)_\si$ to within a value $>H$.


\<{prop} \lbl{lift-val}    Let $L=L^a$ be a transformal valued field of transformal dimension $1$ over an inversive difference
field $F=L_{res}$; with value group $\Q-si$. 
Let $F < K=K^a \leq \widehat{L}$.       Then  
 there exists a difference field $K' \leq \widehat{L}$,   with \begin{enumerate}
  \item $K \leq K'$, and $K'/K$ is purely ramified.
  \item With $e =  rk_{\Qq} \val(L)/\val(K')$, $\si^e \widehat{L} = \widehat{K'}$.
\end{enumerate}
                       \>{prop}

\proof   

Let $K'$ be a maximal   purely ramified difference field extension of $K$ within $\widehat{L}$.  Then $K' = (K')^{a}$.

By  \ref{char-comp}, $\widehat{L} = \widehat { F(t)_\si^a}$ for some $t$.  Let $e = rk_{\Qq} \val(L) / \val(K')$.  
Then there exists $s \in K'$, $\val(s) >0$, $\val(s) \in H_{e+1}$.  Fix $N \geq e+1$ for a moment.  
 Let $t_n = \si^n(t), s_m = \si^m(s)$, 
$\bar{t_n} =\RES_N(t_n)$, $\bar{s_n} = \RES_N (s)$.  

We saw that $tr. deg._F \RES_{N} (L) = N$.   Thus
 $\bar{t_0},\ldots,\bar{t_e},\bar{s_0},\ldots,\bar{s}_{N-e-1}$
are algebraically dependent over $F$.   Since $s \in H_{e+1}$, 
$\bar{s_0},\ldots,\bar{s}_{N-e-1}$ are algebraically independent over $F$, using the valuation $\val^*_e$     .  Let $m$ be maximal such that 
 $ \bar{t_m},\ldots,\bar{t_e},\bar{s_0},\ldots,\bar{s}_{N-e-1}$ are algebraically dependent over $F$.  
 
Write $f( \bar{t_m},\ldots,\bar{t_e},\bar{s_0},\ldots,\bar{s}_{N-e-1})   =0$, $f \in F[ X_m,\ldots,X_e,\bar{s_0},\ldots,\bar{s}_{N-e-1}]$ of minimal $X_m$- degree.

 We now argue that $m=e$. Lift $f$ to $ f( X_m,\ldots,X_e,s_0,\ldots,s_{N-e-1}) \in {\cal O}_L [X_m]_\si$,   viewed  as a difference polynomial in $X_m$ ( with coefficients involving the $s_i$, and $X_{m+j} =  (X_m)^{\si^j}$.)
  In characteristic $0$,
let $f'$ be the (first) derivative of $f$ with respect to $X_m$; then $f' \neq 0$, and so by minimality of $f$,
 $RES_N f'( \bar{t_m},\ldots,\bar{t_e})  \neq 0$, i.e. $\val  f'( t_m,\ldots,t_e) \in H_N$.  
By \ref{hensel-c}, applied to the element $t_m$, there exists $b \in \widehat{L}$,
 $\val(t_m -b) > H_{N}$, $f( b,\ldots,\si^{e-m}(b),s_0,\ldots,s_{N-e-1}) = 0$.  Let $L' = K'(b)_\si$.  
Clearly $tr. deg._K L' \leq e-m$,
while $rk_{\Qq} \val(L') / \val(K') \geq e-m$ 
(since $\val(\si^i(b))=\val(t_{m+i})$, and  $\val(t_m) << \val(t_{m+1}) << \ldots << \val(t_{e-1}) << \val(c)$
for any $c \in K'$, $\val(c) > 0$.)  
However 
$rk_{\Qq} \val(L') / \val(K') \leq tr. deg._{K'} L' $, so 
these numbers are equal, and $L'/K'$ is purely ramified, of transcendence degree $m-e$.
  But  $K'$ is a maximal purely ramified extension; so $K'=L'$, and   $m=e$.

In characteristic $p$, we first replace $\si$ by $\si \circ \phi^{-l}$ where $\phi(x)=x^{p}$,
reduce the equation if it is now a $p$'th power, till $t_{m}$ occurs with an exponent
that is not a power of $p$.   We then use the above argument.  So $m=e$ in any
characteristic.  

 Thus  $\bar{t_{e}} \in \RES_{N}(K')$ (since this residue field is algebraically
closed).  So for some $b_N \in K'$, $\val(t_{e }-b_N ) > H_{N }(L)$.  Now letting $N \to \infty$,
 we see that $t_e \in \widehat{K'}$.   Let $K'' = F(t_e)^a_\si $;  then
 $  \si^e \widehat{L} = \widehat{K''}  \subset \widehat{K'}$,
and it remains to show equality.   Note that $\val(K'')  = \Qq[\si] \val(t_e) = \val (K')$, and
$\res(K'') = F = \res(K')$.   Thus $K'/K''$ is immediate, so for any $c \in K'$, 
$T(c) = \{\val(c-b): b \in K'' \}$ can have no maximum value (\ref{val-lem}).  But by \ref{vt4-e}, 
for any $c \in \widehat{L} \setminus \widehat{K''}$, $T(c)$ does have a maximal element.
Thus $K'=K''$.  \qed

 \<{cor}\lbl{char-comp-cor}  If $L$ satisfies \ref{char-comp} (1)-(3), and $K$ is an algebraically closed, complete   difference
 subfield of $L$, $K_{\res}=L_{\res}$, then either $K$ is trivially valued or $K$ satisfies the same conditions.  \>{cor}

 \proof  By \ref{res-lift}, we may assume $L,K$ are transformally valued fields over $F$, $F=K_{\res}=L_{\res}$.    Let $K'$
be as in \ref{lift-val}.  Then $K'/K$ is purely ramified, of finite transcendence degree $r$; there
exist  therefore $K_0 \leq K$,$K'_0 \leq K'$, $tr. deg._F (K_0)$ finite, $tr. deg. _{K_0} (K'_0) = r = rk_{\Qq} \val(K'_0)/\val(K_0)$,
 $K'=K K'_0$.  We can choose ${K_0}$ algebraically closed, and
with value group equal to $\val(K)$.   So $\widehat{K_0}   \leq \widehat{K'_0} \leq \widehat{K'}$.   $\widehat{K'}$ is isomorphic (by \ref{char-comp} (1 $\to$ 3)) to $L$, so it   satisfies
\ref{char-comp} (1)-(3); by (2), since $\widehat{K'_0} $ is a complete difference subfield with the same
value group and residue field, $ \widehat{K'} =  \widehat{K'_0}$.  Thus
$\widehat{K_0}   \leq K \leq \widehat{K'} = \widehat{K'_0}$.
If $c \in \widehat{K'} \setminus \widehat{K_0}$, as $K'_0/K_0$ is purely ramified, and as in the proof of 
\ref{vt4-e}, $K_0(c)/K_0$ is also purely ramified.  In particular if $c \in K \setminus \widehat{K_0}$ then
$\val (c) \notin K_0$; but we chose $\val(K_0) = \val(K)$.  Thus $K=\widehat{K_0}$. 
Now $\widehat{K_0}$
satisfies  \ref{char-comp} (1). 
 \qed
 

When $E/F$ is an extension of valued fields, let
$$\vrk(E/F) = \dim_{\Qq}(\Qq \tensor \val(E)/\val(F) )+ tr. deg._F \res E$$ 

When $E,F$ are subfields of a valued field $L$, we let $\vrk(E/F) = \vrk(EF/F)$,
$EF$ being the compositum of $E,F$ in $L$.  A similar convention holds
for transcendence degree.

\def\hK{\widehat{K}}
\<{cor} \lbl{vrk} Let $K$ be an $\omega$-increasing transformal valued field of transformal dimension $1$ over an inversive
 difference field $F$, with value group $\Qq_\si$.   Let $L=L^a$ be a difference field extension of $K$ of
 finite transcendence degree.     Then the following are equivalent:  
\begin{enumerate}
  \item $\vrk (L/K) \leq e$
  \item There exist  difference fields  $K \leq K_r \leq \widehat{K}$,
and  $K_r \leq K_1 \leq K_2 \leq   \widehat{L}$,  with: \begin{enumerate}
  \item $tr.deg._K K_r \leq tr. deg._F K_{\res}$  
  \item  $tr. deg._{K_r} K_1  =  tr.deg._{K_{\res}} L_{\res} = e_{\res}$
  \item  $tr. deg._{K_1} K_2 = rk_{\Qq} \val(K_2)/\val(K_1)  = e_1$
  \item $\widehat{K_2} = \widehat{L}^{\si^{e_2}}$
  \item $e_{\res}+e_1 +e_2 \leq e$
\end{enumerate}

  \item There exists a difference field $K'$, $ \widehat{K} \leq K' \leq  \widehat{L}$, 
$tr.deg._{\hK} K' + e_2 \leq e$, $\widehat{K'} = \widehat{L}^{\si^{e_2}}$.

\item  There exist difference fields $K',M $, $\hK \leq K', \widehat{K'} \leq M, \widehat{M} = \widehat{L}$, with
$tr.deg._{\hK} K' + tr. deg._{\widehat{K'}}(M) \leq e$.
\end{enumerate}

Thus, if $\vrk(L/K)< \infty$, there exists a chain 
of difference subfields $K = M_0 \leq M_1 \leq \cdots \leq M_5=\hat{L}$
such that  $M_{i+1} = \widehat{M_i}$ for even $i<5$, while for odd $i<5$
  $tr.deg._{M_i} M_{i+1} = \vrk(M_{i+1}/M_i)$.

  \>{cor}

\proof
 Since $F$ is inversive, and $L_{\res}$ is algebraically closed, by \ref{vt2} and \ref{inv-fin}, $L_{\res}$ is also inversive.  Thus for any $d \geq 1$, $\res ( L^{\si^{d}}) = \res (L)$.
By \ref{vt1}, $\vrk(\widehat{L} / \widehat{L}^{\si^{d}}) = rk_{\Qq} (\Q-si / \si^d(\Q-si)) = d$.   

 To show that (1) implies (2),  we may assume $K=K^a$. 
 By \ref{res-lift}, there exists a field of representatives $  F_r \leq \widehat{K}$ for the residue map of $\hK$.  We have 
$tr.  deg_F F_r = tr. deg._F K_{\res}$.  Let $K_r = F_rK$.

 By \ref{res-lift} again, applied to $\widehat{L}$ over $K_r$, there exists a difference field $F_r \leq L_1$ such that  $\res: L_1 \to L_{\res}$
is an isomorphism.    Let $K_1 = K_r L_1$.
Then $tr. deg._{K_r} K_1 \leq e_{\res}$; equality holds by 
comparing the residue fields.  Thus (b).  

Now over $K_1$, \ref{lift-val} applies, and gives $K_2$ 
with (c,d).  We have $K \leq \widehat{L}^ {\si^{e_2}}$,
and $e_1 \leq rk_{\Qq} \val( \widehat{L}^ {\si^{e_2}} / K)$;
so
$  e_1 + e_2 \leq  rk_{\Qq} (\val(\widehat{L} )/ \val K) =
 rk_{\Qq} (\val L / \val K) $.  Thus (e).

(2) obviously implies (3), with $K'=K_2 \hK$.  To go from (3) to (4), let 
$K'$ be as in (3).  By \ref{char-comp}, $\widehat{L} \iso
\widehat{(L_{\res}(t)_\si)^a}$, $\widehat{K'} = \widehat{L}^{\si^{e_2}} = 
\widehat{(L_{\res}(t^{e_2})_\si)^a}$.  
Let $M  = \widehat{(L_{\res}(t^{e_2})_\si)^a} (t)$.
Then $\widehat{M} = \widehat{L}$ and $tr. deg._{\widehat{K'}}(M) = e_2$.

(4) implies (1) since $\vrk$ is additive in towers, bounded by
transcendence degree, and vanishes for completions ($\vrk_K \hK = 0$.)

The final conclusion follows from (4) by taking $e =\vrk(L/K)$.
We obtain $M_0,\ldots,M_5$ with 
$M_{i+1} = \widehat{M_i}$ for even $i<5$, and 
$tr.deg._{M_1}M_2 + tr. deg._{M_3} M_4 \leq e$.
But $e = \vrk(M_2/M_1) + \vrk(M_4/M_3)$, and 
$ \vrk(M_{i+1}/M_i) \leq
tr.deg._{M_{i+1}}M_i$; so all inequalities must
be equalities.
  \qed

\<{cor} \lbl{vrk-dec}  Let $L$ be a  transformal valued field of transformal dimension $1$ over an inversive
 difference field $F$, with value group $\Q-si$.   Let $K$ be a nontrivially valued subfield of $L$.    
Let $K'$ be an extension of $K $ within some $\omega$-increasing transformal valued field   containing $L$;
assume $tr. deg._K K' < \infty$.   
Then $\vrk (L/K') \leq \vrk(L/K)$.     \>{cor}  

\proof  Note that taking algebraic closure or completion does
not change $\vrk$.  We will thus  assume all these fields are algebraically closed; and we will prove a more general statement, 
allowing $L$ to be the completion of a transformal valued field of transformal dimension $1$.

  All residue fields are therefore also algebraically
closed and, being extensions of $K_{\res}$ of transformal
transcendence degree $0$, are inversive.  Note 
(say by \ref{char-comp}) that   $L' = \widehat{K'L}$
is also  the completion of a transformal valued field of transformal dimension $1$.

If $K \leq M \leq  \widehat{L}$, then $\vrk(L/K) = \vrk(L/M) + \vrk(M/K)$, and $\vrk(L/K') = \vrk(L/MK') + \vrk(M/K')$;
so the lemma for $M/K ; K'$ and for $L/M; MK'$,   implies the statement for $L/K; K'$.  By \ref{vrk}, we may thus
assume one of the following cases holds:  \begin{enumerate}
  \item $L \leq \widehat{K}$
  \item  $tr.deg._K (L) = \vrk (L/K)$

\end{enumerate}

In case (1), $LK' \leq \widehat{K'}$, so $\vrk(L/K') = \vrk(\widehat{K'L} / \widehat{K'}) = 0$.  

In case (2), $\vrk(L/K') \leq tr.deg._{K'}(K'L) \leq tr.deg. _K (L) = \vrk (L/K)$.

\qed

\<{example}  In \ref{vrk-dec}, the hypothesis that the value group of $L$ is   $\Qq_\si$ is necessary.  \>{example}

Let $F$ be an inversive difference field,   $K$
the inversive hull of $F(t)_\si$.  Consider also the field of generalized
power series $F((t^{G}))$, where the coefficient group $G$
is $\union _n \si^{-n} \Qq_\si \subset \Qq[\si,\si^{-1}]$ (ordered
naturally.)  Pick any nonzero $a \in F$
 
Let 
$$c = a^{\si^{-1}} + a^{\si^{-2}} t^{\si^{-1}} +  a^{\si^{-3}}
 t^{\si^{-1}+ \si^{-2}}+ 
a^{\si^{-4}}  t^{\si^{-1}+ \si^{-2}+ \si^{-3}  }     + \cdots  \in F((t^{G}))$$ 
So  $K' = K(c)$ is an immediate extension of $K$, and we have $\si(c)-tc  = a$

Let $L=K(b)$, where $b$ is a solution of the equation
 $\si(x)-tx  = a$, generic   over $F((t^{G}))$.

Then $\vrk(L/K) = \vrk(K'/K) = 0$.  But $\vrk(L/K') = 1$, 
since $d=b-c$ is a nonzero solution of $\si(x)-tx = 0$, and
hence $\val(d) = (\si-1)^{-1}\val(t) \notin \Qq_\si $.

Recall that two valued field extensions $L_1,L_2$ of an algebraically closed valued field $K$ are {\em almost orthogonal} 
if, up to isomorphism, there exists a unique pair of embeddings $f_i: L_i \to M$ into a valued field 
$M$  with $f_1|K=f_2|K$ and $f_1L_1 \union f_2L_2$
generating $M$.  

\<{prop}\lbl{domin} Let $K=K^a$ be a  transformal valued field of  transformal dimension $1$ over $F$,
with value group $\Qq_\si$.  Let 
  $L$ be a   transformal valued field extension  transformal dimension $0$ over $K$.
 Let $\bar{K}$ be an  immediate extension of $K$.  Then $L,\bar{K}$ (and even  $\widehat{L},\bar{K}$) 
 are almost orthogonal over $K$ as valued fields.   
\>{prop}


\proof   The value group $\G(K)=\G(\bar{K})$
is cofinal in $\G(L)$, since if $a \in L$ has valuation $> \G(K)$ then it must be transformally transcendental 
over $K$.

It follows in particular that  the completion $\widehat{L}$ embeds uniquely into $\widehat{M}$.  From this,
once we show that $L,\bar{K}$ are almost orthogonal over $K$,  almost orthogonality of  $ \widehat{L},\bar{K}$ over $K$ follows.  


By \ref{vrk}, there exists a tower 
$K= K_0 \leq \cdots \leq K_6 = \widehat{L}$  of difference field extensions of $K$, 
with $K_{i+1}$ purely ramified or purely inertial over $K_i$, or contained in $\widehat{K_i}$.  (Namely $K \leq K_r \leq K_1 \leq K_2 \leq \widehat{K_2} \leq M \leq \widehat{L}$; with $M$ as in \ref{vrk} (4).) 
The proposition 
follows from  \ref{domin-val}.

{\bf Discussion}   

This can be formulated in terms of domination for  types of $ACVF_K$, or quantifier-free types of the language
of valued fields over $K$;  as in the Remark following \lemref{domin-val}.  
  
This implies the existence of 
canonical base change of $L/K$ to  extensions $K'$ of $K$ as transformally valued fields, relative to a base change in $\res(L)$ 
 and
in $\val(L)$ making $\val(L),\val(K')$ independent over $\val(K)$; in particular this applies when $\val(K') = \val(K)$,
more generally when $(\val(K') \meet \Qq(\si) \val(K))$ and $(\val(L) \meet \Qq(\si) \val(K))$ intersect 
in $\val(K)$ within $\Qq(\si)\val(K)$.      But note that in the example above, restricted to 
finitely generated approximations, this condition will not hold.  In any case, this material will not be used in the present paper.

\end{paragraph}
\end{subsection}

\<{subsection}{Finite generation of the value group}

In this subsection, $L$ is an  $\omega$-increasing transformal valued field of 
transformal dimension $1$ over $F$.   We further assume $L $ is a finitely generated difference field extension of $F$.

The results of the present subsection will not be used in the present paper. 

A valued field $K$ is said to be {\em Abhyankar} over a trivially valued subfield $F$ if the transcendence degree of $K$ equals the sum of the transcendence degree $t$ of
the residue field extension, and the dimension $\dim_{ram}(K/F)$ over $\Qq$ of $\Qq \tensor \G(K) $.  It is {\em strongly} Abhyankar if 
$tr.deg._F(K) = t + rk_{ram}(K/F)$; equivalently, $K$ is Abhyankar, and in addition the value group $\G(K)$ contains a strictly increasing chain of convex subgroups,
with quotients embeddable in $\Qq$>  
   A subfield of a (strongly) Abhyankar valued field is always (strongly) Abhyankar over $F$.  Also if $KF'$ is (strongly) Abhyankar over $F'$,
where $F'/F$ is a trivially valued field extension, then $K/F$ is (strongly) Abhyankar.  A finitely generated Abhyankar extension has finitely generated value group and residue field.

\<{lem}\lbl{abh}   Every finitely generated subfield of $L$ is strongly Abhyankar over $F$.  \>{lem} 
 
\proof  By \lemref{vt1}, the value group of $L$ is contained in $\Qq_\si$.  Then condition (1) of  \lemref{char-comp} holds for $\widehat{L^a}$; hence by (3) of that lemma,
$L \leq \widehat{F'(t)_\si^a}$ for an appropriate trivially valued field $F'$.   By (2) of the same lemma, we may take $t \in L$.   By the remarks above on Abhyankar extensions, it suffices to show (re-denoting $F'$ by $F$):

\claim{}  Let ${K} \leq \widehat{F(t)_\si^a}$ be a finitely generated extension field of $F$.  Then ${K}$ is an 
Abhyankar extension of $F$.  

\prf  Let $t_n = \si^n(t)$.
The value group $v({K}^*)$ of ${K}$ is contained in a finitely generated $\Qq$-subspace of $\Qq[\si]$.  
Hence it is contained, for some $m$, in the convex hull $C$ of $\Qq v(t_n)$.  There is a natural 
place on $\widehat{F(t)_\si^a}$, with residue field $F(t_0,\ldots,t_n)^a$, and value group $\Qq[\si] / C$.
The residue map $r$ of this place annihilates only elements $x$ with $v(x)>C$; thus it is injective
on ${K}$; moreover it defines a valued field isomorphism on ${K}$ into its image.  As every subfield of
$F(t_0,\ldots,t_n)^a$ is Abhyankar over $F$, the claim and the lemma are proved.  

 \eprf
%
 \qed

\<{lem}\lbl{vt1fg}  Let  $L$ be an $\omega$-increasing transformal valued field of 
transformal dimension $1$ over $F$.   Assume    $L $ is a finitely generated difference field extension of $F$.   Then $\val(L)$ is isomorphic to a
$\Zz[\si]$-submodule of $\Z-si$. 
\>{lem}
\proof  
Let  
$L_0$ be a subfield of $L$,
finitely generated over $F$ as a field,   generating $L$ as a difference field.  
Let    $L_n$ be the subfield of $L$ generated by $ \union_{k \leq n   } \si^k(L_0)$.
Then $[L_{n+1}: L_n] \leq [L_1:L_0] < \infty$.   Let $A_n = \val L_n$.
Then $A_0 \subset A_1 \subset \ldots \subset \val L^a \iso \Qq[\si]$, 
$A_0$ is finitely generated (using \lemref{abh}), and $ A_{n+1} / A_n$ is bounded.  By  \lemref{vt1} and  lemma
\ref{vt1-0}  below, $\val L$ is a finitely generated $\Zz[\si]$-module.   And
by \ref{vt1-2}, every finitely generated $\Zz[\si]$ submodule of $\Qq[\si]$ is contained in a free $\Zz[\si]$-module of rank one. 
\qed

\<{lem} \lbl{vt1-0}  Let $B_0 \subseteq B_1 \subseteq \ldots $ 
  be   finitely generated subgroups of
  $\Qq[T]$ such that   $TB_n \subseteq B_{n+1}$, with $B_{n+1} / B_n  $ finite
and bounded.  Then $B = \union_n  B_n$ is a finitely generated $\Zz[T]$-module.
 
\>{lem}

\proof    Let $A_i= (B_i + \Zz[T]) / \Zz[T]$, $A= \union_i A_i$.  We will show $A_i$ is finitely generated as
a  $\Zz[T]$-module.  It follows that $B+ \Zz[T]$ is a finitely generated $\Zz[T]$-module; and by Noetherianity of $\Zz[T]$,
so is $B$.    Note that $A_{n+1} / A_n$ is finite and bounded, and $TA_n \subseteq A_{n+1}$.  

We use induction on $m = \limsup _n \abs { A_{n+1} /  A_n  }$.   If $A \neq 0$,
there exists  $c \in A$, $ c \neq 0$, $lc = 0$, $l$ prime.  Say $c \in A_{n_0}$.

\claim Let $E=E_l = \{x: l x=0 \}$.  Then for $n \geq n_0$, $E \meet A_{n+1} \not \subseteq A_n$.

Suppose the claim is false; then any $x \in E \meet A_n$ satisfies $Tx \in E \meet A_{n+1} = E \meet A_n$.  Thus $\Zz[T] c \subset A_n$.  But $A_n$ is a finitely
generated group, while $\Zz[T] c $ is not:  $c,Tc,\ldots$ are
$\Zz/l \Zz$ - linearly independent.  A contradiction.

Thus the natural surjective map $e_l: A_{n+1}/A_n \to lA_{n+1} / lA_n$, $e_l(x)=lx$,
 is not 
injective, for $n \geq n_0$.  So $\limsup _n \abs { lA_{n+1} /  lA_n  } < m$. 
By induction,  $lA$ is a finitely generated $\Zz[T]$-module. 

Now $E_l  / ( \Zz[T] c)$ is finite; indeed $E_l \iso ( \Zz/l \Zz)[T]$, so $ E_l  / ( \Zz[T] c) \iso ( \Zz/l \Zz)[T] / f$
for some nonzero polynomial $f(T)$.  Thus $(A \meet E_l) / (\Zz[T]c)$ is 
finite.  Since $\Zz[T]c$ is a finitely generated $\Zz[T]$-module, so is
$A \meet E_l = \ker (x \mapsto lx)$.  We saw $lA$ was finitely generated; hence
so is $A$.  \qed

\<{lem}\lbl{vt1-2}  Let $M$ be a finitely generated $\Zz[\si]$-submodule of $\Q-si$.  Let 
$\tM$ be the union of all submodules $N$ of $\Q-si$ with $N/M$ finite.  Then $\tM$ is
1-generated. \>{lem} 

\proof    If $M$ is generated by $a_1,\ldots,a_n$, $m_ia_i \in \Zz_\si$ for some
$0< m_i \in \Nn$; so each $a_i \in \Zz_\si[1/m]$, where $m = \Pi_i m_i$.
  Thus   $M \subset N$ for some principal $N   \leq \Q-si$.

\claim  For any  ideal $J$ of $\Zz[\si]$, there exists a principal ideal $J'$  containing $J$ with $J'/J$ finite.

\proof  $J$ is $k$-generated for some finite $k$.  The claim reduces inductively to the case $k=2$.
So say $J= \Zz[\si](a,b)$.  Since $\Zz[\si]$ is a  unique factorization domain,  we
 may write $a=a'c,b=b'c$ with no prime element of $\Zz[\si]$
dividing both $a'$ and $b'$.  As $\Zz[\si]$ has
Krull dimension $2$, $\Zz[\si]/(a',b')$ must be finite.  

Let $J' = \Zz[\si]c$.  Then $a',b' $ annihilate $J'/J$, so that  $J'/J$ is a quotient
of $\Zz[\si]/(a',b')$, and hence is also  finite.

Actually, $J'$ is unique:

\claim If $N \leq \Q-si$ is principal, then $\Q-si $ has no  $\Zz[\si]$-submodules $N'$ such that 
$N \subseteq N'$, and $N'/N$ is finite and nonzero.

 For otherwise it would have a finitely
generated one, so again one contained in a principal $\Zz[\si]$-module.  Thus it suffices
to show that if $J$ is a submodule of $\Zz[\si]$ containing a principal module $K$, then $J/K$ 
is zero or infinite.  If $J/K$ is nonzero and finite, so is $J'/K$, where $J'/J$ is finite and $J'$ is principal.
So $\Zz[\si]/f\Zz[\si]$ is finite nonzero; but this is clearly absurd.  

Thus $\tM \subset N$.  Let $I = \{r \in \Zz[\si]: rN \subset \tM\}$.  So $\tM = IN$. By 
definition of $\tM$, if $I'/I$ is finite then $I'=I$.  Thus by the Claim $I$ is principal.  So   $\tM$ is 1-generated.

\qed

\<{example} $\val(L)$ need not itself be free.  \>{example}  
Take $F(t^\si,t^2)_\si \leq F(t)_\si$;
the value group is $M = \Zz[2v  , \si(v)] \subset \Zz[v,\si(v),\ldots]$.  We have 
$M/  \si(M) \iso   \Zz \oplus (\Zz/2\Zz)[\si]$. \qed

At all events, when $L$ is finitely generated, of limit degree $d$, $\Zz[1/d!] \val(L) $ is already free on one generator.

This suggests  that $L$ may have an algebraic extension  whose 
value group is a free $\Zz[\si]$-module of rank one.   We prove such a  statement in 
characteristic $0$; in positive characteristic $p>0$, the proof works up to localization at $p$, 
 and we did not investigate it further.

\<{prop}  Let
  $F$ be a difference field of characteristic  $0$, closed under roots.     Let   $L$ be a   
  transformal discrete valuation ring   over $F$, Henselian as a valued field.
   Then $L$   has a finite $\si$-invariant extension $M$, whose value group is
isomorphic to $\Zz[\si]$.  \>{prop}

\proof  

 By \ref{vt1-2}, $\val L \subset \tM$ for some free 1-generated $\Zz[\si]$-module $\tM$ with $\tM / \val L$
finite;
we have $\val L = I \tM$ for some ideal $I$, and necessarily $n \in I$ for some $n >0$. 

 Since $L$ is Henselian and $F$ is closed under roots, the homorphism
 $val: L^* \onto \val L$  is a surjection.  It induces a surjection 
$$(L^*)^n \onto n \val L$$
Let $H$ be the pullback of $n \tM \subset \val L$.  We have $\si(H) \subset H$; and 
  $$H/(L^*)^n \iso n \tM / n \val L \iso \tM / \val(L)$$
In particular,  $H/(L^*)^n$ is finite.    By Kummer theory, there exists a unique Galois extension
$K$ of $L$, such that $(K^*)^n \meet L = H$, and   $[K:L] = [H:  (L^*)^n] = \tM / \val L$.   Now $\si(K)L$
is also a Galois extension of $L$, generated by $n$'th roots of elements of $H$; hence it is contained
in $K$; so $K$ is a difference field.   
 Tracing back the isomorphisms we see that $\tM \subset \val(K) $, so $\tM = \val(K)$. \qed

\<{rem}\lbl{vt1c2z}  Let $M \subset \Zz-si$ be a $\Zz[\si]$-module.      Let $Y$ be a downward-closed subset of $M$.  
  Then for some nonzero convex subgroup
$S$ of $M$, and $0 \leq a \in \Rr[\si]$, 
$$Y  = \{y \in M:  (\exists c \in S) y \leq a+c \}$$
  \>{rem} 
The proof is easy and left to the reader (see \ref{vt1c2}).   
%
%
%

\>{subsection} 

\end{section}
\begin{section}{Transformal valuation rings and analyzability}
\lbl{transvalan}
\<{subsection}{The residue map on difference varieties}

We begin with  some criteria for the residue map to be total, injective and surjective.   We will not use them much, but they clarify the picture.

\<{lem}\lbl{res-dom}  Let $D \subset D'=D[a]$ be difference domains;
$a=(a_1,\ldots,a_n)$ a tuple of generators of $D'$ as a difference $D$-algebra.  
 Assume $\si(a) \in D[a]^{int}$ ( the integral closure of the domain $D[a]$.).  Then for large enough $m$, for 
any morphism $h: D' \to K$,
with $(K,R,v)$ an $m$-increasing transformal valuation field, if $h(D) \subset R$ then $h(D') \in R$.  \>{lem}

\proof  For each $i$, $\si(a_i)$ is a root of a monic  polynomial 
$f_i = \sum_{l=0}^{d_i} c_{i,l}X^l \in D[a][X]$.
The coefficients $c_{i,l}$ are themselves polynomials in $a$ 
over $D$; take $m$ bigger than the total degree of all these polynomials.
Then if $v$ is a valuation, with $v(d) \geq 0$ for $d \in D$ and $v(a_i) \geq  \delta$ for each $i$, (where $\delta < 0$), then  for each $i,l$,
$v(c_{i,l}) > m \delta$.  

  Assume  $(K,R,\val)$ is $m$-increasing, and $h(D) \subset R$.
We may replace $D,D'$ by their images under $h$, and assume 
$h$ is the inclusion.
Say $\val(ha_1) \leq \val(ha_2) \leq \ldots \leq val(ha_n)$.  We have to show that $\delta=\val(a_1) \geq 0$.    
Otherwise,  as the (monic) leading monomial of $f_i(\si(a_i))$ cannot have
valuation less than all other monomials, we have $d_1 val(\si(a_1)) \geq   \val(c_{1,l})+l val(\si(a_1)) $ for some
$l< d_1$.  Thus $(d_1-l) val(\si(a_1)) \geq \val(c_{1,l}) > m \delta$.
So  $\val \si({a_1}^{-1}) \leq m \val {a_1}^{-1}$.  This contradicts the assumption
that $\val$ is $m$-increasing, and proves the lemma.  \qed

\<{definition} \ \lbl{res-inj-not} \>{definition} 
For the sake of the   lemma \ref{res-inj}, define a {\em standard unramified map}
to be a scheme morphism 
 $\spe S \to \spe R$, where $R$ is a (not necessarily Noetherian)
commutative ring, $S=R[a_1,\ldots,a_n]$, and we have $f_i(a_1,\ldots,a_n) = 0$
for some polynomials $f_1,\ldots,f_n $ over $R$, with invertible Jacobian matrix.  
 If the $(f_i)$ give  a presentation of $S$, i.e. $S = R[X_1,\ldots,X_n] / (f_1,\ldots,f_n)$,
  the map is said to be {\em \'etale}.   
Let us say that a map $f: X \to Y$
of schemes is {\em unramified} 
if for any two points of $Y$ over a field, there exists an open subscheme $Y'$
of $Y$ containing the two points, $X' = f^{-1}(Y')$, with $X' \to Y'$ (isomorphic to)
a standard unramified map.   When $Y$
is a difference scheme, we view the map $Y \to B_1Y$ as a map of schemes, and apply the same
terminology.

\<{lem}\lbl{res-inj}  Let $(K,R,\bar{K})$ be a  strictly increasing transformal valued field. 
Let $X$ be a finitely generated quasi-projective difference  scheme over $R$, $X_0 = X \tensor_R {\bar{K}}$.    Assume the  reduction sequence map $X_0 \to  B_1X_0$ is
unramified (\ref{res-inj-not}).  Then the residue map $X(R) \to X({\bar{K}})$ is injective.
\>{lem} 

\proof    
Consider first the special case   $X = \spe R[x]/f(x)$, $f(x) = x  + g(x) + h(x)  \in R[x]_\si$; 
where  $g $ has coefficients in the maximal ideal $M$ of $R$, and $h$ is a sum of monomials $x^{\nu}$, $\nu \in \Nn[\si]$, $\nu > 1$.  Let  $a,b \in X(R)$ with $\res(a)=\res(b)$; we will
  show $a=b$.  We may assume $b=0$.  Then 
   $\res(a)=0$ so $\val(a) > 0$; but then 
   $\val (h(a)) > \val (a)$, and also $\val (g(a)) > \val (a)$, so $f(a)  = 0$ implies
$a=0$.   

In general, we may assume $X = \spe A$;  $X_0 = \spe A_0$, 
$B_1X_0 = \spe {A_0}'$, with   $A_0= A \tensor_R {\bar{K}} $, ${A_0}' = \si(A_0)$,
with   $ X_0 \to  B_1X_0$  a   standard  
unramified map (\ref{res-inj-not}):
  $A_0$ has generators $y_1,\ldots,y_n$ over $ {A_0}'$, admitting relations  $F_1,\ldots,F_n$, 
   where $(F_{i}) \in {A_0}'[y_1,\ldots,y_n] $ has 
invertible Jacobian matrix $\det (\partial_{y_i}F_j) \in G_m({A_0}')$.

Complete $y_1,\ldots, y_n$ to a system   of  generators of $A_0$ over $\bar{K}$; since
the $y_i$ already generate $A$ over $ {A_0}'$, the additional generators $y_j$ may
be chosen from $ {A_0}'$.  Each new $y_j$ solves the inhomogeneous linear polynomial
$Y_j-y_j \in {A_0}'[Y_j]$; adding these generators and relations to the system clearly leaves
the Jacobian invertible.

Now the coefficients of the $F_i$ lie in $\si(A_0)$, so they may themselves be expressed
as difference polynomials in the $y_i$, indeed in the $\si(y_i)$.  Replacing the coefficients
by these polynomials, we may assume that $F_i \in {\bar K}[Y_1,\ldots,Y_n]$ has coefficients
in $\bar{K}$, rather than in $\si(A_0)$.  If we convene that $\partial \si(Y_j) / \partial Y_k = 0$,
the Jacobian matrix remains invertible.

  Lift  the generators $y_i$ to $x_i \in A$.  Then any element of $A$ can be written as a difference 
polynomial over $\si(A)$ in the $x_i$, up to an element of $\M A = \ker A \to A \tensor_R {\bar{K}} $.
(Here $\M$ is the maximal ideal of $R$.  )
So $A$ has generators $x_1,\ldots,x_n,x_{n+1},\ldots,x_N$, where for $i>n$ the image of $x_i $ 
in $A_0$ is $y_i=0$.  We still have $A_0= {A_0}'[y_1,\ldots,y_N]/ (F_1,\ldots,F_N)$, where
$F_i = y_i$ for $i>n$.   Lift the $F_i$ to $f_i \in \si(A)[x_1,\ldots,x_N]$.  Then
in $A$ we have a relation $f_i(x_1,\ldots,x_N) = g_i(x_1,\ldots,x_N)$, where $g_i$ 
 has coefficients in $\M$.  

Consider two elements $e=(e_1,\ldots,e_N), e'=(e'_1,\ldots,e'_N)$ of $X(R)$ specializing to the same point
of $X_0$.  Replacing $x_i$ by $x_i-e'_i$,  we may assume 
$e'=0$.      Suppose $e \neq 0$; let $\min_{i} \val (e_i) =  \beta > 0$.   In these coordinates,
let $L_i$ be the linear (monomials of order $1$) part of $f_i$.  The $L_i$ are then the rows
of an invertible matrix $L$ over $R$, and $L_i x - g_i(x) + h_i(x) = 0$, where 
$h_i$   is a sum of $\si$-monomials of degree $>1$.  Since $e \in X(R)$, we 
have $Le - g(e) + h(e) = 0$.   (Where $g,h$ are the matrices of $\si$-polynomials whose
rows are $g_i,h_i$.)

Now $\min_i \val (\si(e_i)) = \si(\beta) > \beta $, so $\val h_i(e)   > \beta$.
Also $\val g_i(e) > \beta$.  Writing $e = - G^{-1} ( g(e) + h(e))$, we see that $\min_i \val (e_i) > \beta$,
a contradiction.     \qed
%


We differentiate difference polynomials using the rule $\partial ({x_j}^\si) / \partial x_i = 0$.

 \<{lem} \lbl{res-surj}      Let $(K,R,\bar{K})$ be a  strictly increasing transformal valued field,
 with $K$ maximally complete as a valued field.  
Let $X = \spe R[X]/ (F, J \inv)$, where $X=(x_1,\ldots,x_n)$, $F=(f_1,\ldots,f_n)$, where the $f_i$
are difference polynomials over $R$, and $J= \det \partial f_i / \partial x_j $.   
Then $X(R) \to X({\bar{K}})$ is surjective.
\>{lem}

\prf   Let $p \in X({\bar{K}})$; we may assume $p = (0,\ldots,0)$.  
Replacing $F$ by $J^{-1}(p)F$,   we can write:
$F(x)=  x+ \hbox{ (higher terms) + \hbox{terms with coefficients from } $M$ }$.  
We seek $e \in M^n$, $F(e)=0$.    The usual proof of 
Hensel's lemma provides such an $e$, by a transfinite sequence of successive approximations.
Maximal completeness permits jumping over limit steps.

 \eprf
 
 \<{rem}\   \lbl{res-surj-p}\>{rem} In this paper, we will only use \ref{res-surj} when 
$\si$ is the standard Frobenius, $\si(x)=x^q$ on $K$,
  $K$   a complete discrete valuation ring.   Here the convention that $\partial \si(x_j) / \partial x_i = 0$ follows
from Leibniz's rule:  $\partial ({x_j}^q) / \partial x_i = 0$ in characteristic $p$. 
  
\>{subsection}

\<{subsection}{Sets of finite total dimension are residually analyzable}

\lbl{fintdra}

\<{definition} \lbl{scatter-def}  Let $K$ be a valued field.  A subset $X \subset K$ is {\em scattered} if
$\{|x-y|: x,y \in X \}$ is finite.   $X \subset K^n$ is scattered if $pr_i X  \subset K$
is scattered for each $i$.    \>{definition}

The definition can be made more  generally an ultrametric spaces.  We will use it in
contexts where $X$ is definable in some expansion of $K$, and  one can envisage $X$ in arbitrary elementary extensions.  Then the notion appears sufficiently close to the usual topological one to 
permit our choice of the term "scattered".  

A scattered set $X$ may consist of clusters of points of distance $\alpha < 1$.
Each cluster can be enlarged, to reveal a new set of clusters separated by
the residue map.  After finitely many iterations, this process separates points of $X$.

\<{lem}\lbl{scatter-analysis}
 Assume $X$ is scattered.  Then there are a finite number of equivalence relations
$E_0 \subset E_1 \subset \ldots E_n$,  
such that $E_0 = Id$, and for each $E_{i+1}$-class $Y$ of $X$, there exists a map
$f_i^Y$ embedding $Y/E_i$ into the residue field.  The $E_i$ and $f_i^{Y}$ are 
quantifier-free definable in the language of valued fields; $f_i^{Y}$ is defined by a formula with  a parameter
depending on $Y$, while $E_i$ requires no parameters beyond those needed to define $X$. \>{lem}

\proof  First take $X \subset K$. Let $\rho_0 < \rho_1 < \ldots \rho_n$ be the possible values of  $|x-y|$ for $x,y \in X$.
Say $\rho_i = |c_i|$.  
Let $E_i(x,y) \equiv |x-y| \leq \rho_i$.  Given $Y$, pick $b \in Y$, and let 
$f_i(x) = \res {c_i}^{-1} (x-b)$.  

When $X \subset K^2$, we let $E_i(x,y) \equiv ( |pr_0(x)-pr_0(y)| \leq \rho_i \& pr_1(x)=pr_1(y)) $
for $i \leq n$, $E_{j}(x,y) \equiv |pr_1(x)-pr_1(y)| \leq \rho_{j-n}$ for $j>n$; etc.  
\qed

\<{prop} \lbl{scatter-1} Let $K$ be a transformal valued field, with $\Gamma$ a torsion-free
$\Zz[\si]$-module.  Let $F \in K[x]_\si$ be a nonzero transformal polynomial.  Then 
$X = \{x: F(x)=0\}$  is scattered.  
\footnote{In fact if $F(a)=F(b)=0$ where $a \neq b$ lie in any transformal valued field extension $L$ of $F$, 
then $\lambda \val (a-b) \in \val(K)$ for some $0 \neq \lambda \in \Nn[\si]$.  }
 More generally, any $X \subset K^n$ of finite total dimension is
scattered.  \>{prop}

\proof  This reduces to $X \subset K$, using projections.   When  $X \subset K$
has finite total dimension, there exists a nonzero difference polynomial $F$ over $K$ with
$X \subset   \{x: F(x)=0\}$.  Recall the transformal derivatives $F_\nu$.
We have $F(x+y) = \sum_\nu F_\nu(x) y^\nu$.  So if $a,a+b \in X$, 
then $b$ is a root of $\sum_{\nu >0} F_\nu(a) Y^\nu = 0$.  If   $F_\nu(a)=0$ for all $\nu > 0$,
then $F$ is constant, so $X = \emptyset$.  Otherwise, either $b=0$ or 
$\val F_\nu(a) b^\nu = \val F_\mu(a) b^\mu $ for some $\nu < \mu$;
so $ (\mu - \nu)  \val (b) = \val F_\mu(a) -  \val F_\nu(a)$.   By the Claim below,
$\{ \val F_\nu (a) : a \in X \}$ is finite for each $\nu$; so $\{ \val(b): b \neq 0, a,a+b \in X \}$ is 
also finite, and contained in the $\Zz[\si]$-division points of $\G(K)$.  

\claim Let $G \in K[x]_\si$.  Then $\{\val G(a): \ F(a)=0 \}$ is finite.

It suffices to see that $\{\val G(a): \ F(a)=0 \}$ remains bounded in any  elementary extension $L$ of $K$.
In fact if $F(a)=0$ then $\lambda \val (a) \in \val(K)$ for some $0 \neq \lambda \in \Nn[\si]$.  
Let $a \in L$, $F(a)=0$, $c = G(a)$.  Then $K(a)_\si$ has transcendence degree 
$\leq \deg_\si F$ over $K$; hence so does the subfield $K(c)_\si$.  So $H(c)=0$ for some
nonzero $H \in K[x]_\si$, $H= \sum_\nu d_\nu x^\nu$.
  Arguing as above, we see that $(\mu-\nu) \val(c) = \val d_\nu - \val d_\mu$ for some $\nu < \mu$.   \qed 

\<{paragraph}{Remarks on the liaison groups}

 In most  run of the mill cases, if $S$ is the set of zeroes of a difference polynomial $F$ within some 
ball $B$,  some transformal derivative $F'$ of $F$ will   have a unique root  in $B$, and applying
$F'$ will map $S$ into a ball of known radius around $0$; so a direct coordinatization will be
possible (internalization without additional parameters.)  This process will fail in those some
cases where a transformal derivative is constant, but non zero.  The polynomial 
$F(X) = X^\si +X -a$ is an example.  It seems likely that the examples can all be shown to have
an Artin-Schreier aspect, and to become generalized Artin-Schreier extensions upon application of 
  $M_q$.  

Another approach, that does not explicitly look at the form of the polynomial, is
in \cite{HB}.  
The associated groups  of automorphisms of clusters over the residue field
are   all subgroups of the additive group; this has to do with "higher 
ramification groups" (cf. \cite{serre-local-fields}).

  We will not go in these directions here.   However, in the next subsection, we will explain
  how to assign a dimension to scattered sets, given a notion of dimension over the
  residue field.  (We will apply this to the total dimension over the residue field, obtaining something
finer than total dimension of the closure over the valued field.)   It will be convenient to explain the
way that this dimension is induced in a more abstract setting.  Our only application however  will be the above mentioned
one, and  the reader is welcome to use  \ref{val-case},\ref{R-descent-ex} as a dictionary.  
\>{paragraph}

\>{subsection} %
\<{subsection}{Analyzability and inertial dimension}
\lbl{anresdim}

Let $L$ be a language  with a distinguished sort $V$. 
$x,y,v$ will denote   tuples of variables; and we will 
write $a \in M$ to mean:  $a$ is a tuple of elements of $M$.

{\bf Variables $v=(v_1,\ldots,v_m)$ will be reserved
 for   elements of $V$. }  We use the notation $\phi(x;y)$  when  we have in mind 
the formulas $\phi(x;b)$ with $b \in M \models  T$.    ( This corresponds to
the use of relative language, for schemes over a given scheme, in algebraic
geometry.) 

{\bf Data}

A theory $T$  (not necessarily complete).  
 A set $\Phi $   of  quantifier-free $L$-formulas and  a set 
$\fn $   of basic functions of $L$.   We assume $\Phi$ is closed under conjunctions and 
substitutions of  functions from $\fn$.   The set of formulas of $\Phi$ in the  variables $x$
is denoted $\Phi(x)$; similarly  $\fn (x,t)$   refers to domain variables $x$
and range variable $t$.  We allow partial functions, whose domain is 
given by some $P \in \Phi(x)$.   $\Phi(x;y)$ denotes formulas in $\Phi(xy)$,
together with a partition of the variables, as indicated.

A map $d_V: \Phi(v;y) \to \Nn$.

\<{example}\lbl{val-case}  \rm 
 
 $L=$  the language of transformal valued rings, over   $M_F(t)_\si^a$, with a distinguished
sort $V=V_{\res}$ for the residue field.  Here $F$ is a trivially valued inversive difference field,
and $\val(t)>0$.  

 $T=$    the theory of transformally Henselian, $\omega$-increasing transformal valued fields. 
 If we wish to view the theory as two-sorted,  with a sort for the residue field $V$,  we include the surjectivity 
of $\res$ in the axioms.  
\footnote{Thanks to Tom Scanlon and Immi Halupczok for pointing out the need for care in the passage to a two-sorted viewpoint.
The  Henselianity axioms are considerable overkill here, but will be convenient in the proof of Proposition \ref{vrk-vdim} as they imply that surjectivity continues to hold when restricted
to finite transformal dimension.}


 $\Phi(x,v)  = $ all quantifier-free formulas $\phi(x,v )$ implying an ACFA quantifier-free formula of 
  transformal dimension $0$ over $M$, in both $x$ and $v$ variables.

$\fn$ includes  polynomials, and  maps of the form $x \to \res( (x-y_1)/(y_1-y_2))$;
with domain $\{(x,y_1,y_2): y_1 \neq y_2, \val(x-y_1)=\val(y_1-y_2)\}$. 

$d_V(\phi(v;y)) \leq n$ iff for any $M \leq M' \models T$, and $b \in M'$,  $\phi(v;b)$ implies $v \in X$ for some difference scheme $X$ over $M'_{\res}$ of total dimension $\leq n$.   

Later, we will use a rank $Rk(L/K)$ on substructures; it corresponds to $\vrk(L/K)$.  

In this case, the $V$-dimension $\dimV$ defined below will be referred to as the {\em inertial dimension}.
\end{example}

We define $V$-co-analyzability (relative to $T,\Phi $) and  the $V$-dimension of $P \in \Phi$
(relative also to $d_V$),   using recursion on $h \in (1/2)\Nn$.   
(Compare \cite{HHeM}.)
The second half-step
between any two integers  does
not directly relate to $x$, but serves to provide additional parameters for future steps.  
    Recall that the variables $v$ refer to $n$-tuples from $V$.  Since $V$ is fixed, we will
    simply refer to co-analyzability (but will continue to talk of $V$-dimension.)
 
\<{definition}  \lbl{co-an}    
\<{enumerate} 
\item
$P(x;u)$ is $0$-step coanalyzable over $V$ (with $V$-dimension $0$, i.e. $\leq n$ for all $n$)
if  $T \models P(x,u) \wedge P(x',u).  \implies x=x'$.

\item
$P(x;u)$ is  co-analyzable in $h+1/2$ steps    (with $V$-dimension $\leq n$)
 if there exist $Q \in \Phi(x; u,v)$,  co-analyzable in $h$ steps (with $V$-dimension $\leq n_1$ ),
   $R \in \Phi(v;u)$  (with $d_V(R) \leq n_2$), and $g \in \fn (x,u;v)$, such that 
$$T \models  P(x;u)  \implies (Q (x,u,g(x,u))  \wedge  R (g(x,u);u) ) $$ 
(and $n_1+n_2 \leq n$ ).

\item
$P(x;u)$ is  co-analyzable in $h+1 $ steps     (with $V$-dimension $\leq n$)  if
there are  finitely many  
$Q_j \in \Phi(y;u) $  
such that $T \models P(x;u)   \implies \bigvee_j (\exists y) Q_j(y,u)$, and for each $j$,
 $$ (P(x;u)   \wedge  Q_j(y,u)) \in \Phi(x;uy)$$
  is
 coanalyzable in $h+1/2$ steps    (with $V$-dimension $\leq n$).    
\end{enumerate}
\>{definition}

We say $P$ is ($V$-)co-analyzable if it is so in some finite number of steps.  We will
write
$$ \dimV (P) \leq r$$
for:  `` $P$ is ($V$-)co-analyzable, of $V$-dimension $\leq r$. ".  Clearly if $T \models P \implies P'$ and $\dimV (P') \leq r$ then $\dimV (P) \leq r$.

In applications, it is convenient to apply this terminology to $\infty$-definable $P = \wedge_{i \in I} P_i$.  An $\infty$-definable $P(x,u)$ is given, by
definition, by a collection $\{P_i(x;u_i): i \in I \}$, $P_i \in \Phi(x;u_i)$.
 $u$ may be an infinite
list, containing all the finite lists $u_i$.    .  We have in mind
$P(M) =_{def} \meet_{i \in I} P_i(M)$.   We simply define:
$\dimV (P) = \inf  \{ \dimV ( \wedge_{i \in I'} P_i): 
                                    I' \subset I, I' \hbox{ finite.}  \}  $.  

When $K \leq L \leq M \models T$, and $a \in L$, we will
write $\dimV(a/K) \leq r$ if   
there exists $P \in \Phi(x;y)$, $\dimV(P) \leq r$,   and $ b \in K$, such that $M \models P(a,b)$.

 \<{example}     
 Assume $P(x;y)$  is $V$-co-analyzable in $1$ step.  \rm  In this case,$P(x;y)$ is  {\em $V$-internal}, i.e.
if  $b \in M \models T$, there exists 
an $M$-definable   injective map $f: P(x,b) \to V^{eq}$.
 All such injective maps have the same image,
up to an $M$-definable bijection; if $V$ is stably embedded,
and the given dimension on $V$ is an invariant of definable bijections, then
the $V$-dimension of $P$ equals the dimension of the image  of any of these maps $f$. \>{example}

See \cite{HB} (appendix B)  for a general treatment of internality, and associated
definable groups.

\<{example}   Let $P(x) \in \Phi(x)$, 
$E_0 \subset \ldots \subset E_h \in \Phi(x,x')$ equivalence relations
on $P(x)$,   such that for each $i$,  the set of $E_{i+1}$-classes 
is internal relative to the $E_i$-classes.  \rm (I.e. for  every
class $Y$ of $E_i$,  
there exists a function $f_Y = f(y,b):  Y \to V^{eq}$, $f \in \fn$,  
parameters $b$ depending on $Y$, 
such that $f_Y$ is injective modulo $E_{i-1}$.)  
  In this situation $P$ is said to be {\em $V$-analyzable.}      
This is a stronger notion than 
co-analyzability.  If the $V$-dimension of $Y/E_{i+1}$ is $n_i$ for each $E_i$-class $Y$,
then $\sum n_i$ is an obvious upper bound for the $V$-dimension of $P$.  
\>{example}

\<{remark}  In the case of valued fields, \ref{val-case}, if $P$ is scattered, 
then it is in fact analyzable.
\>{remark}

Consider   structures  $K  \models T_{\forall}$; (more precisely, subsets $K$ of models $M$ of $T$, closed under
$\fn$, with the formulas in $\Phi$ interpreted according to $M$.) 
   We will consider pairs $K \leq L$ 
with $L$   generated over $K$ by a tuple $c$;  write $L=K(c)$ in this case.  Let
$$tp_\Phi (c/K) = \{\phi(x;b)  : \phi(x;y) \in \Phi ,  b \in K,  L \models \phi(c,b) \}$$
 
Assume given an  $\Nn \union {\infty}$-valued function  $Rk$  on such $\Phi$-types.
We will assume that the $Rk$ does not depend on the choice of a generator $c$ of $L/K$,
and write $Rk(L/K)$ for $Rk(c/K)$ when $L = K(c)$.

\<{lem}\lbl{str-form}
  Assume:
\begin{enumerate}
  \item If $Rk(L/K ) \leq n$, $a \in V(L)$, then for some $\phi(v;u) \in \Phi(v;u)$ 
       with $d_V(\phi) \leq n$, and $b \in K$, $M \models \phi(a,b)$.

  \item If $K \leq K' \leq L \leq M$, $Rk(L/K') + Rk(K'/K) = Rk(L/K)$.

  \item  If   $K \leq K'  \leq M$, $c \in M$, $Rk(K(c)/K)< \infty$,  then $Rk(K'(c)/K') \leq Rk(K(c)/K)$. 
\end{enumerate}

 Let $K \leq L  \models T_{\forall}$, with $Rk(L/K) \leq n$.   
Let $a \in L$, $a/K$  co-analyzable.  Then  $\dimV(a/K) \leq n$.                  \>{lem}

\proof  We use induction on the number of  steps of co-analyzability (the case of $0$ steps being clear.)
Suppose $a/K$ is  co-analyzable in $h+1/2$ steps, i.e. $L \models P(a,b)$, $b \in K$, 
$P$  co-analyzable in $h+1/2$ steps.    Then there exists $Q \in \Phi(x;u,v)$,   co-analyzable in $h$ steps,  $g \in \fn(x,u;v)$, and 
$d=g(a,b)$, such that $M \models Q(a;b,d)$.       
  By induction, there exists
$Q' \in \Phi(x;y',v)$  and $b' \in K$  with $M \models Q(a,b',d)$, $\dimV(Q') \leq Rk(L/K(d))$.  (Any parameters from $K(d)$ can be written as terms in elements $b'$ of $K$, and $d$.)

By (1), there exists $R'(v;u) \in \Phi(v;u)$, $d_V(R') \leq Rk(K(d)/K)$, 
and $b'' \in K$, with  $R'(d,b'')$.

Let $P'(x;y,y',u) = Q' (x,y',g(x,y))  \wedge  R (g(x,y);y'') $.  Then $P'(a;b,b',b'')$ shows 
that the $V$-dimension of $a/K $ is $\leq Rk(L/K(d)) + Rk(K(d)/K)    = Rk(L/K)$.
 
Now suppose  $a/K$ is  co-analyzable in $h+1$ steps.  Then $L \models P(a,b)$, $b \in K$, 
 $T \models P(x;u)   \implies \bigvee_j (\exists y) Q_j(y,u)$, and for each $j$, (*)
 $ P(x;u)   \wedge  Q_j(y,u)$ is
 coanalyzable in $h+1/2$ steps.  Let $r= Rk(L/K)$, 
$$\Xi = \{S \in \Phi(x;u,u',y): \, \dimV(S) \leq  r \}$$
and let $\Xi' = \{ \neg S(x,b,b',y): S \in \Xi , b' \in K\}$.  

{\bf Claim}  $T_{\forall} \union tp _\Phi (a/K) \union \Xi' \union \{Q_j(y,b) \}$ is inconsistent.

\proof  Suppose otherwise.  Then in some  $M \models T_{\forall}$, $K \leq M$, 
we can find $a' \models  tp _\Phi (a/K)$, and $d$ with $Q_j(d,b)$, such that (**)
for any $S \in \Xi$, $M \models \neg S (a,b,b',d)$.  Let $K'=K(d)$, $L' = K(a',d)$.
   by (3), $Rk(L'/K') \leq r$.  By (*), $a'/K'$ is   
 coanalyzable in $h+1/2$ steps.  Thus by induction, $\dimV (a'/K')  \leq r$.  So there exists 
  $S \in tp_\Phi (a'/K')$  with $\dimV (S)  \leq r$.  We can take $S = S(x,b,b',d)$,
$S \in \Phi (x, u,u',y)$.  So $S \in \Xi$.   But this contradicts (**).  \qed

By compactness, for some finite disjunction $\bigvee_{j'} S_{jj'}(x,u,u',y)$ 
(with $S_{jj'} \in \Xi)$ and some $P' (x,b,b') \in tp_\Phi(a/K)$ implying $P(x,b)$,

$$T_{\forall} \models   P' (x,u,u') \wedge Q_j(y,u). \implies \bigvee_{j'} S_{jj'}(x,u,u',y)$$
  
so $\dimV (P' (x,u,u') \wedge Q_j(y,u)) \leq r$.  By \ref{co-an} (3), $\dimV (P')  \leq r$.  
\qed

%
%
%
%

Specializing this to our example  \ref{val-case}, we have:

\<{prop}\lbl{vrk-vdim}  Let $M$ be an algebraically closed $\omega$-increasing
  transformal valued fields of transformal dimension $1$ over an inversive
 difference field, with value group $\Q-si$.  Let $\phi(x) \in \Phi(x) $ be
a quantifier-free formula  in the language of transformal valued fields
over $M$, cf. \ref{val-case}.   Assume $\phi$- is $V_\res$-analyzable,
and: for any
 $\omega$-increasing
  transformal valued field extension $L=M(c)$ of $M$, 
with $\phi(c)$, $\vrk(L/M) \leq n$.  Then $\phi$   has inertial dimension $\leq n$.  
  \>{prop}

\proof   Note that the proof of  \ref{str-form} uses   $Rk(K'/K'')$  only for   finitely generated extensions $K',K''$  of $K$ of transformal dimension $0$ over $K$.  For such
extensions, take $Rk(K'/K'') = \vrk(K'/K'')$.  Then hypothesis (1) of \ref{str-form} is clear since if $\vrk(K'/K'') \leq n$ then $tr. deg. K'_{\res} / K''_{\res} \leq n$.  (2) follows from additivity of transcendence degree and
vector space dimension.  And the truth of (3) is the content of \ref{vrk-dec}.

Now if $K=(K_{field},K_{res})$ is a two-sorted substructure of some model $M$ of $T$, 
and   $c \in \res(K)$ has  finite transformal dimension over $\res(K_{field})$,  
there exists   $c' \in M$  of  finite transformal dimension over $K$ with $\res(c')=c$;  this follows from \lemref{res-lift-lem-0}.  
  
  Hence by compactness we can find $P(x,y)$ such that (i) $P$ implies $\phi(x,\res(y))$; (ii) If $\phi(a,b')$ then in some extension 
  we can find $b$ with $P(a,b)$ and $\res(b)=b'$.  By definition of $\dimV$, we have $\dimV(\phi) \leq \dimV(P)$.  
  
  For any $M \models T$, if $a,b \in M$, $M \models P(a,b)$, letting $K$ be the 
substructure of $M$ generated by $b$, we see that there exists $Q \in \Phi(x;y)$,
$M \models  Q(a,b)$, $\dimV(Q) \leq r$.  
By compactness, $T \models P \implies \bigvee_j Q_j$, with $\dimV(Q_j) \leq r$.  
It follows that $\dimV(P) \leq r$.  Hence $\dimV(\phi) \leq r$ also.

\qed

\<{remark}  There is a canonical function $Rk$ satisfying \ref{str-form} (1-3). \rm
To define it, let \begin{itemize}
  \item $Rk_0(L/K) = \sup \{ d_V(a/K): a \in V(L) \}$
  \item $Rk_{n+1/2}(L/K) = \sup \{ Rk_n(L/K') + Rk_n(K'/K): K \leq K' \leq L \}$
  \item $Rk_{n+1}(L/K) = \sup \{ Rk_n(LK'/K'): K \leq K',L \leq M \}$
  \item  $ Rk_\infty (L/K) = \sup _n Rk_n(L/K)$ 
\end{itemize}
  \>{remark}

\<{remark}\lbl{R-descent}  Let $P(x;y) \in \Phi(x;y)$ be $V$ co-analyzable in $h$ steps, with $V$-dimension $\leq n$, with respect to $T,\Phi$.
\begin{enumerate}
  \item For some finite $T_0 \subset T$, $\Phi_0 \subset \Phi$, $\fn _0 \subset \fn$,
     $P$ has $V$-dimension $\leq N$ with respect to $T_0, \Phi_0, \fn _0$.
     
   \item  Let $M_q \models T_0$ be a family of models of $T_0$, indexed by an infinite
   set of integers $\{q\}$,  and suppose that for any $P(v,y) \in \Phi_0$ , for some $\beta$,
   for all $q$ and all $b \in M_q$,  $|P(M,b)| \leq \beta q^{d(P)}$. 
   
    Then for any $P(x,y) \in   \Phi_0$ of $V$-dimension $\leq n$, for some $\beta$,
   for all $q$ and all $b \in M_q$,  $|P(M,b)| \leq \beta q^{n}$.
  
\end{enumerate}

\>{remark}

\proof  (1) is obvious from the definition; (2) also follows immediately from the definition, by induction
on the number of steps.  (For each $n$, the half-step from $n$ to $n+1/2$ increases the $V$-dimension
and the exponent; the half-step from $n+1/2$ to $n+1$ increases only the constant $\beta$.)

\<{example}\lbl{R-descent-ex} \rm (continuing \ref{val-case}).  Here
$T_0$ will be, for some $k$, the theory of $k$-increasing transformal
valued fields.  The family of models $M_q$ can be taken to be 
the  fields $K_q(t)^{alg}$, endowed  with a nontrivial  valuation
over $K_q$, and with the $q$-Frobenius automorphism.  The validity
of the assumption will be seen in \ref{si-dim}.  \>{example}
%
%
%
\>{subsection}

\<{subsection}{Direct generation}  \lbl{dirgen}

 Let $F$ be an inversive difference field, $K$ a
 transformal discrete valuation ring
 over $F$.  Then  $K$ is generated as a difference field by a subfield $K(0)$ with  $tr.deg._F K(0) < \infty$, and (letting $K(n)$ be the field generated
by $K(0) \union \ldots \union \si^n K(0)$), $tr. deg._{K(0)}  K(1) = 1$.  We wish to
show that this automatically implies:  $K_{\res} \subset (K(0)_{\res})^a$.  We phrase it a little more generally.

 \<{lem}\lbl{val2.1} Let $F$ be an inversive difference field, and let $(K,R,M,\bar{K})$ be a weakly transformal valued field
extending $(F(t)_\si,\Ft,t \Ft,F)$.  Let $K(0) \subset K(1) \subset \ldots$ be subfields of $K$ with  $K= \union_n K(n)$, 
 $t \in K(1)$; let $R(n)=K(n) \meet R$,
$\bar{K}(0) = res(R(0))$.

  Assume $tr. deg._F K(0)    < \infty$,
  $tr.  deg.  _{K(n)} K(n+1) \leq 1$, and   $\si(R(n)) \subset R(n+1)$.

Then $\bar{K} \subset \bar{K}(0)^a$.  \end{lem}

\proof Let us first reduce to the case that $K$ is a transformal valued field, and
the $val(t_n)$ are
cofinal in $\Gamma$.  Let $\Gamma_\infty$ be the convex subgroup of $\Gamma$ generated by the elements $val(t_n)$.  (If $r \in R = \union_n R(n)$,
$val(r) \leq val(t_n)$ then $\si(r) \neq 0$, so $\Gamma'$ of \ref{tvrl}(3) contains  $\Gamma_\infty$.)  Let
$\hat{K},\hat{R},\hat{M},\pi,R_{\Gamma_\infty}$ be as in \ref{tvrl}(6), and let $\hat{K}(n) = \pi(K(n) \meet \hat{R})$.

Note that $\si$ extends to $\hat{R}$:  an element of $\hat{R}$ has the form $ba^{-1}$, with $a,b \in R$, $0 \leq val(a) \leq
val(t_n)$ for some $n$.  Then $val(\si(a)) \leq val(t_{n+1})$ (in particular, $\si(a) \neq 0$.)  Let $\si(ba^{-1}) =\si(b)
\si(a)^{-1}$.

It follows that $\hat{K}$ is a transformal valued field.  Since all elements of $\Ft$ have
valuation in  $\Gamma_\infty$, $(\hat{K},\hat{R},\hat{M},\bar{K})$
extend $(F(t)_\si,\Ft,t \Ft,F)$.  This effects the reduction.

Thus $\si$ extends to an endomorphism of $K$; and we have $\si(K(n)) \subset K(n+1)$.

Let $\Gamma_n$ denote the divisible hull within $\Gamma (=\Gamma_\infty)$ of $\{val(a):  a \in K(n) \}$.  So $\Gamma_n$ is a group of finite rank.

{\bf Claim 1} Let $\G_I = \{u \in \Gamma:  u << \si(u) \}$.  If one of $ u, \si(u) \in \G_I$ and $0 < u \leq v
\leq \si(u)$ then $v \in \G_I$.  Each $t_m \in \G_I$.

\proof If $\si(u) \leq mu$, then $\si^2(u) \leq m \si(u)$, so neither $u$ nor $\si(u)$ are in 
 $\G_I$.  Thus $u << \si(u)$.  Let $0 < u \leq v \leq \si(u)$.  If $mu \leq v$ for all $m$, applying $\si$ we obtain $mv < m \si(u) \leq \si(v)$ for all $m$, so $v <<
\si(v)$.  If $v< mu $ for some $m$, then $v \leq mu << \si(u) \leq \si(v)$.

{\bf Claim 2} If $u \in \Gamma_n \meet \G_I$ , then $\Gamma_n \neq \Gamma_{n+1}$, and in fact
$\Gamma_{n+1}$ has an element greater than any element of $\Gamma_n$.  

\proof In an ordered Abelian group, if $0< u_1 << u_2 << u_3 << \ldots$ then $u_1,u_2,u_3,\ldots$ are linearly independent.
Thus there is no infinite $<<$-chain of elements of $\Gamma_n$.  So we may take $u$ to be $<<$-maximal in $\Gamma_n \meet \G_I$.  But then $u << \si(u)$, $\si(u) \in \G_I$ so we cannot have $\si(u)  \in \Gamma_n$.  At the same time
$\si(u) \in \Gamma_{n+1}$.  So $\Gamma_n \neq \Gamma_{n+1}$.

{\bf Claim 3} For all $n \geq 0$, $\Gamma _n \neq \Gamma_{n+1}$.  Moreover, 
$\Gamma_{n+1}$ has an element greater than any element of $\Gamma_n$.  

\proof As $val(t) \in \Gamma_1 \meet \G_I $, the previous claim applies for $n \geq 1$.  It applies to $n=0$ too if $\Gamma_0$ 
has an element $v \geq val(t)$.  (Since $v \leq val(t^{\si^m})$ for some $m$, so $v \in \G_I$.)
Otherwise,   $\Gamma_0 < \val(t) \in \Gamma_1$.   This covers all the
cases.

{\bf Claim 4}     $K(n+1)_{\res} \subset  (K(n)_{\res})^a$.

\proof This follows from valuation theory, the previous claim, and the assumption $tr.  deg.  _{K(n)} K(n+1) \leq 1$.

This finishes the proof of the lemma.  \qed

As in definition \ref{tvr-def},let

$$\Delta  = \{v \in val(K):  -val(t) < nv < val(t), \ n=1,2,\ldots \}$$

\<{lem}\lbl{val2.2}  In \ref{val2.1}, we can also conclude:  $\Delta \subset  \G_0$. \>{lem}

\proof   Note $\Delta \subset \Gamma_\infty = \union_n \Gamma_n$.  By Claim 3 of \ref{val2.1},
$\Gamma_{n+1}$ has an element greater than any element of $\Gamma_n$.  But
$tr.deg. _{K(n)} K(n+1) \leq 1$ implies $rk_{\Qq}(\Gamma_{n+1}/\Gamma_n) \leq 1$, and it
it follows that
$\Delta \meet \Gamma_n = \Delta \meet \Gamma_{n+1}$.  By induction,
$\Delta \meet \Gamma_n \subset \Gamma_0$, and the lemma follows.
   \qed

\<{cor}\lbl{val2.4}  In \ref{val2.1}, let $d=tr. deg._F K(0)$.  Then 
 $R/tR$ has a unique minimal prime ideal $P$.  The total dimension of $R/tR$ 
is $\leq d$;  if equality holds, then 
 $tR \meet R(0) = P \meet R(0) = (0)$.   Ditto the total valuative dimension $\vrk(K)$.
 
If $K$ is $\omega$-increasing, then $^{\si} \sqrt{P} =M$.
\end{cor}

\proof Let $\Delta$ be as in \ref{val2.2}, and let $P = \{a \in R: val(a) \notin \Delta \}$.
Then $P  /tR$ is the ideal of nilpotents of $R/tR$, and $P$ is prime.

Let $R' =   \{a/b: a,b \in R, val(b) \in \Delta \}$.  
  $R'$ is a   valuation ring, with maximal ideal $PR'$, and  value group $\Gamma/\Delta$.  
$R'$ is not in general invariant under $\si$, but $R' \meet K_n$ has value group  
$\Gamma_n / \Delta$, and by
\ref{val2.1}, \ref{val2.2}   the rank of this value group grows with $n$.   Thus as in \ref{val2.1},
if $K'$ is the field of fractions $K'$ of $R/P$, then 
$K' \subset K'(0)^a$, where $K'(0) = res R'(0)$.  
 So $tr. deg._F K'  =   tr. deg. _F K'(0) \leq tr. deg. K(0)$;
  the second inequality is strict unless the $R'$-valuation is trivial on $K(0)$, and in this case $R(0) \meet P = (0)$.  

The total dimension of $R/tR$ equals that of $R/PR$ since $P/tR = \sqrt (0)$ in $R/tR$, and so
is bounded by the  transcendence degree of the  field of fractions of this domain.  The total   dimension
of $R/tR$ equals the reduced total dimension, $r.dim(R/tR) = \tdim (R/ M) = tr. deg._F  \bar{K}$,
 plus the  transformal multiplicity, or total dimension of $\spec R/t$ over $\spec K$.
As there is a chain of prime difference ideals  of length $rk_{ram}(K)$ between $M$ and $P$,
this total dimension is $\geq  rk_{ram}(K)$.  Thus 
 
$$\tdim (R/tR) \geq  tr. deg._F  \bar{K} + rk_{ram}(K) = \vrk (K)$$   
So  $\vrk(K)  \leq d$, 
and if equality holds then   $P \meet R(0) = (0)$. 

The last statement is clear from the definition of $P$:  if
$K$ is $\omega$-increasing, and $a \in M$, then $\val(a) >0$;
we have $\val(a) << \val(\si(a) ) << \ldots$, so they are $\Qq$-
linearly independent, and hence  cannot all be in $\Delta$; thus
$\si^m(a) \in P$ for some $m \leq rk_{\Qq} \Delta < \infty$.
 \qed

\end{subsection}  
\end{section}

\<{section}{Transformal specialization} \lbl{transspec}

\<{subsection}{Flatness }
  
Intersection theory leads us to study  difference schemes "moving over a line"; the behavior
of a difference scheme $X_t$ depending on a parameter $t$, as $t \to 0$.  In his Foundations
of Algebraic Geometry, Weil could say that $(a,t)$ {\em specializes} to $(a',0)$ (written
$(a,t) \to (a',0)$ )   if 
$(a',0)$ lies in every Zariski closed set that $(a,t)$ lies in; i.e. the point of $\spe X$
corresponding to $(a',0)$ lies in the closure of the (generic) point $(a,t)$.    By Chevalley's lemma, this
automatically implies the existence of a valuation ring, and the  theory that comes with that.  
We think of that as indicating the existence of a "path" from $(a,t)$ to $(a',0)$.  In 
difference schemes, closure and pathwise closure do not coincide.  We will use transformal 
valuations to define the latter.  

We  can view the valuations (or open blowing up) as revealing new components, that cannot be separated with
difference polynomials, but can be separated using functions involving for instance $x^{\si-1}$. 
(This can also be given responsibility for the failure of the dimension theorem in its
 original formulation, cf. \cite{cohn}.)

The above phenomenon can lead to a special fiber with 
higher total dimension than the generic fiber, or even with infinite total dimension; in this case it
cannot be seen as a smooth movement of a single object.  To prevent this, we need to
blow up the base.  It is 
convenient to replace the affine line once and for all with the
spectrum of a transformal valuation scheme (a posteriori,  a finite blowing-up will suffice to separate off   occult components in any particular instance.
 In general, we   will use the language of valuation theory in
 the present treatment, so that blowing ups occur only in the implicit background.)
Over a valuative base, we will show \ref{flat2} that  total dimension does  not increase.  

The lack of jumps in total dimension is  analogous to the
 preservation of dimension ("flatness") of the classical dynamic theory;  but it is { not}
the Frobenius transpose of the classical statement.   
The latter 
corresponds rather  to the preservation of {\em transformal} dimension, a fact that holds true already over
the usual transformal  affine line, without removing hidden components.  The good behavior of  
  total dimension is related under Frobenius  to another
principle of classical algebraic geometry; what Weil called  the   "preservation of number".

We actually require something more than the flatness of total dimension when measured
globally over a fiber.  Consider a difference subscheme $X$ of an algebraic variety $V$, or more
generally a morphism $j: X \to [\si]_k V$.   Say $X$ is {\em evenly spread} (along $V$, via $j$) if 
 for any proper  subvariety $U$ of $V$, $j^{-1}([\si]_k U)$ has total dimension $<\dim(V)$.
We need to know that if the  the  generic fiber $X_t$ 
is evenly  spread, then the same is true of the special fiber.
To achieve this, we need to replace the naive closure $X_0$ of  with a pathwise closure 
$\Xvs = ''lim_{t \to 0} X_t ''$, and 
the total dimension by a valuative dimension.
 (cf. \ref{key}).

As a matter of convenience, since we are interested in the generic point and in one special point
at a time, we  localize away from the others.

$\Ak$ will be used as a base for moving difference varieties, analogously to the affine line in rational equivalence theory of
algebraic varieties.  Let $X$ be a difference scheme over $\Ak$; we will write $X_t$ for the generic fiber $X \times_{\Ak}
\spe k(t)_\si$, and $X_0$ for the special fiber $X \times_{\Ak} \spe k$ (referring respectively to the inclusion $\kt \to
k(t)_\si$ and the map $\kt \to k$, $t \mapsto 0$.)

  Let $X$ be a difference subscheme of $\Vv = [\si]_k V \times_k \Ak$, $V$ an algebraic variety over $k$.  We
  denote by $X[n] \subset V \times V^\si \times \ldots \times V^{\si^{n}}$ the $n$'th weak Zariski closure of $X$, and similarly ${X_t}[n] $,
  ${X_0}[n]$.

 \<{definition}\lbl{flat1.5}  (cf. \ref{flat1}.)   We will say that a difference scheme $X$ of finite type over $\Ak$ is {\em flat} over $\Ak$ if
in every local ring, $y \mapsto t_n y $ is injective.    \>{definition}
When $X$ is algebraically reduced, this is equivalent to:  $X$ has no
component contained in the special fiber $X_0$.

\
\<{lem}\lbl{flat2} Let $X$ be flat over $\Ak$.  If $X_t$ has total dimension $d$ over $ k(t)_\si$, then $X_0$ has total
dimension $\leq d$ over $k$.  \>{lem}

\proof 
View $X[n]$ as a scheme over $\spec \kt$.  Since ${X_0}[n] \subset X[n]_0$, it suffices to show that
$\dim X[n]_0 \leq d$ (where $X[n]_0$ is the fiber of $X[n]$ above $t_0=t_1=\ldots=0$).  $X[n]$ arises by base extension from a
scheme $Y$ of finite type over $\spec S'$, where $S'$ is a finitely generated $k$-subalgebra of $S=\kt$.
Let $K = k(t)_\si$, $K[m] =   k(t_0,\ldots,t_m)$, 
 $S[m] = \kt[m] = \kt \meet K[m]$; then for some $m$ we can take $S'=\kt[m]$; and
  $\dim_{K[m] } (Y \tensor_{S[m]} K[m]) =d$.     

\claim Let $A$ be a finitely generated $S[m]$-algebra, such that $x \mapsto t_mx$ is injective on $A$.  Suppose
$\dim_{K[m] } (A \tensor_{S[m]} K[m]) = d$.  (This refers to Krull dimension.)  Let $A' = A/JA$,
where $J$ is the ideal generated by $(t_m,t_m t_{m-1}^{-1}, t_m t_{m-1}^{-2}, \ldots) $.  Then 
$A'$ is an $S[m]/J = S[m-1]$-algebra, 
$x \mapsto t_{m-1}x$ is
injective on $A'$ , and $\dim_{K[m-1] } (A' \tensor_{S[m-1]} K[m-1]) \leq d$

\proof Suppose $t_{m-1}c \in JA$; so $t_{m-1}c = t_{m} t_{m-1}^{-r} a$ for some $r \in \Nn$ and $a \in A$.  As $t_{m-1} | t_m
$, $ x \to t_{m-1}x $ is injective on $A$.  So $c= t_mt_{m-1}^{-r-1} a$.  But then $c \in JA$.  This shows that $x \mapsto
t_{m-1}x$ is injective on $A'$.

To see that the dimension remains $\leq d$, let $B = A \tensor _{K[m]} K[m-1][t_m]'$, where 
$K[t_m]'$ is the localization of $K[m-1][t_m]$ at $t_m=0$.  It is easy to see, by looking at the numerator, that $t_m$ is not a
0-divisor in $B$.  We are given that $B \tensor _{K[m-1][t_m]'} K[m-1](t_m)$ has Krull dimension $\leq d$.  Thus $\spec B \to \spec K[m-1][t_m]$ has 
generic fiber of   dimension $\leq d$, and has no component sitting over $t_m=0$, so it has
special fiber of dimension $\leq d$ as well (over $t_m=0$); thus 
 $B/t_m B$ has Krull dimension $\leq d$.  But $B/t_mB = A' \tensor _{S[m-1]} K[m-1]$, proving the claim. 

The lemma follows upon $m+1$ successive applications of the Claim. \qed

\begin{remark} \lbl{move1} \ \end{remark}

 Using the main theorems of this paper concerning Frobenius reduction, and the
translation this affords, one can conclude:         \begin{enumerate}

 \item A statement similar to \ref{flat2} holds for transformal
dimension; this corresponds to \cite{HA} III 9.6 (or, using the dimension theorem in $\Pp^m \times \Pp^1$, where $X
\subset \Pp^m$, note that no component of the special fiber $t=0$ can have smaller codimension than a component of a generic
fiber of an irreducible variety projecting dominantly to $\Pp^1$.)

\item When $X$ is a closed subscheme of $V \times_k \Ak$, $V$ a proper algebraic variety over $k$, one can conclude that $X_0$
has total dimension equal to $d$ (though not necessarily reduced total dimension $d$).

 (The Frobenius specializations give systems of curves over the affine line, having about $ q^d $ points over a generic point
of the affine line, hence the fiber over $0$ cannot be of size $O(q^{d-1})$.)

\item One cannot expect every component of $X_0$ to have total dimension $d$ (even with completeness and irreducibility
assumptions).  Nor can one expect the reduced total dimension to equal $d$.  It suffices to note that the diagonal $\Delta$ of
${\Pp^1}$ can be moved to $(\{pt\} \times \Pp^1) \union (\Pp^1 \times \{pt\})$ (and pull back with the graph of $\Sigma$.)

\end{enumerate}

The conclusion of \ref{flat2} below for ordinary schemes is
usually obtained from flatness hypotheses; cf.  e.g.  \cite{HA} III 9.6.  This suggested the terminology, 
as well as the following remark.  We have
not been able to use it, though, since it is not obvious how to reduce to schemes of finite type while retaining flatness.
The lemma refers to the $k$-algebra structure of $\kt$, disregarding the difference structure.

\<{rem}\lbl{flat1} Let $k$ be a difference field, $M$ be a $\kt$ module, such that $t_my=0$ implies $y=0$.  Then $M$ is a flat
$\kt$- module.  \>{rem}

\proof By \cite{HA} III 9.1A, it suffices to show that for any finitely generated ideal $I$ of $A_k$, the map $I
\tensor _R M \to M$ is injective.  Now $\kt$ is a valuation ring, so any finitely generated ideal is principal, $I = \kt c$.
The map in question is injective if $cy=0$ implies $y=0$.  But for some $m \in \Nn$, $ c | t_m$ in $\kt$, so $t_my=0$, and
thus $y=0$.

\end{subsection}

\<{subsection}{Boolean-valued valued difference fields}

We wish to consider finite products of  transformal valuation domains
(corresponding to a finite disjoint unions of difference schemes.)  
  As we wish to include
them in an axiomatizable class (for purposes of compactness and 
decidability), we will permit arbitrary products, and so will
discuss briefly Boolean-valued transformal valuation domains. 
Compare \cite{LS}, \cite{mac}; but here we will need them 
only at a definitional level.

\<{definition} (cf. \cite{shoenfield}.)  Let $T$ be a theory in a language $L$.  The language $L_{boolean}$
has the same constant and function symbols as $L$; and for each $l$-place formula
$\phi$ of $L$, 
 a new $l$-place function symbol $[\phi]$, taking values in a new sort $B$; also,
functions $\union, \meet ,\neg,0,1$ on $B$.   $T_{boolean}$ is the theory of all pairs
$(M,B)$, where $(B,\union, \meet,\neg ,0,1)$ is a Boolean algebra, and $M$
is a  Boolean-valued model of $T$.  Thus $[R](a_1,\ldots,a_l)$ is viewed as the truth value of $R(a_1,\ldots,a_l)$.

Axioms:  universal closures of $[ \phi \& \psi ]= [\phi ] \meet [\psi]$, 
$[\neg \phi] = \neg [\phi]$; if $\phi = (\exists x) \psi$, $[\phi] \geq [\psi]$, and    $(\exists x)( [\phi]=[\psi]$).  

And:  
$[\phi]=1$, where  $ T \models \phi$.

 \>{definition}

When $B= [\bf 2]$ is the 2-element Boolean algebra, we obtain an ordinary model
of $T$.

If $(M,B) \models T_{boolean}$, and $h: B \to B'$ is a homomorphism of Boolean
algebras, we obtain another model $(M,B')$ of $T_{boolean}$, by 
letting $[\phi]' = h ([\phi])$.

Thus when $(M,B) \models T_{boolean}$,   $X = Hom(B,{\bf 2})$, for each 
$x \in X$ we obtain a model $M_x$ of $T$.

When $T$ is a theory of fields, it is not necessary to have $B$ as a separate
sort.  A model $M$ of $T_{boolean}$ is a ring, and $B$ can be identified with the idempotents of $M$, via the map $[x=1]: M \to B$.  (The sentence
$(\forall y)(\exists x) (y=0 \to x=0  \, \& \, y \neq 0 \to x=1  )$, true in $T$,
must have value $1$ in $T_{boolean}$, and shows that for any $m \in M$
there exists an idempotent $b$ with $[m=0] = [b=0]$.)

When $L$ has an additional unary function symbol $\si$ and $T$ states that 
$\si(0)=0, \si(1)=1$, $T_{boolean}$ implies that $\si$ fixes the idempotents.

We will just need the case when $T$ is the theory of transformal valued fields
$(K,R,M)$.  
A model of $T_{boolean}$ is then a difference ring, called a {\em multiple 
transformal valued field.}   $\{r: [val(r) \geq 0] =1 \}$ is a subring $S$, called
a {\em boolean-valued
transformal valuation domain.}  
Thus a boolean-valued
transformal valuation domain is a certain kind of difference ring $R$,
without nilpotents,  such that for any ultrafilter $U$ on the Boolean algebra 
$B$ of idempotents of $R$, $R/U$ is a transformal valuation domain.  This class of rings is closed under Cartesian products.   

When the Boolean algebra is finite, we say $T$ is {\em finitely valued.}  A finitely valued
transformal valuation domain is just a finite product of transformal valuation domains.

\>{subsection} 

\begin{subsection}{Pathwise specialization}

  We will formulate a notion of
a    ''specialization" $\lim_{t \to 0} X_t$   of a difference variety $X_t$ over $F_t= F(t)_\si$,
as $t \to 0$.    This will be a difference subscheme of $X_0$; we will denote
it more briefly as $\Xvs$.  In addition, there will be  certain data over each point of $\Xvs$ determining 
the "multiplicity" of that point as a point of specialization.    

A {\em valuative} difference scheme over $\Ak$ is a difference scheme $T = \spe R$, $R$ a   
 transformal valuation domain extending $\kt$, the valuation ring of 
 $k(t)_\si$.  

 Note that any family of closed difference subschemes of a difference scheme $X$ has a ''union", a smallest
closed subscheme containing each element of the family.  This 
correspond to the fact that  the intersection of
any family of well-mixed difference ideals is again one.  This makes possible the following definition.  
 
 \<{definition}\lbl{vs}
   Let $X$ be a   difference scheme of finite type over $\Ak$. 
 We define the pathwise specialization 
$\Xvs$   to be the smallest well-mixed difference subscheme
 $Y$ of $X_0$ such that for any valuative difference scheme $T  $ over $\Ak$ , and any morphism 
$ f:T \to X $ over $X \to \Ak$, $f(T_0) \subset Y$ 
(i.e. the restriction of $f$ to the fiber over $0$ factors through $Y$.)
                      \>{definition}

Intuitively, $T$ is a kind of transformal smooth curve, and $f(T)$ marks a path from points on $X_t$ to points on $X_0$.  

\<{remark} \ \lbl{vs-r} \>{remark} 
\begin{description}
   \item[(i) ]   Since $T$ is perfectly reduced, $\Xvs$ depends
only on the perfectly reduced difference scheme underlying $X$.
(Though of course $\Xvs$ need not itself be perfectly or
even transformally reduced.)

  \item[(ii) ] 
One could think of allowing  more singular paths, using weakly 
 transformal valuation rings, with nilpotent elements permitted.   Presumably, when $X_t$ is reduced, and
irreducible even in the valuative sense, the two definitions yield the same object.   

\item[(iii)]   A simple point of $X_t$
may specialize to a singular point of $X_0$; or several simple points may specialize to one.  This multiplicity
of specialization is not faithfully reflected in geometric multiplicity on $\Xvs$.  (For instance,
when $X$ is the union of of several ''transformal curves" $C_i $ meeting at a point $p$ of the special fiber,   each $C_i$
 will map into $X$, but separately;  so that $p$ is   registered just once on $\Xvs$.)   Thus the local algebra of points
on $\Xvs$ is insufficient to capture the data needed over each point. 

One can modify the definition of $\Xvs$  so that the local
algebra will capture the multiplicities; cf. $\Xvs '$ in  \S \ref{pinched-sec}.  
 We will  deal with the issue differently,    replacing the local algebra information by  scattered sets (definable in the language of transformal valued fields, and residually analyzable.)   This approach appears to us much more transparent.

 \end{description}

\<{lem}\lbl{vs-lim} 
 
  Let $Z$ be an algebraic component   of $\Xvs$.
 Then there exists a valuative scheme $T$ over $\Ak$,   and a 
morphism $ T \to X $ over $\Ak$,  such that $Z$ is contained in the image of $T_0$ in $X_0$.
  
    \>{lem}

\proof     Moving to an open affine difference subscheme of $X$, we can assume $X = \spe R$; 
$Z$ corresponds
to an algebraically prime difference ideal $p$; 
 there are homomorphisms $h_j: R \to A_j \subset L_j$, $L_j$
a  transformally valued field   extending $\kt$,  with corresponding 
transformal valuation domain $A_j$,  such that $p$ contains $\meet_{j \in J} {h_j}^{-1}(tA_j)$.
  Let $h: R \to A_*  = \Pi_j A_j \subset L_* = \Pi_j L_j$ 
be the product homomorphism; then $h^{-1}(tA_* ) \subset p$, and $A_* $ is
a Boolean-valued transformal valuation domain, with Boolean algebra $E$ of
idempotents.
 
Let $U_0$ be the filter in  $E$ generated by all $e \in E$ such that for some   $r \in R \setminus p$,
 $h (r) \in tA_* + eA_*   $.  
If $r_1,\ldots,r_n \in R \setminus p$, and $h(r_i) \in tA_* + e_iA_*$, then (since $p$ is prime)
$r=r_1 \cdot \ldots \cdot r_n \in R \setminus p$; and with $e = e_1 \meet \ldots \meet e_n$,  
  $h(r) \in tA_* + eA_*   $; so $e \neq 0$.  Thus $I_0$ is a proper filter, and so extends to an ultrafilter $U$,
with complementary maximal ideal $I$.
Composing $h$ with the map $A_* \to A_*/I =: \bar{A} \subset   L_*/I =: \bar{L} $, we obtain 
$h^*: R \to \bar{A} \subset \bar{L}$, $\bar{A}$ a transformal valuation domain   extending $\Ak$,  such that for $r \in R \setminus p$,
$h(r) \notin t \bar{A}$.  Thus $p$ contains $h^{-1}(t \bar{A})$.   \qed

We come to the main estimate for the ''equivalence of the infinite components".  Given a difference scheme $X$, we have the residue map $X_t \to X_0$
(on the points in a transformal valued field.) If $X$ has
 finite total dimension over $\Ak$,
the pullback of a   
difference subvariety $W \subset X_0$ is a certain scattered set 
$^*W \subset X_t$.    
Under certain conditions of direct presentation, we find an upper
bound on  the inertial dimension of $^*W$; it can be 
smaller
than the total dimension of the closure of $^*W$ within $X_t$.

 Let $V$  be an algebraic variety over  an inversive difference field $k$,
$\dim(V) =d$.  
Let $V_{\Ak} = [\si]_k V \times_{k} \Ak$, $V_t=  V \tensor _k k(t)$. 


 \<{prop}\lbl{key} Let $X $ be a closed difference subscheme of  
$V_{\Ak}$.   Then (1) implies (3), and (1 \& 2) implies (3 \& 4).

 \begin{enumerate}

 \item   Consider transformally valued fields  $L$  
generated over $k(t)_\si$ by  $c \in X_t(L) \subset V(L)$. 
 Whenever $tr. deg. _{k(t)_\si} L \geq d$,   equality holds, and $\si(c) \in k(t,c)^a$.  

  \item  Every  weakly  Zariski dense  (in $V$) algebraically integral  subscheme of $X_0$ 
 is transformally  reduced.

\item    $\Xvs$ is evenly spread out along $V$. 
 
\item 
for any subvariety  $W$  of $V$ with  $\dim(W) < d = \dim(V)$,  
$$^*W = \{ x \in X_t(R):   \res_X (x) \in W  \}$$
 has inertial dimension $< \dim(V)$.  
\end{enumerate}
\>{prop}  

Explanations:

(1e) A point $c$ of $X_t(L)$ corresponds to a morphism
from a local difference ring of $X_t$, into $L$; let $k(c)_\si$
denote the field of fractions of the image.  The inclusion 
$X_t \subset [\si]_{k(t)_\si} V_t$ allows us to restrict $c$ to a point
of $V(L)$, with a morphism from a local ring of $V$ into $L$;
let $k(c)$ denote the field of fractions of the image.  Then (1)
states that if $tr. deg._{k(t)_\si} k(c)_\si = d$, then $k(c)^a = (k(c)_\si)^a$.  It follows also that $k(c)^a = k(V)$.  

 (2) implies that a weakly Zariski dense algebraic component 
is actually a component of $X_0$.  Note (2e) that any weakly
Zariski dense (in $V$) difference  scheme has 
total dimension $\geq \dim(V)$.  

(3) means:  every algebraic component $Y$ of $\Xvs$ of total dimension $\geq d$  is 
weakly Zariski dense in $V$.  Note that (1) is a   weak form of 
the statement that $X_t$ is evenly spread out along $V$.

 (4):  Here $R$ denotes the valuation ring, a quantifier-free
formula in the language of transformal valued rings.   $\res_X$
denotes the map  $X_t(R^L) \to X_0( L_{\res})$ induced by $\res$ 
on a transformally valued field $L$.

\proof  $(1) \implies (3):$ 

 We may  assume $X= \spe A$.  Let $Y$ be an algebraic component of
$\Xvs$, and let $p$ be the corresponding algebraically prime difference ideal of $A$.   By \ref{vs-lim}, 
 there exists a  transformally valued field $(K,R)$ 
  extending $(k(t)_\si, \Ak)$  and a difference ring 
morphism $ h: A \to R $ over $\kt$,  such that $p$   contains   $h^{-1}(t R)$; so $A/p$ embeds
into a quotient of $R/tR$. 

 Note that $h^{-1}(tR)$ may be bigger than $(ker h, t)$, and $h(A)$
cannot be assumed to generate $R$ (cf. \S \ref{pinched-sec}).   But 
we may assume 
  $h(A)$ generates  $K$ as a field over $k(t)_\si$.  
By (1), $tr. deg. _k K \leq d$.

Now $t$ is not  a 0-divisor in $A/\ker(h)$.
 By \ref{flat2}, $A/(\ker(h),t)$ has total dimension $\leq d$,
hence so does $A / h^{-1}(tR)$.    Moreover 
if    $tr. deg. _k K < d$,  or if $p$ is not a minimal prime
of $A/ h^{-1}(tR)$, then 
 then $Y = \spe A/p$ has total dimension $<d$.  So we may assume $tr. deg._k K =d$ (hence by (1), $k(V)$ embeds into $K$ via $h$),
 and $p$ is a   minimal prime
of $A/ h^{-1}(tR)$.

   Let $K[0] = k(V)^a \meet K$.  Let $K[n] = K[0](t_0,\ldots,t_{n-1})^a \meet K$.
Then 
 $\si(K[0]) \subset {K[1]}^a$; so $\si(K[n] ) \subset K[n+1]$.
Now \ref{val2.1}-\ref{val2.4} apply to $\union_n K[n]$. 
  If $R/tR$ has total dimension $<d$, then so does $A/p$, hence $Y$.
In this case (3) holds trivially.   Otherwise by  \ref{val2.4}, the nilradical   of $R/tR$ lifts to 
 an algebraically prime difference
ideal $P$ of $R$, with $P \meet K[0] = (0)$.
If $h(x) \in P$ then $h(x^n) \in tR$ for some $n$, so $x^n \in p$,
and thus $x \in p$.  So $h^{-1}(P) \subset p$.  But $p$ is a minimal 
prime, so $h^{-1}(P)=p$.  Thus  
 $p \meet k(V) = (0)$, i.e. $Y$ is weakly Zariski dense in $V$.

\smallskip
$(1 \& 2) \implies (4):$   We will use \ref{vrk-vdim}.
 Let $M$ be an $\omega$-increasing transformal valued field
extending $K=k(t)_\si$, with valuation ring $R=R^M $; let
$c \in X_t(R^M)$, $\res(c) \in W$.  Let $K = k(t,c)_\si$.
We must
show that $\vrk K/k(t)_\si <  \dim(V)$.

  If $tr. deg. _{k(t)_\si} K < \dim(V)$, this is 
clear.  Otherwise, let $K(0) = k(c)$, $K(n) = k(c,\ldots,\si^n(c),t,\ldots,\si^{n-1}(t))$.  By (1),   $\si(c) \in k(c,t)^a$
so  $\si^{n+1}(c) \in k(\si^n(c),\si^n(t))^a$.  
So the hypotheses of \ref{val2.1} hold.  By \ref{val2.4}, 
$R/tR$ has a unique minimal prime ideal $P$; and
$\vrk(K) < d$, unless $\vrk(K)=\tdim (R/tR) = d$, 
and  $P \meet R \meet K(0) = (0)$. 

 In the  latter case, we will obtain a contradiction.   Let $A$ be the image in $R$ of the local ring
of $X$, corresponding to the point $c$.  Let $A(0)$ be the image
of the corresponding restriction to a point of $V$, as in (1e).  So $A(0) \subset R^M \meet K(0)$. 

 Thus
$A(0) \meet P = (0)$.  So 
$Z=\spe A/(P \meet A)$ is a closed difference subscheme of
$X_0$, weakly Zariski dense in $V$.   By \ref{flat2}, $Z$ has total dimension $\leq d$;
being weakly Zariski  dense in $V$, it must have dimension $d$.  
By (2), $Z$ is a component of $X_0$.   So $P \meet A$ is transformally prime in $A$.

Let  $K_A$ (resp. $K'$ ) be the  field of fractions of $A/P$
(resp. $R/P$).  Then $tr. deg._k K_A = d = tr. deg._k K'$.  So
$K' \subset (K_A)^a$.  
By lemma \ref{alg-h}, $P$ is transformally prime.  But if $a \in R$, $\val(a) >0$,
then $\si^k(\val(a)) > \val(t)$ for some $k$; so $\si^k(a) \in P$;
and thus $a \in P$.  Hence $P$ is the maximal ideal of $R$.  Since
$P \meet K(0)=0$, the residue map is an isomorphism on $K(0)$.
But $c \notin W$, $\res(c) \in W$; a contradiction  
   \qed

 \end{subsection}



\end{section} 

\begin{section}{Frobenius reduction}   \lbl{frobred}

\begin{subsection}{The functors $M_q$}

These functors may be viewed as difference-theoretic
 analogs of ``reduction mod p"; we reduce mod $p$,
and also ``reduce" $\si$ to a Frobenius map. In the case of
subrings of number fields, the situation is what  one ordinarily
describes using the Artin symbol.

Let $R$ be a difference  ring. Let $q$ be a power of a prime number $p$.
We let $J_q(R)$ be the ideal generated by $p = p \cdot 1_R$ together with
 all elements
$r^q-\si(r)$. Let $M_q(R) = (R/J_q(R))$. This is a difference ring,
 on which $\si$ coincides with the Frobenius map $x \mapsto x^q$. In
particular, $J_q$ is a difference ideal.

\begin{lem}\lbl{localizations} Let $T$ be a multiplicative subset of $R$
 with $\si(T) \subset T$, $\bar T$ the image
of $T$ under the quotient map in $M_q(R)$. Let $R[T^{-1}]$ be the localized
ring. Then $M_q(R[T^{-1}]) = M_q(R)[\bar T^{-1}]$. \end{lem}
\proof  
In whatever order it is obtained, the ring is characterized by the universal
property for difference rings $W$ with maps $R \to W$, $p=0$ in $W$,
such that the image of
$T$ is invertible, and $\si(x)=x^q$.   \qed

\begin{lem}  Let $f: S_i \to R$ be surjective maps of difference rings.  Then
\[ M_q(S_1) \times_{M_q(R)}  M_q(S_2) = M_q( S_1 \times_R S_2) \] \end{lem}

\proof  Similar to that of \ref{localizations}

\begin{lem}\lbl{kernels} Let $h: R \to S$ be a surjective difference 
ring homomorphism, with kernel $I$.
Then the kernel of $M_q(h): M_q(R) \to M_q(S)$ is $(I+J_q(R))/J_q(R)$. \end{lem}

\proof The kernel is $h^{-1}(J_q(S))/J_q(R)$. If $h(r) \in J_q(S)$,
$h(r) = \sum t_j(s_j^q - \si(s_j)) + ps'$ for some $t_j,s_j,s' \in S$. Writing
$s_j = h(\bar s_j)$, $t_j = h(\bar t_j)$, $s' = h(\bar s')$,
$\bar r = \sum_j \bar t_j (\bar s_j^q -\si(\bar s_j)) +p \bar s'$, we obtain $r - \bar r \in I$,
$\bar r \in J_q(R)$, so $r \in I+J_q(R)$. \qed

\begin{definition}
Let $X$ be a difference scheme.  Let  ${\cal J}_q$ be the difference-ideal sheaf
on $X$ generated by the presheaf:  $ {\cal J}_q (U)  = J_q ( \cO_X(U) ) $.
$M_q(X)$ is the closed subscheme corresponding to ${\cal J}_q$. \end{definition}

The above lemmas show that for a difference ring $R$ without zero - divisors,
$M_q ( \spe R) = \spe M_q(R) $.

Let us say that a {\em $q$-Frobenius difference scheme} is a difference scheme in which
$\si$ coincides
with $x \mapsto x^q$ on every local ring.   Any ordinary scheme together
with a map into $\spec \Ff_q$ admits a canonical $q$-Frobenius difference scheme
structure.

$M_q$ yields a functor on difference -schemes into $q$-Frobenius difference
schemes.

\begin{rem} \ \end{rem} Let $K$ be a number field, $\si$ an automorphism, $R$
a finitely generated difference subring of $K$. Then $R$ is also finitely generated
as a ring. $\spec M_q(R)$ may be empty, and may be reducible. It is always
a reduced scheme: if $d = [K:\Qq]$, then $\si^d(x) = x$ on $R$,
so $x^{q^d} = x$ on $M_q(R)$, hence $M_q(R)$ has no nilpotents.

\begin{rem} \ \end{rem} Let $R$ be as above,
$S$ a finitely generated domain over $R$ of positive but finite transcendence
degree $d$, and let $\si$ be the identity automorphism.   Then
for all but finitely many $q$,
$M_q(S)$ is reducible. Indeed every homomorphism of $S$ into $K_q$
factors through $M_q(S)$.  We will see later that
the number of homomorphism of $M_q(S)$ into $K_q$ is $O(1)q^d$;
these all fall into the field $GF(q)$. The relevant Galois group has order $d$,
so the number of prime ideals of $M_q(S)$ is at least $O(1)q^d / d$, and in particular
$>1$.  But $M_q(S)$ is finite, hence it is reducible.

On the other hand, if $S = \Zz[X], \si X = X+1$, then $M_p(S)$ is a domain for
all primes $p$.

Even when the fraction field of $S$ is
a ``regular" extension of $K$, in the sense that it is linearly disjoint
over $K$ from any difference field of finite transcendence degree over $K$,
 and $M_q(R)$ is a field,
$M_q(S)$ can be reducible.   Example: $R = \Zz[X,Y: \si(X) = YX]$.

It would be interesting to study this question systematically.

\begin{notation} \lbl{Xqy} Let $X$ be a difference scheme over $Y$,
$q$ a prime power, $y$ a point of $M_q(Y)$ valued in a difference field $L$.
We denote $X_{q,y} = (M_q(X))_y$ (cf. \ref{Xy}).  
Thus if $X= \spe R$, $Y = \spe D$, $R$ a $D$-algebra, then 
$X_{q,y} = \spec M_q(R) \tensor_{M_q(D)} L$.
\end{notation}

\end{subsection}

\begin{subsection}{Dimensions and Frobenius reduction}

We will see that the functors $M_q$ take
transformal dimension to ordinary algebraic dimension.  When the transformal
dimension of $X$ is zero, $M_q(X)$ will be a finite scheme; in this
case the total dimension is a (logarithmic) measure of the rate of growth
of $M_q(X)$ with $q$.   Transformal multiplicity becomes, in a similar
sense, geometric multiplicity, while  transformal degree 
is related  to the
logarithm of the projective degree.

By the {\em size} $|Y|$ of a $0$-dimensional scheme $Y$ over a field, we will mean the
number of points, weighted by their geometric multiplicity.   Thus for a $k$-algebra $R$,
 the size of $\spec R$ over $\spec k$ is just $\dim_k (R)$.

\begin{lem}\lbl{si-dim}
 Let $X$ be a difference scheme of finite type over  a Noetherian
difference scheme $Y$.
\begin{enumerate}

\item
Assume $X$ has transformal
dimension $\leq d$ over $Y$.
Then for all  prime powers $q$, and all points $y$ of $M_q( Y)$,
$X_{q,y}$ has dimension at most $d$ as a scheme over $L_y$.

\item
If   $X$ has reduced total dimension $\leq e$ over $Y$, then
there exists $b \in \Nn$ such that for all large enough prime powers $q$, and all
points $y \in M_q(Y)(L)$, $L$ a difference field,
 the zero-dimensional scheme $X_{q,y}$ over $L_y$
has   at most $bq^e$ points.

\end{enumerate}
\end{lem}

{\bf Remark} The statement of \ref{si-dim}(2) refers to the number of points,
multiplicities ignored.  We will see later that it remains true if multiplicities are
taken into account.  The proof of this refinement is  more delicate in that
when  passing to a proper difference subscheme $X'$ of $X$ one cannot forget $X$;
even if a point is known to lie on $X'$ one must take into account the multiplicity
of $X$ , not of $X'$, at that point.

\proof
The lemma reduces to the case
$X = \spe (R)$, $Y = \spe(D)$,  $R$ a finitely - generated $D$-algebra.   We may assume
$Y$ is irreducible and perfectly reduced, i.e. $D$ is a difference domain.
We wish to
reduce the lemma further to the one-generated case, $R = D[a,\si(a),\ldots]$.
We use Noetherian induction
on $\spe D$ and, with $D$ fixed, on $\spe S$, and thus assume the lemma is
true for proper closed subsets.

\claim  If $D \subset S \subset R$ and
the lemma holds  for $\spe R _{U'}$ over $U'$ and for $\spe S _U$ over $U$, whenever
$U',U$ are Zariski open in $\spe S$, $\spe D$ respectively, then
  it holds for $\spe R$ over $\spe D$.

\proof   We prove this for (1); the proof for (2) is entirely similar. 
  Let $p \in \spe R$,
$p \meet D = p_D$, $p \meet S = p_S$.  Let
$d'$ the relative dimension of $S_{p_S}$ over $D_{p_D}$,
$d''$ the relative dimension of $R_p$ over $S_{p_S}$.  Then
$d' + d'' \leq d$.  There exists an nonempty open  affine  neighborhood $U'$
of $p_S$ in $\spe S$
such that  $\spe R$ has relative dimension $\leq d''$ over any point
of $U'$.  There exists an nonempty open  affine  neighborhood $U$ of
$p$ in $\spe D$ such that $\spe S$ has relative dimension $\leq d'$
over any point of $U$.   Let $q$ be a large prime, $y \in M_q(\spe D)$.
Then as the lemma is assumed to hold in those cases, $(\spe S _U')_{q , y}$
has dimension at most $d''$ over $U'$, while $(\spe R _U)_{q,y}$ has
dimension at most $d'$.  It follows that for any $p \in \spe R$
with $p \meet R \in U'$ and $p \meet D \in U$, the lemma holds. By
additivity of transcendence degree in extensions, the lemma holds for
$\spe R$ over $\spe D$ at a neighborhood of any such $p$.  Outside
of  $U$, or of $U'$, the lemma still holds by Noetherian induction.  \qed

   We continue to use Noetherian induction on $\spe D$.
Let $k$ be the field of fractions of $D$.  We may assume $R$ is perfect,
since factoring out the perfect ideal generated by $0$ changes neither the
transformal or reduced total
 dimension nor the physical size of the 
zero-dimensional schemes $X_{q,y}$.  In this case $O$ is the intersection
of finitely many transformally prime ideals $p_i$ of $R$, and it suffices to prove
the lemma separately for each $R/p_i$.  Thus we may assume $R$ is a difference
domain.

By the Claim, we may assume $R$ 
is generated a single element $a$
as a difference $D$-algebra. 

If $R$ is isomorphic to the difference polynomial ring over $D$ , the transformal
dimension of $\spe R$ over $ \spe D$ equals $1$, and we
must   show  the dimension of $spec M_q(S)$ over
$spec M_q(D)$ is everywhere at most $1$.  This is clear since
$M_q(R) $ is generated over $M_q(D)$ by a single element.

Otherwise, there exists
a relation $G(a,\si(a),\ldots,\si^m(a))=0$ in $F$, $G$ over $D$.
Moreover in case (2), we can choose $m \leq e$.
We may write $G(a) = \sum c_i a^{\nu_i(\si)}$, where $c_i \in S, c_i \neq 0$, the
$\nu_i$ are finitely many integral polynomials,
of degree at most  the total dimension $r$ of $\spe R$ over $\spe S$.

Note that for all sufficiently large $q$, the values $\nu_i(q)$ are
distinct, indeed $\nu_i(q)<\nu_j(q)$ if $i<j$.  The set $Y'$ of  transformally prime
ideals of $D$
containing one of the nonzero coefficients of $G$ is a proper closed subset
of $\spe S$; by induction the lemma  is true over $Y'$.
On the other hand if $y \in \spe S \setminus Y'$, then the monomials
of $M_q(G)$ are distinct and their coefficients are non-zero  in $L_y$;
it follows that $M_q(R) \tensor L_y$ is $0$-dimensional over $L_y$, and moreover
has at most $\max_i \ \nu_i(q)  \leq O(q^e)$ points.  This
finishes the proof of the lemma.    \qed

\end{subsection}

\begin{subsection}{Multiplicities and Frobenius reduction}

The order of magnitude of multiplicity upon Frobenius reduction will be
shown to be bounded by the transformal multiplicity.  We will begin with
 two  elementary and purely algebraic  lemmas regarding multiplicities.

\begin{subsubsection}{An explicit example}

We begin with an example in one variable; it shows explicitly the distribution of 
multiplicities of the points of $M_q(X)$, where $X$ is a difference scheme defined
by $F(X,X^\si,\ldots,X^{\si^n}) = 0$.   If $F$ is irreducible, 
there will be at most  about $q^{n-l-1}$ points with multiplicity
of order $q^{l}$.  The general picture is similar (cf. \ref{si-dim3.2}),
but we have not succeeded in reducing
it to the example (because of tricky  additivity behavior of transformal multiplicity),
and will give an independent proof.

  The Hasse derivatives
$D^{\nu}f$ of a
polynomial $f$ are defined by the Taylor series expansion:

\[ f(X+U) = \sum_{\nu} (D^{\nu} f)(X) U^{\nu} \]
Here $X = (X_0,\dots, X_n)$, $U = (U_0,\ldots,U_n)$,
$X+U = (X_0+U_0,\ldots,X_n+U_n)$,
  $\nu$ is a multi-index
$(\nu_0,\ldots,\nu_k)$, $U^{\nu} = \Pi_i {U_i}^{\nu_i}$.  Write $sup \  supp \ \nu$
for the highest $i$ such that $\nu_i > 0$.

\begin{example} \lbl{mult-1} Let $K$ be a field, $f \in K[X_0,\ldots,X_n]$ an irreducible
polynomial, $f \notin K[X_1,\ldots,X_n]$.  For $k \geq 1$,  let $V_k$ be the
Zariski closed subset of  $\Aa^{n+1}$ defined by the vanishing of all $D^{\nu}f$
with $sup \ supp \  \nu < k$.
\begin{enumerate}
\item $V_k = \Aa^{k} \times U$,
 $U \subset \Aa^{n+1-k}$, $\dim(U) < n-k $.
\item Assume $K$ has characteristic $p >0$, and let $q$ be a power of $p$
with
$q > \deg_{X_i} (f)$ for $i<k$.
Let $F(X) = f(X,X^q,\ldots,X^{q^n})$.  Let $S$ be a subscheme of $\Aa^1$ defined by
an ideal $I$ with $F \in I$, and let $a \in S$, $(a,\ldots,a^{q^n}) \notin V_k$.
Then the
multiplicity of $S$ at $a$ is bounded by 
$$ \mult_a (S) \leq {\sum_{i<k} \deg _{X_i} f}q^{k-1} \leq (\deg f) q^{k-1}$$

\end{enumerate}   \end{example}

\proof  (1) If $(a_0,\ldots,a_n) \in V_k$, then using the Taylor series expansion,
we see that all polynomials $D^{\nu} f$  are constant on $\Aa^k \times \{(a_k,\ldots,a_n)\}$;
so $\Aa^k \times \{(a_k,\ldots,a_n)\} \subset V_k$.  Thus $V_k = \Aa^k \times U$ for
some $U$.  Since $f$ is irreducible, and involves $X_0$, $V_k$ is a proper subset
of the zero set $V(f)$.  Thus $\dim V_k < \dim V(f) = n$.  So $\dim U < n - k $.

(2) We may assume $a=0$.  Let $g$ be the sum of monomials of $f$ involving only
the variables $X_0,\ldots,X_{k-1}$.  Then $D^{\nu}(f-g)(0) = 0$ if $sup \ supp \ \nu < k$.
As $0 \notin V_k$, $D^{\nu}f (0) \neq 0$ for some such $\nu$.  Thus $g \neq 0$.
  Let $G(X) = g(X,\ldots,X^{q^{(k-1)} }  )$.  A monomial $\Pi_i X_i^{a_i}$ of $g$ turns  into
the monomial $X^{\sum a_iq^i}$ of $G$; as by assumption $a_i < q$ for each $i$,
no cancellation occurs; so $G \neq 0$.  Write $G(X) = X^m H(X)$,
$H(0 ) \neq 0$.  Then  for some  $I$, 
$$F(X)    = G(X) + X^{q^k} I(X) =  X^m( H(X) + X^{q^k - m}I(X))$$
Note $m < q^k$.
So $F(X)/X^m$ is a polynomial not vanishing at $0$, thus $m$ bounds the multiplicity of
$F$ at $0$.  \qed

\end{subsubsection}
\begin{subsubsection}{Algebraic lemmas}
 \begin{lem}\lbl{gen-M}  Let $k$ be a field of characteristic $p$, $[k:k^p]=p^e$ (in
 characteristic $0$,
let $e=0$.)  Let $S$ be a
$k$-algebra, generated by $m$ elements, $M$ a maximal ideal of $S$
with $dim_k(S/M) < \infty$.  
   Then $M$   is generated by at most $m+e+1$ elements. \end{lem}

\proof   Let $h: S \to K$ be a surjective homomorphism with kernel $M$, 
$K$ a finite extension field
of $k$.
 Let $s_1,\ldots,s_m$ be generators for $S$ as a $k$-algebra.

Consider first the case $e=0$, i.e. $K/k$
separable.  Using the primitive element theorem, $K$ is generated over $k$ by
one element $a=h(s_0)$; and one relation, the minimal monic polynomial 
$P(a)=0$, $P \in k[X]$.  We have $h(s_i) = Q_i(a) = h (Q_i(s_0))$
 (some $Q_i \in k[X]$), so $s_i-Q_i(s_0) \in M$; and clearly $M$ is generated
by the $m+1$ elements $P(s_0),s_i - Q_i(s_0)$.  

In general, $K/k$ is generated by $e+1$ elements, and $e+1$ relations.  
(Let $K \leq K_s \leq K$, where $K_s/k$   separable, $K/K_s$ purely
inseparable.  Then $K_s/k$ can be presented by one generator $a_0$ and one relation,
as above. $K/K_s$ is generated by $\leq e$ generators $a_1,\ldots,a_e$
and $\leq e$ relations; these can be viewed 
as a relation over $k$ between $a_0$ and the $a_i$.)
As above, the maximal ideal $M$ is generated by the preimages of the above
relations.  \qed

\begin{lem}\lbl{mult-2-0}  Let $k$ be a field, $R$ a finite dimensional local $k$ - algebra,
with maximal ideal $M$.  Let $S$ be a finitely generated $R$-module.  Then
\[ \dimk  S \leq  (\dimk  R) (\dimk  (S/ MS) ) \]
\end{lem}

\proof  Let $A$ be a $k$-subspace of $S$ with $\dimk  (A) = \dimk  (S/ MS)$ and
$A + M S = S$.  Let $T$ be the $k$-span of $RA$.  So $T$ is an $R$-submodule
of $S$, and $\dimk (T) \leq  (\dimk  R) (\dimk  (S/ MS) )$.  Also
$T+MS = S$.  By Nakayama (applied to the $R$-module $S/T$), $T=S$.

\begin{lem} \lbl{mult-2}  Let $f: X \to Y$ be a morphism of finite schemes over $\spec(k)$.
Let $p \in \spec Y$, $X_p$ the fiber of $X$ above $p$, and $q \in   X_p$.  Then
\[ \mult_q X \leq (\mult_q X_p) (\mult_p Y) \]  \end{lem}

\proof  Let $R,S$ be the local rings of $Y,X$ at $p,q$ respectively.  Then the local ring
of $X_p$ at $q$ is $S / MS$, $M$ being the maximal ideal of $R$.  We must show
that $(\dimk  S)   \leq  (\dimk  R ) (\dimk  (S/ MS)   ) $.  This follows from \ref{mult-2-0}.

\begin{lem}\lbl{mult-2-3}  Let 
 $S$ a   local ring, with finitely generated maximal ideal $M$.  
 Suppose $S/M^r$ has length $n < r$.   
 Then  $M^n=0$.  
\end{lem}

\proof   We have $S \supset M \supset M^2 \supset \ldots \supset M^r$.
$M^i/M^{i+1}$ is an $S/M$-space of some finite dimension.  
As the length of $S/M^r$ is $n<r$, we must have 
$M^i = M^{i+1}$ for some $i \leq n$.  By Nakayama, $M^i = 0$.  So $M^n=0$.

{\bf N.B.} Let $I$ be an ideal of a ring $R$ with $p \cdot 1_R = 0$, and let $q$
be a power of $p$.    Let $\phi_q(x)=x^q$.
Then $I^q$
might mean one of three things:  the ideal $RI^q$ generated by $q$-fold products of
elements of  $I$;
the ideal $\phi_q(I)$ of $\phi_q(R)$; or possibly the ideal $R \phi_q(I)$ of $R$.  We will
refrain from writing $I^q$ altogether, and use the above notations.

\begin{remark}\lbl{q-bq}  Let $S$ be a ring, with ideal $M$ generated by
$s_1,\ldots,s_b$.   Assume $p=0$ in $S$, and let $q$ be a power of $p$.
Then $SM^{bq} \subset  S  \phi_q(M) \subset SM^q $ \end{remark}

\proof The ideal   $SM^{bq}$ has generators of the form 
${s_1}^{m_1} \cdot \ldots \cdot {s_b}^{m_b}$
with $m_1+\ldots+m_b = bq$.  Necessarily $m_i \geq q$ for some $i$, so
${s_i}^{m_i}$ and hence the product are in $S \phi_q(M)$. \qed

\begin{lem}\lbl{mult-2-4}  Let $S$ be a local ring of characteristic $p>0$, 
with maximal ideal $M$, and let
$q$ be a power of $p$.  Assume $M$ is finitely generated.  
Then $\phi^q(S)$ has a unique maximal ideal $\phi_q(M)$, with residue field 
 $k = \phi_q(S)/\phi_q(M)$.  Note $S$ is an $S^q$-module,  so $S/S\phi_q(M)$ is
a $k$-space.  Assume $\dimk (S / S\phi_q(M)) < q$, or even just

{\em (\#) $\dimk (S/ SM^q) = n < q$ }

Then $S$ has length $\leq n$ (as an $S$-module.) \end{lem}

\proof If $a = \phi^q(b) \in \phi^q(S)$ then $a$ is a unit if $b$ is, and otherwise $a \in \phi^q(M)$; so $\phi^q(S)$ is local.   Clearly $S/SM^q$ has length $\leq n$ as an $S$-module.  
By Lemma \ref{mult-2-3},
$SM^n=0$, so $SM^q=0$.  Thus $S = S/SM^q$ has length $\leq n$.  \qed

Let $S$ be a local ring of characteristic $p>0$, 
with maximal ideal $M$, generated by $b$ elements.   For $r$ a power of $p$,
write $S_r = \phi_r(S)$, $M_r = \phi_r(M)$; so $S_r$ is local, with maximal
ideal $M_r$ (the elements of $S \setminus M$ are units;
 hence $\phi_r(S \setminus M)$ consists of units of $S_r$, while
$\phi_r(M)$ is a proper ideal of $S_r$.  These sets are thus disjoint, while
their union equals $S_r$.)  

\begin{lem}\lbl{mult-2-5}  Let $S$ be a  local ring of characteristic $p>0$, 
with maximal ideal $M$, generated by $b$ elements.  
Let $r,q$ be  powers of $p$.
Assume $\phi_r(S/SM_{q} )$ has length $n$ as an $S_{r}$-module, and 
$brn<q$.
then $SM^{q} = 0$.   In particular $(S_rM_r)^{{n}} = 0$, and
  $S_r$ has length at most $1 + b + \ldots + b^{n-1} \leq max(n,b^n-1)$.
 \end{lem}

\proof  The assumption that $\phi_r(S/SM_{q} )$ has length $n$ refers
to the ring $S/SM_{q}$, and the function $\phi_r(x)=x^r$ in that ring.
Being a quotient of $S$,  $S/SM_{q}$ is an $S$-module, and 
$\phi_r(S/SM_{q} )$ is an $S_r$-module.  Note 
$$\phi_r(S/SM_{q}) = (\phi_r(S) + SM_{q}) / (SM_{q}) \iso_{S_r} 
   \phi_r(S)   / (SM_{q} \meet \phi_r(S))  $$

Thus $S_r / (SM_{q} \meet S_r) $ has length $n$ as an $S_r$-module.
(It is an $S_r/(M_r)^{q}$-module; this ring has finite length as a module over
itself, so 
all nontrivial finitely generated modules below have nonzero finite length.)
  Consider the $S_r$-
modules $N_i = (SM^{q} \meet S_r) + (S_r M_r)^i$.  They lie between 
$S_r$ and $SM_{q} \meet S_r $, and form a descending chain.  
As the length of $S_r / (SM_{q} \meet S_r) $ is $n$, we have $N_i=N_{i+1}$
for some $i \leq n$.   By Nakayama,  $S_r(M_r)^i \subseteq (SM^q \meet S_r)$,
hence $S_rM_r^n \subseteq SM^q$.


By Remark \ref{q-bq}, $SM^{rb} \subset SM_r$, so 
$SM^{rbn} \subset  S(M_r)^n  \subset SM^q$.  By Nakayama again, using $rbn<q$,
we obtain $SM^{rbn} = (0)$.  

In particular $(S_rM_r)^{n} = 0$.  As $M_r$ is also generated by $b$ elements,
$(S_rM_r)^l$ can be generated by $\leq b^l$ elements, so $S_r$ has length
$\leq  length(S_r/M_r) + length(M_r/(S_rM_r)^2) + \ldots 
\leq 1+b+b^2 + \ldots + b^{n-1}  $.
\qed

\end{subsubsection}
\begin{subsubsection}{Relative  reduction multiplicity}

We first observe a natural  relationship between transformally radicial
extensions and  {\em relative} reduction multiplicity.   The proof is  
less straightforward than it ought to be; one 
reason is that the notion of multiplicity of a point on a scheme is somewhat delicate,
and does not easily permit devissage.  Recall that when
$q$ is a power of the prime $p$, $K_q$
denotes a (large) algebraically closed field of characteristic $p$, endowed with the
map $x \mapsto x^q$.  For any scheme $Y$, $Y(K_q)$ denotes the set
of $K_q$-valued points of $Y$, i.e. difference scheme maps  $y: \spe K_q \to Y$.
Recall also that $X_{q,y}$ denotes $M_q(X) \times_{f,y} \spe(K_q)$.

\begin{definition}  Let $f: X \to Y$ be a morphism of difference schemes.
$X$ {\em has reduction multiplicity $\leq k$ over $Y$} if for some integer $B$,
for all large
 prime powers $q$, and all $y \in Y(K_q)$, and
 $z \in X_{q,y}$,
\[ \mult_z X_{q,y} \leq B q^k \] \end{definition}

\begin{definition}  Let  $X$ be a
difference scheme of finite type over a finitely generated difference field $K$. We will say that $X$
is generically of $k$-bounded reduction multiplicity  over $K$
if there exists a finitely generated difference domain 
$D\subset K$, and a difference scheme $X_0$ of finite type over $D$,
with $X \iso X_0 \times_{\spe D} \spe K$, such that $X_0$ is of
reduction multiplicity $\leq k$ over $\spe D$. \end{definition}

\begin{lem}\lbl{si-dim-3.0}
 Let $X$ be a difference scheme of finite type over  a Noetherian
difference scheme $Y$.  Assume $X$ is transformally
radicial  over $Y$, and of  total dimension $\leq e$ over $Y$.
Then $X$ has reduction multiplicity $\leq e$ over $Y$
\end{lem}

\proof


We assume inductively the lemma
holds for $X_{Y'} $ over $Y'$ for any proper closed
$Y' \subset Y$.   Note that if the lemma holds for $X_{Y'}$ as well as for
$X_{Y \setminus Y'}$, then it holds for $X$.   Thus  it suffices to prove the lemma
for any open subset of $Y$.  We may thus assume $Y $ is irreducible; and indeed
that $Y=  \spe R$,  $R$ a difference domain; and it suffices to prove the lemma for
some difference localization $R'$ of $R$ in place of $R$.

The case of $X$ is more delicate; we may still assume the lemma is true for any
proper closed subset; but if the lemma holds for a closed subset and its complement,
it is not clear that it holds for $X$. Still if $X$ is a union of {\em open} subsets,
it suffices to prove the lemma for each separately.  Thus 
we may take $X = \spe S$, $S$ a finitely
generated transformally radicial $R$-algebra.  We may assume $S$ is
well-mixed.

\claim{1}  Let $R'$ be a difference $R$-subalgebra of $S$,  finitely generated as
an $R$-module.  If   the conclusion of the  lemma holds for $X$ over $\spe R'$,
then it holds for $X$ over $\spe R$.

\proof For any $y \in M_q( \spe R) (L)$
and $y_1 \in (\spe R')_{q,y}$, there exists $d \in \Nn$ with
$$\mult_ {y_1 } (\spe R')_{q,y} \leq d $$
By assumption, for some $b_1$,
$$ \mult_ z X_{q,y_1} \leq b_1 q^e $$
So by \ref{mult-2}, with $b=db_1$,
$ \mult_ z X_{q,y} \leq b q^e $.
 \qed

\smallskip

Find a sequence of subrings $R = S_0 \subset S_1 \subset \ldots
\subset S_n = S$ such that $S_{k+1} = S_k[a_k]$ and $\si(a_k) \in
S_k$.
We will also use induction also on the length $n$ of this chain.

Let $L$ be the field of fractions of $R$,$S_L = L \tensor_R S$.
  Effecting a finite localization of $R$, we may assume $S$ embeds into $S_L$. (The 
ideal of polynomials over $L$  vanishing at $(a_1,\ldots,a_k)$ can be taken to be
generated by polynomials over $R$.) We will further use induction on the Noetherian
rank of $S_L$.

Consider  the ideal $I_L$ of $S_L = L \tensor_R S$ generated by elements
 $a$ with $\si(a)=0, a^2=0$.  $I_L$ is generated by finitely many elements $b_1,\ldots,b_m$.
After a finite localization of $R$, we may assume that the $b_i$ lie in $S$.
  Let $R' = R[b_1,\ldots,b_m]$.
$R'$ is finitely generated over
$R$ as an $R$-module.  By the claim, it suffices to prove the lemma for $S$ over $R'$.
Every prime of $R'$ must contain the $b_i$, so this amounts to proving the lemma
for $S' = S/I_L$.   If $I_L \neq 0$, then $S'_L = S' \tensor_R L$ has smaller Noetherian
rank than $S_L$; so using induction on this ordinal, we have the lemma for $S'$
and hence for $S$.  Thus we may assume $I_L = 0$.  

So $S_L$ has no nonzero elements $a$ with $\si(a)=a^2=0$.  (Since $S$ embeds in $S_L$, $S$
has no such elements either.)

Next, suppose $S_L$ has a zero-divisor $c$ with $\si(c) \in L$.  If $\si(c) \neq  0$,
then passing to a finite localization we may assume $\si(c)$ is a unit in $R$;
thus if $cd=0$ then $\si(d) = 0$; so replacing $c$ by $d$ if necessary
we may assume $\si(c) = 0$.  Let $J = \{s \in S_L: cs = 0 \}$.  Then $S_Lc$,$J$ are nontrivial
difference ideals of $S_L$, and after further localization of $R$, $J \meet S$
is nontrivial. The lemma holds for $S/(J \meet S)$ and for 
$S/Sc$ by Noetherian
induction.  Moreover if $e \in S_Lc \meet J$ then $e^2=0$ and $\si(e) = 0$, so $e = 0$
by the previous paragraph.  So $S_Lc \meet J = 0$.  The lemma follows for $S$.

Thus we may assume there are no zero divisors $c \in S_L$ with $\si(c) \in L$.

Consider first the case that there exists a polynomial $F \in R[X]$ with nonzero leading coefficient
$r$, such that $F(a_1)=0$.  By the remark preceding the claim, it suffices to
prove the lemma for $S[r^{-1}]$ over $R[r^{-1}]$; so we may assume $F$ is monic.
In this case $S_1$ is a finitely generated $R$-module.  By the induction hypothesis,
 the lemma holds true for $S$ over $S_1$.  But now by the Claim it holds for $S$ over $R$.

Assume now there is no such $F$; so $S_1$ as an $R$-algebra is the polynomial
ring in one variable.    Moreover for $f \in R[X]$, $f(a_1)$ is not a zero-divisor
in $S$.  And $S$ has Krull dimension $e$ over $R$.
 It follows algebraically that - after replacing $R$ by a finite localization -
{\em every } fiber of  the map of algebraic schemes $\spec S \to \spec R[a_1]$ has
dimension $\leq e-1$.

Let $Y' = \spe S_1$.   Clearly for any $q$ and any $y \in M_qY$,
and $y' \in Y'_{q,y}$,
\[\mult_{y'} Y'_{q,y} \leq q \]

Using the case $e-1$ assumed inductively, there exists $b$ such that
for any large $q$, and any $y' \in Y'(L)$,  and $z \in X_{q,y'}$,
\[ \mult_ z X_{q,y'} \leq b q^{e-1} \]

Let $y \in M_q Y (L)$, $z \in X_{q,y}$; let $y' \in Y'_{q,y}$ be the intermediate. We then
have $\mult_ z X_{q,y'} \leq b q^{e-1}$ and   $\mult_{y'} Y'_{q,y} \leq q$;
by \ref{mult-2},
\[ \mult_z X_{q,y} \leq  b q^e \]
 the lemma follows.  \qed

 \end{subsubsection} 

 \begin{subsubsection}{Transformal multiplicity and reduction multiplicity}
 
\begin{lem}\lbl{trans-mult-red} Let $X$ be a difference scheme of
finite type over a finitely generated difference domain $D$.  Assume $X$
is of generic transformal multiplicity $\leq n$.  Then there exists a
nonempty open $Y \subset \spe D$ and a bound $b_0$ such that for all
inversive difference fields $L$ and $L$ -valued points $y$ of $Y$ and
$x$ of $B_{n+1}(X_y)$, for every local ring $R$ of $(X_y)_x$,
$\si^n(R)$ is a finite-dimensional $L$-space, of dimension $\leq b_0$.
 \end{lem}

\proof Let $K$ be the inversive closure of the field of fractions of $D$.
We may take $X = \spe S$.    Let $S_K = S \tensor_D K $.

We first consider the case of generic $y \in \spe D$, i.e. of
difference field extensions $L$ of $K$.  
Let $S_L = S \tensor_D L$ (so $\spe S_{L} = X_{y}$). 
Let $x$ be an $L$-valued point of $B_{n+1}(X_{y})$, i.e. $x:
\si^{n+1}(S_{L}) \to L$ a difference ring homomorphism.  Let $x_{k } = x | {\si^{n+1}(S_K)}$.

We have a natural homomorphism 

$$ S_K \tensor_{x_{k}} L \to  S_L \tensor _{x} L $$

and it is easily seen to be surjective.  By \ref{trans-mult-n},
$\si^{n}(S_K \tensor_{x_{k}} L)$ has total dimension $0$.  Thus the
homomorphic image $\si^{n}( S_L \tensor _{x} L )$ also has total dimension $0$.
By Lemma \ref{rad-fg}, $S_L \tensor _{x} L$ is a finitely generated
$L$-algebra, hence so is $\si^{n}(S_L \tensor _{x} L)$.  By 
\ref{total-krull-0}, $\dim_{L}(\si^{n}(S_L \tensor _{x} L)) < \infty$.   

It follows by compactness that for some $d_{0} \in D$ and $b_{0}$, for
all inversive difference fields $L$ and difference ring homomorphisms $y: D \to L$
with $y(d_{0)} \neq 0$, and all $L$-difference algebra homomorphisms $x:
\si^{n+1}(S_{y}) \to L$,  $\dim_{L}(\si^{n}(S_y \tensor _{x} L)) < b_{0}$.
To see this, let $F_{0}$ be a finite set of generators of $S$
as a difference $D$ -algebra, $F= F_{0} \union \ldots \union \si^{n}(F_{0})$.  
Then the image of $F$ generates $S_y \tensor _{x} L$ as an
$L$-algebra, for any $x,y$.  Let $F_{1} = F, \ldots , F_{m+1} = FF_{m}$.
Then $\dim_{L} (\si^{n}(S_y \tensor _{x} L)) < \infty$ iff for some $m$,

(*)   for any
$c,d \in \si^{n } F_{m}$, the image of $cd$ in $S_y \tensor _{x} L$ is in the
$L$-span of the image of $\si^{n } F_{m}$.

Let $y_{i},x_{i},L_{i}$ be a sequence of such triples, with $y_{i}$
approaching the generic point of $\spe D$; let $S_{i} =S_{y_{i}} \tensor _{x_{i}} L_{i}$; we have to show that
$\dim_{L_{i}}(\si^{n}(S_{i})) $ is bounded.  Let $(y_{\infty},x_{\infty},L_{\infty},S_{\infty})$
be an nonprincipal ultraproduct of the $(x_{i},y_{i},L_{i})$.  Then
$\dim_{L}(\si^{n}(S_{y_{\infty}} \tensor _{x_{\infty}} L_{\infty})) = < \infty$,
so a fortiori the image of this ring in $S_{\infty}$ has finite
dimension over $L_{\infty}$.  Thus for some $m$, (*) above holds 
for the image of $F_{m}$ in $S_{\infty}$, hence also for the image
of $F_{m}$ in $S_{i}$, for almost all $i$.  This proves that $\dim_{L_{i}}(\si^{n}(S_{i})) $ is bounded,
as required.

The existence of the open difference subscheme $Y$ and the bound
$b_{0}$ now follow by a standard compactness argument.
 
\qed

\begin{cor}\lbl{trans-mult-red-q}   Let $f: X \to Y$ be a  
a morphism of difference schemes of finite type, of relative
transformal multiplicity $\leq n$.  Then there exists $b \in \Nn$
such that  for all prime powers $q>b$, for any perfect field $L$ 
and any
$y \in (M_qY)(L)$,   every closed point $x  \in M_q(X_y)$ has
multiplicity $\leq bq^n$ on $M_q(X_y)$.  \end{cor} 

\proof We may assume $X=\spe R$, $Y= \spe D$, $R$ a finitely generated
$D$-algebra.   By lemma \ref{gen-M}, the local rings of $M_q(X_y)$,
$y : M_q(D) \to L$, $L$ a perfect field, have maximal ideals generated by a 
bounded number of elements.  ($D$ itself is generated by a finite number $e$
of elements, so that quotient domains of $D$ have fields of fractions
$k$ with $[k^q:k] \leq q^e$, but we just consider perfect $L$, so
actually \ref{gen-M} applies with $e=0$.)  Thus Lemma \ref{mult-2-5}
applies, and translates Lemma \ref{trans-mult-red} to say that $x'= x |
\si^{n+1}(R \tensor_{y } {L})$
has bounded multiplicity on
$M_q(B_{n+1}(X_y))$.  By Lemma \ref{si-dim-3.0}, $X/(Y \times
B_{n+1}(X))$ has reduction multiplicity $\leq n$.  So $Mult_x M_{q}(X)_{y,x'}
\leq O(q^{n})$.  By Lemma \ref{mult-2}, $x$ has multiplicity $O(q^n)$ on
$M_q(X_y)$.    \qed

\begin{cor}\lbl{si-dim+}
 Let $X$ be a difference scheme of finite type over  a Noetherian
difference scheme $Y$, of total dimension $d$ over $Y$. Then
there exists $b$ such that for all large enough prime powers $q$, and all
$y \in M_q( Y)$, the zero-dimensional scheme $X_{q,y}$ over $L_y$
has  size at most $bq^d$.
\end{cor}

\proof  By the usual Noetherian induction on $Y$ we may assume $Y = \spe D$,
$D$ a difference domain with field of fractions $K$; and  we may pass
to a localization of $D$ by a finite set.   Moreover base change will
not change the total relative dimension or the size, so we may
replace $D$ by a finitely generated extension   within the inversive hull
$K'$ of $K^{alg}$.  Let $X' = X \tensor_Y \spe K'$.  We enlarge $D$ within
$K$ so that the $Mlt_k X'$ and their components $X_{k,j}$ are defined
over $D$.  The components $X_{k,j}$ have reduced total dimension $\leq d-k$;
so by \ref{si-dim}, for all $y \in Y$ and all sufficiently large $q$,
$(X_{k,j})_{y,q}$ has $O(q^{d-k})$ points, and by \ref{si-dim3.2},
away from $Mlt_{k+1}(X)$, these points have multiplicity $\leq O(q^k)$
on $M_qX$.  Thus $X_{q,y}$ is divided into $d$ groups $(Z_kX)_{q,y}$,
 the $k$'th
having at most $O(q^{d-k})$ of multiplicity (on $M_qX$) at most
$O(q^k)$; so altogether there are $O(q^d)$ points, multiplicity
counted.

The corollary \ref{si-dim+} can also be obtained using Lemmas
\ref{direct-ritt-0} and  Bezout theorem methods (see \ref{fin}).  One could
take things up from that point with \ref{si-dim-3.0}.  However, the   methods of this section     
give information about the distribution of multiplicities (\ref{si-dim++}), and in particular allow
 estimating the multiplicity on $X$ of points on smaller-dimensional $X' \subset X$, 
(\ref{si-dim+++}).  This does not seem apparent from
the Bezout approach (though one can use Bezout methods to estimate the multplicities of these points on $X'$.)

\begin{cor}\lbl{si-dim++}   
 Let $X_0$ be a difference scheme of finite type over a difference
 field $K$.  Assume $X_0$ has total dimension $d$, and 
that no algebraic component of $X_0$ of total dimension $d$ has
transformal multiplicity $>0$.  
Then there exist
a finitely generated subdomain $D$ of $K$, a difference scheme $X$ over $D$
with $X_{0}= X \times _{\spe D} \spe K$,  and an
integer $b$, so that   for all large
enough prime powers $q$, and all $y \in   M_q( \spe D)$, the
zero-dimensional scheme $X_{q,y}$ over $L_y$ has size at most $bq^d$;
moreover all but $bq^{d-1}$ points of $X_{q,y}$ (counting
multiplicities) have multiplicity $<b$.
\end{cor}

\proof As in the proof of \ref{si-dim+}, we consider reduced subschemes $Y$ of $Z_kX = Mlt_k(X) \setminus Mlt_{k+1}(X)$.
Here however, for $k \geq 1$, we use \ref{si-dim3.2}(3) to conclude
that 
the reduced total dimension of $Y$ is $\leq d-k-1$.  
 This   gives an 
$O(q^{d-k-1})$-bound on multiplicities, by \ref{si-dim3.2}.  As we are considering $Y$ away from $Mlt_{k+1}$, 
we obtain $O(q^{d-1})$ points counted with their multiplicities on $X$.     \qed

\begin{cor}\lbl{si-dim+++}  Let $X_0$ be as in Corollary \ref{si-dim++}.  
For any proper subscheme $X'$
of $X$, of total dimension $d'<d$,  
the number of points of $X'$, counted with their
multiplicities on $X$,  is $O(q^{d-1})$. \end{cor}
\proof  By Corollary \ref{si-dim++}, only $O(q^{d-1})$ points have high multiplicity on $X_{q,y}$.
These can therefore be ignored.  The rest have multiplicity $O(1)$ on $X_{q,y}$, and by
Corollary \ref{si-dim+}, their number is $O(q^{d'}) $.  Thus even with multiplicity,
the number is $O(q^{d-1})$.
\smallskip

\end{subsubsection}

\end{subsection}  


 \end{section}
\end{part}
\begin{part}{Intersections with Frobenius}

\begin{section}{Geometric Preliminaries}    \lbl{geopre}
 
We show here that the main theorem is invariant under
birational changes (\ref{birational}) and in the appropriate sense
under taking finite covers (\ref{birat2}). We use this together with
de Jong's version of resolution of singularities to reduce to the smooth 
case,
where intersection theory applies. We deal with the possible
inseparability, and note another proof of a  crude initial bound
(\ref{fin}) on the set whose size we are trying
to estimate.

Our basic reference, here and in later sections, is \cite{fulton}.
We will begin with some basic lemmas on proper and improper
intersections, degrees and correspondences.

\begin{subsection}{Proper intersections and moving lemmas}
\lbl{moving}

Let $X$ be a nonsingular variety over an algebraically closed field.
We write $\codim{X}{U}$ for $\dim(X)-\dim(U)$. A variety (or cycle) is 
said to have
pure (co)dimension $k$ if each irreducible component has (co)dimension 
$k$.
Two  subvarieties $U$,$V$ of $X$ of pure codimension $k$,$l$ 
respectively are said to meet
{\it properly} if either $U \meet V = \emptyset$, or $\codim{X}{U \meet 
V} \geq k+l$.
By the dimension theorem for smooth varieties,
each   component of $U \meet V$ must have codimension at most $k+l$;  
thus the intersection
is proper iff each (nonempty) component of $U \meet V$ has codimension 
$k+l$ in $X$.  The properness
of the intersection is equivalent to the properness of intersection of 
all irreducible
components of $U$ with those of $V$.

\begin{lem}\lbl{proper1}  Let $U,V,W$ be pure-dimensional varieties.
 Assume $U,V$ meet properly, and $U \meet V, W$ meet properly.  Then $U$ 
meets
$V \meet W$ properly.
\end{lem}

\proof
Let $U,V,W$ have codimensions $k,l,m$.  Then $U \meet V = \emptyset$ or
$\codim{X}{ U \meet V} = k+l$; so
$U \meet V \meet W = \emptyset$ or $\codim{X} {U \meet V \meet W} = 
k+l+m$; this last
condition is equivalent to the assumptions, and is symmetric in $U,V,W$.

\begin{lem}\lbl{proper2}  Assume $X,Y$ are smooth varieties.  Let $T$ be 
an irreducible subvariety of $X \times Y$,
and let $\pi[k] T$ be the subvariety of $X$, whose points are
$\{p \in X:  \dim( T \meet (\{p\} \times Y) ) \geq k \}$.  If $U$ is a 
subvariety of
$X$ meeting each component of each $\pi[k] T$ properly, then $U \times 
Y$ meets $T$ properly.
If $U \times Y$ meets $T$ properly, then $U$ meets $\pi T$ properly.
\end{lem}
\proof
Note that $\codim {X \times Y} {U \times Y} = \codim{X}{U}$; let $c$ be 
the common value.
Let $W$ be an irreducible component of $(U \times Y) \meet T$.  Say
$\dim(W) = \dim( \pi W) + k$.  Then $\pi W \subset \pi[k] T$.  So
$\pi W   \subset \pi[k] T \meet U$,and
$  \dim(\pi W  ) \leq \dim(  \pi[k] T \meet U) \leq \dim ( \pi[k] T ) - 
c
\leq \dim(T) -k -c $
Thus $\dim(W) \leq \dim(T) -c$.  This shows that $U \times Y$ meets $T$ 
properly.

For the remaining statement, let $k = \dim(T) - \dim( \pi T)$.  Then 
$\pi T = \pi[k] T$.
So  $\dim(T) - \codim{X}{U}  \geq \dim( (U \times Y) \meet T)
\geq k + \dim( U \meet \pi T)$; hence $\dim( \pi T) - \codim{X}{U} \geq 
\dim( U \meet \pi T)$.

\begin{lem}\lbl{proper3}  Let $U \subset (X \times Y)$, $V \subset (Y 
\times Z)$
be complete varieties, all of pure dimension $d$.  Let $W = pr_{XZ} ( (U 
\times Z) \meet (X \times V) )$.
Let $T$ be a subvariety of $X \times Z$.  Assume

$U$ meets each component of $X \times pr_Y[l] V$ and of $pr_X[l] T 
\times Y$ properly for each $l$;

$V$ meets each component of $pr_{YZ} [l] (   \left( T \times Y) \meet ( 
U \times Z)   \right)$ properly;

Then $U \times Z$  meets $X \times V$ properly, and  $W$ meets $T$ 
properly.
\end{lem}

\proof  \begin{enumerate}
\item  By the first assumption and \ref{proper2}, $U \times Z$ meets $X 
\times V$
properly.
\item   Similarly, $T \times Y$ meets $U \times Z$ properly.
\item   By the second assumption, $X \times V$ meets $\left( ( T \times Y) 
\meet ( U \times Z)    \right)$
properly.
\item   By \ref{proper1}, and (2),(3), $T \times Y$ meets $ (U \times Z) 
\meet (X \times V)$
properly.
\item  By the last statement
of \ref{proper2}, $T$ meets $W$ properly.
\end{enumerate}

Let $X$ be a smooth algebraic variety over
an algebraically closed field.  A {\em cycle} is a formal sum,
with integer coefficients, of irreducible subvarieties of $X$.

To each subscheme $U$  of $X$, one associates
a cycle $[U]$;  it is the sum of the irreducible components of $U$,
with certain nonnegative integer coefficients   (\cite{fulton}, 1.5).

\begin{notation} 
 For any cycle $I$, write
$I = \Sum_k \{I\}_k$, where $\{I\}_k$ is a $k$-dimensional cycle.

$B^i(X)$ denotes the group of  cycles whose components have
codimension $i$, up to algebraic equivalence (\cite{fulton} 10.3).
 An operation $\cdot$
can be defined on $B^*(X) =  \oplus_i B^i(X)$, making it into a commutative graded
 ring with unit.  This operation is determined by
the fact that $[U] \cdot [V] = [U \meet V]$ when $U,V$ are irreducible 
subvarieties
of $X$ meeting properly.
\end{notation}

Intersection numbers depend only on the algebraic equivalence class of a cycle.  
In particular,
a $0$-cycle $\sum n_i p_i$ is determined up to algebraic equivalence
by  $\sum n_i \in \Zz$; we will identify the group of $0$-cycles
with $\Zz$.  

 The moving lemmas are available for a   rational equivalence, where the variety parametrizing the equivalence is $\Pp^1$.
 It would be possible to work with  rational equivalence, but this would require a delicate  treatment of multiplicities later on;   see Remark \ref{macaulay}.
To avoid this, we use a   still stricter equivalence, where only reduced cycles appear.  Many thanks to Yves Laszlo for clarifying these issues.

 \subsection{Reduced rational equivalence}

   Let $Y$ be a smooth variety over the perfect field $k$.  A subscheme is {\em generically smooth} if
 it has dimension $n$, and is smooth away from a closed subset of dimension $<n$.  For varieties over a perfect field, this is equivalent to {\em generically
 reduced.}    For zero-dimensional subschemes, this implies that each point has multiplicity one.  In particular if $X$ is zero-dimensional,  irreducible
 and generically smooth over $k$, then $k(X)$ is a finite separable extension of $k$.
  
   Let $A_n(Y)$ be the group of formal linear combinations of $n$-dimensional subvarieties of $Y$.  Let $\B_n(Y)$ be the subgroup generated
by all cycles of the form $[T(\infty)]-[T(0)]$, where $T$ is a closed, irreducible subvariety of $Y \times \Pp^1$, projecting onto $\Pp^1$, {\em such that the schemes $T(\infty),T(0)$ are generically smooth.}
Here  $T(a)$ denotes the scheme $\pi \inv(a)$, where  $\pi: Y \times \Pp^1 \to \Pp^1 $ is the second projection.  An element of $\B_n(Y)$ will be said to be
{\em reduced- rationally equivalent to zero}. 

In this definition, $\Pp^1$ may be replaced by any Zariski open subset of  $\Pp^1$ containing $0,\infty$; one can return to $\Pp^1$ by  taking Zariski closure.   

We also observe that if $Z$ is an smooth subvariety of $Y$, in the above situation, and if 
 $T(\infty) \meet Z$  and $T(0) \meet Z$ are generically transverse intersections of dimension $n$, then  $[T(\infty) \meet Z] - [T(0) \meet Z] \in \B_n(Z)$
To see this let $T' = T \meet (\Pp^1 \times Z)$.  Then any component of $T'$ has dimension $\geq n+1$; if there are any vertical components,
we may safely remove them, since they cannot be at $0$ or $\infty$.  Let  $C_1,\ldots,C_j$ be the components of $T'$ 
 $T'$ projecting onto $\Pp^1$.  Then $\dim (C_i(0) ) \geq n$
so equality holds, and $C_i(0)$ must be one of the components of $T(0) \meet Z$.  On the other hand any component of $T(0) \meet Z$ must be
contained in some component of $T \meet Z$ of dimension $\geq n+1$; this must be one of the $C_i$.  It follows easily that $[C_i(0)]-[C_i(\infty)]  \in \B_n(Z)$
and that $[T(\infty) \meet Z] - [T(0) \meet Z] = \sum_i [C_i(0)]-[C_i(\infty)] $.  The same holds over $f$, if we assume that $T(\infty) \meet Z, T(0) \meet Z$
are $f$-admissible.  


$\B_n(Y)$  is a subgroup of the group of cycles rationally equivalent to $0$; where only generically smooth subschemes are allowed.   (It may well be equal
to the group of cycles rationally equivalent to zero.) 

We will also need a relative statement.  Assume $f: Y \to {\bY}$ is a given morphism, whose generic fibers are generically smooth. 
Say a subvariety $W$ of $Y$ is $f$-admissible if $W$ is generically smooth, and moreover, for a generic $v \in {\bY}$,  $T(\infty) \meet f \inv(v)$ is nonempty and generically smooth.   

In practice the fibers $f \inv (y)$ will be smooth; in this case at least, the condition on intersections with generic fibers
of $f$ implies generic smoothness of $Y$ itself.  Indeed if $[Y] = m [Y_{red}]$ and $Z$ is a local complete intersection then $[Y \meet Z ] = m [(Y \meet Z)_{red}]$(this follows inductively when $Z$ is the intersection of $n$ hypersurfaces from \cite{fulton} 2.4); we can take $Z= f\inv(v)$, so generic smoothness
and nonemptyness of $ Y \meet Z$ implies $m=1$, hence $Y=Y_{red}$ and (being over a perfect field) $Y$ is generically smooth.

  Let $\B_n^f(Y)$ be the subgroup generated
by all cycles of the form $[T(\infty)]-[T(0)]$, where $T$ is a subvariety of $Y \times \Pp^1$  such that every component of $T(\infty),T(0)$ 
is $f$-admissible.   Two $f$-admissible  cycles are said to be   reduced- rationally equivalent over $f$ if they are $f$-admissible, and their difference is in $\B_n^f(Y)$.

    As above, when the fibers of $f$ have complementary dimension to the components $C$ of $T(0)$ in $Y$,
this implies that the morphism $f|C$ is separable.

\<{lem}  \lbl{moveall}  Let $Y$ be a smooth projective variety over an 
algebraically closed field $k$,
and let $U$ be an irreducible subvariety of $Y$.  
There exists a purely transcendental extension $K$ of $k$, and a cycle
$U'$ on $Y$ defined over $K$, such that $U,U'$ are reduced- rationally equivalent,
and for any subvariety $W$ of $Y$ defined over $k$, each component $C$ of $U'$
meets $W$ properly.   

If  $f: Y \to {\bY}$ a   morphism, and $U$ is $f$-admissible, then $U'$ can be taken to be $f$-admissible, and reduced- rationally equivalent to $U$ over
$f$.  
\>{lem}

\proof In this proof, all cycles are homogeneous of  dimension $\dim(U)$.   If no $f$ is given, the ``over $f$" part of the statements should be ignored.


\claim{}  
There exists a cycle $U[1]$  defined over a purely transcendental extension of $k$, reduced-rationally equivalent to $U$ over $f$,
such that for any subvariety $W$ of $Y$ defined over $k$,
each component $C$ of $U[1]$ either meets $W$ properly, 
or satisfies $\dim(C \meet W ) < \dim (C'  \meet W)$ for some component
$C'$ of $U$.    

\proof 
 This is a version of Chow's moving lemma proved in  \cite{hoyt} (see also
\cite{fulton} 11.4.1.)  We briefly recall the proof.   

Say $Y \subset \Pp^{N_0}$; for the transversality and surjectivity arguments later,
we take a $d$-uple embedding of $\Pp^{N_0} \to \Pp^N$ for large enough $d$.  \footnote{$d \geq 3$ is required by \cite{hoyt}; we also need  any four distinct points of $\Pp^{N_0}$   linearly  independent in $\Pp^N$, in order to apply Lemma \ref{moveall-tl}.}  We use the composed embedding $Y \subset \Pp^N$.  

Let $G$ be the Grassmanian
of linear subspaces of $\Pp^N$ of codimension one more than the codimension of $U$ in $Y$.
Pick a point  $l$ of $G$ generic over $k$, representing a
linear subspace $L$.  
Given any $k$-subvariety $U$ of
$Y$, form the cone $C(L,U )$.  Away from $L$, this cone can be described as follows.  Let $H$ be a generic hyperplane containing $L$.
We view the complement $\Pp^N \m H$ as an affine space, with corresponding vector space $V_L$.   
 The union of lines through any fixed point $0 \in \Pp^N \m H$ and a point of $L$, restricted to $\Pp^N \m H$, is a vector subspace $E_L \leq V_L$.
 Let $\pi_L: (\Pp^N \m H) \to (\Pp^N \m H)/E_L$ be the natural linear map.  Then $C(L,U) \m H = \pi_L \inv \pi_L (U \m H) $ (the union of all lines through $U \m H$ whose closure meets $H$ in a point of $L$.)  Also $L$ is a generic linear subspace of $H$
 of codimension =$codim_Y U$.   Thus $E_L$ is a generic  subspace of $V_L$ of   dimension $\dim(L)$.  By \cite{hoyt}, Lemma 4, $C(L,U)$ meets $Y$
 transversally; we have $[C(L,U) \meet Y] = [U] + \sum_i  [U'']_i$ for some generically smooth subvarieties $U''_i$ of $Y$.

 Each $U''_i$ is $f$-admissible:   The argument that $f$ maps each component
  onto ${\bY}$ is given separately, in \ref{moveall-tl}.      The
   generic transversality of the intersection of $U''_i$ with a generic fiber $f \inv(\by)$ of $f$ can be seen as follows.  Let ${\bY}_0$ be the set of all $\by \in \bY$
  such that $f \inv(\by)$ meets $U$ generically transversally.    By admissibility of $U$, ${\bY}_0$ is Zariski dense.  By   \cite{hoyt}, Proposition 6, 
  for any $k$-rational point $\by$ of 
  ${\bY}_0$, the intersection $f \inv(\by) \meet U''_i$ is also proper and generically transverse.   (In the notation there, $I'(U''_i,f \inv(\by)) = \emptyset$.)   Since ${\bY}_0 (k)$ is Zariski dense in ${\bY}_0$, 
  it follows that for a generic $\by \in {\bY}_0$, $f \inv(\by) \meet U''_i$ is generically transverse.  
  

Now pick also  $g \in PGL_N$, generic over $k(l)$,
and consider   the translate $gC(L,U)$.  
  Then $gC(L,U) \meet Y$ is generically transverse, and $gC(L,U)$ has a transverse, nonempty intersection
with a generic fiber of $f$.  These facts   are easier versions of the corresponding statements about generic cone intersections above; they also 
   follow from Bertini's theorem, \cite{kleiman-transversality}, Corollaries 11 and 12, applied to the nonsingular part of $C(L,U)$.  \footnote{Kleiman speaks about regularity
   rather than smoothness, but applying the theorem over the perfect closure gives the same statement with smoothness.}  Surjectivity of $f$ on $gC(L,U)$ 
   is again part of the next lemma.  Thus   Let  the components $\{U'_j \} $ of $gC(L,U) \meet Y$ are $f$-admissible, and  $[gC(L,U) \meet Y= \sum   U'_j$.
   
 Let $U[1] = \sum m'_j U'_j - \sum m_i U''_i$.

It remains only to exhibit the rational equivalence between $C \meet Y, gC  \meet Y$ \footnote{in a way involving no cancellation of additional  subvarieties of $Y$},
where $C$ is a component of $C(L,U)$.  
It is easy to find a morphism $\g: R_0 \to PGL_N$ with $R_0$ a Zariski open subset of $\Pp^1$, $0,\infty \in R_0$, $\g(0)= 1$, $\g(\infty) = g$.  
Let $T  \leq R_0 \times \Pp^N$ be defined by:  $(t,y) \in T  \iff \g(t) \inv y \in C$.  Then $T$ is closed and irreducible, $T(0) \meet Y, T(\infty) \meet Y$
are $f$-admissible, so the remark below the definition of reduced-rational equivalence applies.   
  \qed

   If $U$ is not a variety but a cycle, a formal sum
$\sum U = n_i U_i$ of varieties of the same dimension, define
$U[1]  = \sum n_i U_i[1]$.  

Let $U[2]=U[1][1],U[3]=U[2][1]$, etc.   Then 
 $U[m]$   is  $f$-  reduced- rationally equivalent to $U$, is defined over a purely transcendental extension of $k$, and
for any subvariety $W$ of $Y$ defined over $k$,
each component $C$ of $U[m]$ either meets $W$ properly,
or satisfies $\dim(C \meet W ) \leq \dim (C'  \meet W)-m$ for some component
$C'$ of $U$.   It follows that $U'=U[m+1]$ satisfies our requirements.

\qed

\<{lem}\lbl{moveall-tl}  Let $U \subset Y \subset \Pp^N$ be   projective varieties, with no three or four points of $Y$ linearly
dependent in $\Pp^N$.  Let $f: Y \to V$ a flat morphism to an irreducible variety.
Let $\dim(Y)=n $, $dim(U)=d = \dim(V)$, $n+1<N$. 

\begin{enumerate}
  \item  Let $G$ be the Grassmanian of linear subspaces of $\Pp^N$ of codimension $n+1$.
For $L \in G$, let $C(L,U)$ be the cone on $U$ with center $L$.  Then for generic $L \in G$, and any component $U'$ of
$C(L,U) \meet Y$ other than $U$, $f(U') = V$.  

  \item  Let $W$ be any subvariety of $\Pp^N$   of codimension $n-d$, and let $g \in PGL_N$ be generic.  Then
for any component $U'$ of $gW \meet Y$, $f(U')=V$.

\end{enumerate}
\>{lem}  

\proof   
(1)   Let

\[   M_2=\{(L,y,y',p,p') \in G \times (Y \setminus U)^2 \times (\Pp^N)^2:   p \neq p' \in L, y \in C(p,U),y' \in C( p',U), fy=fy'\} \]

Then $M_2$ projects to $G \times Y^2$, and the image of the projection contains 

\[M =  \{(L,y,y' ) \in G \times (Y \setminus U)^2 :  y \neq y' \in C(L,U), fy=fy' \}    \]

Indeed if $(L,y,y') \in M$, then there exist $p,p' \in L$ and $u,u' \in U$ with $y,u,p$ and $y',u', p'$ colinear.  Since
no three or four distinct points of $Y$ are linearly dependent in $\Pp^N$, we must have $p \neq p'$.

So $\dim(M) \leq   \dim(M_2)$. 

Let $G_{p,p'} = \{L \in G: p,p' \in L \}$; we have $\dim G_{p,p'} + 2N \leq \dim(G) + 2 (N-(n+1)) $.

  We will also  use: 
 if $y \in C(p,U) \m U$, then $p \in C(y,U)$, and $\dim C(y,U) = \dim(U)+1$.   

 Now compute,   using the maps $(L,y,y',p,p') \mapsto (y,y',p,p) \mapsto (y,y') \mapsto fy$:

\[ \dim(M_2) \leq  \dim(V) + 2 (\dim(Y)-\dim(V))  + 2(d+1) + \dim G_{p,p'}  = \]
$$= d+2(n-d)+2(d+1)+ (\dim(G) - 2(n+1) ) = d + \dim(G).   $$

Thus we see that for generic $L \in G$, $\dim  \{(y,y') :  (L,y,y') \in M \} \leq d$.  

Now let $U'$ be  a component  of $C(L,U)$ not containing $U$, and  let $\tilde{U} = U' \setminus U$.
 If $\dim f(U') < d$, then $\dim \{(y,y') \in \tilde{U}^2: fu=fu' \} \geq (\dim f(U')) + 2( \dim(U') - \dim f(U')) > d$,
a contradiction.  So $\dim   f(U') = d$.  Since $V$ is irreducible and $U'$ is projective,   $f(U')=V$.  

 For (2), a similar computation shows that $\dim \{(g,y,y'): y \neq y' \in Y \meet gW, fy=fy'  \} \leq d$, and we conclude as above.


 \qed

  If $f: S \to V$ is a morphism of varieties, and $U$ a subvariety of $S$, recall the notation:
$f[1](U) := \{v \in V: \dim f ^{-1} (v) \geq 1 \}$.  

\<{lem}\lbl{move2}   Let $V$ be a smooth projective variety over a
difference field $k$, and let   $S$ be a subvariety of $Y = (V \times V^\si)$,
$\dim(S) = \dim(V) = d$, with $ \pr_1: S \to V^\si$   dominant (respectively dominant and separable.)

 There exists  a cycle
$S' = \sum_i m_i [U_i]$ on $Y$ defined over $k$, such that 
$S,S'$ are   reduced- rationally equivalent, and 
    $$[\si_k] pr_V [1] (U_i) \meet (U_i \star \Sigma)  = \emptyset$$

Moreover, each $U_i$ as well as the varieties involved in the rational equivalence of $S$ with $S'$  
all have dominant (resp. dominant and separable) projections to  $V^\si$.     

  \>{lem} 
\proof  Using \ref{moveall}, find a purely transcendental extension $k(s)$
of $k$ and a cycle $S'$ on $Y$, reduced-rationally equivalent to $S$, such that
any subvariety of $Y$ defined over $k^{alg}$ meets $S'$ properly, and such that
each component of $S'$ maps dominantly to $V$.
We will show that $S'$ has the required properties.   It then follows
easily, by specializing $s$ into $k$ (avoiding finitely many proper 
Zariski closed sets), that such a cycle and such a rational equivalence
exist over $k$ too.

 Let $s_0=s$, and
$K = k(s_0,s_1,s_2,\ldots)$, $\si^n(s_0)=s_{n}$ 
(with $s_0,s_1,\ldots$ algebraically independent over $k$.)  Let  $V_n = \si^n(V)$; note $V_n$ is defined over $k$.

 Let $S_0$ be any component of $S'$, $S_n = \si^n(S_0)$.  
 By \ref{moveall} applied to the inversive hull of $k^{alg}$,   any 
 subvariety of $\si^n(Y)$ defined over $ k^{alg}  $  
meets $S_n$ properly .   Since $\si^n(Y)$ and $S_n$ are defined over $k(s_{n})$, and
$k(s_0,\ldots,s_{n-1} ,s_{n+1},s_{n+2},\ldots)^{alg}$ is linearly free from $k^{alg}(s_{n})$ over $k^{alg}$,
in fact any  subvariety of $\si^n(Y)$ defined over $k(s_0,\ldots,s_{n-1} ,s_{n+1},s_{n+2},\ldots)^{alg}$
 meets $S_n$ properly.  

{\bf Claim}  Let $W_1$ be any proper subvariety of $S_0$, defined over $k(s_0)$.  Then 
$ W_1 \star \Sigma = \emptyset$. 
To prove the claim,  let $c = (d-\dim(W_1))/2 > 0$.  We will inductively define $W_n$   satisfying:

\begin{enumerate}
  \item $W_n \subset V_0 \times \ldots \times V_{2^n-1}$
  \item  $\dim(W_n) \leq d-2^n c$ (or $W_n = \emptyset$).
  \item  $W_n$ is defined over $k(s_0,\ldots,s_{{2^n-2}})$
  \item  $( W_1 \star \Sigma)[{2^n-1}] \subset W_n$
\end{enumerate}

For $n=1$,
(1-4) hold by assumption.  Assume they hold for $n$.  Then $(W_n)^{\si^{2^n}} \subset V_{2^n} \times \ldots \times V_{2^{n+1}-1}$,
has dimension $\leq d-2^nc$, and is defined over $k(s_{2^{n}},\ldots,s_{2^{n+1}-2})$.

Thus   $W^*=W_n \times (W_n)^{\si^{2^n}}$    has dimension $\leq 2d-2^{n+1}c$, and is defined
over $k(s_0,\ldots,s_{{2^n-2}},s_{2^{n}},\ldots,s_{2^{n+1}-2})$.  Let $\pi: (V_0 \times \ldots \times V_{2^{n+1}-1}) \to  (V_{2^n} \times V_{2^n+1})$
be the projection.  Then $\pi[k]W$ is defined over $k(s_0,s_{{2^n-2}},s_{2^{n}},s_{2^{n+1}-2})$, hence meets 
 $S_{2^n-1}$ properly; by  \ref{proper2}, $\pi^{-1}S_{2^n-1} $ meets $W^*$ properly.   Let $W_{n+1} = \pi^{-1}S_{2^n-1} \meet W^*$.
Then $\dim W_{n+1} \leq \dim(W^*)-d = d-2^{n+1}c$, and (1)-(4) hold.  This proves the claim.  

When $2^nc >d$, (2) forces $W_n=\emptyset$, so by (4), $ W_1 \star \Sigma  = \emptyset$.  

In particular, as $pr_0: S_0 \to V$   dominantly, and $\dim(S_0)=\dim(V)$,  
  $  pr_0[1](S_0) \neq V$; so letting  $W_1 = S_0 \meet {pr_0}^{-1}(pr_0[1](S_0))$, we have $\dim(W_1)<d$, and we can
apply the claim. 
 
\qed

\<{remark} \  \lbl{move-rem} \>{remark}   In particular, in \ref{move2} one can conclude:
$[\si]_K U_i \star \Sigma$  has transformal dimension
$0$.  (See \ref{zd}).

An alternative treatment of a moving lemma for difference schemes can be given by  
transposing the  proof of \cite{fulton},11.4.1 to difference algebra, almost verbatim.
 The methods used there - projective cones,   moving via
 the group of automorphisms of projective space, ''counting constants'' - work very well for difference varieties and transformal dimension in place of  varieties and dimension.  One can also use the dimension growth sequence to get finer results.

\end{subsection}

\begin{subsection}{Degrees of cycles}

In general, the intersection product gives no direct information about
improper intersections.  In  products of projective spaces, we can
obtain such information,   as in \cite{fulton}, Example 8.4.6.

\begin{notation}
Let $H,H'$ be hyperplane
divisors on $\Pp^n$, $\Pp^m$, respectively, $s=pr_1^* H$, $t=pr_2^*H'$.
For a subvariety $U$ of $\Pp^n \times \Pp^m$, $U = \Sum a_{ij}s^it^j$ as
 a cycle up to rational equivalence; the $a_{ij}$ are
called the bidegrees. (\cite{fulton}, Example 8.4.4).   If $U$ is of 
pure
dimension $k$ and $i+j=k$, we have
$a_{ij} = (U \cdot s^{n-i}t^{m-j})$. By taking a representative in 
general
position, it follows that the bidegrees $a_{ij}$ are non-negative 
integers.

 In a product
of more than two projective spaces, multi-degrees are defined 
analogously.

A divisor $H$ on a variety $X$ is said to be {\em very ample} if it is 
the pullback
of a hyperplane divisor of projective space, under some projective 
embedding
of $X$.   The projective degree of  $U \in B^i(X)$ under this embedding 
can then be expressed
as $U \cdot_X H^i$, and written $\deg_H (U)$.

We will denote projection from a product $X \times X' \times \ldots$
to some of its factors $X \times X'$ by $pr_{X,X'}$.
\end{notation}

\begin{notation}\lbl{nota3}
 If $A,B$ are cycles on $P$,  a multi-projective spaces, write
$A \leq B$ if $B-A$ is effective.
Equivalently, if $\{V\}^{ijk} = V \cdot h_1^{i}h_2^{j}h_3^{k}$
are the multi-degrees of a cycle $V$ on $P$, $$\{A\}^{ijk} \leq 
\{B\}^{ijk}$$
for each $i,j,k$. Note that this partial order is preserved
by the intersection product.
\end{notation}

Let $H_i$ be the hyperplane divisor on 
$\Pp^n$,
$H = pr_1^* H_1 + pr_2^* H_2$.   Then $H$ is very ample.
   (It corresponds to the Segre embedding of $\Pp^n \times \Pp^n$ in 
$\Pp^N$.)

The following lemmas will be used in \S \ref{virtual}.

\begin{lem}\lbl{proj0} Let $D_i$ be a very ample divisor on a projective
variety $X_i$, $Y = X_1 \times X_2$, $\dim(X_i) = d_i$,
$D = {\pr_1}^* D_1 + {pr_2}^* D_2$.  Let $\pi = pr_1: Y \to X_1$ be the 
projection.
\begin{enumerate}
\item  Let $U$ be a $k$-dimensional subvariety of $X_1$.  Then
$\deg_D  (U \times X_2) \leq   {k+d_2 \choose k} \deg_{D_2}(X_2)  
\deg_{D_1}(U)$
\item  Let $W$ be an $l$-dimensional subvariety of $Y$.  Then
$\deg_{D_1} (\pi_* W) \leq  \deg_D (W)$
\end{enumerate}
\end{lem}

\proof  \begin{enumerate}
\item  $\dim(U \times X_2) = k+d_2$,  and
 $\deg_D(U \times X_2) = ({\pr_1}^* D_1 + {pr_2}^* D_2)^{k+d_2} (U 
\times X_2) =
= \sum_{i+j=k+d_2} {i+j \choose i} ({D_1}^i {D_2}^j \cdot U \times X_2)$ 
.  The
only nonzero factor is $i=k,j=d_2$.
\item  $\deg_D(W) =\sum_{i+j=l} {l \choose i}  ({{pr_1}^*D_1}^i 
{pr_2}^*{D_2}^j) \cdot W \geq
{pr_1}^*{D_1}^l \cdot W = {D_1}^l \cdot pr_* W$, using  the  projection 
formula for
 divisors (\cite{fulton} 2.3 ).  The last quantity equals $\deg _{D_1} 
(\pi_* W)$.
\end{enumerate}
\qed

We can obtain some information about a
proper intersection  in a smooth variety $Y$  by viewing it
as an improper intersection in projective (or multiprojective) space.
Note that two rationally
equivalent cycles on $Y$ are a fortiori rationally equivalent on
the  ambient projective  space, so have the same degrees.

\begin{lem}\lbl{trans1}  Let $Y$ be a smooth  subvariety of projective 
space $\Pp^m$.
Let $U,V$ be properly intersecting
subvarieties of $Y$,  $W = U \cdot_Y V$.
 Then $\deg (W) \leq  \deg(U) \deg(V) $
\end{lem}

\proof  We may represent $W$ by an effective cycle $\sum m_i W_i$, where
$W_i$ are the components of $U \meet V$ and
$m_i = i(W_i,U \cdot V; Y) $ are the intersection
multiplicities.

 Let $L$ be a generic linear subspace of $\Pp^m$
of dimension $m - \dim(Y) -1$, and let $C$ be the cone over $U$
with vertex $L$, cf. \cite{fulton} Example 11.4.1.  $C$ is a subvariety
of $\Pp^m$, ``union of all lines meeting $U$ and $L$".

Note that $\deg(C) = \deg(U)$:  take
a generic linear space $J$ of dimension complementary to $C$ in $\Pp^m$,
so that  $J$  meets $C$ transversally in $ \deg(C)$ points.
The cone $E$ on $J$ with center $L$ is (over the original base field)
a generic linear space, and so meets $U$
 transversally,
in $\deg(U)$ points.  But both numbers equal the number of pairs
$(p,q) \in J \times L$ such that $p,q,r$ are colinear for some $r \in 
U$.

 By
\cite{fulton} Example 11.4.3, each $W_i$ is a proper component
of $C \meet W$, and $m_i = i(W_i,C \cdot V; \Pp^m)$.

By the refined Bezout theorem
\cite{fulton} Example 12.3.1, $\sum_i m_i \deg(W_i) \leq 
\deg(C)\deg(V)$.
Thus $\deg(W) \leq \deg(U) \deg(V)$.

\begin{notation} \lbl{not2.0}  Let $D$ be a very ample divisor on $Y$,
and let $U$ be a cycle, or a rational equivalence class of cycles. 
Define
$$|U|_D = \sup_S \inf_{U',U''} \deg_D(U') + \deg_D(U'') $$
where $S$ ranges over all cycles of $Y$ , and $U',U''$ range over all 
pairs
of effective cycles such that $U$ is rationally equivalent to $U'-U''$, 
and $U',U''$
meet   $S$ properly.
  \end{notation}

\begin{cor}\lbl{mult1} Let $Y$ be a smooth variety, $D$ a very ample
divisor on $Y$.
Let $|X|=|X|_D$.   Let $U$,$V$ be cycles on $Y$.
  Then $|U \cdot V| \leq |U||V|$
\end{cor}

\proof
Let $W$ be an arbitrary cycle on $Y$, of dimension
 $\dim(Y)-\dim(U)-\dim(V)$.
Write $V = V' - V''$, with $V,V'$ effective and
meeting $W$ properly, and
$|V| = \deg(V')+\deg(V'')$.  Write $U = U' - U''$ similarly, with
$U',U''$ meeting $V' \meet W,V'' \meet W$ properly.
Let $X_1 $ be an effective cycle representing $U' \cdot V'$, supported
on the components of $U' \meet V'$; and similarly $X_2,X_3,X_4$
for $U' \meet V'', \ldots, U'' \meet V''$.
Then by \ref{proper1}, the $X_i$
meet $W$ properly.
By \ref{trans1}, $\deg(X_1) \leq \deg(U')\deg(V')$, etc.
Thus
 $  \sum_i | \deg(X_i) |    \leq |U||V|$.
Since $W$ was arbitrary, $|U \cdot V| \leq |U||V|$.

\begin{lem}\lbl{cyfin}  For any cycle $R$ on $Y$, and very ample $H$ on 
$Y$,
$|R|_H $ is finite.
\end{lem}
\proof Say $D,Y,U$ are defined over an algebraically closed field $k$. 
Since $k$ is an elementary submodel of an algebraically closed field
$K$ with $tr. deg._k K $ infinite, in the
definition \ref{not2.0} of $|R|_H$, one can restrict the $sup_S$ to range
over cycles $S$ defined over $k$, while allowing the $ U',U''$ to be defined
over $K$.  The finiteness is now immediate from Lemma \ref{moveall}.   \qed

\end{subsection}  

\begin{subsection}{Correspondences}
\lbl{correspondences}

Let us recall the language of correspondences, (cf. \cite{fulton}, 
16.1).
Let $X$,$X'$ be smooth, complete varieties
of dimension $d$ over an algebraically closed field.
A correspondence $R$ on $X \times X'$ is a $d$-cycle on
$X \times X'$; we will write somewhat incorrectly $R \subset (X \times 
X')$.
 The transpose of $R$ is the corresponding cycle
on $X' \times X$; it is denoted $R^t$. If $T$ is a correspondence
of $X' \times X''$, then one defines the composition by:
 \[T \circ R = pr_{X,X'' \ *}( {pr_{X,X'}}^* R \cdot {pr_{X',X''}}^* T ) 
\]

When the context does
not make clear which product is in question, we will use $X \circ Y$ and 
$X^{\circ n}$
for the composition product and power, and $X \cdot Y$, $X^{\cdot n}$ 
for the
intersection product and power.

\begin{notation} \lbl{not2}
Suppose $X$ and $X'$ are smooth complete varieties of dimension $d$,
given together with very ample divisors $H$,$H'$ on $X$,$X'$. Let
$R \subset (X \times X')$ be a correspondence.
Write
\[  \deg_{cor}(R) = [p \times X'] \cdot R \]
where $p$ is a point of $X$

\end{notation}

$\deg_{cor}(R)$ is intended to denote the  degree of $R$ as a 
correspondence.
If $R$ is an irreducible subvariety $R$ and $pr_X$ restricts to a 
dominant morphism
$\pi: R \to X$, then  $\deg_{cor}(R)$ is the degree of $\pi$.  In 
particular that
if $R$ is the graph of a function, then $\deg_{cor}(\Phi)=1$.
By contrast, we will write $\deg(R)$ or $\deg_{cy}(R)$ for 
$\deg_{{pr_1}^*H + {pr_2}^*H'} R$.
If $R = R_1 -R_2$ with $R_1,R_2$ effective (irredundantly), we let
 $\deg_{|cy|}(R) =  \deg_{cy}(R_1) -\deg_{cy}(R_2)$.

Note that $\deg_{cor}(\Phi)=1$, $\deg_{cor}(\Phi^t) = q^d$,
$\deg_{cor}(S) = \delta$.

\begin{lem}\lbl{cy} Let $R \subset X \times Y$ and $T \subset Y \times 
Z$
be correspondences. Computing their degrees with respect to very ample 
divisors $H_X$,$H_Y$,$H_Z$
on $X$,$Y$,$Z$ respectively: \begin{enumerate}
\item $\deg_{cor}(T \circ R) = \deg_{cor}(T)\deg_{cor}(R)$.
\item $\deg_{cor}(R) \leq \deg_{|cy|}(R)$
\item $|T \circ R| \leq {2d \choose d}^2 |T| |R|$
\item Suppose $Z=X$ and $H_X = H_Z$. Then
$$\left( (T \circ R) \cdot_{X^2} \Delta_X  \right)
 = (R \cdot T^{t}) \leq  |T^t| |R|$$
\end{enumerate}
\end{lem}

\proof \begin{enumerate}
\item   Let $P_{XY}$ be the correspondence ${p} \times Y$, where $p$ is 
a point of $X$.
Then $\deg_{cor}(R) = (P_{XY} \cdot R) = \Delta_X \cdot (P_{XY}^t \circ 
R)$.  (See (4) below
for the last equality.)   For $R$ the cycle of an irreducible variety, 
and hence
in general, $P_{XY}^t \circ R$ is a multiple of $P_{XX}^t$.  So
$$P_{XY}^t \circ R  =   \deg_{cor}(R) P_{XX}^t$$
Similarly $P_{YZ}^t  \circ T  =   \deg_{cor}(T) P_{YY}^t$, so
  $P_{XZ}^t \circ T   =  P_{XY}^t \circ P_{YZ}^t \circ T =  
\deg_{cor}(T) P_{XY}^t$.
Thus:
\[ (P_{XZ}^t \circ T \circ R) =  \deg_{cor}(T) (P_{XY}^{t} \circ R)  \]
 The desired formula  follows. upon intersecting this with $\Delta_X$.

\item   Here we may assume $R \geq 0$, so that $\deg_{|cy|}(R) = 
\deg_{cy}(R)$.
 On $X$, $H^{\cdot n} = c [pt]$ for some $c \geq 1$; so
$$\deg_{cor}(T) = (P_{XY} \cdot R) =
c^{-1} ({pr_{X}}^*H)^n \cdot R \leq \deg_{cy}(R)$$

\item  Both $R \circ T$ and $|R|$,$|T|$ depend on $R$,$T$ only up to
rational equivalence,
so we may change $R$ and $T$ within their rational equivalence class.
By \ref{proper3}, $R$,$T$ may be replaced so that $\deg_{|cy|}(R) \leq 
|R|$, $\deg_{|cy|}(T) \leq |T|$,
$(R \times Z) \meet (X \times T) $
is a proper intersection, and $pr_{X,Z} ( (R \times Z) \meet (X \times 
T) )$
meets properly a given cycle $S$.

 Then
 $R \circ T = {pr_{X,Z}}_* ( (R \times Z) \meet (X \times T) )$, and
 by   \ref{proj0}(1),\ref{trans1}, and \ref{proj0}(2),
$\deg_{|cy|}(W) \leq  ({2d \choose d} \deg_{|cy|}(R) {2d \choose d} 
\deg_{|cy|}(T)| \leq
{2d \choose d}^2 |R ||T|$;
and $W$ meets $S$ properly.
Taking supremum over $S$, the statement follows.

\item The equality $(T\circ R) \cdot \Delta_X = (R \cdot T^{t})$ follows
from the projection formula, and the definition of composition of 
correspondences;
both are equal to the triple intersection product
$(R \times X) \cdot (X \times T) \cdot \Delta_{13}$.

The inequality is clear from \ref{trans1}.
\ \ \ \ \ $\Box$ \end{enumerate}

\end{subsection} 

\begin{subsection}{Geometric statement of the uniformity}
\lbl{geostat}

To formulate the theorem geometrically, we will need to consider
$V$, $q$, and $S \subset V \times V^{\phi_q}$
all varying separately.  Below, $B$ will be the base over which $V$ 
varies. $B'$ will be the base for $S$;
it need not be a subscheme of $B \times B$, though one loses little or 
nothing by thinking of that case.
A scheme {\em over} $B$ is a scheme $S$ together 
with a morphism $\alpha: S \to B$.
We will not always have a notation for $\alpha$; Instead, if $L$ is a 
field and
$b \in B(L)$, we will write $S_b$ for the fiber of $S$ over the point 
$b$.
If $U \to V$ is a dominant map of irreducible varieties over $k$, then
the (purely inseparable) degree of $U/V$ is
 the (purely inseparable) degree of $[k(U):k(V)]$.

\begin{notation}\lbl{not} \begin{enumerate}

\item    $B$ and $B'$ are reduced, irreducible, separated schemes
 over $\Zz[1/m]$ or $\Ff_p$.
    (Not necessarily absolutely irreducible.)
    A (base change) map $\beta: B' \to B^2$ is also given.
 We will refer to the
     cases as ``characteristic 0" and ``characteristic p", respectively.

 Let $V$ be a scheme over $B$. View $V^2$ as
a scheme over $B^2$. Let $S$ be a ($B' -$) subscheme of $V^2 
\times_{B^2} B'$.

For $L$ a field and $b \in B'(L)$, denote $\beta(b) = (b_1,b_2) \in 
V^2(L)$. We will
assume that $V_{b_1}$,$V_{b_2}$ and $S_b$ are varieties, and view
$S_b$ as a subvariety of $V_{b_1} \times V_{b_2}$

\item We further assume that for $b \in B'$, $S_b$,$V_{b_1}$,$V_{b_2}$ 
are absolutely irreducible;
 the projection map $S_b \to V_{b_2}$ is a quasi-finite map;
 $\dim(V_{b_1}) = \dim(V_{b_2}) = \dim(S_b) =d$,
$\deg(S_b/V_{b_1}) = \delta < \infty$,
and if $B$ is over $\Ff_p$, the purely inseparable degree of $S_b$ over 
$V_{b_2}$
is $\delta'_p$.

\item Let $q$ be a prime power (power of $p$), $a \in B(K_q)$, $b \in 
B'(K_q)$ such that
 $\beta(b) = (a,\phi_q(a))$. In this situation, we denote:

$$ V_b(S,q) = \{c \in V_a(K_q): (c,\phi_q(c)) \in S_{b}(K_q) \}$$

\end{enumerate}
\end{notation}

In this language, the asymptotic version of ACFA 1 is the following:

\smallskip \noindent {\bf Theorem 1B'}  {\em 
Let $B$ be a reduced, separated scheme
of finite type over $\Zz[1/m]$ or $\Ff_p$.
Let $B',V,S,\beta$ be as in \ref{not}. For any sufficiently large $q$, 
if $b \in B'(K_q)$,
and $b_2 = \phi_q(b_1)$, then there exists $c \in V_{b_1}(K_q)$ with 
$(c, \phi_q(c)) \in S_b(K_q)$.
} \smallskip

While we need a mere existence statement, we see no way to prove it 
without
going through a quantitative estimate. We formulate this estimate as 
follows.
Note the similarity to the Lang-Weil estimates; these are the special 
case when
the $S_b$ are the diagonals.

\smallskip \noindent {\bf Theorem 1B}  {\em 
 
Let assumptions be as in Theorem 1B'.
Then there exists a nonempty open  
subscheme $B''$ of $B'$ and constants $\rho$ and $\delta^* > 0$
such that
if $b \in B''(K_q)$, $\beta(b) = (a,\phi_q(a))$, then $V_b(S,q)$ is 
finite, of cardinality
$$ \#(V_b(S,q)) = \delta^* q^d + e \hbox{ \ \ \ \ with \ } |e| \leq \rho
q^{d-\half} $$
with  $\delta/\delta'_p$} \smallskip

\smallskip

Theorem 1B' follows from Theorem 1B by
restriction to an open subset of $V$, and using
 Noetherian induction on $B'$.  
 
 The statement 1B refers to number of points, without multiplicities.  In fact 
 each point (excepting $O(q^{d-1})$) occurs with 
multiplicity $\delta'_p$,  so that
the number of points counted with multiplicity is $\delta q^d + O(q^{d-
\half})$.

\end{subsection}

\begin{subsection}{Separability}

\begin{lem}\lbl{sep} In Theorem 1B, we may assume that 
 for $b \in B'$, $S_b/V_{b_2}$ is separable.    \end{lem}

\proof If $B \tensor \spec \Qq \neq \emptyset$,   this may simply be achieved 
by replacing
$B'$ by an open subscheme $B''$ so that for $b \in B''$,
$deg(S_b/V_{b_2})$ is a constant $\delta'$, and then further by the
open subscheme $B''' = B'' \tensor_{\Zz[1/m]} \Zz[1/(m \delta')]$.  Othewise we may assume
$B$ is over $\Ff_p$, and $\delta'_p= p^l$.
We will modify $B'$ and $S$ so as to obtain a similar situation
with $\delta'_p= 1$; the modified objects will be denoted by a $\tilde{\ 
}$ but
$B$ and $V$ are left the same.

 Let $F_B: B^2 \to B^2$ be the map $(Id,\phi_{p^l})$.
Let $\tilde B' = B' \times _{B^2} B^2$, where the implicit map $B^2 \to 
B^2$
is $F_B$. (Thus if $B'$ is a subscheme of $B^2$,
 $\tilde B' = F_B^*(B')_{red}$)
is the reduced scheme underlying the pullback of $B'$ by $F_B^*)$.
Now $F_B$ induces a map $F_2: \tilde B' \to B'$. Also let $F_1$ denote
the map $(Id,\phi_{p^l}) : V^2 \to V^2$, and let
$$F = (F_1 \times_{F_B} F_2) \ : \ \ \ V^2 \times_{B^2} \tilde B' \to 
V^2 \times_{B^2} B' $$
Finally let $\tilde S = F^*(S)_{red}$ be the reduced pullback of $S$ via 
$F$. Also if $q$ is a power
of $p$, let $\tilde q = q/p^l$. We now claim that $\tilde \delta'_p= 1$,
$\tilde \delta = \delta p^{ld}/\delta'_p $,
and that if $\tilde b \in \tilde B'(K_q)$, $b = F(\tilde b)$, then
$V_b(S,q)$ and $V_{\tilde b}(\tilde S, \tilde q)$ coincide as sets. This 
is
easily seen by going to fibers: if $\tilde \beta ( \tilde b) = 
(b_1,\tilde b_2)$
then $\beta(b) = (b_1,b_2)$ with $b_2 = \phi_{p^l}( \tilde b_2)$. If 
further
$(a_1,\tilde a_2)$ is a generic point of $\tilde S_{\tilde b}$, and
$a_2 = \phi_{p^l}( \tilde a_2)$, then $(a_1,a_2)$ is a generic
 point of $S_b$.
We have $\delta = ( K_q(a_1,a_2):K_q(a_1) )$,
 $\tilde \delta =( K_q(a_1,\tilde a_2):K_q(a_1))$,
and the degree computations are elementary. The truth of Theorem 1B for 
$S$ now follows
from the same for $\tilde S$.

\end{subsection}

\begin{subsection}{Rough bound (Bezout methods)}

We give here an alternative proof of the rough
upper bound on the number of points on the intersection of a variety
with Frobenius, \ref{si-dim+}, using   Bezout's theorem.  We require Fulton's "refined" version
since the intersection can have   components of various dimensions (if only at infinity.)
We need to reduce to an intersection of divisors, partly in order to arrive at a
situation where   intersection multiplicities agree with geometric multiplicities.

\<{lem}\lbl{bezout-1}
Let $l,k \in \Nn$.  Let $U_i$ be a $d_i$-dimensional irreducible subvariety of $\Pp^n$ ($i=1,\ldots,r$).  Let $V$ be a subscheme of $\Pp^n$ cut out by   hypersurfaces of degree $b$.  Then there exists
a subscheme $V' \supseteq V$ such that for any $i$,  any $d_i-k$-dimensional 
irreducible component  of $U_i \meet  V$ is also an irreducible component of $U_i \meet V'$; and $V'$ is the intersection of $k$ hypersurfaces of degree $b$.    \>{lem}

\<{proof}
  Say $V$ is   cut out by $P_j = 0  \ (j=1,\ldots,l)$,
with $P_j$ a homogeneous polynomial of degree $b$.
 Let
$M$ be a $k \times l$ matrix, with generic coefficients, i.e. avoiding a finite number of proper Zariski closed sets of matrices, to be specifed below.
Let $(Q_1,\ldots,Q_k) = M (P_1,\ldots,P_l)$.  Let $V'[m]$ be the zero set of $\{Q_1,\ldots,Q_m\}$, $V'=V'[k]$.  Evidently  $V' \supseteq V$.  Fix $U=U_i$, $d=d_i$, and a   $d-k$-dimensional 
irreducible component $C$ of $U \meet  V$.  Then $C \subseteq (U \meet V')$, and we have
to show that $C$ is an irreducible component of $U \meet V'$,
i.e. that $C$ is contained in no larger irreducible subvariety of $U \meet V'$.
  Let $U'$ be a Zariski open neighborhood in $U$  of the generic point of $C$, not meeting any of the other irreducible components
of $U \meet V$.  Let $C' = (C \meet U')$.  Then $\dim(C') = d-k$.
 
\claim{}  For $m \leq k$,  $\dim(U' \meet V'[m]) \leq d-m$.  

 For $m=0$ this reads $\dim(U')=d$, which is clear.  Assume $m<k$ and the inequality is true for $m$;
 to show that it holds for $m+1$,  let $E$ be an irreducible component of $U' \meet V'[m]$.
Then $\dim(E) \leq d-m$, and it suffices to show that $\dim(E \meet V'[m+1]) < d-m$.
If $\dim(E) < d-m$ there is nothing to show.  Otherwise, $\dim(E) = d-m > d-k$; 
as $\dim(U' \meet V) = \dim(C') = d-k$, we have 
 $E \not \subseteq V$.  Thus not all $P_j$ vanish on $E$.  By genericity of the $m+1$'st row of $M$  over the previous 
rows,
$Q_{m+1}$ does not vanish on $E$.  So $E \not \subseteq V'[m+1]$.  Thus
$\dim(E \meet V'[m+1]) < dim(E) \leq d-m$, proving the claim.  

In particular, $\dim (U' \meet V') \leq d-k$; so $C'$ is an irreducible component of $U' \meet V'$.
It follows that $C$ is an irreducible component of $U \meet V'$.
        \eprf

If $U$ is a variety and $C$ is an irreducible subvariety of $U$, let $\dim_C(U)$ be the dimension of $U$ at a neighborhood of  the generic point of $C$; i.e. $\dim_C(U) = \min \dim(U \meet W)$ where $W$ ranges over Zariski open subsets of $\Pp^n$ containing a dense subset of $C$.  
Equivalently, $\dim_C(U)$ is the maximal dimension of a component of $U$ containing $C$.
Let $\cdim_C(U) =  \dim_C(U) - \dim(C)$.  

 Let $V$ be a subscheme of $\Pp^n$ cut out by   hypersurfaces of degree $b$,  and let $V'[k]$ be defined as in the proof of Lemma \ref{bezout-1} above.  
\<{cor}  \lbl{bezout-1c}  Let 
$U$ be a variety,
and let $C$ be an irreducible component of $U \meet V$.   Then $C$ is a component
of $U \meet V'[k]$, with $k  = \cdim_C(U)$. 

The geometric multiplicity of $U \meet V'[k]$ along $C$ is at least that of $U \meet V$ along $C$.   
 \>{cor}  

\<{proof}  Let $\{U_i\}$ be the irreducible components of $U$ containing $C$, $d_i = \dim(U_i)$.
Then $C$ is an irreducible component of each $U_i$.   By Lemma \ref{bezout-1}, $C$ is an irreducible component of $U_i \meet V'[k_i]$ with $k_i=d_i-\dim(C)$.  Now $k_i \leq k$,
so $V  \subseteq V'[k] \subseteq V'[k_i]$; hence $C$ is also an irreducible component of $U_i \meet V'[k]$.  Since this is true for each $U_i$, $C$ is an irreducible component of $U \meet V'[k]$.

The statement regarding geometric multiplicity, is immediate since $V \subseteq V'[k]$
as schemes.  
\eprf

Here are a few special cases, with $U = \Pp^n$.   Let $C$ be an irreducible subvariety of a variety
$P$, ${\cal F}$ a set of  hypersurfaces of $P$.  We say $C$ is {\em 
weakly cut out
by ${\cal F}$} if it is a component of
the intersection scheme $\meet {\cal F}$.

More generally, if the cycle $[\meet {\cal F}] = \sum n_i C_i$ with
$C_i$ the cycle of an irreducible variety, and $0 \leq m_i \leq n_i$,
we say that
the cycle $\sum m_i C_i$ is  {\em weakly cut out
by ${\cal F}$}.

 Let $S$ be a $n-k$-dimensional subvariety of  $\Pp^n$,
weakly cut out by hypersurfaces of projective degree $b$.
Then $S$ is weakly cut out by $k$ hypersurfaces of  degree $b$.

More generally, if $C = \sum m_i C_i$ is weakly cut out by $l$ 
hypersurfaces
of degree $b$, and $\dim(C_i) \geq n-k$ for each $i$, then
$C$ is weakly cut out by $k$ hypersurfaces $D_i$ of  degree $b$.   



The proof of the next lemma will use positivity properties of
intersection theory in multi-projective space.  Let
$\Pp = \Pp^{n_1} \times \Pp^{n_2} \times \Pp^{n_3}$
be a product of three projective spaces.  Let $T_i$
be divisors on $\Pp$, with multi-degrees $(a_i,b_i,c_i)$
with respect to the standard divisors $H_i$, pullbacks of the
hyperplane divisors on $\Pp^{n_i}$.  Then

$$ \{ [T_1 \meet \ldots \meet T_m] \}_{\dim(\Pp) - m}
                                  \leq T_1 \cdot \ldots \cdot T_m $$

Here on the left side we have the proper part of the cycle corresponding
to the scheme-theoretic intersection, and on the right-hand side the
intersection cycle in the sense of intersection theory.
The inequality $\leq$ has the sense of \ref{nota3}. More explicitly,
one considers the intersection of $T_1 \times \ldots \times T_m$
with the   diagonal embedding of $\Pp$ in $\Pp^m$. This intersection has
proper components $W_1,\ldots,W_l$, and
other distinguished varieties $Z_j$.  By \cite{fulton} 6.1, we have
the canonical decomposition

$$ T_1 \cdot \ldots \cdot T_m = \sum_i m_i [W_i] + \sum_j m_j \alpha_j$$
where $\alpha_j$ are certain cycles on the $W_\nu$.  We push forward
to $\Pp$ to obtain the same equality for cycles on $\Pp$.
Now  the tangent bundle to multi-projective
space is generated by global sections  (\cite{fulton} 12.2.1(a,c)),
hence ( \cite{fulton} 12.2(a)) each $\alpha_j$ is represented by a non-
negative
cycle.   So
$$T_1 \cdot \ldots \cdot T_m \geq \sum_i m_i [W_i]$$
On the other hand, by \cite{fulton} 7.1.10, $m_i$ is precisely the 
geometric
multiplicity of $W_i$ in the intersection scheme $T_1 \meet \ldots \meet 
T_m$;
so $\{[T_1 \meet \ldots \meet T_m]\}_{\dim(\Pp) - m} = \sum m_i [W_i]$.  
The claimed Bezout
inequality follows.

We now consider an analogue of Lemma \ref{bezout-1} when two degrees $b,B$  are involved.
Assume $U$ is a subscheme of $\Pp^n$ cut out by hypersurfaces of degree $\leq b$,
while $V$ is a subscheme of $\Pp^n$ cut out by hypersurfaces of degree $\leq B$.  

  Let
$\Xi(U,V,l,k)$ be the set of components $C$ of $U \meet V$ with $\dim(C) = l, \dim_C(U) = k$;
write $\deg ( \Xi(U,V,l,k)) = \sum _{C \in \Xi(U,V,l,k)} m_C \deg(C) $ where
$m_C$ is the geometric multiplicity of $C$ in $U \meet V$, and $\deg(C)$ is the projective degree.

\<{lem} \lbl{bezout-2}  $ \deg(\Xi(U,V,l,k)) \leq b^{n-k} B^{k-l} $ \>{lem}

\<{proof}  
Say $U$   is cut out by homogeneous polynomials $P^U_i$ 
of degree $b$, while $V$
 is cut out by homogeneous polynomials 
  $P^V_j$ of  degree $B$.  As in Lemma \ref{bezout-1}, 
 let $U'[m]$ (respectively $V'[m]$)  be the scheme cut out by  $m$ generic linear combinations of the $P^U_i$ (respectively $P^U_j$.)  
 
 Let $C \in \Xi(U,V,l,k)$.  
By the dimension theorem, $ \dim(C) \geq \dim_C V + \dim_C U - n$, so
$\dim_C V \leq n-k+l$.   It follows using the genericity in
the definition of $V'[j]$ that for $j \leq k-l$ we have $\dim_C V'[j] = n- j$.  (The inequality
$\geq $ is true by the dimension theorem, since $V'[j]$ is the intersection of $j$ hypersurfaces.)  In particular, $\dim_C V'[k-l] = n-(k-l)$.

 By Corollary \ref{bezout-1c}, $C$ is an irreducible component of $U \meet V'[k-l]$.  We have   
 $\cdim_C V'[k-l] = n-(k-l)-l = n-k$.  Hence again by \ref{bezout-1c}, 
$C$  is an irreducible component of $U'[n-k] \meet V'[k-l]$; with geometric multiplicity
  at least $m_C$.  Now this is an intersection of $n-k$ hypersurfaces of degree $b$ and
  $k-l$ hypersurfaces of degree $B$, so 
  by the refined Bezout theorem as discussed above, $\sum m_C \deg(C) \leq b^{n-k} B^{k-l}$.  \eprf

\<{cor}  \lbl{bezout-3}  Let $L$ be an algebraically closed field.   $f_1,\ldots,f_r,g_1,\ldots,g_s$ be polynomials in $n$-variables, with $\deg(f_i) \leq b, \deg(g_j) \leq B$,   where $b \leq B \in \Nn$.  Let $J$ be the set of dimensions of
components of the zero set $Z(f_1,\ldots,f_r)$, and assume 
$W=Z(f_1,\ldots,f_r,g_1,\ldots,g_s)$ is finite.  Then $|W| \leq \sum_{j \in J} b^{n-j}B^j$.   \>{cor}

\begin{lem} \lbl{fin-proj}  Let $H$ be the hyperplane divisor on 
$\Pp^n$,
 $H_{12} = {pr_1}^*H + {pr_ 2}^*H$ on $\Pp^n \times \Pp^n$.
 Let $U$ be a purely $k$-dimensional cycle
of $\Pp^n \times \Pp^n$,   cut out by hypersurfaces with
$H_{12}$-degrees $\leq b$.
Then $$\deg_{H_{12}} [U \meet \Phi_q]_l  \leq  b^{2n-k} (q+1)^{k-l}$$ 
 \end{lem}

\prf
The subvariety $\Phi_q \subset \Pp^n \times \Pp^n$ is cut out by ${n 
\choose 2}$
hypersurfaces  $F_{ij}$
 of bidegree $(q,1)$ (namely, those defined by ${X_i}^qY_j = 
{X_j}^qY_i$). The statement thus follows from Lemma \ref{bezout-2}.    \eprf

\begin{lem}\lbl{fin} Let assumptions be as in \ref{not}  (1),(3).
Assume further that
for any $b \in B'$, the projection maps $S_b \to V_{b_2}$ has finite 
fibers.
Then
for all large enough $q$, and any $b \in B'$, $V_b(S,q)$ is finite.    
Moreover,
 for some constant $\rho$, $\#V_b(S,q) \leq \rho q^d$. \end{lem}

\proof 
Note that the assumptions remain  true if $V$ is replaced by a proper 
subscheme $U$
and $S$ by any component of
$S \meet U^2 \times_{B^2} B'$. Hence using Noetherian induction, we may 
assume
the lemma holds for proper closed subschemes of $V$.

 As in \ref{sep}, we may assume the
projection $S_b \to V_{b_2}$ is generically separable.
Let $b \in B'(K_q)$; note that $V_b(S,q)$ is the set of points of a 
scheme over $K_q$,
and let $C$ be an irreducible component over $K_q$. Let $a_1$ be a 
generic point of $C$,
$a_2 = \phi_q(a_1)$.
Then $(a_1,a_2) \in S_b$ so the field extension 
$[K_q(a_1,a_2):K_q(a_2)]$ is
finite and separable. Since it is also purely inseparable, 
$K_q(a_1)=K_q(a_1^q)$,
so $K_q(a_1)$ is a perfect, finitely generated field extension of $K_q$. 
It follows
that $K_q(a_1) = K_q$, so that $C$ is finite. Since $C$ was an arbitrary 
component,
$V_b(S,q)$ must be finite.

For the quantitative bound, we may (again using Noetherian induction
and stratifying if needed)
embed $V_a$ in a projective space $\Pp_n$; the dimension and degree of the
embedding are bounded irrespective of $a \in B$, as well as the 
bidegrees $(f_1,f_2)$ of
the resulting embedding of $S_b$ in $\Pp_n \times \Pp_n$.    Now 
$V_b(S,q)$ is the
projection to $V_{a_1}$ of the intersection of $S_b$ with the graph 
$\Phi_q$ of
the Frobenius map $\phi_q: \Pp_n \to \Pp_n$. By \ref{si-dim}  or \ref{fin-proj} the 
intersection
has  $O(q^d)$ points, and the lemma follows. \qed

\begin{rem} \lbl{transverse} \   \end{rem} 

(1)  If  $S  \subset V \times V^\si$
are varieties over $k$, and $S \to V^\si$ is separable in the sense that $k(a) \in k(b)^{sep}$
for any $L \supset k$ and $(a,b) \in S(L)$, then $S, \Phi_q$ meet properly.  Indeed 
if $(a,b) \in (S \meet \Phi_q)(L)$ then $k(a) \subset k(a^{1/q})^{sep}$, so $k(a)^{sep}$
is perfect, and it follows that $a \in k^{sep}$.   See   \ref{towers}.

(2)  (\cite{Lang-Abelian}, \cite{pink})  if $S \to V^{\si}$ is 
\'etale,
then $S$ and $\Phi_q$ intersect transversally. (Look at the tangent 
spaces:
 $T\Phi_q$ is horizontal at each point, while $TS \to TV^{\si}$ is
 injective.) 

\end{subsection}

\begin{subsection}{Finite covers}
\begin{lem}\lbl{birational}(Birational invariance). Let $B$,$V$,$B'$,$S$ 
satisfy the hypotheses of
Theorem 1B. Suppose $\tilde V$ is an open subscheme of $V$, so that for 
$a$ the generic
point of $B$, ${\tilde V} _a$ is an open subvariety of $V_a$. Let
$\tilde S = S \meet {\tilde V}^2 \times_{B^2} B'$. Let $\tilde B'$ be an 
open
subscheme of $B'$ such that the hypotheses of \ref{not} hold.
Then Theorem 1B  is true for $B$,$V$,$B'$,$S$ iff it holds for
$B$,$\tilde V$,$\tilde B'$,$\tilde S$. \end{lem}

\proof The invariants $\delta$,$\delta'_p$ are the same for the two 
families. The sets
$V_b(S,q)$
and $\tilde V_b(\tilde S,q)$ are also the same, except for points
in $F_b(S^*,q)$, where $F$ is a closed subscheme of $V$ complementary to 
$\tilde V$
at the generic fiber, and $S^*$ is the restriction of $S$ to $F$. By 
\ref{fin},
the size of $F_b(S^*,q)$ is of the order of the error term in the 
expression
for $V_b(S,q)$ in 1B. \qed

The first item of the next lemma will be used to reduce to smooth varieties.
The second will not be used, but originally gave some indication of
the correctness of the form of Theorem 1B, since it shows that the 
desired lower bound in this theorem follows from the upper bound.  (Once
once one knows the theorem for the easy case of single difference equations),

\begin{sta}\lbl{2qin}
$$ \#(V_b(S,q)) \leq \delta^* q^d + e \hbox{ \ \ \ \ with \ } e \leq 
\rho q
^{d-\half} $$
\end{sta}

We will write such statements as
$$ \#(V_b(S,q)) \leq \delta^* q^d + O(q^{d-\half})$$
The point is that the implicit coefficient depends only on the system 
$(B',V',S,
\beta)$ and not on the choice of $b$ or of $q$.

\begin{lem}\lbl{birat2} Let $B$,$V$,$B'$,$S$ satisfy the hypotheses of
1B. Let $h_B: \tilde B \to B$ and $h_V: \tilde V \to V$ be compatible
quasi-finite maps (base extension for $B$, and finite cover of $V$).
Let $\tilde B' = B' \times _{B^2} {\tilde B}^2$. We get an induced map
$h: (\tilde V)^2 \times_{{\tilde B}^2} {\tilde B'} \to V^2 \times_B^2 
B'$.
Let $\tilde S(j)$ be the various irreducible components of $h^{-1}S$, 
that are
nonempty in the generic fiber. \begin{description}
\item{(a)} Suppose the conclusion of 1B  holds for each
$\tilde B$,$\tilde V$,$\tilde B'$,$\tilde S(j)$. Then it holds also for 
the original
system.
\item{(b)} Conversely, suppose the conclusion of 1B hold for 
$(B,V,B',S)$,
and the conclusion of \ref{2qin} holds for each $\tilde B$,$\tilde 
V$,$\tilde B'$,$\tilde S(j)$.
Then the conclusion of 1B holds for each $\tilde B$,$\tilde V$,$\tilde 
B'$,$\tilde S(j)$.
\end{description}
\end{lem}

\proof We can first apply the base change from $B$ to $\tilde B$ without
changing $V$ (i.e. replacing $V$ by $V \times_B \tilde B$); this clearly
makes no difference to 1B. Thus we may assume $\tilde B = B$, and 
$\tilde B' = B'$.
Let $a$ be the generic point of $B$. Let $h_a: \tilde V_a \to V_a$ be
the map induced by $h$ on the fibers. Outside of a proper subvariety of 
$V_a$,
all fibers of $h_a$ have the same size $s$. This proper subvariety may 
be removed
using \ref{birational}. Similarly we may assume the various $\tilde 
S(j)$ are disjoint.
Let $b$ be a generic point of $B'$. Let
$$H = h_{b_1} \times h_{b_2} : \
{\tilde V}_{b_1} \times {\tilde V}_{b_2} \to V_{b_1} \times V_{b_2}$$
We have $H^{-1} S_b = \union _j S(j)_b$. Thus
$$H^{-1} V_b(S,q) = \union _j V_b(\tilde S(j),q)$$
By the constant size of the fibers of $h_a$,
$$s \#V_b(S,q) = \Sum_j \#V_b(\tilde S(j),q)$$
A simple exercise in degrees of field extensions shows:
$$\delta^* s = \sum_j \delta^*_j$$
Where $\delta^*_j$ refers to $\tilde S(j)$, cf. 1B.
(Because of \ref{sep}, only the case $\delta'_p= 1$ need be considered 
in the applications).

Now (a) if 1B holds $\tilde S(j)$, we have
$$\#V_b(\tilde S(j),q) = \delta^*_j q^d + O(q^{d-\half})$$
   Summing over $j$ , using the previous two
displayed equalities, we obtain
$$s \#V_b(S,q) = \delta^* s q^d + O(q^{d-\half})$$
and division by $s$ yields (a). (b) follows similarly, using the 
principle that
if $\sum_i a_i = \sum_i a_i' $
 and each $a_i \leq a_i'$, then each $a_i=a_i'$
 (here all up to $O(q^{d-\half})$.) \qed

\end{subsection}

\begin{subsection}{Reduction to smooth varieties}
\begin{lem}\lbl{smooth} In Theorem 1B, we may assume that $V$ is an open
subscheme of $\bar V$, $S = \bar S \meet V^2$, with ${\bar V}$ smooth 
and projective over $B$. \end{lem}

\proof In characteristic zero, we can use Hironaka's resolution of 
singularities;
the generic fiber $V_b$ of $V$ is birational to an open subvariety of a 
smooth, proper
variety, and by \ref{birational} the theorem is invariant under 
birational changes
of this type. Similarly, in general, we can use \ref{birat2} and
the following theorem of de Jong (\cite{dejong}):

{\em For any variety $V$ over a field $k$ there exists a finite 
extension
$\tilde k$ of $k$, and a smooth projective variety $\tilde V$ over 
$\tilde k$,
and a finite, dominant morphism from an open subvariety of $\tilde V$ to 
$V$.}
\qed

\begin{remark} \ \end{remark}
For
theorem 1B' the reduction to the smooth case is not needed on this 
(geometrical
and numerical)
level; a similar but easier reduction shows directly (using
 \cite{dejong}) that the instances of the axiom scheme ACFA
referring to open subvarieties of smooth complete varieties
 imply the rest of the axioms.

\end{subsection}

\end{section}

\begin{section}{The virtual intersection number}

\lbl{virtual}

   This section is devoted to a proof of theorem \ref{2qc}. It is
   formally similar to
Theorem 1B;  but it refers to the virtual intersection number, in the 
sense
of intersection theory, and not to the actual intersection.  By virtue 
of
\ref{smooth},\ref{sep}
 we can now assume that
the varieties
$V_a$ are smooth
and complete, and $S_b/V_{b_2}$ separable.   But as a result of the 
completion
process, the intersection may have infinite components.  The effect
of these components on the intersection number will have to be estimated
in the next section.

\begin{lemnot}\lbl{sm} Let $\mu$, $\nu$ be bounds for the absolute 
degree
and the sum of Betti numbers of $V_b$, so that:  \begin{itemize}
\item For $b \in B$,
there exists a very ample divisor $H_b$ on $V_b$,
 such that $|S_b| = |S_b|_{H_b} \leq \mu$
\item For $b \in B$, letting $X$ be the variety $\bar V_b$ base-extended 
to
an algebraically closed field, we have
 $\Sum_{i \leq 2 \dim(X)} \dim H^i(X,\Qq_l) \leq \nu$.
\end{itemize}
 \end{lemnot}

\proof  To show the existence of either $\mu$ or $\nu$, it suffices to 
show
that the numbers in question are bounded on a Zariski open subset of 
$V$;
then one may use Noetherian induction.  Let $b$ be a generic point of 
$V$,
and pick a very ample $H_b$ on  $V_b$.  The construction of
\ref{cyfin} can be described by a
first-order formula  concerning $b$.  This  formula must remain
valid for a Zariski neighborhood of $b$.   Thus $|S_{b'}|$ is
uniformly bounded  on this neighborhood. The same argument, using
the constructibility of the higher derived images of constant sheaves, 
shows
that $\Sum_{i \leq 2 \dim(X)} \dim_{Z/lZ} H^i(X,Z/lZ)  $ is bounded on a 
Zariski
neighborhood; this number bounds also the sum of the Betti numbers.

\begin{thm}\lbl{2qc} Let notation be as in 1B and \ref{smooth}.
Let $b \in B'(K_q)$ ,
$\beta(b) = (a,\phi_q(a))$, $X = \bar V_a$, $X' = \bar V_{\phi_q(a)}$.
Let $\Phi_b \subset (X \times X')$ be the graph of Frobenius.
Then
\[ S_b \cdot_{X \times X'}  \Phi_b = \delta q^d + e
   \hbox{ \ \ \ \ with \ } |e| \leq {2d \choose d}^{2} \mu \nu q^{d-
\half} \]
\end{thm}
\

\begin{example}  \ \lbl{cycle-t} \end{example} In the case of projective space,
i.e. ${\bar V}_{b_1} = \Pp_d$,
 we can immediately prove Theorem \ref{2qc} cycle-theoretically
(without transcendental cohomology.) 

Let $\bar S_b \subset {\Pp_d}^2$ be the closure of $S_b$.
Consider the intersection of $\bar S_b$ with $\Phi_q$,
the graph of Frobenius on ${\Pp^d}$. In the intersection theory sense, 
it can
be computed in terms of the bidegrees:
$$ [\Phi_q] = \sum_i q^i h_1^ih_2^{d-i}$$
$$ [S] = \sum_i {a_i} {h_1}^i{h_2}^{d-i}$$ \ with $a_0 = \delta $.
So
$$ [\Phi_q \cdot S] = \sum_i a_i q^{d-i} = \delta q^d + O(q^{d-1})$$
\qed
\

For curves, Weil's ``positivity'' proof of the Riemann Hypothesis
(\cite{Weil})  works here too.  ( The notation follows \S \ref{correspondences}.  )

\<{example}\lbl{weil}  Let $C$ be a smooth, complete curve of genus $g$ over
a field $K$ of characteristic $p>0$.  Let $q=p^m$.  Let $\phi_q$
be the $q$-Frobenius, and let $C' = C^{\phi_q}$.  Let  $S \leq C \times C'$ be an irreducible subvariety.  Then
 $$S \cdot \Phi_q  =  q \deg_{cor} (S) + deg_{cor}(S^t) + e$$
with 
$$|e| \leq ((2g)  (2 \deg_{cor}(S) \deg_{cor}(S^t) - |S \cdot S^t|))^{1/2} q^{1/2} (\leq O(1)q^{1/2})$$
\end{example}

\proof  Let $A_0$ be the group of divisors on $C  \times C'$,
up to rational equivalence, $A = \Qq \tensor A_0$.   
 Define a symmetric bilinear form:
\def\dc{ { \deg_{cor}} }
$$\beta(X,Y) = \deg_{cor} (X) \deg_{cor}(Y^t) + \deg_{cor}(Y)\dc(X^t) - X \cdot Y$$

If $X$,$Y$ are divisors on $C  \times C^{'}$, let $X^{t}Y$ denote
the intersection-theoretic composition $X^t \circ Y$.  It is a
correspondence on $C \times C$; and we have:
$$ \Delta.(X^{t}Y)  = X.Y, \ \ \ \ deg_{cor}(X^{t}Y) = deg_{cor}(Y)deg_{cor}(X^{t})$$

Thus $\beta(X,Y) = \beta_C(X^{t}Y,\Delta)$ where $\beta_C$
is defined like $\beta$, but on $C \times C$.

Let $Z= \Phi_q$.  
We have $Z^{t}Z = ZZ^{t} = q\Delta$ (where $\Delta$
denotes the diagonal of $C$, respectively $C'$.)
By associativity of composition, for any $X$ and $Y$,
$$(Z^{t}X)^{t}(Z^{t}Y) = X^{t}(ZZ^{t})Y = qX^{t}Y$$
Hence $\beta(X,Y) = \beta_C(X^{t}Y, \Delta) = (1/q)\beta_C((Z^{t}X)^{t},(Z^{t}Y))$. 
Now according to Weil,  $\beta_C(W^t,W) \geq 0$; hence 
$\beta(X,X) \geq 0$.    This non-negativity extends to $A \tensor \Rr$.
 
  By Cauchy-
Schwartz, we have $\beta(X,S) \leq \beta(X,X)^{1/2} \beta(S,S)^{1/2}$.
Now
$$\beta(Z,Z) = \beta_C(q \Delta, \Delta) = q \beta_C(\Delta,\Delta) =
q(2-  \Delta \cdot \Delta ) = q(2+(g-2)) =  2gq$$
$$(Z \cdot S) = qdeg(S) + deg (S^t)  - \beta(S,\Phi_q)$$
So  $e^2 \leq \beta(S,S) (2gq)$, and the estimate follows.  \qed

But in general, we will have to decompose $[S]$ and $\Phi_q$ 
cohomologically
and not cycle-theoretically.

It will be convenient to renormalize the norm, and
write $\doub{S} = {2d \choose d}^{2}\abs{S}$.   Then (by \ref{cy})
$\doub{S \circ T} \leq   {2d \choose d}^{4}  \abs{S} \abs{T} =  
\doub{S}\doub{T}$

When dealing with the \'etale cohomology groups, we
will always work over an algebraically closed field $k$, fix a prime
$l \neq char(k)$,
and fix an isomorphism of $\Zz_l(1)$ with $\Zz_l$. See the
section of \cite{deligne1} on ``orientations", and \cite{kleiman}. (For 
the groundwork,
see \cite{sga4h} and the references there, or \cite{freitag} or 
\cite{milne}).

\begin{notation}\lbl{not3}
If $R$ is a correspondence on $X \times X'$, we
 let $\eta_i(R): H^i(X,\Qq_l) \to H^i(X',\Qq_l)$
 be the endomorphism of $H^i(X,\Qq_l)$ induced
by $R$ ; cf \cite{kleiman}, 1.3.
 When $X=X'$, we let $\tau_i(R)$ the trace of this
 endomorphism.
\end{notation}

When $R$ is a correspondence on $X$ as in \ref{not3},
the Lefschetz trace formula
(\cite{kleiman}, 1.3.6) relates the intersection number of $R$ with the 
diagonal
to the endomorphisms $\eta_i(R)$ as follows:
\begin{sta} \lbl{lef}

\[ R \cdot \Delta_X = \Sum_{i=1,\ldots,2n} (-1)^i \tau_i(R) \]
\end{sta}

For the top cohomology, we have

\begin{sta}\lbl{top} Let $n = \dim(X)$.   Then $\dim H^{2n}(X) = 1$,
and $\tau_{2n}(R) = \deg_{cor}(R)$ is the degree of $S$ as a 
correspondence (\ref{not2}).
\end{sta}

(The first statement is in \cite{kleiman} 1.2; hence $\tau_{2n}$ acts as
multiplication by a scalar, and this scalar can be seen to be 
$\deg_{cor}(R)$
by computing the effect of $\eta_{2n}(R)$ on the image of a point under 
the cycle map.)

We will make weak use of Deligne's theorem (\cite{deligne1}).

\begin{sta}\lbl{deligne} Let $X$ be a smooth projective variety over
$k$, the algebraic closure of a finite field. Suppose $X$ descends to 
$\Ff_q$.
Let $\Phi$ be the graph of the Frobenius correspondence $x \mapsto x^q$ 
on $X \times X$.
Then all the eigenvalues of $\eta_i(\Phi^{t})$ are algebraic numbers, of
complex absolute value $q^{i/2}$.
\end{sta}

We wish to study eigenvalues of compositions of two correspondences 
$S$,$T$.
 This can probably be done using the remark \ref{2vars} below.
However we will take another route, starting with a slight 
generalization
 of \cite{Lang-Abelian},   V.3     lemma   2.

\begin{lem}\lbl{cp0} Let     $b , t > 0$ and $\omega_j$, $c_j$ be
 complex numbers
($j=1,\ldots,s$). Suppose
$$\limsup t^{-n} \abs{\Sum_j c_j \omega_j^n } \leq b$$
Then one may partition $\set{1,\ldots,s} $ into disjoint sets 
$J_1$,$J_2$
such that $$\Sum_{j \in J_2} c_j \omega_j^n = 0$$ for all $n$,
and
$$\abs{\omega_j} \leq t$$ for each $ j \in J_1$.        \end{lem}

\proof Let $M$ be the maximal norm of an $\omega_j$. Dividing the
$\omega_j$ and $t$ by $M$, we may assume $M=1$.

Consider first the $\omega_j$ lying on $T$,
the group of complex numbers of norm one.
Say they are $\omega_1,\ldots, \omega_{s'}$;
note that summing only over $j \leq s'$ does not effect the $\limsup$ in 
the hypothesis.
Let $S$ be the
closed subgroup of $T^{s'}$ generated by   the point
$\alpha = (\omega_1,\ldots,\omega_{s'})$.
Now $S$ is a compact group. In every neighborhood of $S$, there are 
infinitely
many powers $\alpha^m$; for these we have $(c_1,\ldots,c_{s'}) \cdot 
\alpha^m$ arbitrarily
small; hence there exists a point $\gamma$ in the closure of this 
neighborhood,
with $(c_1,\ldots,c_{s'}) \cdot \gamma= 0$.
So such points are dense in $S$, and thus
$(c_1,\ldots,c_{s'}) \cdot \gamma = 0$ for all $\gamma \in S$. We may 
thus put all the
$\omega_i$, $i \leq s'$,  into $J_2$.   Now removing them from the sum
has no effect, so the hypothesis applies to $\omega_{s'+1},\ldots , 
\omega_s$, and
we may proceed by induction.
 \ \ \ \ \ $\Box$

\begin{lem}\lbl{deg1} Let $R$
and $T$ be commuting correspondences on a smooth projective
variety $X$ of dimension $d$ over an algebraically closed field $K$,
$T\circ R = R \circ T$.
Fix an embedding of $\Qq_l$ into $\Cc$.
Let $c_{ij}, e_{ij}$ be the corresponding
eigenvalues of the endomorphisms $\eta_i(R)$, $\eta_i(T)$, considered
as complex numbers.
Let $J_1$ be the set of $(i,j)$ such that $e_{ij}$ has complex
norm at most $\doub{T}$, and $J_2$ the rest.
 Then for any $m$, $k$
\[   \Sum _{ (i,j) \in J_2}   (-1)^i c_{ij}^k e_{ij}^m = 0 \]
\end{lem}

\proof
  The following claim is immediate
from \ref{cy} (3) and (4):

\claim For all $m$, $ (R \circ T^{\circ m} \cdot \Delta_X) \leq \abs{R} 
\doub{T}^m$

Let $\gamma = \sum_{(i,j) \in J_1} \abs{c_{ij}} + |R|$.
By definition of $J_1$,
\[ \abs{ \Sum _{(i,j) \in J_1} (-1)^i c_{ij}e_{ij}^m } \leq  (\gamma - 
|R|) {\doub{T}}^m \]
By \ref{lef},
\[ (R \circ T^{\circ m}) \cdot \Delta_X = \Sum_i (-1)^i\tau_i(RT^m) =
         \Sum_{ij} (-1)^i c_{ij} e_{ij}^m \]
Thus

\[ \abs{\Sum_{(i,j) \in J_2}(-1)^i c_{ij} e_{ij}^m }
 \leq \gamma {\doub{T}}^m \]
By \ref{cp0}, one can partition $J_2$
into disjoint sets $J_3$,$J_4$
such that $\Sum_{(i,j) \in J_3}(-1)^i c_{ij} e_{ij}^m = 0$
for all positive integers $m$,
and
$\abs{e_{ij}} \leq {\doub{T}}$ for each $(i,j) \in J_4$.
But then $J_4 \subseteq J_1$ by definition of $J_1$, so
$J_4 = \emptyset$. Thus $J_2 = J_3$, and we have shown what
we wanted for $k=1$. Applying this to $R^k$ in place of $R$
yields the lemma. \ \ \ \ \ $\Box$

The following Proposition does not mention Frobenius, and could be
stated over arbitrary fields, but the proof we give uses 
\cite{deligne1}.

\begin{prop}\lbl{bound}
 Let $R$ be a correspondence on a smooth projective
variety $X$ over $k$, the algebraic closure of a finite field.
Then every eigenvalue of $\eta_i(R)$ is an algebraic number,
of complex absolute value at most $\doub{R}$.
\end{prop}

It suffices to show that in every complex embedding, every
eigenvalue has absolute value at most ${\doub{r}}$. Fix therefore
a complex embedding, and view the eigenvalues of $\eta_i(R)$ as complex 
numbers
$e_{ij}$. Let $T = \Phi_q^{t}$ be a Frobenius correspondence on $X$,
with $q$ large enough so that $T \circ R = R \circ T$ as 
correspondences.
Let $J_2$ index the set of $e_{ij}$ of absolute value $>|R|$.
Let $c_{ij}$ be the eigenvalues of $T$, with corresponding indices.
By \ref{deg1}, for any $m$ and $k$,

\[       \Sum _{ (i,j) \in J_2}   (-1)^i c_{ij}^m e_{ij}^k = 0 \]

Suppose for contradiction $J_2$ is nonempty. Let $i_{\max}$ be the 
highest
index represented in $J_2$, and let $J_3 = \{j: (i_{\max},j) \in J_2\}$.
Let $d_{ij} = (-1)^ic_{ij}q^{-i_{\max}/2} $, and let $d_j = 
d_{i_{\max},j}$,
$e_j= e_{i_{\max},j}$.
Then $|d_j|=1$, $|d_{ij}|<1$ for $i< i_{\max}$, and
$\Sum _{ (i,j) \in J_2}   d_{ij}^m e_{ij}^k = 0 $ for any $k$,$m$.
For fixed $m$,$k$ we have
$\Sum _{ (i,j) \in J_2}   d_{ij}^{m+m'} e_{ij}^k = 0 $
Let $m'$ approach infinity, becoming highly divisible;
then   $d_{ij}^{m'}$ approaches $0$
for $i < i_{\max}$, while $d_j^{m'}$ approaches $1$.   We obtain
\[ \Sum _{ j \in J_3}   d_j^m e_{j}^k = 0 \] for any $k$,$m$.
Dually, let $J_4 = \{j \in J_3: |e_j| = L \}$, where $L$ is
the highest value attained by an $e_j$. Then as above we get
\[ \Sum _{ j \in J_4}   d_j^m (e_{j}/L)^k = 0 \] for any $k$,$m$.
But now letting both $m$ and $k$ approach infinity, each $d_j^m$
and each $(e_{j}/L)^k$ approach $1$, and we obtain $|J_4| = 0$,
a contradiction. Thus $J_2 = \emptyset$ and we are done.

\begin{cor}\lbl{bound2} Let $X$ be a smooth complete variety over
the algebraic closure $k$ of a finite field.
Let $X'$ be the image of $X$
under the $q$-Frobenius automorphism of $K$, and let
$T = \Phi_q^{t} \subset X' \times X$ be the transposed Frobenius
correspondence. Let $R \subset X \times X'$ be a correspondence.
Then every eigenvalue of $\eta_i(T \circ R)$ is an algebraic number,
of complex absolute value at most $q^{i/2}\doub{R}$.
\end{cor}

\proof For some $m$, $X$ and $S$ are defined over $\Ff_{q^m}$.
Taking $m$'th powers, it suffices to prove that every eigenvalue of
$\eta_i((T \circ R)^m)$ is an algebraic number,
of absolute value at most $q^{im/2}{\doub{R}}^m$. Let $X_i$
be the image of $X$ under $\phi_q^i$, and $R_i \subset X_i \times 
X_{i+1}$
the image of $R$. Then
\[ (T \circ R)^m = \Phi_{q^m}^t \circ R_{m-1}\circ \ldots \circ R \]
Let $R^* = R^{m-1}\circ \ldots \circ R$. Then  (\ref{cy})
$\doub{R^*}
\leq \Pi_{i=0,\ldots,m-1} \doub{R_i} = {\doub{R}}^m$.
$R^*$ is a correspondence on $X$.
By \ref{bound}, every eigenvalue of $\eta_i(R^*)$ is algebraic of 
absolute
value at most $|R|^m$. Now note that $ \Phi_{q^m}^t$ is
a correspondence on $X$, commuting with $R^*$, and (\ref{deligne}) with
algebraic eigenvalues of absolute value $q^{im/2}$.
Thus $\eta_i( (T \circ R)^m) = \eta_i( \Phi_{q^m}^t \circ R^* )$
has algebraic eigenvalues of absolute value $({\doub{R}}q^{i /2})^m$.
The conclusion follows.

\begin{lem}\lbl{2qcr} In Theorem \ref{2qc} we may assume $K_q$ is the algebraic
closure of the prime field; in other words it suffices
to prove the result for $b$ chosen from a finite field. \end{lem}
\proof   We may assume $V$ (and $V \times_B V$) are flat
over $B$.
Let $T$ be an irreducible component of $B' \meet \beta^{-1} \Phi_{B  }$, 
where
$\Phi_B \subset B^2$ is the graph of Frobenius on $B$.
By \cite{fulton} (10.2), the number
$${\bar S}_b \cdot_{\bar V_{b_1} \times \bar V_{b_2}} \Phi_b$$
is constant for $(b_1,b_2) \in T$. Thus one
can replace $X,X'$ by $V_{b_1},V_{b_2}$, where
$b$ is a point of $T$ rational over a finite field.  \qed

\

{\bf proof of  Theorem \ref{2qc} }

By \ref{2qcr}, we may assume $K$ is the algebraic
closure of a finite field.
By \ref{lef},
\[ S_b \cdot \Phi_b = ((\Phi_b)^t \circ S_b) \cdot \Delta_X
    = \Sum_{i \leq 2d} \tau_i((\Phi_b)^t \circ S_b)\]
By \ref{bound2},
$$|\Sum_{i < 2d} \tau_i((\Phi_b)^t \circ S_b)| \leq {2d \choose d}^{2} 
\mu \nu q^{d-\half} $$
By \ref{top}, $\tau_{2d} ( (\Phi_b)^t \circ S_b) =
\deg_{cor}((\Phi_b)^t \circ S_b) = \deg_{cor}(S_b)q^d$. \qed

We conclude the section with the rationality lemma mentioned but not 
used
above.

\begin{lem}\lbl{2vars}  Let $X$ be a smooth projective variety over a 
field $K$,
and let $S$, $T$ be correspondences on $X$.  Write a power series in two
variables:
$$ F(s,t) = \Sum_{m,n \geq 0} S^m \cdot T^n s^mt^n$$
Then $F$ is a rational function, with denominator of the form
$f(s)g(t)$.
\end{lem}

\proof  Using Grothendieck's cohomological representation, it suffices
to prove more generally that if $V$ is a complex vector space, and
$S,U$ two linear transformations of $V$, then
$ \Sum_{m,n \geq 0} tr(S^m U^n) s^mt^n$
is a rational function of two variables.

We will actually show that $$\Sum_{m,n} S^mU^n s^mt^n \in 
End(V)[[s,t]]$$ lies
in $End(V)\tensor_{\Cc}C(s,t)$, so that all the matrix coefficients
are rational functions.  (With denominators as specified.)
Now
$$ \Sum_{m,n \geq 0} S^m U^n s^mt^n = (\Sum_m S^m s^m)(\Sum_n U^nt^n 
) $$
So it suffices to consider $\Sum_m S^m s^m$.    We can
sum the geometric series:
$$ \Sum_{m \geq 0} S^m s^m = (1-Ss)^{-1} = (det(1-Ss)^{-1})S'$$
For a certain matrix $S'$ with coefficients polynomial in $s$.  
\qed

\end{section}
\begin{section}{Equivalence of the infinite components:  asymptotic estimate}
 
\lbl{equivalence}

While Theorem \ref{1} promises actual points of intersection, Theorem 
 \ref{2qc} 
refers rather to an intersection number, i.e. to the intersection with a 
better placed $S'$ rationally equivalent to $S$.  
We will now bridge the gap by giving an asymptotic estimate of the 
difference between the number of isolated points of these intersections.
Here we do so  numerically and roughly (up to lower order of magnitude); 
this requires only tying together some previous threads.  This suffices for the
proof of Theorem \ref{1}.  However our methods are essentially precise and
motivic, and later we will describe the language needed to bring this out.

We will actually bound  at once both the 
equivalence of the infinite components, and the contribution of isolated points
of high multiplicity.

We work with a single fiber of the data \ref{not}:
for some $b \in B'(K_q)$ with $\phi_q(b_1)=b_2$, we let
$\tilde X_1 = V_{b_1}$, $\tilde X_2 = V_{b_2}$, $\tilde S = S_b$.
Thus we have a $d$-dimensional variety $\tilde X_1$ over $K_q$,
$\tilde X_2 = \tilde X_1^{\phi_q}$, and an irreducible subvariety
$\tilde S \subset \tilde X_1 \times \tilde X_2$.
We will bound a certain quantity associated
with the intersection of $\tilde S$ with the graph of Frobenius.
We will write $O(q^l)$ for a quantity bounded by $Cq^l$,
where $C$ is a constant independent of $q$, and dependent on $\tilde 
X_1$,
$\tilde X_2$,
$\tilde S$ in a way that remains bounded when the varieties are
obtained as above from the data \ref{not}, with $b$ varying.

    We assume (cf. \ref{sep},\ref{smooth}) that
the projection from $\tilde S$ to $\tilde X_2$ is \'etale, and that 
$\tilde X_1$ is an open
subvariety of a smooth, complete variety $ X_1$. (And thus $\tilde X_2$ 
of
$X_2 = X_1^{\phi_q}$). Let $\Phi = \Phi_q \subset Y=  X_1 \times X_2$ denote 
the
graph of $\phi_q: X_1 \to X_2$, and let $S$ be the Zariski closure 
of $\tilde S$
in $X_1 \times X_2$.

\begin{prop}\lbl{.1}
 Let $e = (\Phi_q \cdot S) - |\{\Phi_q \meet \tilde{S}\}|$. Then
$|e| \leq O(q^{d-1})$ \end{prop}

Observe that the points of $\Phi_q \meet \tilde{S}$ are simple, by
 by \ref{transverse}; so  it makes no difference whether 
$|\{\Phi_q \meet \tilde{S}\}|$ is counted with multiplicities. 
 

 \<{paragraph}{Almost all Frobenius specializations}

Recall that a Frobenius field $K_q$ is an  algebraically closed difference field of characteristic $p>0$ 
 made into a difference field via   $x \mapsto x^q$ ($q=p^m, m \geq 1$).   
 $\size{Y}$ denotes the   size  of a   0-dimensional scheme, i.e.  the  number of points  with multiplicities. 
 
We will consider   difference schemes $Y$ of finite type over a
  difference field  $k$.  $Y$ arises by base extension from a scheme $Y_D$ over a finitely generated
difference domain $D$ with fraction field $K$.  We will say:  "for almost all $(q,h), \, \cdots$ " to mean:  for some such $D$ and $Y_D$, for
all sufficiently large prime powers $q$, and all $ h:  M_q(D) \to K_q \, \ldots$. 

When $K$ is infinite, by enlarging $D$, one sees that a simpler formulation is equivalent:

"for some $D$, $Y_D$, for all prime powers $q$ and all $h: h:  M_q(D) \to K_q \ldots$."

When concerned with only almost all  $(q,h)$, we will not mind increasing $D$ within $K$, and so can assume that some such   $Y_D$ has been fixed, and write $Y_{q,h}$ for $M_q(Y_D) \tensor _h K_q$.   
Different choices of $Y_D$ will give the same $Y_{q,h}$ for almost all $(q,h)$.

    Let $D \subset k$ be a finitely generated difference domain, and let $D[t]'$ be a finitely generated
sub-$D$-difference algebra of $\kt$.  By enlarging $D[t]'$, and $D$, we may obtain the form $D[t]' = D[t^{\si-n},f^{-1}]_\si$,
where $f \in D[t]_\si$ and $f(0) $ is invertible in $D$.  Consider morphisms $h:  M_q(D) \to K_q$ into Frobenius fields $K_q$
with $q \geq n$.  For any such $h$, let $\htt:  D[t]' \to K_q(t)$ be the unique difference ring morphism extending $h$ and
with $\htt(t^{\si-n}) = t^{q-n}$.  Sometimes we will consider it as a map 
$\ D[t]' \to K_q[t,f^{-1}(t)]$ or
$\ D[t]' \to K_q[\breve{t}] :=K_q[t,g^{-1}: g(0) \neq 0]$; in this case we will
write $\hT$.   

When $Y$ is a difference scheme of finite type over $\kt$, 
we fix $D[t]'$ as above; given $q$ and $h: M_q(D) \to K_q$,
we define $Y_{q,h} = Y_{q,\htt}$.

 
 \>{paragraph}

Let $k$ be a perfect, inversive difference field, $V $ an absolutely irreducible, projective (or complete) algebraic variety over $k$,   $\dim(V)=d$.

Let $\CV   $ be the set of  subvarieties $S$ of $V \times V^\si$, such that  any component $C$ of $S$ (over $k^{alg}$) has dimension $d$, 
and the projections   $\pr_0^{S}: S \to V, \pr_1^{S}: S \to V^\si$ 
are dominant.  Let $\ZCV$ denote the free Abelian group generated by $\CV$.
For any subscheme $W$ of $V \times V^\si$, let $[W]$ be the formal sum of the components of 
$W$ of dimension $d$, weighted by their geometric multiplicities.  (cf. \cite{fulton}.)

Write $S \in \CVsep$ if in addition, 
  $k(S)$ is a separable finite extension of $k(V^\si)$. 
In this case,  the projection $\pr_1^{S}$ is \'etale above
a  nonempty Zariski open subset of $V^\si$.  Let $ V_{\etl}(S)$
be the largest smooth open subvariety $U$ of $V$ such that $\pr_1^{S}$ is \'etale above $U^\si$.

In general, for   $S \in \CV$, let $ \pr_1[1](S )  = \{a \in V:  \dim {\pr_1^S}^{-1}(a) >0 \}$.  Let $V_{fin}(S)  = (V \setminus \pr_1[1](S)^{\si^{-1}}$, and 
$\tS = S \meet (V_{fin}(S) \times V_{fin}(S)^\si)$.  Thus $V_{fin}(S)$
is   the 
largest open subvariety of $V$ such that, with this definition,
$\pr_1: \tS \to \tV^\si$ is quasi-finite.

Let $X(S) = (\tS \star \Sigma) $.  So $X(S)$ has total dimension $d$.

 Let  $\CVf = \{S \in \CV:  \pr_0[1](S ) \meet (S \star \Sigma) = \emptyset \}= \{S \in \CV:  X(S) = S \star \Sigma \}$.  
  (cf. \reff{zd}.)  $\ZCVf$ is the free Abelian group generated by $\CVf$.

The $X(S)$ with $S \in \CVf$ are  directly presented,
and of total dimension $d$; and moreover of   transformal 
multiplicity $0$.  

	{\bf Asymptotic equivalence} 
Given $S \in \CV$, defined over a finitely
generated difference domain $D$, let $X=X(S)$,   let $\eta_0(S;q,h)$ be
the number of points in $X_{q,h}(K_q)$, and let $\eta(S;q,h) = \delta'_p \eta_0(S;q,h)$
(cf. \ref{not}).  


We say $S,S' \in \CV$ are {\em asymptotically equinumerous} ($S \sim_a S'$) if
if there exists a finitely
generated difference domain $D \subset k$ (containing any given finite set), such that for some $b$, for almost all $q$, for all $h: M_q(D) \to K_q$, 
$|\eta(S;q,h) - \eta(S';q,h)| \leq bq^{d-1}$.

We extend the terminology  to $\ZCV$ by additivity.

\<{remark} \rm We could instead count the points of $X(S)_{q,h}$ with their geometric
multiplicities, or with appropriate intersection multiplicities; asymptotically
this makes no difference.    However, it is natural to extend the definition to non-reduced schemes $S$.  In this case one should 
count points with geometric multiplicities.  
\>{remark}

Let $\NCV$ be the group of cycles in $\ZCV$ that are reduced- rationally equivalent
to $0$ over $\pr_1$.

Recall
\reff{move2}:

\<{lem}\lbl{ae0}  For any $S \in \CVsep$ there exists
  $S' \in \Zz \CVf$ with $S-S' \in \NCV$.   
\>{lem}

\<{lem} \lbl{ae1}  Let $S \in \CV$, and let $Z$ be a proper Zariski closed subset of $V$.
Then $\size{X_{q,h} \meet Z} =O(q^{d-1})$.  \>{lem}
\proof By \reff{zd}, \reff{si-dim+++}.   \qed

\<{lem} \lbl{ae2}  Let $\uS$ be an  irreducible $d+1$-dimensional
variety of $(V \times V^\si) \times \Pp^1$.  Let $t$ denote a generic point
of $\Pp^1$ over $k$.   Assume  $\uS_0$ is generically reduced and generically smooth over $\pr_1$.
  Then $S_t \sim_a S_0$.  \>{lem}

(Here  we identify $S_0$, a variety over $k$, with the same  variety
viewed by base change as a variety over $k(t)$.   

\proof  Write $S_t = \uS_t $,   $S_0 = \uS_0$, and let
Let 
  $V_ {\Ak} = V \times_{\spe k} {\Ak}$. 
  $V_t = V \times_{  k} {  {k(t)_{\si}}} $, $V_0=V$ denote the fibers of this map over the generic point of $\Ak$ and the point $t=0$ respectively.
 
  Let $X_t $ be the closure of $X(S_t)$ in $V_t$ (i.e. the smallest $k(t)_\si$-
definable closed difference subscheme containing $X(S_t)$.)   
  Let $X$ be  the closure in $V_{\Ak}$ of $X_t$; it is a  $k$- closed difference subscheme of $\Vv = [\si]_k V \times_k \Ak$, flat over $\Ak$, whose fiber over $k(t)_\si$ is $X_t$. 


Let $W = V_0 \setminus V_{\etl}(S_0)$; and let $^*W$ be as in Proposition \reff{key}
( a formula in the language of transformal valued fields over $k(t)_\si$.)
By definition of $X(S_t)$, and \reff{zd}, condition (1) of \reff{key} holds; 
hence it holds also for the closure $X_t$.  (Cf. \reff{dim-cl}).    Since
$S_0 \in \CV$, condition (2) of \reff{key} holds too; see \ref{zd} (2). Thus (4) holds,
i.e. $^*W$ has inertial dimension $< \dim(V)$.  

Consider a (large) prime power $q$.  Fix a discrete valuation of $L=L_q=K_q(t)^a$ over $K_q$ with $\val(t)=1$; let $R_q$ denote the
valuation ring.  We will show that $\size{ X(S_t)_{q,h}(L)  } =  \size{ X(S_0)(L)} \ \mod O(q^{d-1})$.
Since $q$ is large enough, $X(S_0)_{q,h}$ is a finite scheme;
so $X(S_0)_{q,h} (L) =  X(S_0)_{q,h}(K_q)$.  

{\bf  Since $V_t$ is projective}, $V_t(R_q) = V_t(L)$, and
the residue homomorphism induces a map $r_V: V_t(L) \to V_0(L)$;
restricting to $r_0: (X_t)_{q,h} (L) \to  (X_0)_{q,h}(K_q)$.  
Let  $Y =  (X_t)_{q,h} \setminus ^*W$.  
Let $r_1= r_0| Y $.  By definition of $^*W$, $r_1: Y \to V_{\etl}(S_0)$.

{\bf Injectivity of $r_1$} follows from \ref{res-inj}.    An alternative argument:  by 
\reff{transverse}, $X(S_0)_{q,h}$ has only simple points on $V_{\etl}(S_0)$;
by 
\reff{gauss},  $r_1$ must be injective.

{\bf Surjectivity.}  We can use a Hensel lemma approach
(analogous to \reff{res-surj};   here only algebraic and not transformal
valuation fields are in question,   \reff{res-surj-p}.)  

But we prefer to 
vary the visualization.  Let $\phi$ denote the $q$-Frobenius,
$\Phi \subset V \times V^\phi$ the graph of $\phi$, 
$\uPhi =  \Phi  \times \Pp^1$.  Then $\uS \meet \uPhi$ is an algebraic
set of relative dimension $0$ over $\Pp^1$.  By the dimension theorem, it has no isolated points.
So it consists of horizontal curves,  together with other curves and higher-dimensional
components above the non-quasi-finite locus
Let $p \in (X_0)_{q,h}(K_q)$,  $p \notin W$;
then
$(p,0)$  must lie on one of the 
horizontal curves $C$ of  $\uS \meet \uPhi$.  Let
$(p',t) \in C$ with $t$ a generic point of $\Pp^1$.   
Then  $p' \in X(S_t)$. 
(If $p' \in  \pr_1[1](S_t)$ then $p \in \pr_1[1](S_0) \subseteq W$.)
 Thus there exists a place of $K_q(t)^a$ into $K_q$
with $p' \mapsto p$.   For the valuation corresponding to this place,
$p$ lies in the image of $X_t(L)$ under the residue map; but all
such valuations are isomorphic.  This  proves the surjectivity of $r_1$.  \qed

We thus exhibited a bijection between
$X(S_0)_{q,h} (K_q) \setminus W = X(S_0)_{q,h} (K_q) \meet V_{\etl}(S)$ and 
 $X(S_t)_{q,h}(L) \setminus \ {^*W} (L)$.  Now by  \reff{ae1}, 
at most $O(q^{d-1})$ of $X(S_0)_{q,h} (K_q)$
lies on $W$ (even when counted with multiplicities on $X(S_0)_{q,h} (K_q)$.)   By \reff{R-descent} (cf. \reff{R-descent-ex}), $^*W(L)$ has $O(q^{d-1})$ points.  This proves the lemma.  \qed 

\<{remark} \label{macaulay}   \rm With the natural definition of asymptotic equivalence for non-reduced schemes, counting with geometric multiplicities, \lemref{ae2} remains valid when $S_0$ is not   reduced. This is due to Yves Laszlo, who   noted that the  weaker moving lemma used in a previous version requires  it.    The  complication in the non-reduced case is that the dimension theorem, as used in the proof of \lemref{ae2} above, is required also for embedded components and not only the apparent ones; this implies 
flatness in \lemref{gauss}, allowing to improve the inequality there to an equality.  The required strong version of the dimension theorem      is the unmixedness theorem for Cohen-Macaulay rings.  Essentially, the   correspondence $S$ is generically Cohen-Macaulay, and this is maintained under intersection with Frobenius.   \>{remark}   

\<{lem} \lbl{ae3}  Let $S \in \NCV$.  Then $S \sim_a 0$.
\>{lem}

\proof
We may take $S = [\uS_0] - [\uS_\infty]$, $ \uS$ an irreducible $d+1$-dimensional
variety of $(V \times V^\si) \times \Pp^1$, with $\uS_0,\uS_\infty$ generically reduced and generically smooth over $\pr_1$.  Let $t$ denote a generic point
of $\Pp^1$ over $k$.  Then by \reff{ae2}, 
$[\uS_t] \sim_a [\uS_0]$ and $[\uS_t] \sim_a [\uS_\infty]$.  So
$[\uS_0] - [\uS_\infty] \sim_a 0$ , as required.  \qed

We repeat the statement of \ref{.1} in this language, and prove it.
\<{cor}  Let $S \in \CV$.  
For almost all $(q,h)$, the numbers
$  \size{X(S)_{q,h}}$, $ (S_{q,h} \cdot \Phi_q) $, $\Phi_q \meet \tilde{S}$
all differ by at most $O(q^{d-1})$.  \>{cor}

\proof  We will use \cite{fulton} Prop. 7.1 and the notions of  multiplicity there.  
By \ref{sep}, 
we may assume  in addition that $S \in CVsep$. 
By \ref{ae0}, there exists $S' \in \Zz \CVf$
with $S-S' \in \NCV$.   By \ref{ae3}, $S \sim_a S'$.  Since rational
equivalence preserves intersection numbers, $S'_{q,h} \cdot \Phi_q = 
S _{q,h} \cdot \Phi_q$.  Thus we may assume $S=S'$, i.e.
we may assume $S \in \CVf$.   By \ref{transverse},
the intersection $S_{q,h} \meet \Phi_q$ on $(V  \times V^\si)_{q,h}$
is proper; so $ (S_{q,h} \cdot \Phi_q) $ is the number of points
of $S_{q,h} \meet \Phi_q$, counted with intersection multiplicities.
On $V_{\etl}(S)_{q,h}$ the geometric multiplicity of the intersection
is $1$; while on $Z = V \setminus V_{\etl}(S)$, by \ref{ae1}, the 
total number of points of intersection (with geometric multiplicities)
is $O(q^{d-1})$.  Thus  counting the points
of  $X(S)_{q,h}$ with
intersection or with geometric multiplicities, or without multiplicities agree up to $O(q^{d-1})$.  Since the intersection multiplicities are sandwiched
in between, and have value $1$ whenever the geometric multiplicity is $1$,
 they    yield the same result.  \qed

\end{section} 

\begin{section}{Proofs and applications}

\lbl{proofs}
   

\begin{subsection}{Uniformity and Decidability}

Throughout this paper we have been dealing with a variable Frobenius
on an essentially fixed variety. When the variety is not itself
defined over the fixed field of the Frobenius, this is slightly strange;
given the algebraic variety $Y = X \times X^\phi$ there is (at most) one
Frobenius on it. We have explained the uniformity in three different ways:

\begin{enumerate}

\item   By fixing, not $Y$ and $S$, but only the {\em degrees} of their
projective completions. (As in Theorem \ref{1})

\item  Using the formalism of \ref{not} (As in 1B)

\item  Using difference schemes. In this language, one could state Theorem 1A thus:

\end{enumerate}

\noindent{\bf Theorem {\ref{1}C}}  {\em  Let $D$ be a finitely generated difference
domain, $X$ an absolutely transformally integral
 difference scheme of total dimension $d$ over $\spe D$.
 Then there exist
$b \in \Nn$ and $c \in D$ such that  for  prime power $q > b$,and any 
$y \in M_q(D[c^{-1}])$,
$X_{q,y}$ is nonempty, indeed has $\delta^* q^d+O(q^{d-1/2})$ distinct points.}

Here if $K$ is the field of fractions of $D$ and $L$ is the transformal function field of $X$, 
then $\delta^* = \deg_{lim}(L/K) / \iota ' (L/K)$, where $\deg_{lim}$ is the limit degree, 
and $\iota'$ is the purely inseparable dual degree, cf. \ref{towers}.
\smallskip

\smallskip

Uniformity statements such as Theorem \ref{1} amount to the same
in any of the formulations.   To demonstrate this, we
show (2) implies (3) implies (1), and (2).

\smallskip
\noindent{\em proof that Theorem \ref{1}B is equivalent to Theorem \ref{1}C}

In effect Theorems \ref{1}B and \ref{1}C make the same statement about points of
a difference
scheme ${\cal X}$;  but Theorem \ref{1}B assumes the generic fiber  
$X$ of ${\cal X}$ is directly presented, while Theorem \ref{1}C
assumes instead that $X$ is irreducible.   One may   pass from
the directly presented to the general case using   \ref{direct-ritt-0};
and from the irreducible case to the general case using \ref{direct2}.

\smallskip

\noindent{\em proof of Theorem 1B}

We make the assumption stated in \ref{smooth}, and
use the notation there.  (We also accept the conclusion of
\ref{sep}.) By Theorem \ref{2qc},
  ${\bar S} \cdot \Phi_q$ has degree as stated in 1B.  Theorem 1B
now follows from Proposition \ref{.1}.
 
\qed

\smallskip

\noindent{\em proof  that Theorem \ref{1}B implies Theorem 1}

 Suppose     given a family $X_i$ of affine varieties over algebraically closed
fields $k_i$
of characteristic $p_i$,
and $S_i \subset (X \times X^{\phi_{q_i}})$, with the assumptions
of Theorem \ref{1} holding, and the degrees of the projective completions of the $X_i$ and $S_i$
bounded, but with  
 $(q^{-(d-\half)}(|S(k) \meet \Phi_q(k)|)- aq^d$ unbounded.
Let $(L,\si)$ be an ultraproduct of the $(k_i,\phi_{q_i})$. Since the degrees and dimensions
are
bounded, the corresponding ultraproducts of the $X_i$ and $S_i$ are ordinary
subvarieties $X$,$S \subset X \times X^\si$
of finite dimensional affine space over $L$. Let $D$ be a finitely generated
difference sub-domain of $L$, with $X$,$S$ defined over $L$; replacing
$D$ by a localization, we may assume $X$,$S$ remain absolutely
irreducible varieties of dimension $d$ when reduced modulo
any prime ideal of $D$.  This puts us in the situation of Theorem \ref{1}B.
 For almost
every $i$, one has a difference  ring homomorphism $D \to k_i$, mapping $X$ to $X_i$.
By Theorem \ref{1}B, $X_i  \meet \Phi_q$ has about $aq^d$ points.     A contradiction.
\qed

\begin{cor}\lbl{13.1} Let $k$ be an algebraically closed difference field, $X$ 
a difference variety over $K$ of transformal dimension $l$.  Then there
exists a finitely generated difference domain $D \subset k$, such that $X$
descends to $D$, and for almost all $q$, for all $h: D \to K_q$, 
$\dim(X_{q,h}) = l$.  \end{cor}

\proof  $X$ admits a definable map to $[\si]_K \Aa^l$, whose image
contains the complement of a proper difference subvariety of $[\si]_K \Aa^l$.
Thus for almost all $q,h$,$X_{q,h}$ admits a dominant map to $\Aa^l$,
and so has dimension $\geq l$.   The converse is \ref{si-dim}. \qed

We note one more form of Theorem \ref{1}; where in effect
quasi-finiteness is replaced by the weaker assumption of finite total dimension.

\begin{thm}\lbl{1D} Let $V$ be a quasi-projective variety over $K_q$, and let
 $S \subset (V \times V^{\phi_q})$
be an irreducible subvariety. Assume $\dim(S)=\dim(V)=d$, and
$\rm{pr}_{V} |S$ is a generically quasi-finite map of degree $a$.   If the scheme
$S \meet \Phi_q$ is 0-dimensional, then it has
size $\leq aq^d + O(q^{d-\half})$ \end{thm}

\proof  By taking ultraproducts, we find $(V,S)$ over a difference
domain $D$, with $X= [\si]S \meet \Sigma$ of  finite total dimension;
and we must show that $|X_{q,h} | \leq aq^d + O(q^{d-\half})$.
By Theorem \ref{1}, it suffices to show that $X$ has total dimension $d$, 
and that every algebraic component of $X$ of total dimension $d$ is Zariski
dense in $V$.  This follows from \ref{www}.
\qed

\begin{paragraph}{ACFA}
\proof  of Theorem \ref{dec}  That $T_\infty$ contains ACFA follows
from  Theorem \ref{1}.   Now it is shown in \cite{CH} that ACFA 
is nearly complete:   the full elementary theory of a model $(K,\sigma)$
of ACFA is determined by the isomorphism type of $(K_0,\sigma | K_0)$,
where $K_0$ is the subfield of $K$
consisting of points algebraic over the prime field.  By the Cebotarev density theorem,
every  isomorphism type of $(K_0,\sigma | K_0)$ can occur with $\si_0$ 
an ultrapower of Frobenius maps.  So every completion of ACFA is consistent
with  $T_\infty = ACFA$.  Thus every sentence of $T_\infty$ is a consequence of ACFA.
\end{paragraph}
\end{subsection}

\begin{paragraph}{Decidability}

The decidability referred to in this paper, in particular in Theorem \ref{4}, is
in the sense of G\"odel; it corresponds to the dichotomy: finite/infinite, and not
to: finite/bounded, or to distinctions between different degrees of boundedness.
{\em The proof of Proposition \ref{.1} can routinely be seen to be  effective in this sense.}
We mention two instances.   Suppose one has a smooth variety $Y$ over a finite field, and
two subvarieties $V$,$X$ on $Y$, meeting improperly. One knows there
exists $Y'$ rationally equivalent to $Y$, intersecting $X$ transversally. Then
without further work, this is effective: it is merely necessary to search
(over some finite field extension) for a $Y'$ and for data demonstrating the
rational equivalence.   To take another example, suppose a sentence is shown to be true in 
every   (boolean-valued) $\omega$-increasing transformal valuation field, by whatever methods.
  Since the
class of such transformal valuation fields is an elementary class, one can
search for an elementary formal proof, with assured success; the proof
will only use that the valuation is $m$-increasing, for some $m$. Thus   the statement will be true in every   valuation field with the automorphism $\si(x)=x^q$,  as soon as $q \geq m$; and we have found $m$ effectively.

Theorems I  or 1B  follow  from Propositions \ref{2qc} and \ref{.1}.  To use
Theorem \ref{2qc}, we require an
effective bound on
the Betti numbers  $\dim H^i(X,\Qq_l)$ of a
smooth projective variety $X$ over an algebraically closed field.
In characteristic zero, by Artin's comparison theorem, one can
compute the Betti numbers using singular cohomology.  For this it suffices
to search and find a triangulation.  In positive characteristic, the situation
is somewhat less clear a priori, but the proof in  SGA4 or in \cite{freitag}
of the finiteness of these numbers is everywhere effective.  An explicit 
bound is given in \cite{katz} (thanks to Emmanuel Kowalski for the
reference.): 
 
\begin{fact} [\cite{katz}]  Let $b(n,m)$
be the maximum possible Betti number of a smooth subvariety of $\Pp^n$ of degree $m$.
Then $b(n,m)$ is bounded by a recursive function.
\end{fact}

\noindent{\it Proof of Theorem \ref{4}. } \lb
Let $T$ be the theory of all $K_p$. The completions of $T$ are:
\begin{itemize}
\item The theories $T_p$ of the individual $K_p$; these are just the theory
of algebraically closed fields of characteristic $p$, with an axiom stating
$(\forall x)(\si(x) = x^p)$. Call this axiom $\alpha_p$.

\item The completions of ${T^\infty}_0$, the extension of ACFA stating that
in addition, the field has characteristic zero.
\end{itemize}

Let $\beta_n$ be the disjunction of $\alpha_p$ for $p \leq n$.

Thus $T$ is axiomatized by sentences of the form:
$$\beta_n \ \ \vee \ \   \phi$$
with $\phi$ an axiom of ACFA;
the problem is only to know, given $\phi$, a value of $n$ for which this
disjunction is true. Now the proof of \ref{2qc} estimates an intersection number
$\bar S \cdot \Phi_q$ as $aq^d$, with an error of $\leq b(q^{d-\half})$, with $b$ bounded
effectively and independently of $q$. \ref{.1} gives effectively another constant
$c$, such that $|S \meet \Phi_q | - (\bar S \cdot \Phi_q) \leq cq^{d-1}$. Thus
it suffices to choose $n$ such that $cq^{d-1} + b(q^{d-\half}) < aq^d$
for $q>n$.

This shows that the axioms above are recursively enumerable. On the other hand,
every sentence is equivalent, modulo these axioms, to a bounded
Boolean combination of existential sentences
describing a finite extension of the prime field. The decidability problem for
the theory is thus reduced to the same problem for such sentences, and this is settled
by Cebotarev. \qed

\end{paragraph}

\begin{subsection}{Finite simple groups}  

In the case of {\em algebraic}
simple groups $G$ of dimension $n$, over an algebraically
closed field, the proof is classical (cf. e.g.
Humphreys, Linear Algebraic Groups):  any conjugacy class $X$ of $G$
generates $G$ in at most $\dim(G)+1$ steps.  One considers
$X_n$, the set of $n$-fold products of  elements of $X$.  Then
$\dim(X_n)$ is nondecreasing, and must eventually stabilize:
$\dim(X_{n-1}) = \dim(X_{n})$, $n \leq \dim(G)$.  At this point,
one considers the {\em stabilizer} of the Zariski closure of $X_n$, and
concludes it is a group of dimension equal to $\dim(X_n)$, and containing
a conjugacy class.  Simplicity  of $G$ implies that the stabilizer equals
$G$.

Boris Zilber realized that the proof generalizes to groups of finite Morley
rank; the key to his proof was a different definition of the stabilizer,
using the dimension theory directly without reference to closed sets.

Difference equations of finite total dimension do not have finite Morley rank,
but they do have "finite S1-rank", cf \cite{CH}.
More precisely, such a dimension theory applies to their solution sets
in a universal domain for difference fields; not necessarily to smaller fields.

In the presence of finite S1-rank,
the stabilizer must be defined in a significantly different way than
in finite Morley rank; it can however be defined, and enjoys
similar properties, sufficient to make the proof go through.  See
\cite{HP}.

The ultraproducts of the groups $G_n(q)$
are a priori a group defined in the same way over the ultraproduct $U$ of the
difference fields $GF(p)^{alg},\phi_q)$.  The proof presented in
\cite{HP} would then go through, provided one knows that $U$
has an appropriate dimension theory.  The present paper completes the proof
by showing that $U$ is a universal domain for difference fields.

\end{subsection}

\begin{subsection}{Nonstandard powers and a suggestion of Voloch's}

We first restate and prove \ref{1cor}.
Let $F = GF(p)^{\alg}$ be the algebraic closure of a finite
field, and let $G$  be the automorphism group.  $G$ is isomorphic
to the profinite completion of $Z$, and we consider it with this topology.

\

{\em   Let $F = GF(p)^{\alg}$.  Let
$\Upsilon$ be the set of automorphisms $\si$ of $F$
such that $(F,\si)$ models ACFA.  Then $\Upsilon$ is co-meager in the
sense of Baire.
}

\

\proof   Fix a variety $V$ over $F$.  Then $V$ is defined over
$GF(q)$ for some power $q=p^l$ of $p$.
There are  $l$
possibilities for $\si | GF(q)$; fix one of them, say $\phi_{p^j}$;
let $V' = V^{\phi_{p^j}}$; and then
fix $S \subset V \times V'$ as in the statement of ACFA.  There are
countably many possible $V,j,V',S$ altogether.
It suffices therefore to
show that for each such set of data, the set of $\si$ meeting
the particular instance of ACFA relative to $V,V',S$ is comeager
among all $\si$ with $\si | GF(q) = \phi_{p^j}$.
This set is an open set of automorphisms, so it suffices to show
that it is dense.  For this we must show that for any $l'$ and $j'$
with $l | l'$  and $j' = j (mod l)$, there exists $\si \in \Upsilon$
with $\si | GF(p^{l'} )  = \phi_{p^{j'}}$, and such that
the relevant instance of ACFA holds.  However, by Theorem \ref{1},
any given instance of ACFA holds for any sufficiently large (standard)
power of the Frobenius automorphism; in particular one can pick
a Frobenius power $p^{j''}$ large enough, and such that
$j'' = j' \ (\mod l')$.    Thus the open set in question is dense,
and $\Upsilon $ is comeager. \qed

\noindent{\bf The examples of Cherlin-Jarden}  

Let $l $ be a positive rational.
The equation $E_l: \  \si(x)=lx, x \neq 0$ implies $\si^n(x)=l^nx$.
Thus $E_l$ has no solutions fixed by $\si^n$, hence no algebraic solutions.
So $(\bar{\Qq},\si)$ is not a model of ACFA.  

It can also be shown that
if $a_l \in E_l$, with $l$ varying over the positive rational primes,
then the $a_l$ are algebraically independent over $\Qq$.  This justifies
the assertion made after the statement of \ref{1cor}:  any model of ACFA 
of characteristic $0$ must have infinite transcendence degree over $\Qq$.

{\it However},as a rule, an axiom of ACFA {\em is} true in 
  $(\bar{\Qq},\si)$ for a generic $\si$.  For instance, among  
 order-one difference varieties, this rule admits only 
a few concrete families of exceptions, of the form  
$\si(x)=f(x)$ with $f$ a fractional linear transformation, and $\si(x)=ax+b$,
with $x$ ranging over an elliptic curve.

\

\noindent{\bf Voloch's question on two primes}

\

Let $p$ be a rational prime, and let
$\Omega_p$ be the roots of unity of order prime to $p$.
$L_p = \Qq(\Omega_p)$ be the maximal prime-to-$p$ cyclotomic extension
of $\Qq$.
Let $p_1,p_2$ be two distinct primes on the integers of $L_p$, lying above
the rational prime $p$.
Let $V$ be a variety, say given by
linear equations $L$ over $\Zz$.
Voloch suggested investigating the solutions to $L$
modulo
{\em both} primes, with coordinates on $\Omega_p$:

\[ { (1)} \ \ \ \ \ \ \ \ \ \ \ \ V(v_1,v_2) = \{x=(x_1,\ldots,x_k) \in \Omega_p^k: L(x) \in p_1 \meet p_2 \}  \]

He noted
that $x_2$ is conjugate to $x_1$ by an automorphism $\theta$ of $L_p$,
and that on $\Omega_p$, one may write:
$\theta(x)=x^n$, with $n \in (\Pi_{l \neq p} \Zz_l)^*$.  The question
is therefore equivalent to solving, in $GF(p)^{alg}$, the equations:

\[ (2)  \ \ \ \ \ \ \ \ \ \ \ \  (x_1,\ldots,x_k) \in V, (x_1^n,\ldots,x_k^n) \in V  \]

Under what conditions
is the number of solutions finite?

In a very special case, the present results are relevant.  Assume
$n=p^m+1$, where $m  \in (\Pi_{l \neq p} \Zz_l)$; or more generally, that 
$n=f(p^m)$,
 $f$ is a fixed nonconstant polynomial over $\Zz$.
Choose $m \in  (\Pi_{l \neq p} \Zz_l)$ {\em generically}.
Then   the structure $(GF(p)^{alg}, x \mapsto x^n)$ is interpretable
in $(GF(p)^{alg}, x \mapsto x^{p^m})$.  The latter is a model of ACFA by
\ref{1cor}, hence we have a precise criterion for the answer to (2).
In particular, it is easy to see that that if $V$ is defined by a {\em single}
linear form with at least $3$ variables, then (1) has infinitely many
solutions.

\end{subsection}

\begin{subsection}{Consequences of the trichotomy}

The Lefschetz principle suggests transferring the results of \cite{CH}
to statements about the Frobenius.  The propositions below were
found in this way;  though in retrospect they can be proved using only
a very small part of the results of the present paper. Nonetheless
the general Lefschetz principle gives the propositions a general context
and makes the proofs immediate.

The first and third propositions are direct translations of Theorem 1.10 
in \cite{H-MM},
using Theorem \ref{dec}.   Note that the numerical bounds obtained
there can be construed as first-order statements.

\begin{prop}  Let $f \in \Zz[T]$ be an integral polynomial in one 
variable,
with no cyclotomic factors.  Let $L(X_1,\ldots,X_k)$ be a
linear form over $\Qq$, $k \geq 3$.
Let $X(q)$ be the roots of unity of order
$f(q)$, in $\tilde{\Ff_q}$.  Let
$Y(q)$ be the set of $k$-tuples
$a = (a_1,\ldots,a_k)$
from $X(q)$ satisfying $L(a) = 0$, but
such that no proper sub-sum is zero.
Then there exists an absolute bound $B(f,L)$ and $p_0$ such
that for all primes $p \geq p_0$ and all powers $q$ of $p$,
$|Y(q)| \leq B(f,L)$.  The bounds  $p_0$ and $B(f,L)$ are effectively 
computable
from $f$ and $L$.
\end{prop}

\proof  If $\si$ is an automorphism of a field $K$, let
$X(\si) = \{a \in K:  a^{f(\si)} = 1 \}$, and let $Y(\si)$ be
be defined analogously to $Y(q)$.
By Theorem 1.10 in \cite{H-MM}, the axioms of ACFA together with the
axiom  scheme of fields of characteristic $0$, imply that $X(\si)$
is finite; say $|X(\si)| \leq B$ where
$B = B(f,L)$ is finite.  Hence a finite set of axioms of ACFA imply that
if the field characteristic is larger than some $p_0$, then
$Y(\si)$ is finite.  By Theorem \ref{dec}, the said axioms
hold in all algebraically closed fields of large enough
positive characteristic, when $\si(x) = x^q$, $q$ any power
of the characteristic.   \qed

\begin{rem} \ \end{rem} \begin{enumerate}
\item
An upper bound for $B(f,L)$ can be given explicitly; it is doubly exponential
in $k$ and in the degree  of $f$, with coefficients using   the absolute values
of the coefficients. 
%
This follows from 
 \cite{H-MM}, Proposition
1.11, by plugging in the
parameters.
\item  Using a positive - characteristic version of \cite{H-MM}, Theorem 
1.10,
one can determine $p_0$; the proposition holds for all primes $p$
such that $f(p^{a/b}) \neq 0$ for all rational $a/b$.
\item  One could take $L$ over a number field instead of over $\Qq$, 
still with
a uniform bound, independent of $L$.
\item One can also compute the set of numbers that actually occur
as $|Y(q)|$, or  and describe the set of $q$ for which a particular
number occurs.
\end{enumerate}

By taking $f(T) = T-2$, we obtain the following corollary; in hindsight 
it is easy to find a
direct proof, but we keep it as an example.   For a
prime $p$, let $E_2(p)$ be the set of natural numbers $e$ such
that $2$ is in the multiplicative subgroup of $(Z/eZ)^*$ generated
by $p$.  (This includes the $e$ such that $p$ is a primitive root mod 
$e$.)
Let

$$R_2(p) = \{   x \in GF(p)^{alg}: (\exists e \in E_2(p) ) (x^e = 1 )  
\}$$

\begin{cor} Let $p$ be an odd prime.  Let $L=  \sum c_i X_i$ be a linear 
form
in $k \geq 3$ variables.  Let $Y$ be the set $k$ - tuples
$(a_1,\ldots,a_k) \in R_2(p)^k$ with $\sum c_i X_i =0$ , but no proper 
sub-sum
equals $0$.  Then  $|Y| \leq 2^{k 2^{(k-1)}}$
\end{cor}

\proof  This follows from the proposition (and the remarks), since if $x 
\in R_2(p)$
then $x^{q-2} = 1$ for some power $q$ of $p$.

Similar remarks will apply to the propositions below.
First, a non-linear analog.   For simplicity we consider
one-variable $q$-nomials over $\Zz$, but
the results of \cite{CH}
allow an analysis of the
behavior of any  system of $q$-polynomials, in
several variables, with respect to the question raised below; and
the bounds obtained are independent of parameters, if the polynomial
is over a larger field.

Let $h(X,Y)\in \Zz[X,Y]$.  Let us say that $h$ is {\em special \ } if
there exists a curve $C$, either $C=G_m$ or $C$ an elliptic curve,
a coset $S$ of a subgroup of $C \times C$,
correspondences $U,V \subset (C \times \Pp^1)$  ($U$,$V$ are irreducible
curves, projecting dominantly to $C$ and to $\Pp^1$), and
points $(a,b,c,d)$ such that $(a,b)$,$(a,c)$,$(b,d)$, are
generic points of $S$, $U$,$V$ respectively, and $h(c,d) = 0$.

(If $h$ is special, then the roots of the $q$-nomial $h(X^q,X)$
are in uniform algebraic correspondence with either roots of unity as
in the previous proposition, or an elliptic analog, or the points of 
$\Ff_q$.)

Let us say that an affine variety $U \subset \Aa^n$ is {\em special \ }
if it is explicitly a product of curves and points; i.e. there exists a 
partition
of the coordinates of $\Aa^n$, corresponding to an isomorphism
$j: \Aa^n \to \Aa^{k_1} \times \ldots \Aa^{k_l}$, such that
$j(U) = C_1 \times \ldots \times C_l$, with $C_i$ a point or a curve.

\begin{thm}\lbl{tricho1}  Let $h_i(X,Y)\in \Zz[X,Y]$ be non-special
polynomials, 
and let $U \subset \Aa^n$ be an affine variety over $\Zz$.  Then
either $h$ is special, or there exists a finite union $W$ of special
subvarieties of $U$ (with $W$ defined over $\Qq$) such that for all
sufficiently large primes $p$, and all powers $q$ of $p$,
if $h_i(a_i,(a_i)^q)=0$ and  $(a_1,\ldots,a_n) \in U$, then $(a_1,\ldots,a_n) \in W$.
\end{thm}

 \proof  Assume $h$ is not special.  By the trichotomy theorem of 
\cite{CH}, in any model $\Qq \leq M \models ACVF$, 
$$X = \{x:  h(x,\si(x)) \}$$
is a stable, stably embedded definable set, whose induced structure
is superstable of rank one, and {\em distintegrated}:  the algebraic
closure relation on $X(M) \setminus \Qq^{a}$ is an equivalence relation.
 
Let $U \subset \Aa^n$ be an affine variety. If $(a_1,\ldots,a_n) \in U$
and $a_i \in X$, let $w_o(a) = \{i \leq n: a_i \in \Qq^a \}$, 
$w_+(a)= \{i \leq n:  a_i \notin \Qq^a \}$ and let
$\{w_k(a) : k=1,\ldots,m \}$ be the classes of  the equivalence relation on
$w_+(a)$ defined by:  $i \sim j$ iff $\Qq(a_i)^a = \Qq(a_j)^a$.  
Let   $a(k) = (a_i: i \in w(k))$, so that we can write
$a = (a(0),a(1),\ldots,a(m))$.  Let $C(k)$ be the locus over $\Qq^a$
of $a(k)$; so $C(0)$ is finite, and $C(k)$ is a  curve for $k \geq 1$.
Let $S(a) =\{a(0)\} \times C(1) \times \ldots \times C(m) \}$; it is
a special variety, defined over $\Qq^a$. 
By disintegration of $X$, the fields $\Qq(a(0)),\ldots,\Qq(a(m))$
are linearly free.  Thus $a$ is a generic point of
$S(a)$.
Since $a \in U$, we have:  
$ S_a   \subset U $.  
Now $a$ was an arbitrary point of $U(M) \meet X(M)^n$, in a model
$M$, and so by compactness, there are finitely many special $S_i$,
defined over $\Qq^a$, such that $U(M) \subset   \union_i S_i$.
Further taking the union of all $\Qq $-conjugates, we find a $\Qq$-
definable finite union $S$ of special varieties, with 
$U \meet X^n \subset S$ (in any difference field extension of $\Qq$.)
The statement on the large primes follows immediately by compactness.
\qed

Similar results exist for $q$-nomials of higher order, and for
commutative  algebraic groups.  In particular, for 
 semi-Abelian varieties, we have: 

\begin{prop} Let $A$ be a commutative algebraic group scheme over
a scheme $Y$ of of finite type over $\Zz$, with
generic fiber $A_Y$ a semi-Abelian
variety.  Let $f$ be an integral polynomial in one variable,
with no cyclotomic factors.
Let $V$ be a subscheme of $A$.  For any closed point $y \in Y$,
let $A_y$
be the fiber of $A$ over $y$.
 Let $\phi$ denote the Frobenius endomorphism of $A_y$, relative
 to the residue field $k_y$, and let $\psi$ be
 any finite power of $\phi$.  Let
$X(y,\psi,f)$ be the kernel of the endomorphism $f(\psi)$ on 
$A(\tilde{k_y})$.

Then there
exists a Zariski  open $Y_0 \subset Y$ and finitely many group
subschemes $B_i$ of $A$ over $Y_0$, with $B_i \subset V$,
and an effectively computable integer $B$,   such that for
any closed point $y \in Y$,
$ X(y,\psi,f) \meet V$ is contained in at most $B$ translates of
some of the $B_i$.  

\end{prop}

The proof is analogous to that of \ref{tricho1}.  \qed

\end{subsection}

\>{section}

\<{section}{Jacobi's bound for difference equations}
\def\jj{{\bf j}}
\def\JJ{{\bf J}}
\label{jacobiS}
Consider a system $u$ of $n$ difference equations $u_1,\ldots,u_n$  in $n$ variables
$x_1,\ldots,x_n$; assume $u_i$ has order $h_k^i$ with respect to the variable $x_i$.  
Assume that the system defines a difference scheme of finite total dimension $d$.
What bound can one impose on $d$, based on the data $(h_k^i)$?

Jacobi considered this problem for differential equations, and gave as a bound:

$$ \jj = \max _{\theta \in Sym(n)} \Sum_{k=1}^n h_k^{\theta(k)}$$

Jacobi's proof was criticized by Ritt \cite{ritt}, and the problem is considered open.
See \cite{ritt}, \cite{KMP}, \cite{cohn-gr} for partial results.  

The transposition to difference equations was made by Cohn, with analogous partial
results; cf. \cite{lando}, \cite{KMP}.     Jacobi's original idea was reduction to 
the linear case (the "equations of variation"), and then the the case of constant
coefficients; and present results still follow this line.  However when the equations have
hidden singularities, the reduction to the linear case requires additional assumptions,
or new ideas.

We show here that Jacobi's bound is correct for difference equations by a quite
different method, Frobenius reduction.  We also get an ``explanation" for the curious
determinant-like formula.

Cohn  refined the problem by setting $h_k^i =  - \infty$ if $x_i$ does not occur in $u_k$;
we accept this refinement.  (It not only gives a better bound in many cases, but
 allows for greater flexibility in manipulating the equations.)

Consider first the analogous problem for {\em algebraic} equations.  Let $U_1,\ldots,U_n$
be polynomials in $n$ variables $x_1,\ldots,x_n$,  over a field $k$; suppose
$U_i$ has degree $H_i^j$ in the variable $x_i$.  Let $Z(U)$ be the scheme
cut out by $\{U_i\}$, and let $Z_0$ be the 0-dimensional part of $Z(U)$.

\<{lem} \lbl{jacobi-alg}  $|Z_0| \leq \JJ = \sum_{\theta \in Sym(n) } \prod _{k=1}^n H_k^{\theta(k)} $ \>{lem}

\proof   Let $ V= (\Pp^1)^n$ be the $n$'th Cartesian powers of the  projective line.  
The tangent bundle of $V$ is generated by its global sections, since this is true
for each factor $\Pp^1$.  The normal bundle to $V$ embedded diagonally 
in $V^n$ is hence also generated by its global sections.  Thus by \cite{fulton}
Theorem 12.2, each contribution to the intersection of $V$ with a subvariety of
$V^n$ of codimension $n$ is represented by a non-negative $0$-cycle.  The number of components of intersection, in particular the number isolated points,   is therefore 
bounded by the intersection number (just as in \cite{fulton}, Example 12.3.)
 
 Each equation $U_i$ determines a hypersurface $[U_i]$ in $V$ (the closure
 in $V$ of the subscheme of $(\Aa^1)^n$ cut out by $U_i$.)   $Z_0$
 is the 0-dimensional part of 
   $ [U_1] \meet \cdots \meet [U_n] \simeq ([U_1] \times \cdots \times [U_n] ) \meet V$.  By the above discussion, $|Z_0| \leq [U_1] \cdot \ldots, \cdot [U_n]$.  It remains
 to compute this intersection number.  
 
 Now $[U_i]$ is rationally equivalent to $\sum_{j=1}^n h_i^j D_j$, where 
 $D_j$ is the divisor defined by the equation $x_j=0$, the pullback of 
 a point on the $j$'th copy of $\Pp^1$.    We have $D_j^2=0$, since
 $\{pt\}^2=0$ in $\Pp^1$.   Thus $D_\theta(1) \cdot \ldots \cdot D_\theta(n) = 0$
 if $\theta$ is not injective, while $D_\theta(1) \cdot \ldots \cdot D_\theta(n) =
  D_1 \cdot \ldots \cdot D_n =
  \{pt\}$
 if $\theta$ is a permutation.  
 Thus
 $$[U_1] \cdot \ldots, \cdot [U_n] = \prod_{i=1}^n \sum_{j=1}^n h_i^j D_j   = \JJ $$ 

\qed

\<{thm}  Let $k$ be a difference field $u_i \in k[x_1,\ldots,x_n]_\si$ a differential
polynomial of order $h_i^j$ in $x_j$.  Let $W$ be a component of 
the difference scheme $Y$ cut out by $h_1=\cdots =h_n =0$.  Then the   total
dimension of $W$, if finite, is no larger than  Jacobi's bound 
   $ \jj$.  \>{thm}

\proof  By assumption, $W$ is transformally reduced, of   total dimension $w$.  Let $v$ be a differential polynomial that does not vanish on $W$, but vanishes on every other component of $Y$.
Adding a variable $y$ and the equation $vy=1$ does not change the Jacobi bound
(using Cohn's convention, the new equation has order $0$ in the variable $y$,
and all other equations have order $- \infty$ in this variable.)  Thus we may localize
away from the other components; so we may assume that $Y$ itself has finite reduced total dimension $w$.  We must show that $w \leq \jj$.

 Let $D \subset k$ be a finitely generated difference domain, with
$u_i \in D[x_1,\ldots,x_n]_\si$.  Given a homomorphisms $\phi: D \to K_q$, let
$U_i$ be the polynomial obtained by replacing $\si(x)$ by $x^q$ and $d$ by
$\phi(d)$ in $u_i$.   

Let $H_i^j = \deg_{x_j} H_i$.   Writing $u_j$ as a sum of $\si$-monomials,
 they each involve  $u_j^{\nu}$ for some $\nu \in \Nn[\si]$, and we have
 $H_i^j = \nu_j(q)$ for the highest such $\nu$. This $\nu$ has the form 
 $d \si^{h_i^j}+ \cdots $ (lower terms), so $H_i^j = dq^{h_i^j} + O(q^{h_i^j-1})$.   Thus 
 $\lim_{q \to \infty} \log_q H_i^j = h_i^j$.  

By \ref{si-dim},
 for some $D$, for almost all $q$ and $\phi$,  $Y_q = \{x: U_1(x) = \cdots = U_n(x) = 0 \}$
 is {\em finite}; and by Theorem \ref{1}C,   for some $a,b \in \Qq$ with $a > 0$, 
 $|Y_q| = aq^w + e_q, |e_q| \leq bq^{w-1/2}$.    On the other hand , by \ref{jacobi-alg},
 $|Y_q| \leq \JJ$.   So
 $$aq^w + e_q \leq  \sum_{\theta \in Sym(n) } \prod _{k=1}^n H_k^{\theta(k)} $$
  
  Taking $\log_q$ and letting $q \to \infty$ we obtain
 $$ w \leq \max _{\theta \in Sym(n) } \sum _{k=1}^n h_k^{\theta(k)} $$  
 \qed

\<{remark}  {\bf Ritt's conjecture for difference equations}.  \rm
   
As Cohn observed (\cite{cohn-od} for differential equations, personal communication
for difference equations), the validity of the refined Jacobi bound implies the Ritt
dimension conjecture:  given a system of $m$ difference equations in $n$ 
variables, if $m<n$ then every component of the solution set has transformal
dimension $\geq n-m$.  This can also be deduced directly from Theorem \ref{1},
as in \ref{illustrate-dirdim}.  \end{remark}

\>{section}

\end{part}

\begin{section}{Complements}
\lbl{complements}

\<{subsection}{Transformal dimension and degree of directly presented difference schemes}
\lbl{transdimdeg}

Here is a variant of \ref{direct2}, valid for all components. 
 The proof  {\em  assumes} Theorem \ref{dec}, and illustrates the way the
 main result of this paper may be used in difference algebra.   

 \begin{prop}\lbl{illustrate-dirdim}    Assume theorem \ref{dec}.
 Let $X$ be a smooth algebraic variety over a difference field $K$.
Let $S \subset X \times X^\si$ be an absolutely irreducible subvariety,
$\dim(X) = d$, and assume $\dim  S = e+d$.
Let
\[ Z = \{x \in X: (x, \si(x) ) \in S  \}   \]
Then any
component of $Z$
 has transformal dimension $\geq e$.
\end{prop}

         \proof  Let $Z(1),\ldots, Z(m)$ be the components.
Suppose for contradiction that $\hbox{ trans. dim. } (Z(1)) < e$.  We may assume $K$
is the field of fractions of a finitely generated difference domain $D$.
By the main theorem, for almost all $q$, and almost $y \in \spec M_q(D)$,
$M_q(Z)_y = \union _j M_q(Z(j)_y)$,
$M_q(Z(1))_y \not \subset \union _{j>1} M_q(Z(j))_y$,
and $\dim(M_q(Z(j))_y) \leq \hbox{ trans. dim. } (Z(1)) < e$.  It follows that $M_q(Z)_y$
has a component of dimension $<e$.  However
$M_q(Z)_y = S_y \star \Phi_q$.  Since (for almost all $y$)
$S_y \subset X_y \times (X^\si)_y$ is irreducible of dimension $d+e$, and
$X_y \times (X^\si)_y$ is smooth, the dimension theorem implies that
every component of this intersection has dimension $\geq (d+e)+d - 2d = e$;
a contradiction.  \qed
\medskip

\<{remark}\lbl{ww+}   Moreover, in \ref{illustrate-dirdim}, if $h: Z \to W$ is a morphism
of difference schemes, and $W $ has transformal dimension $0$, then
each component of each fiber of $h$ has  transformal dimension $\geq e$.
\>{remark}

The proof is similar, using the fact that $M_q(W)_y$ is finite, and thus the 
components of the fibers of $M_q(h)_y$ are also components of $M_q(Z)_y$.

\

Note that \ref{illustrate-dirdim} is equivalent to to following purely geometric statement:

 \<{lem}\lbl{ww}  Let $U,V,S$ be quasi-projective varieties over an algebraically closed difference field $k$, with $V \subset U^\si$, and $S \subset (U \times V)$.  
  Let ${\cal T} $ be the set of absolutely irreducible varieties $T \subset S$
  such that $(pr_1 T)^\si = pr_2 T$.   (Where $pr_i T$ is the Zariski closure
 of the $i$'th projection.)   Let $e \geq 1$.
 Assume $\dim(S) \geq \dim(U)+e$.  Then   any maximal $T \in {\cal T}$ 
 satisfies  $\dim(T) \geq \dim(pr_1 T)+e$.  \>{lem}

 \<{problem} Is there a simple direct proof of \ref{ww}?  Is the smoothness assumption necessary, even if one
just wants one component of transformal dimension $e$?\>{problem}

The difficulty here is associated with singularities.  If 
 all varieties encountered were smooth, it would be easy to prove \ref{ww}:
Let $T_0  \in {\cal T}$.  
Let  $U'$ be a component of $pr_1[1](S)$ (cf. \ref{zd})
containing $pr_1(T_0)$.    Let
 $S' = (U' \times U) \meet S$; then $\dim(S') \geq \dim(S) - (\dim(U)-\dim(U'))
\geq \dim(U') +e$.   Thus one can proceed by induction.

In attempting to find a direct proof of \ref{ww}, I considered projecting
at the right moment to $\Pp^m$ by a  finite map, and using  the dimension theorem there.  To do this, one needs to know Lemma \ref{pos-stat} below; and in order
to reduce to the case of purely positive transformal dimension, one finds oneself requiring
a stronger statement than  \ref{ww}, namely   \ref{ww+} above.
I did not carry through this proof to the end, but the lemmas that came
up seem   sufficiently suggestive in themselves to be stated here.

\<{definition} \lbl{pure-pos} An irreducible difference variety $U$ over $K$
is {\em purely positive transformal dimensional} if 
every differential rational morphism
on $X$ into a difference variety of transformal dimension $0$ is 
constant. \>{definition}

\<{prop}\lbl{pos-stat} Let $K$  be a difference field, $U,V$ algebraic varieties over $K$,
$f: U \to V$ a quasi-finite map.  Let $X$ be an irreducible difference subvariety of $V$, $X(K)$
Zariski dense in $V$.   Assume $X$ is purely positive transformal dimensional.
Then for some Zariski open $\tilde{U} \subset U$, $f^{-1}(X) \meet U$
 is an irreducible difference variety (and also  purely positive transformal dimensional.)
\>{prop}
\proof \cite{CH}.

\<{cor}\lbl{www}  Let $U, S$ be affine varieties over an algebraically closed difference field $k$,  with $S \subset (U \times U^\si)$.  Let $X = S \star \Sigma \subset [\si]_k U$.
Then either $X$ has positive transformal dimension, or else it has total dimension
$\leq \dim(U)$. 

 In the latter case, if   $U$ is absolutely irreducible, every algebraic component of $X$ of total dimension
equal to $\dim(U)$ is Zariski dense in $U$.   \>{cor}

\proof  Let $(a_0,a_1,\ldots)$ be a generic point of an algebraic component $C$
of $U$.  (I.e $(a_0,a_1,\ldots,a_n) \to (a_1,\ldots,a_{n+1})$ is a specialization
over $k$, but not necessarily an isomorphism; and $(a_0,a_1) \in S$.)

Let $d_n = tr. deg._k k(a_0,\ldots,a_n)$.  If, for some $n$, $d_n < d_{n+1}$,
let $U_n,V_n,S_n$ be the $k$-loci of $(a_0,\ldots,a_n),(a_1,\ldots,a_{n+1}),(a_0,\ldots,a_{n+1})$
respectively.  Then the hypothesis of \ref{ww} holds.  So there exists a point $c$
in a difference field extension of $k$, with $c_n = \si^n(c)$ satisfying:
$(c_0,\ldots,c_{n+1}) \in S_n$, and of positive transformal dimension over $k$.  In particular,
$(c_0,c_1) \in S$, so $c \in X$ and $X$ has positive transformal dimension.

Otherwise, $d_{n+1} \leq d_n$ for each $n$.  So $d_n \leq d_0 \leq \dim(U)$ for each $n$.  This
shows that $C$ has total dimension $\leq \dim(U)$.  Moreover if  $c_0$ is not a generic point
of $U$, then $C$ has total dimension $< \dim(U)$.  \qed
\end{subsection}

\begin{subsection}{Projective difference  schemes}
\lbl{projective}

\begin{paragraph}{ $\Nsi$- Graded rings}

Let
$\Nsi$ denote the set of polynomials  over $\Nn$, with indeterminate
labeled $\si$.
Consider difference rings $R$ graded by the ring $\Nsi$.
In other words, we are provided with a decomposition of $R$ as an Abelian
group:
 $$R= \oplus_{n \in \Nsi} H_n$$
Multiplication in $R$ induces maps $H_n \times H_m \to H_{m+n}$.
The action of $\si$ is assumed to carry $H_n$ to $H_{n \si}$. Such a
structure will be called a $\Nsi$-graded ring.

We will assume $R$ is generated as a difference ring by $H_0 \union H_1$.
A difference  ideal $P$ is called {\it homogeneous} if
it is generated by the union of the sets $P \meet H_n$.

\end{paragraph}  

\
\begin{paragraph}{Projective difference  schemes}
 Let $R$ be a $\Nsi$-graded ring, with homogeneous
components $H_n$. We define an associated difference scheme, $\projs  (R)$, as
follows.
The underlying space is the set of homogeneous transformally prime ideals, not containing $H_1$.
The topology is generated by open sets of the form:
$W_a = \{p: a \notin p\}$, with $a \in H_k$ for some $k \in \Nsi$.
Such sets are called affine open.    One assigns to $W_a$ the difference ring:
$$R_a = \{ R[a^{-n} : n \in \Nsi] \}_0 =
 \{\frac{b}{a^{n}}: n \in \Nsi, b \in H_{nk} \} $$
and glues.

Explanation: the localized difference  ring $R[a^{-n}: n \in \Nsi]$
has a natural grading, with an element of the form
 $b/a^n$ ($b \in R, n \in \Nsi$) in the homogeneous component
of degree $\deg(b) - n\deg(a)$.   $\{ R[a^{-1}]\}_0$ is the
degree-zero homogeneous component of this $\Nsi$- graded ring. To verify that this
yields, uniquely, a difference scheme structure, observe
that $\{ R[ (a_1a_2)^{-n}: n \in \Nsi] \}_0$, the affine ring corresponding
to $W_{a_1a_2}$, is the difference ring localization of the ring corresponding
to $W_{a_1}$ by the element $\frac { a_2^{ \deg{a_1} } } { {a_1}^{\deg{a_2}}   } $.
Note also that the prime ideals of $ R_a$
can be identified with the elements of $W_a$ (a homogeneous prime ideal of $R$
is the kernel of a homomorphism $h: R \to L$, $L$ a difference field,
with $h(a) \neq 0$. Then $h$ restricts and extends to $R_a$. Conversely,
let $g: R_a \to L'$ be a homomorphism. Then $g$ extends to a
graded homomorphism
$\bar{g}: R \to L'[t^{\Nsi}]$, with $\bar{g} (a) = t$.)

We could alternatively directly describe the structure sheaf in terms of section,
as in the definition of $\spe$; the local ring at $p$ is defined as
$$\{\frac{b}{a}: ( \exists k) \ a,b \in H_k, \ b \notin p \} $$

\end{paragraph}

\begin{paragraph}{ $\Nsi$-graded ring associated to a graded ring}

Let $D$ be a difference
domain, and let $R$ be a graded $D$-algebra in the usual sense,
$R = \oplus_{i \in \Nn} R_i$.   Assume the homogeneous component of degree 0 is a domain.

\begin{lem}$[\si]_D R$ has a unique $\Nsi$-grading, compatible with the
grading of $R$.
$\{[\si]_D R\}_0 = [\si]_D({R}_0)$.
\end{lem}

\proof For $n \in \Nsi$,
 $n = \sum_i {k_i}\si^i$, let $H_n$ be the subgroup of $[\si]_D R$
generated by the products $\Pi_i \si^i(a_i)$ with $a_i \in R_{k_i}$.
(We write here $a_i$ also for the image of $a_i$ in $[\si]_D R$.)
Clearly the $H_n$ generate $[\si]_D R$ between them, and any $\Nsi$-grading
must make $H_n$ the homogeneous component of degree $n$. Thus it is only
a question of showing that the $H_n$ are in direct sum. Let $F$ be a free
$R_0$-algebra, graded so that the free generators span $F_1$ as a $R_0$-module,
and let $h: F \to R$ be a surjective homomorphism of graded
$R_0$-algebras. Let ${\bar F} = [\si] F$;
it is clear that the lemma holds for $\bar F$. Let $J$ be the $\Nsi$-ideal
generated by $ker(h)$. Also $J$ is generated by $J \meet {\bar F}_1$,
so $J$ is homogeneous, and $J \meet {\bar F}_0 = (0)$.
The $\Zz$-graded part of $(\bar F)/J$ is $F/ker(h) \cong R$.
Thus we obtain a map $R \to (\bar F)/J$. This map extends to a difference ring
map $[\si] _D R \to (\bar F)/J$. We also have a map $\bar F \to [\si]_D R$.
By the universal properties, these are isomorphisms,
and $[\si] _D R$ is $\Nsi$-graded.

\end{paragraph} 

\begin{paragraph}{ Difference structure on projective schemes}
Let $D$ be a difference domain.

\begin{lem}\lbl{5.13} Let $R$ be a graded $D$-algebra.
Then $\projs [\si]_D R$ is the difference scheme obtained by gluing together
$\spe [\si]_D R_a$, where $R_a$ runs over the various localized rings
$R_a = \{R[a^{-1}]\}_0$, $a \in R_1$. \end{lem}

\proof First, for any $a \in R_n, n \in \Nn$,
$[\si]_D (R[a^{-1}]) = ([\si]_D R)[a^{-1}]$,
the second localization being taken in the sense of difference rings. This
is clear, since both are the universal answers to the following problem:
a difference ring $\bar R$, a homomorphism $h: R \to \bar R$, with $h(a)$
invertible in $\bar R$. Next, one verifies that when $R$ is graded,
the two induced gradings on these rings are the same; in particular,
$[\si]_D (R_a) = ([\si]_D R)_a$, using the notation of the definition
of $\projs$, and the previous lemma. \qed

\

\begin{lem} There are canonical isomorphisms: \begin{itemize}
\item $M_q([\si]_D R) = R \tensor_D M_q(D)$
\item $M_q(\spe [\si]_D R) = \spec(R) \times _{\spec D} M_q(\spec D)$
\item When $R$ is a graded $D$-algebra,
$M_q( \projs ([\si]_D R) ) = \proj (R \tensor_{D} M_q(D))$
\end{itemize} \end{lem}

\proof \begin{itemize}
\item By the universal properties.
\item Apply $\spec$ to the previous.
\item By \ref{5.13} and the previous item.
\end{itemize}

\end{paragraph} 

\begin{paragraph}{Multi-projective varieties}

Suppose $R$ is graded by $\Nsi^k$, in place of $\Nsi$.
We can define $\proj ^{\Nsi^k} R$ analogously to the case $k=1$.
Sometimes if the intended structure of $R$ is clear we will just
write $\projs  R$ or $\proj R$.
The underlying space is the set of homogeneous difference ideals,
not containing a homogeneous component $R_j$ (for any   $j \in \Nsi^k$). The basic affine open sets are the sets of primes
not containing $b=a_1a_2 \ldots a_k$, where $a_i$ has
degree $(0,\ldots,0,1,0,\ldots, 0)$.
 The localization $R[b^{-1} ,\si(b^{-1} ),\ldots]$
has a   graded structure,
 and the associated  ring is defined as
the graded component of degree $0$.

\
\begin{note}\lbl{htd} \ \end{note}

   Let $p \in \proj R$. Then
there exists $a_i $ of degree $(0,\ldots,0,1,0,\ldots, 0)$
with $a_i \notin p$. Since $p$ is homogeneous, and prime, the
$a_i$ are algebraically independent over $R_0/(p \meet R_0)$. Thus
if $R$ is a domain, and a $D$-algebra,
the transformal  transcendence degree of
$R$ over $D$ is $k$ more than the transformal dimension of $\proj \ R$ over $D$.

\

The ideals $J_q$ used in defining $M_q$
are not homogeneous. However, a ring graded by
$\Nsi^k$ can be (forgetfully) viewed as graded by
$\Nn^k$, in many ways; any ring
homomorphism from $\Nsi$ to $\Nn$ will give such a way.
Given $q$, let $h_q: \Nsi \to \Nn$ be the homomorphism
with $\si \mapsto q$, and also let $h_q:  \Nsi^k \to \Nn^k$ denote the product
homomorphism. Use $h_q$ to view
$R$ as $\Nn^k$ graded: the graded component
of degree $n$ is by definition the sum of the graded components
of degree $\nu$, over all $\nu$ with $h_q(\nu)=n$.
Then $J_q$ is homogeneous for this grading; the generators
can be taken to be $\si(r)-r^q$ with $r$ in some homogeneous
component, and these elements are homogeneous in the
$h_q$-induced grading. Therefore $M_q(R) = R/J_q(R)$ has
a natural $\Nn^k$-graded structure. We can thus view $M_q$ as a
functor on $\Nsi^k$-graded difference  rings into $\Nn^k$-graded rings.

\

\begin{notation}
The difference polynomial ring in one variable over a difference ring
$R$ is denoted
$R[t^{\Nsi}]$.   The localization of this ring by $t$ is denoted
$R[t^{\Zsi}]$, and is called the $\Zsi$-polynomial ring over $R$.
The $k$-variable version $R[t_1^{\Zsi},\ldots,t_k^{\Zsi}]$ is called
the  $\Zsi^k$-polynomial ring over $R$.
\end{notation}

Note:
\begin{lem}\lbl{22}   Let $R = \oplus \{R_a: a \in {\Nsi}^k \} $
be an $\Nsi^k$-graded ring. Let $a_i \in R$ have degree $(0,\ldots,0,1_{(i)},
0,\ldots,0)$. Let $S = R[{a_i}^{\si^n}: i = 1,\ldots,k, \ n \in \Zz ] $
be the localization by $a_1,\ldots,a_k$, and let  $S_0$ be the degree-$0$
component.  Then
$S$ is the  $\Zsi^k$- polynomial ring
 over $S_0$:
$$S = S_0[{a_1}^{\Zsi},\ldots,{a_k}^{\Zsi}]$$ \end{lem}

\proof An element $c$ of degree $(n_1,\ldots,n_k)$ of this ring can
be written as $ b{a_1}^{n_1} \ldots {a_k}^{n_k}$, where
$b = c{a_1}^{-n_1} \ldots {a_k}^{-n_k} \in S_0$. By homogeneity, any difference -
algebraic relation among the $a_i$ over $S_0$ implies a monomial relation
among them.  But $ba_1^{n_1} \ldots {a_k}^{n_k} = 0 $ implies $b=0 \in S_0$,
 since the $a_i$ are invertible.

\begin{lem}\lbl{Mq-proj} Let $R$ be a $\Nsi^k$-graded ring. Then
$M_q(\projs (R)) \cong \proj (M_q(R))$ \end{lem}

\proof   To simplify notation we treat the case $k=1$.
Let $a \in R_1$. By   \ref{22},
$ R[a^{-1},a^{-\si},\ldots] = R_a [ a^{\Zz[\si]} ]$
is a $\Zz$-polynomial ring over $R_a$.  Thus
$M_q(R_a) = \{ M_q(R_a [ a^{\Zz[\si]} ] )\}_0$.

Write ${\bar a}$   for the image of $a$
 modulo $J_q$. Then
$$ M_q(R[a^{-1},a^{-\si},\ldots]) = M_q(R_a)[{\bar a}^{-1 }]$$
Taking degree-$0$ components,
$$M_q(R_a) = \{M_q(R[a^{-1},a^{-\si},\ldots])\}_0 = \{M_q(R)[{\bar a}^{-1}]\}_0
= (M_q(R))_{\bar a} $$
Thus
$$M_q(\spe (R_a) ) = \spec (M_q(R)_{\bar a})$$
Now the difference schemes $\spe (R_a)$ glue together to give
$\proj ^{\Nsi^k} R$, while the schemes
$ \spec (M_q(R)_a)$ glue together to give $\proj (M_q(R))$. (As
$a $ runs through $R_1$, ${\bar a}$ runs through a generating set
for ${M_q(R)}_1$.) The lemma follows after verifying that the gluing
maps agree.

\begin{remark} \ \lbl{biprojring} \end{remark}If $D$ is a difference domain,
$ D[X,Y] = D[{X_1}^{\Nsi},\ldots,X_n^{\Nsi},{Y_1}^{\Nsi},\ldots,{Y_m}^{\Nsi}]$,
bi-graded so that $X_1^{a_1} \ldots X_n^{a_n}Y_1^{b_1},\ldots,Y_m^{b_m}$ is
homogeneous of degree $(\sum a_i, \sum b_i)$, then $\proj D[X,Y]$ is the product
over $D$ of the $\si$- projective spaces of dimensions $m$,$n$.
Conversely, $\proj$ of any bi- or multi-graded ring can be viewed as a fiber product of simple
$\projs $ of some component
rings.
\end{paragraph}  

\begin{paragraph}{Morphisms between graded rings}

Let $R$ be a $\Nsi^k$-graded ring, and $S$ a $\Nsi^l$-graded
ring.   Let
$h: R \to S$ be a surjective $\si$ - ring homomorphism. Assume
$h(R_c) \subset S_{\lambda(c)}$, where $\lambda: \Nsi^k \to \Nsi^l$
is an injective $\Nsi$-linear map. Then one obtains a map
$h^*: \proj  S \to \proj  R$, as follows. If $q$ is a homogeneous
transformally prime ideal on $S$, not containing a homogeneous
component, then so is $h^{-1}(q)$ on $R$.
Let $h^*(q) = h^{-1}(q)$. $h^*$ is continuous: the open set $W^R_a$
defined by an element $a \in R$ pulls back to the set
$W^S_{h(a)}$. Finally $h$ induces a difference ring homomorphism
on the graded rings,
 $R[a^{-1},a^{-\si},\ldots] \to S[h(a)^{-1},h(a)^{-\si},\ldots]$,
respecting the grading in the
same sense as $h$; and in particular induces a map $R_a \to S_{h(a)}$.

\end{paragraph}

\end{subsection} 

\begin{subsection}{Blowing up}   \lbl{blowup}

We give two constructions of blowing up.   We show 
that the transformal blowing-up reduces under Frobenius to a scheme containing the usual blowing-up,
and of the same dimension, without resolving the interesting questions regarding their exact relation. One construction 
can be viewed as the 
result of a deformation along the affine $t$ line (of transformal dimension $1$) while the other is a limit
of a deformation along $t=\si(t)$ (total dimension $1$.)  The latter
seems to have no analog in the Frobenius picture.  Nevertheless
they are shown to give (almost) the same result.

\begin{paragraph}{Blowing up difference ideals }

Let $R$ be a difference ring, $X = \spe(R)$,
 and $J$ a finitely generated difference ideal.
Consider the ordinary polynomial ring $R[t]$ over $R$,
and set $\si(t) = t$.
Let ${R[Jt]} $ be the subring of $R[t]$:
$${R[Jt]} = \sum_{n \in \Nn} (Jt)^n$$

Observe that $R[Jt]$ is finitely generated as a difference ring, if $R$ is.
Indeed if $V$ is a set of generators for $J$ as a difference ideal, and
$Y$ a set of generators for $R$ as a difference ring, then $Y \union
Vt$ is a set of generators for $R[Jt]$.

Moreover, the degree $0$ difference subring of the difference ring
localization $R[Jt][(at)^{-1}]$ is finitely generated as a difference  ring.  If
$a \in V$, it is
generated by $Y \union \frac{Vt}{at} \union \left\{ {\frac{\si(a)t}{at}},
{\frac{at}{\si(a)t}} \right\}$.

View ${R[Jt]}$ as a (non-Noetherian) graded ring with an endomorphism;
form the ordinary scheme-theoretic $\proj$, \cite{HA} II 2;
 and consider the difference 
subscheme described in \ref{fixed}  above.
 Let
$$  \widetilde{X}_J = \fixs  \proj ({R[Jt]})$$

A basic open affine of $\proj ({R[Jt]})$ has the form
$W_a = \spec ( R[\frac{J}{a}])$ , with $a \in J$.  The intersection of
$W_a$ with the set of difference ideals is contained in
$\spec ( R[\frac{J}{\si^n(a)}] )$ for each $n$.  Thus one sees that
$\widetilde{X}_J$ has an open covering by open subschemes of the form
$\spe R[\frac{J}{a}]_\si$, $a \in J$, where the square brackets
$[ \ ]_\si$ here refer to
localization of difference rings.

There is a natural map $\pi: \widetilde{X}_{J} \to X$, considered as part
of the structure.

{\em Compatibility  with localization}  There are two (closely related)
points here.

First, suppose $\pi: \widetilde{X}_{J} \to X$ is the blowing up at $J$,
$X' = \spe R'$ is an open subscheme of $X$,
$R' = R[a^{-1},a^{-\si},\ldots]$,
with inclusion map $i: X' \to X$,
and $J' = R'J$.   Then $\widetilde{X'}_{J'} = \pi^{-1} ( X')$.  This is immediately
verified.

Secondly, suppose $J'$ is another finitely generated difference ideal, $J \subset J'$, and
for every  $p \in \spe R$, if $R_p$ is the local ring, then $JR_p = J'R_p$.
Then  $\widetilde{X}_J =  \widetilde{X}_{J'}$.  Indeed,  every prime ideal containing
$J$ must contain $J'$, so $\widetilde{X}_{J'}$
is covered by the open affines $\spe R[\frac{J'}{a}]$ with $a \in J$.
Next, one may find an open covering of $\spe R$, such that the restrictions
of $J$ and $J'$ to each open set agree.  It follows that
$\spe R[\frac{J'}{a}] =  \spe R[\frac{J}{a}]$ since they are equal locally (on $R$).


\

Now define $\widetilde{X}_{\cal J}$ and
$$\pi: \widetilde{X}_{\cal J} \to X$$
for a general difference scheme $X$
and quasi-coherent $\si$ - ideal presheaf $\cal J$ on ${X}$,
by gluing.

If $i:W \to X$ is a closed subscheme of a scheme $X$,
with corresponding difference ideal sheaf
${\cal J} = ker({\cal O}_X \to i_* {\cal O}_W)$,
 we write
 $$\widetilde{X}_W = \widetilde{X}_{\cal J}$$
We define the {\it exceptional
divisor} $E_W$ to be the difference scheme inverse image of $W$ under
this map.

\end{paragraph}

\begin{paragraph}{A second construction: the closed blowing up
of ideals}

In order to obtain an embedding of the blowing-up in $\si$-projective
space, we use a different construction.
For this construction, we will blow up an ideal (or quasi-coherent ideal presheaf) $I$
instead of a difference ideal (or difference ideal sheaf) $J$.
The two constructions coincide when $I=J$ is both a finitely
generated ideal and a difference ideal, so there is no confusion
in denoting both by $\widetilde{R}_I$ or $\widetilde{R}_J$.  However,
to emphasize the difference we will temporarily use a
superscript $\ ^c$  for the closed blowing-up of ideals.

\

Let $R$ be a difference ring, and let $I$ be a finitely generated ideal.
We define the closed blowing-up ring  $\widetilde{R}_I^c$ as a
 subring of $R[t^{\Nsi}]$:
$$\widetilde{R}_I^c = \sum_{n \in \Nsi} I^n t^n $$
where if $n = \sum_i \si^{k(i)}$, $I^n $ is the subgroup of $R$
generated by elements of the form $\Pi_i f_i$, with $f_i \in \si^{k(i)}(I)$.

 $\widetilde{R}_I$ inherits a $\Nsi$-grading from $R[t^{\Nsi}]$,
with $R$ the homogeneous component of degree $0$.

If $X = \spe R$, let
$$\widetilde{X}_{I}^c = \projs ( \widetilde{R}_I^c)$$
and let $\pi:  \widetilde{X}_{I}^c \to X$ be the natural map,
corresponding to the inclusion $R \to \widetilde{R}_I^c$.

\

If $X$ is any difference scheme, and ${\cal I}$ a quasi-coherent
ideal presheaf of ${\cal O}_X$ (considered forgetfully as a sheaf
of rings), one can check as above that the blowing ups of $\spe U$ at
the ideals ${\cal I}(U)$, and their canonical maps to $U$, glue together
to give a difference scheme over $X$, denoted $\widetilde{X}_{\cal I}^c$.

\begin{definition}  Let $X$ be a difference scheme, $I$
a quasi-coherent ideal presheaf.  Let $Y = \widetilde{X}_I^c$,
$\pi: Y \to X$ the structure map.    Then we obtain an
ideal presheaf $\pi^*I$ on $Y$.  We define the {\it exceptional divisor}
 to be the subscheme
of $Y$ defined locally by $\pi^*I$ (or equivalently by the difference ideal
sheaf generated by $\pi^*I$.)  \end{definition}

\end{paragraph} 

\begin{paragraph}{Comparing the blowing-ups}

We include a result comparing the two blowing-ups as the finitely generated ideal $I$
approaches $J$.     

When $J$ is already generated as a difference ideal by $I \meet \si^{-1}(I)$,
the
difference between the two blowing-ups is a
proper closed subscheme $Z$
of the exceptional divisor.  It
can be interpreted as follows. In either version of the blowing up,
a point of the exceptional divisor corresponds to a ``direction of
approach" to the difference subscheme defined by $J$.
 However in the closed blowing up, directions
in which $\si(y)/y$ approaches infinity as $y \to 0$  are allowed;
in the open blowing up they are not.

 If one blows up after applying
$M_q$, $I$ and $J$ become identified, and one obtains only ``directions of approach"
in which $\si(y)/y$ approaches $0$; thus these points are bounded away from $Z$.

There is more to be said here:

\begin{enumerate}
\item
As $I \to J$, the closed blowing up $ \widetilde{X}_I^c$ appears to change 
fairly gently.  Perhaps
the different $\widetilde{X}_I^c$ can be compared by  transformally birational radicial
 morphisms.  (cf. Definition
\ref{radicial}.) 

\item  Outside $Z$,  the closed blowing ups contains points representing
directions in which $\si(y)/y$ is finite but nonzero; upon applying $M_q$,
these points form a detachable divisor.  (Is it possible to get rid
of them before?)

\item  On the other hand it may be interesting precisely to study blowing
up  one-generated difference ideals.  Note that after $M_q$, such ideals
become principal, so blowing them up has no effect on smooth varieties.
For instance, it appears possible that one obtains a better definition
of "irreducible difference variety" by demanding that the irreducibility persist
to the strict transform under such principal blowing ups.
\end{enumerate}

\begin{prop} Let $R$ be a difference ring, $J$ a finitely generated difference ideal,

$V$ a set of generators of $J$, and $I$ the ideal generated by
$V \union \si(V)$.  Let $X = \spe R$.
Then $\widetilde{X}_J$ is isomorphic to an open subscheme
of $\widetilde{X}_I^c$. Specifically, let
$$Z(V) = \{p \in \projs (\widetilde{R}_I): Vt \subset p \}$$
 Then
$$ \widetilde{X}_J = \widetilde{X}_I^c   \setminus Z(V)             $$
\end{prop}

\proof Both difference schemes are obtained by gluing together
affine schemes $\spe S$, where for some $a \in V$,
$S $ is a subring of $R[1/a,1/\si(a),\ldots]$. Namely:
$$S = S_1 = \{\frac{c}{a^n}: n \in \Nn[\si], c \in I^n \}$$
according to the closed construction;
for the other, one checks that
$$S = S_2 = \{\frac{c}{a^n}: n \in \Nn[\si], c \in J^n \}$$

Note that $\spe (S_1)$ (respectively
$\spe(S_2)$)  embeds naturally as an open subscheme of $\widetilde{X}_I^c$ (resp.
$\widetilde{X}_J$).  In the first case,
by definition of the closed set $Z(V)$,
the union of these open subschemes over $a \in V$ is precisely
$\widetilde{X}_I^c \setminus Z(V)$.
in the second case, it is  $\widetilde{X}_J$, since $V$ generates
$J$ as a difference ideal.  We will now fix $a \in V$ and show that $S_1=S_2$.
The resulting isomorphisms of open subschemes are natural and glue together
to give the required isomorphism.

Clearly $S_1 \subset S_2$. For the other direction, it suffices
to check that $\frac{c}{a} \in S_1$ when $c \in J$.    
$\{c: \frac{c}{a} \in S_1$ forms an ideal of $R$.
$J$ is generated as an ideal by $\union_k I^{\si^k}$, so we may
take $c \in I^{\si^k}$ for some $k$, so that
 $\frac{c}{a^{\si^k}} \in S_1$.
But also $a^\si \in \si(V) \subset I $, so $\frac{a^\si}{a} \in S_1$,
hence by applying $\si$, $\frac{a^{\si^{i+1}}}{a^{\si^i}} \in S_1$
for $i < k$. Multiplying these $k+1$ elements, we obtain $\frac{c}{a} \in S_1$.

\end{paragraph} 

\begin{paragraph}{Blowing up and $M_q$-reduction}

\begin{lem} \lbl{5.25}
Let $R$ be a difference ring, $I$ an ideal.
Let $\widetilde{R}_{I}^c$ be the closed
blowing up ring. Let $S = M_q(R)$,
and $\bar I = M_q(I) = (I+J_q(R))/J_q(R) $. Let
$\widetilde{S}_{\bar I}^c = \sum_{n \in \Nn} ({\bar I}t)^n \subset S[t]$.

There exists a natural surjective homomorphism of graded rings
 $$   j_q: M_q( \widetilde{R}_I^c ) \to \widetilde{S}_{\bar I}^c $$
 \end{lem}

\proof Clearly $M_q(R[t^{\Nsi}]) \cong S[t]$; by restriction,
we get an isomorphism
$$\widetilde{R}_I^c \hsp /( J_q(R[t^{\Nsi}]) \meet \widetilde{R}_I^c) \cong
\widetilde{S}_{\bar I}^c$$
 Now clearly
$J_q(\widetilde{R}_{I}^c) \subset ( J_q(R[t^{\Nsi}]) \meet \widetilde{R}_{I}^c) $,
yielding the surjective homomorphism
$$M_q(\widetilde{R}_{I}^c ) = \widetilde{R}_{I}^c /   J_q(\widetilde{R}_{I}^c) \to
\widetilde{S}_{\bar I}^c $$
The naturality can be seen via the
universal property of $M_q$. \qed

\begin{lem}\lbl{5.26} Let $X=\spe R$ be an affine difference scheme, ${ I}$ an ideal of $R$,
${\bar   I} \subset M_q(R)$  the image of ${ I}$ under $M_q$.
 There exists a natural embedding of
$ \widetilde{M_q(X)}_{\bar   I}$ as a closed subscheme of $ M_q(\widetilde{X}_{I}^c) $.

Moreover, if ${V}$ is a sub-ideal of $ I$, generating the same difference
ideal, and $Z(V)= \{p \in \projs (\widetilde{R}_{I}^c): Vt \subset p \}$,
then the image of $ \widetilde{M_q(X)}_{\bar   I}$ is disjoint from
$M_q(Z(V))$.
 \end{lem}

\proof From   \ref{5.25} we obtain a map
$$j_q^*: \proj \widetilde{M_q(R)}_{\bar I}^c \to \proj M_q( \widetilde{R}_{I}^c )$$
This can be composed with  \ref{Mq-proj}.  The "moreover" is clear,
since if $Vt \subset j_q^{-1}(p)$ then $j_q(Vt) \subset p$.  However
$V$ and $I$ generate the same difference ideal, so if $\bar V$ denotes
the image of $V$ in $M_q(R)$, then $\bar V$ and $\bar I$ generate the same ideal.
Thus ${\bar V}t \not \subseteq p$ for a homogeneous ideal $p \in \widetilde{M_q(R)}_{\bar I}^c$.
\qed

\begin{cor}\lbl{5.27}  Let $X$ be a difference scheme,
$I$  a quasi-coherent ideal presheaf.
Let $\bar I$ be the $M_q$-image of $I$.
There exists a natural embedding of
$ \widetilde{M_q(X)}_{\bar   I}$ as a closed subscheme of
$M_q(\widetilde{X}_{I}) $.
\end{cor}

\proof This reduces to the local case, \ref{5.26}, by gluing. \qed

\begin{note} \ \end{note} \lbl{5.28} In \ref{5.27}, assume $I$ is a sub-presheaf
of a difference ideal sheaf $J$, and
$I \meet \si(I)$ generates $J$.
Then the image of $\widetilde{M_q(X)}_{\bar   I}$ in
$M_q(\widetilde{X}_{I}) $
is contained
in  the $M_q$-image of the open blowing up of $X$ at $J$.

\end{paragraph}

\begin{paragraph}{A geometric view of blowing up difference schemes}

When $X = \proj (R)$ is itself a projective or
multi-projective difference variety,
blowing-ups of $X$ have a natural multi-projective structure.
Suppose $R$ is graded by $\Nsi^k$.   Then $R[t^{\Nsi}]$ is
naturally graded by $\Nsi^{k+1}$; this induces a grading
of $\widetilde{R}_{I}^c$.

\   Let $I$ be a homogeneous  
ideal
of the ${\Nsi^{k}}$-graded difference ring $R$.
Let $\cal I$ be the corresponding ideal sheaf on $X = \proj ^{{\Nsi}^{k}} R$.
Then
$$\widetilde{X}_{\cal I} \cong \proj ^{\Nsi^{k+1}} \widetilde{R}_{I}^c $$
naturally.

Let $R$ be a difference ring; view it as a scheme $X = \spec R$, endowed
with an endomorphism.  Let $I=I(0)$ be a finitely generated ideal of $R$.
Let $B_1(R,I)$ be the affine coordinate ring of the blow-up of $\spec R$
at $I$, in the sense of algebraic geometry (\cite{fulton}).  Thus
$B_1(R,I)$ may be viewed as a subring of $R[t]$:

$$B_1(R,I) = R + \sum_{i \geq 1} I^n t^n \subset R[t]$$

We can also view $B_1(R)$ as an extension of $R$.
Let $I(1)$ be the ideal of $B_1(R,I)$
generated by $\si(I) \subset R \subset B_1(R,I)$.  Proceeding inductively, let

$$B_{n+1}(R,I) = B_1( B_n(R,I), B_n(R,I) \si^{n+1} I)$$

$B_1(R,I)$ is naturally $\Zz$-graded, with $R$ the homogeneous component
of degree $0$.  This grading inductively builds up to a $\Zz^n$-grading
on $B_n(R,I)$.  We let

$$X_n = Proj^{\Zz^n} (B_n(R,I))$$

On $X_n$ we have the exceptional divisor $E_n$ corresponding to the ideal
$\si^n(I) B_n(R,I)$.
Let $B_\infty(R,I)$ be the direct limit of the rings $B_n(R,I)$.  Each
of these rings is naturally embedded in $R[t,t^\si,\ldots]$, and $B_\infty(R,I)$
can be taken to be the union of the rings $B_n(R,I)$ within
$R[t,t^\si,\ldots]$.  While the individual $B_n(R,I)$ do not in general admit a
difference ring structure, $B_\infty(R,I)$ is clearly a difference subring
of $R[t,t^\si,\ldots]$, and we thus view it as a difference ring.
 
\begin{lem}  $\widetilde{R}_{I}^c \iso B_\infty(R,I)$ as difference ring extensions
of $R$.
\end{lem}

\proof  They actually coincide as subrings of $R[t,t^\si,\ldots]$.  \qed

$\widetilde{X}_{I}^c$ can thus be viewed as a limit of
the projective system of schemes $X \leftarrow X_1 \leftarrow X_2 \leftarrow \ldots$.
The exceptional divisor $E$ on $\widetilde{X}_{I}^c$ corresponds to the intersection
of the pullbacks of all the $E_n$.  Observe that the image of $E$ on
$X_n$ will usually  have unbounded codimension.

\begin{example} \  \end{example}Let $R = [\si] k [X_1,\ldots,X_n]$ be the affine
difference ring of affine $n$-space, with $k$ a difference field.
Then $\spe R$ is affine $n$-space over $k$.  However $\spec R$ is
an infinite product of schemes $Y_i$, each isomorphic to affine $n$-space
over $k$.  Suppose $I$ is an ideal of $R$ generated by an ideal
$I_0$ of $k[X_1,\ldots,X_n]$, corresponding to a subscheme $S_0$ of affine
space, and let $Z_i$ be the result of blowing up $Y_i$ at $I_0$, $EZ_n$
the exceptional divisor.  Then
$X_n $ can be identified with
$$ Z_1 \times Z_2 \times \ldots \times Z_n \times Y_{n+1} \times Y_{n+2} \times
   \ldots$$  \qed

In the above example, we used the following observation
 concerning ordinary blowing-up of
varieties:  when $S,T$ are schemes,
and $C$ a subscheme of $S$,  $\widetilde{(S \times T)}_{C \times T}$
can be identified with $\widetilde{S}_C \times T$.  

\end{paragraph}


\end{subsection}
\<{subsection}{Semi-valuative difference rings  }
\lbl{pinched-sec}

Valuative schemes of finite total dimension over $\Ak$ can be viewed as
  transformal analogs of smooth curves.  But   moving a finite scheme over
$\Pp^1$, we encounter singular curves as well.    Their ultraproducts
lead to considering semi-valuative  schemes.   We take a look at these
in the present subsection; it will not be required in the applications.  
Lemma \ref{semival} explains their potential role in a purely 
schematic treatment of transformal specializations.

An $m$-semi-valuative $A_0$-algebra $A$ is by definition a local $\kt$-algebra contained in
a Boolean -valued valuation ring $R$ over $A_0$, such that $\si^m(t)R \subset A$.  We will be
interested in this when $A_0 = \kt$ or when $A_0 = k[T]$ in characteristic $p$.  An ultraproduct
of the latter will be an instance of the former.   When $R$ is a transformal valuation domain, 
we say $A$ is an $m$-pinched transformal  valuation ring.

\<{definition} \lbl{vmult}  Let $R$ be a transformal valuation domain, $A$ a difference subring of $R$, $t \in R$.  If 
$A \meet t^mR \subset A$, we will say that $A$ is an $m$-pinched valuative domain (with respect to $(R,t)$.)

Let $X_0 = \spe R/tR$; we have a map $f: X_0 \to \spe A/tA$.  If $Y$ is a component of $\spe A/tA$,
write  $r.dim(X_0/Y)$ for the
relative reduced total dimension of $f^{-1}(Y) \subset X_0$ over $Y$.  
   The   {\em valuative multiplicity} of $Y$ is 
$v.mult(Y) = rk_{ram}(K) + r.dim(X_0/Y)$.  Finally, the valuative   dimension of $Y$ (viewed as 
a part of the specialization of $A$) is 
$vt.dim(Y) = r.dim(f^{-1}(Y)) + rk_{ram}(K) $.  
\>{definition}

\<{example}  A plane curve with a point pinched to order $m$, inside formal ring of normalization. \>{example}  More precisely, the local ring of the curve has the required property,
inside the local ring of the normalization, over $k[t]$ rather than $\kt$; an ultraproduct
of such rings leads to $\kt$.  We consider this class of examples 
in more detail.  
 
Consider curves $C \subset \Pp^n$   over an algebraically closed field $k$, $P \in C$.
$C$ may be reducible and singular.   (For simplicity, assume $C$ is reduced.)
 The normalization $\tC$ of $C$ is defined to be
the disjoint sum of the normalizations of the irreducible components of $C$.   Let
$A$ $A$ be the local ring of $C$ at $P$, $B$ the
product of the local rings of $\tC$ at the points above $P$. 

Let $T \in A$ be a non-zero-divisor, corresponding to the restriction of a linear birational map $\Pp^n \to \Pp^1$.

\<{example} \lbl{pinch} Let $C \subset \Pp^n$ be a curve  
of degree $\dd$.       Let $P \in C$, $A,B,\tC,T$ as above.  
Then  $T^{{\frac{3}{4}}\dd^2}  {B} \subset  {A}$.
 \>{example}   

(For our purposes here, we could equally well replace $A,B$ by their completions $\hat{A} ,\hat{B}$
for the $T$-adic topology; i.e. the conclusion  $T^{{\frac{3}{4}}\dd^2} \hat{B} \subset \hat{A}$ 
 is what we really need.)

\proof  We may assume $P = 0 \in \Aa^n \subset \Pp^n $,  $T$ a linear map on $\Aa^n$.  Let
$Y$ be a generic linear map on $\Aa^n$.  Projecting $C$ via $(T,S)$ to $\Aa^2$ we obtain a plane
curve of degree $d$, whose local ring $k[S,T] \meet A$   is if
anything smaller that $A$, while the local ring of $\tC$ does not change since the map
$(T,S)$ is birational on each component of $C$.      So we can assume
$C \subset \spec k[S,T]$ is a plane curve.  

Assume first that $C$ is irreducible.
 By \cite{HA} I ex. 7.2, the arithmetic genus $p_a(C)$ satisfies  $p_a(C) = {{\dd} \choose {2}} \leq
\dd ^2 / 2$.  
By \cite{HA} IV Ex. 1.8 (a),   since the genus of the normalization $\tC$ is non-negative, 
the $A$-module $B/A$ has length $ l \leq p_a(C)$.    Thus $A+T^i B = A+T^{i+1} B$ for some $i  < l$.
In $B/A$ , we have $T (T^i B/A ) = (T^i B/A ) $, so by Nakayama, $T^i B/A = 0$, and hence
$T^l B \subset A$.  

If $C$ is reducible (at $0$), it is a union of irreducible curves $C_i$, of degrees $\dd _i$,
 with $\sum \dd _i = \dd$.   More precisely,  $C$ is  cut out   by $\Pi_i f_i$,
$f_i \in k[S,T]$ relatively prime polynomials
of degrees $\dd_i$; $C_i$ is the curve cut out by $f_i$.  

Let $B_i$ be the local ring of the normalization $\tC_i$ of $C_i$.  Then $B = \Pi_i B_i$.  Let
$A_i$ be the local ring of $C_i$ itself.  We can view $A$ as a subring of $\Pi_i A_i$.  
We will show that $(T^{\dd ^2 / 4} \Pi_i A_i) \subset A$.  By the irreducible case,
 $T^{\dd^2 /2 } \Pi_i B_i \subset \Pi_i A_i$.  So $T^ {3 \dd^2 /4} B \subset A$.   It suffices therefore
to show that $T^{\dd ^2 / 4}  A_i \subset A $ for each $i$; i.e. that there exists $a_i \in A$
whose image in $A_j$ is $T^{\dd ^2 / 4}$, and whose image in every other $A_j$ is $0$.  Take
for instance $i=1$.

Let ${{f_1}'} = f_2 f_3 \cdot \ldots \cdot f_n  $, ${{\dd_1}'} = \deg({{f_1}'})$.  
  Note that $\dd_1+{{\dd_1}'} = \dd$, so
 $\dd_1 {{\dd_1}'} \leq {\frac{1}{4}} \dd ^2$.  Now the curves $(f_1),({{f_1}'})$ intersect in a subscheme of 
size $\leq \dd_1 {{\dd_1}'}$, by Bezout's theorem.  So $T^{{\dd ^2 / 4}} = rf_1 + r' f_1'$,
$r,r' \in k[S,T]$.  The element $r' f_1'$ of $A$ is the one we looked for:  it equals $T^{{\dd ^2 / 4}}$
module $f_1$, and equals $0$ modulo $f_j$ for $j \neq 1$.  \qed



Here is a view of some of the combinatorics of the above example. 
\<{lem} \lbl{semigroup} Let $E$ be a sub-semi-group of $\Nn$, generated by a finite $X$ with greatest element $b$.
Then for some $Y \subset \Zz/b\Zz$, for all $n \geq b(b-1)+1$, $n \in E$ iff $n \mod b \in Y$.  \>{lem}

\proof  Let $Y(i)$ be the image in $\Zz/b\Zz$ of $E \meet [ib+1,(i+1)b]$.  Clearly $\emptyset \subset Y(0) \subseteq Y(1) \subseteq \ldots$, so for some $i<b$ we have $ Y(k)=Y(k+1)$.  Now $X+Y(k) \subset Y(k+1) (\mod b)$, so $X + Y(k) = Y(k) (\mod b)$, and it follows that $E + Y(k)=Y(k)  (\mod b)$.  The claim
follows, with $Y=Y(k)$.  \qed

\<{rem}\lbl{val2.5}  In \ref{val2.4}, let $A$ be an $m$-pinched valuative subring of $R$, with respect to $t$, $t \in A$.   
Then $A / tA$ has total dimension $\leq d$ over $F$;  if equality holds, then $tA \meet K[0] = (0)$
(in fact $ tR \meet K[0]  = (0)$.   )
  \>{rem}

\proof  The inequality is immediate from \ref{flat2}; so we will concentrate on the case:
$tA \meet K[0] \neq (0)$, and show that $A/tA$ has total dimension $<d$.  (But the argument
also serves in general to show the weak inequality.)

Assume   $  tR \meet K[0] \neq (0)$.   The proof of \ref{val2.1}-\ref{val2.4} shows
that $R/t_mR$ has total dimension $< m+d$. (In general, it has total dimension $\leq m+d$.)
 
Since $t_mR   \subset A$, we have $(t_m R \meet A)^2 \subset t_m A$:
if $t_m r_i \in A$, then $(t_mr_1)(t_mr_2) = t_m (t_mr_1r_2) \in t_m (t_m R \meet A)$.
 So the difference
rings $A/t_m A$, $ A/ (t_m R \meet A)$ differ only by nilpotents, and thus have the same total dimension.
Now $A/(t_mR \meet A)$ embeds into $R/t_mR$, and so is well-mixed, and has total
dimension $<m+d$ (resp. $\leq m+d$.)  

Let $J_n$ be the smallest well-mixed algebraically radical ideal of $A$ containing $t_n$.  
So $J_m = t_mR \meet A$, by the above argument; thus $^A \sqrt J_m = (\sqrt J_m) \meet A$
is an algebraically prime ideal.

Let $A''$ be a finitely generated $F$-subalgebra of $A/tA$.  Lift the generators to $A/ t_mA$, and
let $A'$ be the $\Ft[m-1]$-subalgebra of $A/t_mA$ generated by them.  We have
seen that $A/t_mA$ has total dimension $<m+d$, hence so does $A'$, and thus also
$A' / \sqrt{0}_{A'} \leq A / J_m$.  In this ring, $t_m$ is not a $0$-divisor.  Thus the proof of 
\ref{flat2} takes over, and shows that $A'/tA'$ has total dimension $<d$. Hence so does $A''$.   
  \qed

The usefulness of \ref{val2.5} is limited by the fact that it does not apply to an
 arbitrary algebraic component of $A/tA$.  We can use valuation-theoretic rather
  than schematic multiplicities:

\<{lem}  \lbl{val2.6} In \ref{val2.4}, let $A$ be an $m$-pinched valuative subring of $R$,
 with respect to $t$, $t \in A$.   Let $Y$ be a component of $\spe A/t$.  Then 
$$ v.mult(Y) + r.dim(Y) \leq d$$
If equality holds, then $P=0$.  If every weakly Zariski dense algebraic component of $R/tR$ is Zariski dense, then
$Y$ is Zariski dense in $\spec R[0]$.                                         \>{lem}

\proof By the argument in \ref{val2.5}, $A/t_mA$, $A/ (t_mR \meet A)$ differ only by nilpotents;
they thus have the same reduced total dimension.  Similarly, the maps $A/t_mA \to A/tA$,
$A/(t_mR \meet A) \to A/ (tR \meet A)$ are isomorphisms on points, and also respect the 
reduced total dimension.  Thus   $ r.dim(A/tA) = r. dim(A/ (tR \meet A)) \leq r. dim( R / tR)$;
$r. dim(Y) + r.dim(X_0/Y) \leq r. dim(f^{-1}(Y))$.  So $vt.dim(Y) \leq vt.dim( f^{-1}(Y)) \leq
vt.dim(R/tR)$.  The lemma follows from \ref{val2.4}.  \qed

By an {\em $m$- semi-valuative } difference scheme (over $\Ak$), we will mean one one of the
 form $\spe A$, where $A$ is an $m$-semi-valuative $\kt$-algebra.  We will also say that the singularity of
$\spe A$ is of order $\leq m$.

(We think of $\spe A$ as a disjoint union of generalized curves over $\Ak$, with a pinching of order $m$
over $t=0$.)  

 $\Xvs '$ is the variant of $\Xvs$ defined using semi-valuative difference schemes,
rather than valuative ones.    ($2d+1$ -semi-valuative 
difference schemes will do; any larger number will lead to substantially the same scheme.  cf. \ref{pinch}.)
$\Xvs '$ has the same
points as $\Xvs$, but may be thicker at the edges. $\Xvs '$ can probably also be defined 
as the intersection of 
the closed projections of  of $({\hat{X}})_0$ as ${\hat{X}}$ ranges over all  open blow-ups of $X$,
with vertical components removed.    

\<{lem}\lbl{vs-lim-pinch} 
 Let $Z$ be an algebraic component   of $\Xvs '$.
 Then there exists a  $T = \spe  \bar{A}$, $ \bar{A}$ a pinched valuative domain over $\kt$,   and a 
morphism $ T \to X $ over $\Ak$,  such that $Z$ is contained in the image of $T_0$ in $X_0$.
      \>{lem}

\proof  Same as the proof of \ref{vs-lim},  except that we consider $h_j: R \to B_j \subset A_j \subset L_j$, with $T^{\si^{2d+1}}A_j \subset B_j$.  The condition is that $p$ contains $\meet_{j \in J} {h_j}^{-1}(tB_j)$, and in
the conclusion $p$ contains $h^{-1}(t \bar{B})$.   \qed

\<{lem}  \lbl{semival}    

Let $D  \subset k$ be a difference domain such that $X$ descends to $D$
($X' \times _D k = X$, $X'$ a difference  scheme over $D$.)   Let 
$Y^0$ (resp. $Y$) be a finitely presented 
difference scheme with   $\Xvs \subset Y \times _D k $.  Then for
some $a_0 \in K$, with $F = D [1/a_0]$, 
\begin{description}
  
\item[(i)] For all difference fields $L$ and all homomorphisms $h: F \to L$, if $X_h = X' \tensor_{\Ft} \Lt$,  

${(X_h)_{{ \rightarrow 0}}}' \subset  Y_h$ (as schemes.)
  \item[(ii)] For all sufficiently large $q$, all $h: F \to L=K_q$, if ${\bar X}= X_{\hT}$, for any reduced (but possibly reducible) curve $C$ over $L$ and map $C \to {\bar X}$,
$C \meet {\bar X}_0 \subset Y^0_h$ (as varieties.)

\item[ (iii) ]  For all sufficiently large $q$, all $h: F \to L=K_q$, if ${\bar X}= X_{\hT}$, for any reduced (but possibly reducible) curve $C$ over $L$ and map $C \to {\bar X}$,
$C \meet {\bar X}_0 \subset Y_h$ (as schemes.)
 
\end{description}

\end{lem}

\proof  The assumption that  $Y^0$ is finitely presented over $D$ means
that locally, on $\spe A \subset X$, $Y^0$ is defined by a finitely generated ideal $I$ (among  well-mixed ideals.) Let $r_1,\ldots,r_j$ be generators for $I$.  Then (i) and (ii)
amount  to showing that each $r_i$ lies in a certain ideal.  In (i), the 
condition is that $r_i$ should vanish on the image   $f(T_0)$ for any
  valuative $T$ and $f: T \to (Y^0)_h$.  In other words, given a map
$g: A  \to R$ with $R$ an $2d+1$-semi-valuative ring  over $\Lt$,  where $d$ is the total dimension of $X_t$,
that $g(r_i) \in t R$. If (i) fails, there are $F_\nu$ approaching $K$
   and $h_\nu: F \to L_\nu$, and $g_\nu: A \to R_\nu$, with $g(r_i) \notin t R_\nu$.
Taking an ultraproduct, we obtain $g_*: S \to R_*$ with $g(r_i) \notin t R_*$.
But this contradicts the definition of $\Xvs '$, and the assumption $\Xvs ' \subset Y^0$.

As for (ii) and (iii), let $A$ be the 
product of the local rings at $0$ of $C$, and let $R$ be the corresponding product
of the local rings of the the normalizations
of the components of $C$.  Each such local ring is a discrete valuation domain,
and can be viewed as a transformal valuation domain using the Frobenius
map $x \mapsto x^q$.  By Lemma \ref{pinch},  if $\dd$ is a projective degree of $C$,
then $t^{ {\frac{3}{4}} \dd^2 } R \subset A$.
Now by \ref{fin-proj},  
$\dd \leq O(1) q^d$.  So $\dd^2 \leq O(1)^2 q^{2d} \leq q^{2d+1}$ (for
large enough $q$.)  
 so $t^{\si^{2d+1}}R \subset A$. 
Continue as in (i).  

 \qed

\>{subsection} 

\end{section} 

\end{document}